\documentclass[reqno,english]{amsart}
\usepackage[T1]{fontenc}
\usepackage[utf8]{inputenc}
\usepackage[active]{srcltx}
\usepackage{color}
\usepackage{babel}
\usepackage{mathtools}
\usepackage{amstext}
\usepackage{amsthm}
\usepackage{amssymb}
\usepackage[unicode=true,
 bookmarks=false,
 breaklinks=false,pdfborder={0 0 1},backref=false,colorlinks=false]
 {hyperref}

\makeatletter

\numberwithin{equation}{section}
\numberwithin{figure}{section}
\theoremstyle{plain}
\newtheorem{thm}{\protect\theoremname}[section]
\theoremstyle{plain}
\newtheorem{lem}[thm]{\protect\lemmaname}
\theoremstyle{plain}
\newtheorem{prop}[thm]{\protect\propositionname}
\theoremstyle{plain}
\newtheorem{cor}[thm]{\protect\corollaryname}
\theoremstyle{remark}
\newtheorem{rem}[thm]{\protect\remarkname}

\usepackage{babel}
\usepackage{enumitem}
\usepackage{graphicx}
\usepackage{comment}
\usepackage{xcolor}

\providecommand{\theoremname}{Theorem}
\providecommand{\corollaryname}{Corollary}

\providecommand{\lemmaname}{Lemma}
\providecommand{\remarkname}{Remark}
\providecommand{\propositionname}{Proposition}

\makeatother

\providecommand{\lemmaname}{Lemma}
\providecommand{\propositionname}{Proposition}
\providecommand{\remarkname}{Remark}
\providecommand{\theoremname}{Theorem}

\begin{document}
\global\long\def\R{\mathbf{\mathbb{R}}}%
\global\long\def\C{\mathbf{\mathbb{C}}}%
\global\long\def\Z{\mathbf{\mathbb{Z}}}%
\global\long\def\N{\mathbf{\mathbb{N}}}%
\global\long\def\T{\mathbb{T}}%
\global\long\def\Im{\mathrm{Im}}%
\global\long\def\Re{\mathrm{Re}}%
\global\long\def\Hc{\mathcal{H}}%
\global\long\def\M{\mathbb{M}}%
\global\long\def\P{\mathbb{P}}%
\global\long\def\L{\mathcal{L}}%
\global\long\def\F{\mathcal{\mathcal{F}}}%
\global\long\def\s{\sigma}%
\global\long\def\Rc{\mathcal{R}}%
\global\long\def\W{\tilde{W}}%
\global\long\def\G{\mathcal{G}}%
\global\long\def\d{\partial}%
\global\long\def\mc#1{\mathcal{\mathcal{#1}}}%
\global\long\def\Right{\Rightarrow}%
\global\long\def\Left{\Leftarrow}%
\global\long\def\les{\lesssim}%
\global\long\def\hook{\hookrightarrow}%
\global\long\def\D{\mathbf{D}}%
\global\long\def\rad{\mathrm{rad}}%
\global\long\def\d{\partial}%
\global\long\def\jp#1{\langle#1\rangle}%
\global\long\def\norm#1{\|#1\|}%
\global\long\def\ol#1{\overline{#1}}%
\global\long\def\wt#1{\widehat{#1}}%
\global\long\def\tilde#1{\widetilde{#1}}%
\global\long\def\br#1{(#1)}%
\global\long\def\Bb#1{\Big(#1\Big)}%
\global\long\def\bb#1{\big(#1\big)}%
\global\long\def\lr#1{\left(#1\right)}%
\global\long\def\la{\lambda}%
\global\long\def\al{\alpha}%
\global\long\def\be{\beta}%
\global\long\def\ga{\gamma}%
\global\long\def\La{\Lambda}%
\global\long\def\De{\Delta}%
\global\long\def\na{\nabla}%
\global\long\def\fl{\flat}%
\global\long\def\sh{\sharp}%
\global\long\def\calN{\mathcal{N}}%
\global\long\def\avg{\mathrm{avg}}%
\global\long\def\bbR{\mathbf{\mathbb{R}}}%
\global\long\def\bbC{\mathbf{\mathbb{C}}}%
\global\long\def\bbZ{\mathbf{\mathbb{Z}}}%
\global\long\def\bbN{\mathbf{\mathbb{N}}}%
\global\long\def\bbT{\mathbb{T}}%
\global\long\def\bfD{\mathbf{D}}%
\global\long\def\calF{\mathcal{\mathcal{F}}}%
\global\long\def\calH{\mathcal{H}}%
\global\long\def\calL{\mathcal{L}}%
\global\long\def\calO{\mathcal{O}}%
\global\long\def\calR{\mathcal{R}}%
\global\long\def\Lmb{\Lambda}%
\global\long\def\eps{\varepsilon}%
\global\long\def\lmb{\lambda}%
\global\long\def\gmm{\gamma}%
\global\long\def\rd{\partial}%
\global\long\def\aleq{\lesssim}%
\global\long\def\chf{\mathbf{1}}%
\global\long\def\td#1{\widetilde{#1}}%
\global\long\def\sgn{\mathrm{sgn}}%

\title[Blow-up solutions to CM--DNLS]{Construction of smooth chiral finite-time blow-up solutions to Calogero--Moser
derivative nonlinear Schrödinger equation}
\begin{abstract}
We consider the Calogero--Moser derivative nonlinear Schrödinger
equation (CM-DNLS), which is a $L^{2}$-critical nonlinear Schrödinger
equation with explicit solitons, self-duality, and pseudo-conformal
symmetry. More importantly, this equation is known to be completely
integrable in the Hardy space $L_{+}^{2}$ and the solutions in this
class are referred to as \emph{chiral} solutions. A rigorous PDE analysis
of this equation with complete integrability was recently initiated
by Gérard and Lenzmann \cite{GerardLenzmann2022}.

Our main result constructs smooth, chiral, and finite energy finite-time
blow-up solutions with mass arbitrarily close to that of a soliton,
answering the global regularity question for chiral solutions raised
by Gérard and Lenzmann. The blow-up rate obtained for these solutions
is different from the pseudo-conformal rate. Our proof also gives
a construction of a codimension one set of smooth finite energy initial
data (but without addressing chirality) leading to the same blow-up
dynamics. Our blow-up construction in the Hardy space might also be
contrasted with the global well-posedness of the derivative nonlinear
Schrödinger equation (DNLS), which is another integrable $L^{2}$-critical
Schrödinger equation.

The overall scheme of our proof is the forward construction of blow-up
dynamics with modulation analysis.
We begin with developing a linear theory for the near-soliton dynamics.
We discover a nontrivial conjugation identity, which unveils a surprising
connection from the linearized (CM-DNLS) to the 1D free Schrödinger
equation, which is a crucial ingredient for overcoming the difficulties
from the nonlocal nonlinearity. Another principal challenge in this
work, the slow decay of the soliton, is overcome by introducing a trick
of decomposing solutions depending on topologies, which we believe
is of independent interest.
\end{abstract}

\author{Kihyun Kim}
\email{(current) kihyun.kim@snu.ac.kr; (previous) khyun@ihes.fr}
\address{(current) Department of Mathematical Sciences and Research Institute of Mathematics, Seoul National University, 1 Gwanak-ro, Gwanak-gu, Seoul 08826, South Korea; (previous) IHES, 35 route de Chartres, Bures-sur-Yvette 91440, France}
\author{Taegyu Kim}
\email{k1216300@kaist.ac.kr}
\address{Department of Mathematical Sciences, Korea Advanced Institute of Science and Technology, 291 Daehak-ro, Yuseong-gu, Daejeon 34141, Korea}
\author{Soonsik Kwon}
\email{soonsikk@kaist.edu}
\address{Department of Mathematical Sciences, Korea Advanced Institute of Science and Technology, 291 Daehak-ro, Yuseong-gu, Daejeon 34141, Korea}

\subjclass[2020]{35B44 (primary), 35Q55, 37K10, 37K40}
\keywords{Calogero-Moser derivative nonlinear Schrödinger equation, continuum
Calogero-Moser model, blow-up, modulation analysis, conjugation identity}

\maketitle
\tableofcontents{}

\section{Introduction}

We consider the \textit{Calogero--Moser derivative nonlinear Schrödinger
equation} \textit{\emph{(also known as a }}\textit{continuum Calogero--Moser
model}\textit{\emph{)}} 
\begin{align}
\begin{cases}
i\partial_{t}u+\partial_{xx}u+2D_{+}(|u|^{2})u=0,\quad(t,x)\in\mathbb{R}\times\mathbb{R}\\
u(0)=u_{0}.
\end{cases}\tag{CM-DNLS}\label{CMdnls}
\end{align}
Here, $u:I\times\bbR\to\bbC$, $D_{+}=D\Pi_{+}$, $D=-i\partial_{x}$,
and $\Pi_{+}$ denotes the Cauchy--Szeg\H{o} projection defined by
the Fourier multiplier with symbol $\chf_{\xi>0}$. This equation
has a non-local nonlinearity. Originally introduced in \cite{AbanovBettelheimWiegmann2009FormalContinuum},
it is only recently that the rigorous PDE analysis of \eqref{CMdnls}
has begun by Gérard and Lenzmann in \cite{GerardLenzmann2022}. The
principal goal of the present paper (Theorem~\ref{TheoremChiralBlowup})
is to provide the first rigorous construction of \emph{smooth} and
\emph{chiral} finite-time blow-up solutions to \eqref{CMdnls} of
finite energy, together with a sharp description of the blow-up dynamics.

\subsection{Calogero--Moser DNLS and first main result}

Equation \eqref{CMdnls} is derived in \cite{AbanovBettelheimWiegmann2009FormalContinuum}
as a formal continuum limit of \emph{classical Calogero--Moser
system} 
\begin{equation}
\frac{d^{2}x_{j}}{dt^{2}}=\sum_{k\neq j}^{N}\frac{8}{(x_{j}-x_{k})^{3}},\qquad j\in\{1,\dots,N\}.\label{eq:classicalCM}
\end{equation}
This discrete system is a Hamiltonian system that describes $N$ (pairwise)
interacting particles on a line with an inverse square potential,
and it is known to be completely integrable \cite{Calogero1971ClassicCalogerMoser,CalogeroMarchioro1974ClassicCalogerMoser,Moser1975ClassicCalogerMoser,OlshanetskyPerelomov1976Invent}.
The defocusing version of \eqref{CMdnls} was introduced earlier in
\cite{Pelinovsky1995IntermediateNLSInto} as a special case of the
so-called \textit{\emph{intermediate nonlinear Schrödinger equation}}
(intermediate NLS). A similar periodic version of \eqref{eq:classicalCM}
was investigated by Sutherland \cite{Sutherland1971ClassicCalogerMoser,Sutherland1972ClassicCalogerMoser}
and its periodic continuum limit equation has been recently studied
in \cite{Badreddine2024PAA,Badreddine2023CMDNLSTorusDefocusingTravelingwavearxiv}
under the name \textit{\emph{Calogero}}\emph{--}\textit{\emph{Sutherland
derivative NLS}}\emph{.}

Apart from solutions with data in usual Sobolev spaces $H^{s}(\mathbb{R})$,
a distinguished function space of data for \eqref{CMdnls} is the
\emph{Hardy--Sobolev space} $H_{+}^{s}(\mathbb{R})$: 
\[
H_{+}^{s}(\mathbb{R})\coloneqq\{f\in H^{s}(\mathbb{R}):\text{supp}\widehat{f}\subset[0,\infty)\}\qquad(s\in\bbR).
\]
$L_{+}^{2}(\mathbb{R})\coloneqq H_{+}^{0}(\mathbb{R})$ corresponds
to the Hardy space of holomorphic functions on the complex upper half-plane.
We will denote $H_{+}^{\infty}(\bbR)\coloneqq\bigcap_{k\in\bbN}H_{+}^{k}(\bbR)$.
The positive frequency condition $\text{supp}\widehat{f}\subset[0,\infty)$,
interpreted as a \emph{chiral condition} \cite{AbanovBettelheimWiegmann2009FormalContinuum},
is invariant under the flow of \eqref{CMdnls}. Moreover, the equation
admits an integrable structure in the Hardy space \cite{GerardLenzmann2022}.
We note that \eqref{CMdnls} is locally well posed both in $H^{s}(\bbR)$
and in $H_{+}^{s}(\bbR)$ for all $s>\frac{1}{2}$ \cite{MouraPilod2010CMDNLSLocal,GerardLenzmann2022}.

In the Hardy space $L_{+}^{2}(\bbR)$, complete integrability of \eqref{CMdnls}
is first manifested in the Lax pair structure 
\begin{equation}
\frac{d}{dt}\mathcal{L}_{\textnormal{Lax}}=[\mathcal{P}_{\textnormal{Lax}},\mathcal{L}_{\textnormal{Lax}}]\label{eq:LaxPair-evol}
\end{equation}
with the $u$-dependent operators $\mathcal{L}_{\textnormal{Lax}}$
and $\mathcal{P}_{\textnormal{Lax}}$ defined by \eqref{CMdnls}
\begin{align}
\mathcal{L}_{\textnormal{Lax}}=-i\partial_{x}-u\Pi_{+}\overline{u},\quad\text{and}\quad\mathcal{P}_{\textnormal{Lax}}=i\partial_{xx}+2iuD_{+}\overline{u}.\label{eq:LaxPair}
\end{align}
Second, there exists an explicit solution formula in the sense of
Gérard \cite{KLV2025CAMS}.
A Lax pair was first discovered in \cite{PelinovskyGrimshaw1995IntermediateDefocusingCMDNLSLaxPair}
for the intermediate NLS (which is a defocusing version of \eqref{CMdnls}),
and Gérard and Lenzmann \cite{GerardLenzmann2022} found the Lax pair
structure \eqref{eq:LaxPair-evol}--\eqref{eq:LaxPair} for \eqref{CMdnls}
in the Hardy--Sobolev space.\footnote{More precisely, \eqref{eq:LaxPair} above is a small modification
of what was presented in \cite{GerardLenzmann2022}, including a commuting
part; This formulation was presented in (the introduction of) \cite{KLV2025CAMS}.} Explicit solution formulas were previously found in the works of
Gérard and collaborators \cite{GerardGrellier2015ExplictFormulaCubicSzego1,Gerard2023BOequexplicitFormula,GerardPushnitski2024CMP}
for the cubic Szeg\H{o} equation and the Benjamin--Ono equation,
and the explicit formula for \eqref{CMdnls} was derived in \cite{KLV2025CAMS}
(for data in $u_{0}\in H_{+}^{\infty}(\bbR)\cap\langle x\rangle^{-1}L^{2}(\bbR)$
or $u_{0}\in L_{+}^{2}(\mathbb{R})$ with $M(u_{0})<2\pi$). See also
\cite{Badreddine2024PAA} for the periodic setting.


We briefly discuss symmetries and conservation laws of \eqref{CMdnls},
which are similar to those of the 1D quintic nonlinear Schrödinger
equation 
\begin{equation}
i\rd_{t}u+\rd_{xx}u+|u|^{4}u=0.\tag{NLS}\label{eq:NLS}
\end{equation}
We have \textit{time translation symmetry} 
\begin{align*}
u(t,x)\mapsto u(t+t_{0},x),\quad(t_{0}\in\mathbb{R})
\end{align*}
\textit{phase rotation symmetry} 
\begin{align*}
u(t,x)\mapsto e^{i\gamma}u(t,x),\quad(\gamma\in\mathbb{R})
\end{align*}
and \textit{translation symmetry} 
\begin{align*}
u(t,x)\mapsto u(t,x+x_{0}).\quad(x_{0}\in\mathbb{R})
\end{align*}
Associated to these are the conservation of energy, mass, and momentum:
\begin{gather}
\widetilde{E}(u)=\frac{1}{2}\int_{\mathbb{R}}\left|\partial_{x}u-i\Pi_{+}(|u|^{2})u\right|^{2}dx,\quad\text{(Energy)}\label{eq:intro-E-tilde}\\
M(u)=\int_{\mathbb{R}}|u|^{2}dx,\quad\text{(Mass)}\quad\widetilde{P}(u)=\Re\int_{\mathbb{R}}(\overline{u}Du-\frac{1}{2}|u|^{4})dx.\quad\text{(Momentum)}\nonumber 
\end{gather}
Moreover, it has \textit{Galilean invariance} 
\begin{align*}
u(t,x)\mapsto e^{icx-ic^{2}t}u(t,x-2ct),\quad(c\in\mathbb{R}).
\end{align*}
Of particular importance are \textit{$L^{2}$-scaling symmetry} 
\begin{align*}
u(t,x)\mapsto\lambda^{-\frac{1}{2}}u(\lambda^{-2}t,\lambda^{-1}x),\quad(\lambda>0).
\end{align*}
and \textit{pseudo-conformal symmetry} 
\begin{equation}
u(t,x)\mapsto\frac{1}{t^{1/2}}e^{i\frac{x^{2}}{4t}}u\left(-\frac{1}{t},\frac{x}{t}\right),\quad t>0,\label{eq:pseudo-conf-transf}
\end{equation}
or, in continuous form, 
\begin{align}
u(t,x)\mapsto\frac{1}{(1+at)^{1/2}}e^{i\frac{ax^{2}}{4(1+at)}}u\left(\frac{t}{1+at},\frac{x}{1+at}\right),\quad(a\in\mathbb{R}).\label{eq:pseudo-conf-conti}
\end{align}
Associated to these are the virial identities 
\begin{align*}
\frac{d}{dt}\int_{\bbR}|x|^{2}|u(t,x)|^{2}dx & =4\int_{\bbR}x\cdot\left(\Im(\ol u\rd_{x}u)-\frac{1}{2}|u|^{4}\right)dx,\\
\frac{d}{dt}\int_{\bbR}x\cdot\left(\Im(\ol u\rd_{x}u)-\frac{1}{2}|u|^{4}\right)dx & =4\td E(u).
\end{align*}
Equation \eqref{CMdnls} has two time-reversal symmetries: $u(t,x)\mapsto\ol u(-t,x)$
and $u(t,x)\mapsto\ol u(-t,-x)$. We call \eqref{CMdnls} \textit{$L^{2}$-critical}
(or \textit{mass-critical}) since the above scaling symmetry preserves
the $L^{2}$ norm (or $M(u)$). The presence of pseudo-conformal symmetry
is also special, as it typically exists only for mass-critical nonlinear
Schrödinger equations. Finally, if one considers chiral solutions,
some restrictions are necessary to preserve the chirality of solutions.
For example, $c\geq0$ is required for Galilean invariance, only the
second form $u(t,x)\mapsto\ol u(-t,-x)$ can be used as time-reversal
symmetry, and more importantly, pseudo-conformal symmetry \emph{does
not} preserve chirality.

A central role in the global dynamics of solutions is played by \emph{ground state,}, which is the nonzero minimizer of energy $\td E(u)$.
It is known from \cite{GerardLenzmann2022} that the ground state
is given by 
\begin{equation}
\mathcal{R}(x)=\frac{\sqrt{2}}{x+i}\in H_{+}^{1}(\mathbb{R})\quad\text{with}\quad M(\mathcal{R})=2\pi\quad\text{and}\quad\widetilde{E}(\mathcal{R})=0,\label{eq:solition CMDNLS}
\end{equation}
and it is unique up to scaling, phase rotation, and translation symmetries.
This $\calR(x)$ is a static (or time-independent) solution of \eqref{CMdnls}
and it is also chiral. We also call $\calR$ a \emph{soliton}. More
generally, all nonzero $H^{1}(\bbR)$ traveling wave solutions (i.e.,
the solutions of the form $u(t,x)=e^{i\omega t}\calR_{c,\omega}(x-ct)$
for some $\omega,c\in\bbR$) are given by the Galilean transformed
$\calR$ 
\[
u_{c}(t,x)=e^{i\frac{c}{2}x-i\frac{c^{2}}{4}t}\calR(x-ct),\qquad c\in\bbR,
\]
up to scaling, phase rotation, and translation symmetries \cite{GerardLenzmann2022}.
This solution is chiral if and only if $c\ge0$. Applying the pseudo-conformal
transform \eqref{eq:pseudo-conf-transf} to $\calR$, one obtains
an explicit finite-time blow-up solution 
\[
S(t,x)\coloneqq\frac{1}{t^{1/2}}e^{ix^{2}/4t}\mathcal{R}\left(\frac{x}{t}\right)\in L^{2}(\bbR),\qquad\forall t>0.
\]
It is important to note that $S(t)$ does \emph{not} have finite energy
(i.e., $S(t)\notin H^{1}(\bbR)$) and moreover, $S(t)$ is \emph{not
chiral}.

The ground state $\calR$ is believed to be a natural threshold for \emph{global
regularity}, i.e. the global existence of strong solutions. Gérard
and Lenzmann proved in \cite{GerardLenzmann2022} that any $H_{+}^{1}(\bbR)$-solutions
$u(t)$ with $M(u)\leq M(\calR)$ exist globally in time. On the other
hand, the solution $S(t)$ above provides an example of finite-time
blow-up for the data in $L^{2}(\bbR)$ with mass $M(\calR)$. However,
$S(t)$ does not belong to $H_{+}^{1}(\bbR)$ (so it is consistent
with the above global existence result) and determining the actual
threshold for global existence, say in $H_{+}^{1}(\bbR)$, has remained an intriguing open question \cite{GerardLenzmann2022}.

In this paper, we resolve this problem by showing that $M(\calR)$
is the actual threshold for global regularity, considered in $H_{+}^{k}(\bbR)$
for any $k\geq1$. More specifically, we construct finite-time blow-up
solutions with mass \emph{arbitrarily close} to $M(\calR)$ and initial
data belonging to $\mathcal{S}_{+}(\bbR)\coloneqq\mathcal{S}(\bbR)\cap L_{+}^{2}(\bbR)$
with $\mathcal{S}(\bbR)$ the Schwartz space. In particular, these
blow-up solutions are smooth, chiral and have finite energy. Moreover,
the blow-up rate of our solutions (that is, the spatial scale $\lmb(t)\sim(T-t)^{2}$)
is \emph{different} from the pseudo-conformal rate $T-t$. 
\begin{thm}[Smooth chiral finite-time blow-up solutions]
\label{TheoremChiralBlowup}For any $\epsilon>0$, there exists an
initial data $u_{0}\in\mathcal{S}_{+}(\bbR)$ with $M(\mathcal{R})<M(u_{0})<M(\mathcal{R})+\epsilon$
such that the corresponding forward-in-time solution $u(t,x)$ to
\eqref{CMdnls} satisfies the following property. 
\begin{enumerate}
\item (Finite-time blow-up) $u$ blows up in finite time $T=T(u_{0})\in(0,\infty)$. 
\item (Sharp description of the blow-up) There exist $\ell=\ell(u_{0})\in(0,\infty)$,
$\gamma^{*}=\gamma^{*}(u_{0})\in\mathbb{R}/2\pi\bbZ$, $x^{*}=x^{*}(u_{0})\in\mathbb{R}$,
and $u^{*}=u^{*}(u_{0})\in L_{+}^{2}$ such that 
\begin{align*}
u(t,x)-\frac{e^{i\gamma^{*}}}{\sqrt{\ell(T-t)^{2}}}\mathcal{R}\left(\frac{x-x^{*}}{\ell(T-t)^{2}}\right)\to u^{*}\quad\text{in }L_{+}^{2}\quad\text{as }t\to T.
\end{align*}
\end{enumerate}
\end{thm}

We will discuss Theorem~\ref{TheoremChiralBlowup} in more detail
in Section~\ref{subsec:Discussions}.

Our proof also gives a construction of a \emph{codimension one} set
of (smooth) initial data, but without addressing chirality, leading
to the same finite-time blow-up dynamics. We will state this result
after gauge transform; see Theorem~\ref{TheoremCodimension2Blowup}
below.

\subsection{Gauge transform, self-duality, and second main result}

Most of our analysis in the sequel will be performed after taking \emph{gauge transform} 
\[
v(t,x)=-\mathcal{G}(u)(t,x)\coloneqq-u(t,x)e^{-\frac{i}{2}\int_{-\infty}^{x}|u(t,y)|^{2}dy},
\]
which is the same gauge transform as the well-known derivative nonlinear
Schrödinger equation (\eqref{eq:DNLS} below). Note that it is a diffeomorphism
in $H^{s}(\mathbb{R})$ for every $s\geq0$. This transform loses
chirality, but it is particularly convenient, since it transforms \eqref{CMdnls}
into a Hamiltonian PDE with respect to \emph{the standard symplectic
form} $\omega(v,w)=\Im\int_{\bbR}v\ol w\,dx$. (The original equation
\eqref{CMdnls} is also Hamiltonian, but with respect to a non-standard
one.) The gauge transformed equation now reads 
\begin{align}
\begin{cases}
i\partial_{t}v+\partial_{xx}v+|D|(|v|^{2})v-\frac{1}{4}|v|^{4}v=0,\quad(t,x)\in\mathbb{R}\times\mathbb{R}\\
v(0)=v_{0}.
\end{cases} & \tag{\ensuremath{\mathcal{G}}-CM}\label{CMdnls-gauged}
\end{align}
Here $|D|$ is the derivative with the Fourier multiplier $|\xi|$. In
other words, one can write $|D|=\mathcal{H}\partial_{x}=\partial_{x}\mathcal{H}$
for nice functions, where $\mathcal{H}$ is the Hilbert transform,
given by \eqref{eq:HilbertFourierFormula}. This equation also enjoys
the same symmetries as \eqref{CMdnls}, including time translation,
phase rotation, translation, scaling, Galilean invariance, and pseudo-conformal
symmetries, all in identical forms. Similarly, the conserved energy,
mass, and momentum are given as 
\begin{align}
\begin{gathered}E(v)=\frac{1}{2}\int_{\mathbb{R}}|\partial_{x}v|^{2}dx-\frac{1}{2}|D|(|v(x)|^{2})|v(x)|^{2}+\frac{1}{12}|v(x)|^{6}dx,\\
M(v)=\int_{\mathbb{R}}|v|^{2}dx,\quad P(v)=\int_{\mathbb{R}}\Im(\overline{v}\partial_{x}v)dx.
\end{gathered}
\label{eq:intro energy gauge}
\end{align}
The virial identities read 
\begin{align*}
\frac{d}{dt}\int_{\bbR}|x|^{2}|v(t,x)|^{2}dx & =4\int_{\bbR}x\cdot\Im(\ol v\rd_{x}v)\,dx,\\
\frac{d}{dt}\int_{\bbR}x\cdot\Im(\ol v\rd_{x}v)\,dx & =4E(v).
\end{align*}
The analogue of the complete square form \eqref{eq:intro-E-tilde}
of energy is 
\begin{align}
E(v)=\frac{1}{2}\int_{\mathbb{R}}\left|\partial_{x}v+\frac{1}{2}\mathcal{H}(|v|^{2})v\right|^{2}dx.\label{eq:intro Energy1}
\end{align}
The static solution $\calR$ of \eqref{CMdnls} is transformed as
a static solution to \eqref{CMdnls-gauged} 
\begin{align}
Q(x)\coloneqq-\mathcal{G}(\mathcal{R})(x)=\frac{\sqrt{2}}{\sqrt{1+x^{2}}}\in H^{1}(\mathbb{R}),\quad M(Q)=2\pi,\quad E(Q)=0.\label{eq:solition gauged equation}
\end{align}
Note that we chose the minus sign in the transform $v=-\mathcal{G}(u)$
to make $Q$ positive.

The form \eqref{eq:intro Energy1} represents \emph{self-duality}.
Introducing the operator 
\begin{align}
\mathbf{D}_{v}f\coloneqq\partial_{x}f+\frac{1}{2}\mathcal{H}(|v|^{2})f,\label{eq:intro Dv}
\end{align}
the energy \eqref{eq:intro Energy1} can be rewritten as 
\begin{align}
E(v)=\frac{1}{2}\int_{\mathbb{R}}\left|\mathbf{D}_{v}v\right|^{2}dx.\label{eq:intro Energy2 self dual}
\end{align}
In analogy with \cite{Bogomolnyi1976}, we call the nonlinear operator
$v\mapsto\mathbf{D}_{v}v$ the \textit{Bogomol'nyi operator}, and
the soliton $Q$ solves the \emph{Bogomol'nyi equation} $\bfD_{Q}Q=0$.
Similarly to $\mathcal{R}$, $Q$ is a unique solution to the Bogomol'nyi
equation up to symmetries. The Hamiltonian structure of \eqref{CMdnls-gauged}
and the self-dual form of energy \eqref{eq:intro Energy2 self dual}
yield a factorization at the nonlinear PDE level \eqref{CMdnls-gauged}
(we also call it a \emph{self-dual form}): 
\begin{align}
\partial_{t}v=-iL_{v}^{*}\mathbf{D}_{v}v,\label{eq:CMdnls gauged SelfdualForm sec1}
\end{align}
where $L_{v}^{*}$ is the formal $L^{2}$-adjoint of the linearized
operator $L_{v}$ (at $v$) of the Bogomol'nyi operator $v\mapsto\mathbf{D}_{v}v$.

We end this subsection with our second main result, formally addressing
\emph{codimension one stability} of the blow-up dynamics in Theorem~\ref{TheoremChiralBlowup},
in the gauge-transformed side \eqref{CMdnls-gauged}. This result
will be proved by the same analysis as in Theorem~\ref{TheoremChiralBlowup}. 
\begin{thm}[Codimension one stability of blow-up dynamics (informal statement)]
\label{TheoremCodimension2Blowup}There exists a `codimension one'
set $\mathcal{O}$ in $H^2(\bbR)$ such that for any $v_{0}\in\mathcal{O}$,
the corresponding forward-in-time solution $v(t,x)$ to \eqref{CMdnls-gauged}
satisfies the following property. 
\begin{enumerate}
\item (Finite-time blow-up) $v$ blows up in finite time $T=T(v_{0})\in(0,\infty)$. 
\item (Sharp description of the blow-up) There exist $\ell=\ell(v_{0})\in(0,\infty)$,
$\gamma^{*}=\gamma^{*}(v_{0})\in\mathbb{R}$, $x^{*}=x^{*}(v_{0})\in\mathbb{R}$,
and $v^{*}=v^{*}(v_{0})\in H^{1}$ such that 
\begin{align}
v(t,x)-\frac{e^{i\gamma^{*}}}{\sqrt{\ell(T-t)^{2}}}Q\left(\frac{x-x^{*}}{\ell(T-t)^{2}}\right)\to v^{*}\quad\text{in }L^{2}\quad\text{as }t\to T.\label{eq:thm-v-decomp}
\end{align}
\end{enumerate}
\end{thm}

The statement of Theorem~\ref{TheoremCodimension2Blowup} is rather
informal as the description of codimension one set is not rigorously
given here. We will provide a rigorous restatement of Theorem~\ref{TheoremCodimension2Blowup}
in Theorem~\ref{thm:precise statement codimensionone blowup}. The
terminology `codimension one' can be justified when the set $\calO$
has a certain regularity, but we do not address any regularity issues
of the data set $\calO$ in this work.

\subsection{Previous results on Calogero--Moser DNLS}

Here we discuss earlier results of the PDE analysis on \eqref{CMdnls}.

We begin with local well-posedness. De Moura and Pilod \cite{MouraPilod2010CMDNLSLocal}
proved that the intermediate NLS (the defocusing version of \eqref{CMdnls})
and its infinite depth limit model are locally well posed in $H^{s}(\mathbb{R})$
for all $s>\frac{1}{2}$. It can be verified that \eqref{CMdnls}
is also locally well-posed in $H^{s}(\mathbb{R})$ for all $s>\frac{1}{2}$.
Gérard and Lenzmann \cite{GerardLenzmann2022} proved the local well-posedness
of \eqref{CMdnls} in Hardy--Sobolev spaces $H_{+}^{s}(\mathbb{R})$
for all $s>\frac{1}{2}$ using a method similar to that used in \cite{MouraPilod2010CMDNLSLocal}.

The global well-posedness for subthreshold and threshold data (i.e.,
$M(u_{0})\leq M(\mathcal{R})=2\pi$) in $H_{+}^{k}(\mathbb{R})$ for
all $k\in\mathbb{N}_{\geq1}$ was also established by Gérard and Lenzmann
\cite{GerardLenzmann2022}. For subthreshold data ($M(u_{0})<M(\mathcal{R})=2\pi$),
this follows from conservation laws and a sharp inequality related
to the Lax operator. The result for threshold data follows the argument
of \cite{Merle1993Duke} (which yields that any $H^{1}(\bbR)$ finite-time
blow-up solution with threshold mass must be $S(t)$) and $S(t)\notin H^{1}(\bbR)$
due to the slow spatial decay of $\calR(x)\sim\frac{1}{x}$. Subsequently,
Killip, Laurens, and Vi\c{s}an \cite{KLV2025CAMS}
proved local well-posedness in $L_{+}^{2}(\mathbb{R})$ for the subthreshold
data.

The dynamics beyond the threshold was also studied. Gérard and Lenzmann
\cite{GerardLenzmann2022} employed the Lax pair structure to construct
$N$-soliton solutions of the form 
\begin{align*}
u(t,x)=\sum_{k=1}^{N}\frac{a_{k}(t)}{x-z_{k}(t)}\in H_{+}^{1}(\mathbb{R}),\quad M(u_{0})=2\pi N,\quad\forall N\geq2,
\end{align*}
where the residues $a_{k}(t)\in\mathbb{C}$ and the pairwise distinct
poles $z_{k}(t)\in\mathbb{C}_{-}=\{z\in\bbC:\Im(z)<0\}$ for $1\leq k\leq N$
solve a complexified version of the classical Calogero--Moser system.
These $N$-soliton solutions blow up in infinite time in the sense
that $\|u(t)\|_{H^{s}}\sim|t|^{2s}$ as $|t|\to\infty$ for any $s>0$.
Recently, Hogan and Kowalski \cite{HoganKowalski2024PAA},
by adopting the explicit formula, showed the existence of \emph{possibly
infinite-time} blow-up solutions with mass arbitrarily close to the
threshold $M(\calR)$, but it is not clear whether their solutions blow
up in finite time or in infinite time.

In addition, Badreddine investigated the zero dispersion limit of \eqref{CMdnls} \cite{Badreddine2024ZeroDispersionLimit}.

\subsection{\label{subsec:Discussions}Discussions on main theorems}

\ \vspace{5bp}

\emph{1. Method and novelties.} We use the forward construction with
modulation analysis. The overall strategy of our finite-time blow-up
construction is based on the robust road map: \emph{the energy method
with repulsivity }\cite{RodnianskiSterbenz2010,RaphaelRodnianski2012,MerleRaphaelRodnianski2013Invention,RaphaelSchweyer2013CPAMHeat,RaphaelSchweyer2014AnalPDEHeatQuantized,MerleRaphaelRodnianski2015CambJMath}\emph{.} In our setting of the completely integrable equation, the repulsivity is observed from an inherited nonlinear structure, the higher-order conservation laws, as exploited in \cite{JeongKim2024arXiv}.
We will discuss this method in more detail in Section~\ref{subsec:Strategy}.
In practice, most of our analysis will be performed on the gauge transformed
equation \eqref{CMdnls-gauged}, which we find more convenient to
manipulate, thanks to the standard Hamiltonian structure. In a more
technical level, our analysis is inspired from the works \cite{KimKwon2020blowup,KimKwonOh2020blowup,Kim2022CSSrigidityArxiv}
(in the use of nonlinearly conjugated variables and conjugation identities)
for the self-dual Chern--Simons--Schrödinger equation (CSS).

One of our novelties lies in the discovery of (purely algebraic) new
conjugation identities related to the linearized operator around $Q$.
We find that this non-local linearized operator can be transformed,
by conjugating \emph{non-local operators}, into the \emph{free Schrödinger
operator $-i\rd_{xx}$}. This was very unexpected to us. This transform
is one of the crucial ingredients to overcome \emph{non-locality } of the problem, which is a main challenge of \eqref{CMdnls}.
Similar (but more apparent) conjugation identities were found in \cite{KimKwon2020blowup,KimKwonOh2020blowup}.

This transform also allows us to identify the \emph{repulsive structure}
(or monotonicity estimates) in the linearized dynamics. Although
it is well known that monotonicity for the 1D free Schrödinger flow
$\rd_{t}-i\rd_{xx}$ is delicate due to the zero resonance, our further
crucial trick is to make use of the higher-order conservation laws following from the integrability. For this purpose, we examine the equation of higher-order nonlinear variables that intrinsically arise from the integrable structure. Such nonlinear adapted derivatives were first introduced in the work on (CSS) \cite{KimKwonOh2020blowup}. The use of these conservation laws significantly simplifies the energy estimates required for the analysis. In \cite{JeongKim2024arXiv}, this idea was employed to construct quantized blow-up rates.

Our other novelty is related to another principal challenge of the
equation \eqref{CMdnls} (or \eqref{CMdnls-gauged}): the \emph{slow
spatial decay} of the soliton $Q\sim\frac{1}{x}$. Indeed, we find that
the formal dynamical picture (of blow-up) is similar to the 3D energy-critical
NLS, and it is well known that the radial ground state for energy-critical
NLS has \emph{the slowest} spatial decay in 3D, generating resonances
at zero of the associated linearized operator. The difficulties in
the linearized dynamics are resolved by the aformentioned conjugation
identities.

The slow decay also creates serious issues when using cut-offs on
the profiles. The nonlocal structure would make these issues even
worse. Our strategy to overcome this difficulty is to change the decompositions
of solution $v(t)$ \emph{depending on topologies} (or, on norms);
this enables us to perform the dynamical analysis without imposing
any cut-offs on the profiles. This strategy is explained in more detail
in Section~\ref{subsec:Strategy}. We think this strategy to be useful
when dealing with problems involving slow spatial decays.

The final challenge in our problem is to construct blow-up solutions
among \emph{chiral} solutions. Using the forward construction (and
a topological argument), we reduce the matter to construct a finite-dimensional set of chiral initial data, which is almost orthogonal
to the stable manifold. We construct this finite-dimensional set by
introducing special truncations that preserve chirality.

We conclude this comment with our use of non-standard refined modulation
parameters. At the technical level, this allows us to significantly
simplify profile construction (or expansion) in a similar spirit
to \cite{JendrejLawrieRodriguez2022ASENS}; we expand the blow-up
profiles only up to the linear order.

\vspace{5bp}
\emph{2. Comparison with (NLS).} Consider the mass-critical nonlinear Schr\"odinger equation
\begin{align}
	i\partial_tu+\partial_{xx}u+|u|^4u=0, \tag{NLS}
\end{align}
which shares the same symmetries as \eqref{CMdnls}. In particular, both models possess a ground state soliton $\mathcal{R}$ that plays a critical role in the global dynamics, and the threshold mass $M(\mathcal{R})$ serves as a natural boundary between scattering and blow-up. In the sub-threshold regime $M(u_0) < M(\calR)$, solutions to \eqref{eq:NLS} scatter \cite{KenigMerle2006Invent,Dodson2012JAMSd3,Dodson2015Adv}, and similar behavior is expected for \eqref{CMdnls}, though scattering has not been established yet. 

However, there are several structural differences between the two equations. One of these is that \eqref{CMdnls} is completely integrable, whereas \eqref{eq:NLS} is not. This integrable structure endows \eqref{CMdnls} with additional features, such as a Lax pair and an explicit solution formula, which have no analogue in \eqref{eq:NLS}.

The two equations also differ significantly in the structure of their ground state solitons. Unlike \eqref{eq:NLS}, where the ground state is not static and generates time-dependent traveling waves, the soliton $\calR$ in \eqref{CMdnls} is a static solution. Moreover, the linearized operator around $\calR$ in \eqref{CMdnls} admits a self-dual structure, which is absent in the case of \eqref{eq:NLS}.

This difference also affects the associated pseudo-conformal blow-up solution $S(t)$, which exists in both models with mass equal to the ground state. In \eqref{eq:NLS}, the soliton decays sufficiently fast so that $S(t)$ lies in the energy space $H^1(\bbR)$ and exhibits a finite-time blow-up at the threshold. In contrast, the slow spatial decay $\calR(x) \sim \frac{1}{x}$ in \eqref{CMdnls} implies that the corresponding $S(t)$ does not belong to $H^1(\bbR)$, and the global well-posedness holds at the threshold mass in $H^1(\bbR)$.

One of the notable differences arises in the nature of blow-up. The equation \eqref{eq:NLS} admits stable finite-time blow-up solutions that follow the so-called \emph{log-log} blow-up rate, which is nearly self-similar with a log-log correction \cite{MerleRaphael2005AnnMath,MerleRaphael2003GAFA,Raphael2005MathAnnalen,MerleRaphael2004Invent,MerleRaphael2006JAMS,MerleRaphael2005CMP}. This is possible since the energy is sign-indefinite.
In contrast, the finite-time blow-up solutions constructed for \eqref{CMdnls} are unstable and exhibit a different rate, $\lambda(t) \sim (T - t)^2$, which deviates substantially from the pseudo-conformal rate $\lambda(t) \sim T - t$.

\vspace{5bp}

\emph{3. Comparison with (DNLS).} Our finite-time blow-up construction
draws a strong contrast between the global dynamics of \eqref{CMdnls}
and that of \emph{derivative nonlinear Schrödinger equation} 
\begin{equation}
i\rd_{t}u+\rd_{xx}u+i\rd_{x}(|u|^{2}u)=0,\tag{DNLS}\label{eq:DNLS}
\end{equation}
which is another $L^{2}$-critical Schrödinger equation with an integrable
structure (the presence of Lax pairs) and has been extensively studied.
Some earlier results \cite{Ozawa1996,Wu2013APDE,Wu2015APDE,GuoWu2017DCDS}
addressed the global regularity of \eqref{eq:DNLS} for data with
$L^{2}$-mass below some thresholds. However, Bahouri and Perelman
recently proved in \cite{BahouriPerelman2022Invent} that \eqref{eq:DNLS}
is globally well-posed without any restrictions on $L^{2}$-mass
and their proof heavily relies on complete integrability. The contrast
between the global existence for this integrable equation \eqref{eq:DNLS}
and the finite-time blow-up for a non-integrable one \eqref{eq:NLS},
might have made the global regularity problem of \eqref{CMdnls} for
finite-energy chiral solutions (i.e., in $H_{+}^{1}(\bbR)$) more
dubious. Our Theorem~\ref{TheoremChiralBlowup} shows that finite-time
blow-up is indeed possible for the integrable equation \eqref{CMdnls}
(in Hardy space) as opposed to \eqref{eq:DNLS}.\vspace{5bp}

\emph{4. Comparison with (CSS).} The results and techniques in this
work are inspired from the works \cite{KimKwon2020blowup,KimKwonOh2020blowup,Kim2022CSSrigidityArxiv}
on the so-called self-dual Chern--Simons--Schrödinger equation (CSS).
This equation is similar to the mass-critical NLS in $\bbR^{2}$,
as they share analogous symmetries and conservation laws including
pseudo-conformal symmetry. (CSS) is not known to be completely integrable.

Distinguished common features of (CSS) and \eqref{CMdnls-gauged}
are \emph{self-duality} (cf. \eqref{eq:intro Energy2 self dual})
and \emph{non-locality} of the problem. Phenomenologically, these
two equations also share the so-called \emph{rotational instability}
of finite-time blow-up solutions (see the discussion below). The general
scheme of the proof in this paper is much inspired by \cite{KimKwon2020blowup,KimKwonOh2020blowup,Kim2022CSSrigidityArxiv},
such as the use of nonlinear conjugation identities, nonlinearly conjugated
adapted variable $w_{1}$, and so on.

However, there are notable distinctions in our proof. (i) Our discovery
of the conjugation identity \eqref{eq:intro-first-conj-idty} is a
lot more nontrivial than that of (CSS). (ii) We deal with slower spatial
decay of solitons; this slow decay problem is resolved by introducing
a trick of decomposing solutions depending on topologies. (iii) We
also introduce more efficient refined modulation parameters that enable
us to significantly simplify the profile ansatz. (iv) Finally, we
deal with nonradial solutions. There is also an interesting distinction
in the blow-up rates; this will be discussed in the next item.

Remarkably, it turns out that the method of nonlinear conjugation
identities from (CSS) unveils an unconditional Lax pair structure
for \eqref{CMdnls-gauged}; see Section~\ref{SubsectionLaxStructure}\vspace{5bp}

\emph{5. Finite energy blow-up solutions.} The blow-up rate $(T-t)^{2}$
of our constructed solutions is \emph{different} from the pseudo-conformal
rate $T-t$. This polynomial gap is also related to the slow spatial
decay of $Q$. For (CSS), the finite-time blow-up solutions constructed
in \cite{KimKwon2020blowup} (using fast decaying solitons) still
have the same pseudo-conformal rate, and those constructed in \cite{KimKwonOh2020blowup}
have a logarithmically corrected pseudo-conformal rate due to the
logarithmic failure of $S_{\mathrm{CSS}}(t)$ (the CSS analog of
$S(t)$) to belong to the energy space. In our case (for $S(t)$),
this failure is even stronger and leads to \emph{polynomial
gap} in the blow-up rates. We believe that pseudo-conformal (one-bubble)
blow-up is impossible for finite energy solutions, in analogy with
the radial case of \cite{KimKwonOh2022arXiv1}.\vspace{5bp}

\emph{6. Quantized rates and similarity to energy-critical equations.}
In view of the formal dynamical picture of blow-up (or the modulation
laws), we find that \eqref{CMdnls-gauged} is similar to the 3D energy-critical
NLS, or to the 3D critical semilinear heat equation in a more simplified
picture. The spatial scale $(T-t)^{2}$ in this paper formally corresponds
to the first rate in \cite{FilippasHerreroVelazquez2000,delPinoMussoWeiZhang2020HeatQuantized}
for the heat equation. We believe that (more unstable) higher blow-up
rates $(T-t)^{2k}$, $k=2,3,\dots$, called \emph{quantized rates},
can also be constructed for smooth chiral (say $H_{+}^{\infty}(\bbR)$)
initial data. We note that such rates have been constructed under even symmetry in \cite{JeongKim2024arXiv}.

For non-smooth data (say $u_{0}\in H^{s}(\bbR)$ in a limited range
of $s$), we believe that there is a continuum of admissible blow-up
rates; see \cite{OrtolevaPerelman2013,Schmid2023arXiv} for works
in 3D energy-critical NLS, which we believe are similar to \eqref{CMdnls-gauged}.
See also \cite{KriegerSchlagTataru2008Invent,KriegerSchlagTataru2009Duke,Perelman2014CMP}.
We do not know if \emph{chiral} finite-time blow-up solutions can
be constructed by the methods therein.\vspace{5bp}

\emph{7. Codimension 1 stability and rotational instability.} We explain
why the blow-up dynamics in Theorem~\ref{TheoremCodimension2Blowup}
is formally codimension 1 and discuss its instability mechanism. By
modulation analysis, we approximate the blow-up dynamics by a formal
dynamical system \eqref{eq:adiabatic-ansatz} and \eqref{eq:1. modulation first equ}
of six parameters that describe the motion of solitons. This system can
be integrated explicitly and the formulas of the scale $\lmb(t)$
and the phase rotation parameters $\gmm(t)$ are given by \eqref{eq:lmb-gmm-formal-law};
the typical solutions are given as follows. For a fixed time $T\in(0,+\infty)$
and a small free parameter $\eta_{0}\in\bbR$ that varies near $0$,
we have solutions 
\begin{equation}
\begin{aligned}\lmb(t) & =(t-T)^{2}+(\eta_{0}^{2}),\\
\gmm(t) & =\begin{cases}
0 & \text{if }\eta_{0}=0,\\
{\displaystyle \sgn(\eta_{0})\Big\{\arctan\Big(\frac{t-T}{|\eta_{0}|}\Big)+\frac{\pi}{2}\Big\}} & \text{if }\eta_{0}\neq0,
\end{cases}
\end{aligned}
\label{eq:intro-lmb-gmm-formula}
\end{equation}
where $\sgn(\eta_{0})$ denotes the sign of $\eta_{0}$. Now we see
that the blow-up scenario ($\lmb(t)\to0$) is unstable as it can happen
if and only if $\eta_{0}=0$. The analysis of the formal dynamical
system also says that \emph{non-vanishing of $\eta_{0}$} is the only
way of avoiding blow-up, so the blow-up scenario is codimension 1
stable. Moreover, the formula \eqref{eq:intro-lmb-gmm-formula} also
demonstrates the mechanism of avoiding blow-up; the spatial scale
$\lmb(t)$ contracts until it reaches the scale $(\eta_{0})^{2}$,
then \emph{the phase rotation parameter $\gmm(t)$ abruptly changes
in a short time scale $|t-T|\lesssim|\eta_{0}|$ by the fixed amount
of angle $\text{sgn}(\eta_{0})\pi$}, and finally $\lmb(t)$ starts
to grow (the spreading) as a backward evolution of the blow-up. This
instability mechanism is referred to as \emph{rotational instability}
and it was observed in \cite{KimKwon2019,Kim2022CSSrigidityArxiv}
as the instability mechanism of pseudo-conformal blow-up solutions
for the self-dual Chern--Simons--Schrödinger equation. A formal
observation of rotational instability goes back to \cite{vandenBergWilliams2013,MerleRaphaelRodnianski2013Invention}
for the Landau--Lifshitz--Gilbert equation.

\vspace{5bp}

\emph{8. Relation to explicit formula.} Killip, Laurens, and Vi\c{s}an
\cite{KLV2025CAMS} derived
a global-in-time explicit solution formula for \eqref{CMdnls} in
the sense of Gérard, for initial data $u_{0}\in H_{+}^{\infty}(\R)\cap\langle x\rangle^{-1}L^{2}(\R)$.
However, that solution $u(t,z)$ is in fact defined for all $t\in\bbR$
and $z\in\bbC_{+}=\{z\in\bbC:\Im(z)>0\}$ and the solution $u(t,x)$
for $t,x\in\bbR$ are obtained as the limit $z=x+iy\to x$. This $u(t,x)$
might be singular in time and may not be a \emph{strong} solution
(say, in $C_{t}H_{x}^{s}$) to \eqref{CMdnls}. Hence, the formula says that if a strong solution $u(t,x)$ exists in time interval $I$, then $u(t,x)$ must follow the explicit formula mentioned above. The finite-time blow-up construction in this
paper is consistent with the explicit formula, providing a negative
answer to the global \emph{regularity} question for large data.

\subsection{\label{subsec:Strategy}Strategy of the proof}

We use the notation collected in Section~\ref{subsec:Notation}.

We use modulation analysis. We start from the gauge-transformed equation
\eqref{CMdnls-gauged}. We decompose our solution $v(t)$ to \eqref{CMdnls-gauged}
as 
\begin{align}
v(t,x)=\frac{e^{i\gamma(t)}}{\lambda(t)^{\frac{1}{2}}}[Q+P(b(t),\eta(t),\nu(t))+\eps(t)]\left(\frac{x-x(t)}{\lambda(t)}\right),\label{eq:intro v decompose}
\end{align}
where $\lambda,\gamma,x,b,\eta,\nu$ are time-dependent parameters
(called \emph{modulation parameters}) such that $b,\eta,\nu$ small,
$Q+P(\cdot;b,\eta,\nu)$ is some modified profile with $P(\cdot;0,0,0)=0$,
and $\eps$ is the remainder (or, \emph{radiation}). The main steps
of the proof are as follows: 
\begin{enumerate}
\item Derive an effective dynamical system (called \emph{modulation laws})
of the modulation parameters $\lambda,\gamma,x,b,\eta,\nu$ and construct
the corrector profiles $P(\cdot;b,\eta,\nu)$, 
\item Prove sufficient decay estimates for $\eps$ to justify the modulation
laws. 
\end{enumerate}
In the first part, we follow the method of \emph{tail computations}
\cite{RaphaelRodnianski2012,MerleRaphaelRodnianski2013Invention,MerleRaphaelRodnianski2015CambJMath,KimKwonOh2020blowup,Kim2022CSSrigidityArxiv}.
In the second part, we will obtain decay estimates using monotonicity
estimates (or \emph{repulsivity}) and performing energy estimates
for certain higher-order derivatives of $\eps$.

This strategy was employed to address the forward construction of
blow-up dynamics in various contexts. To name a few relevant examples,
we refer to Rodnianski--Sterbenz \cite{RodnianskiSterbenz2010} and
Raphaël--Rodnianski \cite{RaphaelRodnianski2012} for critical wave
maps (and Yang--Mills), Raphaël--Schweyer \cite{RaphaelSchweyer2013CPAMHeat,RaphaelSchweyer2014AnalPDEHeatQuantized}
for critical harmonic map heat flows, and Merle--Raphaël--Rodnianski
\cite{MerleRaphaelRodnianski2013Invention,MerleRaphaelRodnianski2015CambJMath}
for critical Schrödinger maps and energy-supercritical NLS. The
strategy has also been implemented in a wider range of equations;
see, for example, \cite{HillairetRaphael2012AnalPDE,MartelMerleRaphael2014Acta,RaphaelSchweyer2014MathAnn,Collot2017AnalPDE,Collot2018MemAmer,KimKwon2020blowup,KimKwonOh2020blowup,Kim2022CSSrigidityArxiv}.
This list is not exhaustive. Among these works, the works \cite{KimKwon2020blowup,KimKwonOh2020blowup,Kim2022CSSrigidityArxiv}
on the self-dual Chern--Simons--Schrödinger equation are particularly
relevant to the present work.\vspace{5bp}

\emph{1. Linearized dynamics, linear conjugation identities, and repulsivity.}
Since the blow-up in this paper will be formed by the concentration
of soliton $Q$, understanding \emph{the linearized dynamics} of \eqref{CMdnls-gauged}
around $Q$, say 
\begin{equation}
\partial_{t}\eps+i\calL_{Q}\eps=0,\label{eq:intro-linearized-eqn}
\end{equation}
is crucial. As this linearized operator $i\calL_{Q}$ has not been
studied yet, we need to begin with the analysis of this operator.
The goal of this step is twofold: (i) to formally decompose a function
space (one may imagine $L^{2}$ at this moment although it is not
a rigorous choice) into \emph{two invariant subspaces} $N_{g}(i\calL_{Q})$
and $N_{g}(\calL_{Q}i)^{\perp}$, where $N_{g}(i\calL_{Q})$ is the
generalized kernel of $i\calL_{Q}$ and $N_{g}(\calL_{Q}i)^{\perp}$
is its symplectic orthogonal, and (ii) to observe \emph{repulsivity}
(or monotonicity estimates) in the infinite-dimensional space $N_{g}(\calL_{Q}i)^{\perp}$. 

Our analysis in Section~\ref{Section 3 Linearization} will show
that $\Lmb Q,iQ,Q_{x}$ and $ixQ,ix^{2}Q,\rho=Q^{-1}$ satisfy generalized
kernel relations \eqref{eq:GeneralizedKernelRelation}. Set 
\begin{align*}
N_{g}(i\calL_{Q}) & =\mathrm{span}_{\bbR}\{\Lmb Q,iQ,Q_{x},ixQ,ix^{2}Q,\rho\},\\
N_{g}(\calL_{Q}i)^{\perp} & =\{i\Lmb Q,Q,iQ_{x},xQ,x^{2}Q,i\rho\}^{\perp},
\end{align*}
where $\perp$ is defined with respect to the real $L^{2}$-inner
product $(f,g)_{r}=\Re\int_{\bbR}f\ol gdx$. These two spaces are
\emph{formally} invariant under the linearized flow. Moreover, $N_{g}(i\calL_{Q})$ and $N_{g}(\calL_{Q}i)$ are formally
transversal (the $6\times6$ matrix of real inner products formed
by the defining objects of $N_{g}(i\calL_{Q})$ and $N_{g}(\calL_{Q}i)$,
after suitable truncations, is invertible. See Lemma~\ref{LemmaZkTransversality}.) so that their direct sum $N_{g}(i\calL_{Q}) \oplus N_{g}(\calL_{Q}i)$
formally spans the whole function space. The six base elements of
$N_{g}(i\calL_{Q})$ explain why we have six modulation parameters
in the decomposition \eqref{eq:intro v decompose}, and we hope to
put $\eps$ in (a truncated version of) $N_{g}(\calL_{Q}i)^{\perp}$.
The invariance of these subspaces says that the dynamics of the modulation
parameters and $\eps$ are formally decoupled at the linear level.
The arguments in this paragraph are only valid at the formal level
due to the slow decay of $Q\sim\frac{1}{x}$.

In the infinite-dimensional space $N_{g}(\calL_{Q}i)^{\perp}$ where
$\eps$ formally belongs, we hope to observe \emph{repulsivity}
(or monotonicity estimates). As the operator $i\calL_{Q}$ itself
is nonlocal, only $\bbR$-linear, and also has generalized kernels,
it would be difficult in general to analyze $i\calL_{Q}$ directly.
In this work, motivated by the analysis in \cite{KimKwon2020blowup},
we will uncover a repulsive structure of $i\calL_{Q}$ in the space
$N_{g}(\calL_{Q}i)^{\perp}$ through \emph{conjugation identities},
which are purely algebraic operator identities.

The Hamiltonian structure of \eqref{CMdnls-gauged}, the self-dual
form of energy \eqref{eq:intro Energy2 self dual}, and $\bfD_{Q}Q=0$
yield the factorization of $\calL_{Q}$: 
\[
i\calL_{Q}=iL_{Q}^{\ast}L_{Q},
\]
where $L_{Q}$ is the linearized operator of the Bogomol'nyi operator
$v\mapsto\bfD_{v}v$ around $Q$. As in \cite{RodnianskiSterbenz2010},
we conjugate $L_{Q}$ to \eqref{eq:intro-linearized-eqn} to obtain
\[
(\rd_{t}+L_{Q}^{\ast}iL_{Q})(L_{Q}\eps)=0.
\]
A direct computation shows that the operator $H_{Q}$ defined by the
relation 
\[
iH_{Q}=L_{Q}^{\ast}iL_{Q}=i\Big(-\rd_{xx}+\frac{1}{4}Q^{4}-Q|D|(Q\cdot)\Big)
\]
becomes \emph{$\bbC$-linear} and self-adjoint. More remarkably, it
turns out that $iH_{Q}$ has another new factorization (which we call \emph{first linear conjugation identity}) 
\begin{equation}
L_{Q}^{\ast}iL_{Q}=iA_{Q}^{\ast}A_{Q}\qquad\text{with}\quad A_{Q}=\partial_{x}(x-\mathcal{H})\langle x\rangle^{-1},\label{eq:intro-first-conj-idty}
\end{equation}
and its supersymmetric conjugate satisfies the following identity
(which we call \emph{second linear conjugation identity}): 
\[
iA_{Q}A_{Q}^{\ast}=-i\rd_{xx}.
\]
As a result, we see that $A_{Q}L_{Q}\eps$ solves \emph{1D free
Schrödinger equation} 
\[
(\rd_{t}-i\rd_{xx})(A_{Q}L_{Q}\eps)=0.
\]
($L_{Q}\eps$ and $A_{Q}L_{Q}\eps$ are often called \emph{linear
adapted derivatives}.) This $A_{Q}L_{Q}$ is a correct conjugation
in the sense that the kernel of $A_{Q}L_{Q}$ (in an appropriate space)
is equal to $N_{g}(i\calL_{Q})$, and therefore $\eps\mapsto A_{Q}L_{Q}\eps$
is one-to-one on (a truncated version of) $N_{g}(\calL_{Q}i)^{\perp}$.
A similar process was performed in (CSS) \cite{KimKwon2020blowup}. However,
the discovery of \eqref{eq:intro-first-conj-idty} in our case seems
highly nontrivial compared to that of (CSS).

We need one more observation to identify the repulsivity of the flow
$\rd_{t}+i\calL_{Q}$ in $N_{g}(\calL_{Q}i)^{\perp}$. The flow $\rd_{t}-i\rd_{xx}$
already seems nice, but monotonicity estimates for this 1D flow is
delicate due to the resonance at zero, as is well known. Our additional input is the higher-order conservation laws that naturally arise from the integrability.

\vspace{5bp}

\emph{2. Nonlinear conjugation identity.} In our actual nonlinear
analysis, we view \eqref{CMdnls-gauged} as a system of equations
for $v(t)$ and its canonical higher-order derivatives. This idea
was first introduced in \cite{KimKwonOh2020blowup}, inspired by
the use of the Hasimoto transform \cite{ChangShatahUhlenbeck2000CPAM}
in the study of near-soliton dynamics for Schrödinger maps \cite{GustafsonKangTsai2007CPAM,GustafsonKangTsai2008Duke,GustafsonNakanishiTsai2010CMP}.

First, we introduce a new (nonlinearly conjugated) variable $v_{1}=\textbf{D}_{v}v$.
This nonlinear transform is motivated by $\bfD_{Q}Q=0$ and the
self-dual form of \eqref{eq:CMdnls gauged SelfdualForm sec1}; it
causes the modulated part $Q$ of \eqref{eq:intro v decompose} to
disappear. This induces a degeneracy in the $v_{1}$-variable in the
sense that the profile for $v_{1}$, approximately $\bfD_{(Q+P)}(Q+P)$,
has no zeroth-order term. This degeneracy will simplify the profile
for $v_{1}$. Moreover, thanks to the factorization in \eqref{eq:CMdnls gauged SelfdualForm sec1},
the evolution equation for $v_{1}$ also takes the simple form 
\[
\partial_{t}v_{1}+iH_{v}v_{1}=0
\]
with the self-adjoint operator $H_{v}=-\partial_{xx}+\frac{1}{4}|v|^{4}-v|D|(\overline{v}\cdot)$.
Notice that if $v=Q+\eps$, then we have $v_{1}\approx L_{Q}\eps$
at the linear level, but the evolution equation for $v_{1}$ is significantly
simpler than that for $L_{Q}\eps$.

In the following analysis, we will decompose the nonlinear solution
$v(t)$ and its nonlinearly conjugated variable (also called adapted nonlinear derivative) $v_{1}=\bfD_{v}v$ into the form 
\begin{align}
v(t,x)= & \frac{e^{i\gamma(t)}}{\lambda(t)^{\frac{1}{2}}}[Q+\widehat{\eps}(t,\cdot)]\Big(\frac{x-x(t)}{\lmb(t)}\Big)\label{eq:intro-v-decom-1}\\
= & \frac{e^{i\gamma(t)}}{\lambda(t)^{\frac{1}{2}}}[Q+P(b(t),\eta(t),\nu(t);\cdot)+\eps(t,\cdot)]\Big(\frac{x-x(t)}{\lmb(t)}\Big)\label{eq:intro-v-decom-2}
\end{align}
and 
\begin{equation}
v_{1}(t,x)=\frac{e^{i\gamma(t)}}{\lambda(t)^{\frac{3}{2}}}[P_{1}(b(t),\eta(t),\nu(t);\cdot)+\eps_{1}(t,\cdot)]\Big(\frac{x-x(t)}{\lmb(t)}\Big).\label{eq:intro-v1-decom}
\end{equation}
The profiles $P$ and $P_{1}$ will be constructed in the next step,
together with the derivation of formal modulation laws.\vspace{5bp}

\emph{3. Formal modulation laws and the profiles $P(\cdot;b,\eta,\nu)$
and $P_{1}(\cdot;b,\eta,\nu)$.} In this step, we derive an effective
dynamical system of modulation parameters leading to the blow-up dynamics
in our main theorems.

First, we renormalize the flow by introducing the renormalized spacetime
variables $(s,y)$ defined by 
\[
\frac{ds}{dt}=\frac{1}{\lambda^{2}(t)}\quad\text{and}\quad y=\frac{x-x(t)}{\lambda(t)}
\]
and the renormalized functions $w$ and $w_{1}$ defined by 
\[
w(s,y)=\lambda^{\frac{1}{2}}(t)e^{-i\gamma(t)}v(t,\lambda(t)y+x(t))|_{t=t(s)}\quad\text{and}\quad w_{1}=\D_{w}w.
\]
In view of \eqref{eq:intro v decompose} and the smallness of $b,\eta,\nu$,
the renormalized solution $w(s)$ stays in the vicinity of $Q$. The
evolution equations of $w$ and $w_{1}$ are given by \eqref{eq:w-equ}
and \eqref{eq:w1-equ}, respectively.

Now we assume (cf. \eqref{eq:intro-v-decom-2} and \eqref{eq:intro-v1-decom})
\[
w(s)\approx Q+P(b(s),\eta(s),\nu(s);\cdot)\quad\text{and}\quad w_{1}(s)\approx P_{1}(b(s),\eta(s),\nu(s);\cdot).
\]
Note that $P_{1}$ is not completely independent of $P$ due to the
compatibility condition $P_{1}\approx w_{1}=\bfD_{w}w\approx\bfD_{Q+P}(Q+P)$.
Next, we assume \emph{slow adiabatic ansatz} 
\begin{equation}
\frac{\lmb_{s}}{\lmb}+b=0,\qquad\gmm_{s}-\frac{\eta}{2}=0,\qquad\frac{x_{s}}{\lmb}-\nu=0,\label{eq:adiabatic-ansatz}
\end{equation}
and that $b_{s},\eta_{s},\nu_{s}$ are at least quadratic in $b,\eta,\nu$.
We substitute our ansatz in \eqref{eq:w-equ} and take only the linear
terms in $b,\eta,\nu$ from \eqref{eq:w-equ}; it turns out that the
linear terms of $P$ should satisfy certain generalized kernel relations
related to $i\calL_{Q}$. This motivates the choice 
\begin{equation}
P=-ib\frac{y^{2}}{4}Q-\eta\frac{1+y^{2}}{4}Q+i\nu\frac{y}{2}Q.\label{eq:intro-P}
\end{equation}
In fact, \eqref{eq:adiabatic-ansatz} can be understood as a nonlinear
adaptation of the structure of the generalized kernel $N_{g}(i\calL_{Q})$.
Then, the compatibility condition $P_{1}\approx\bfD_{Q+P}(Q+P)\approx L_{Q}P$
determines the linear terms of $P_{1}$ as 
\begin{equation}
P_{1}\approx-(ib+\eta)\frac{y}{2}Q+i\nu\frac{1}{2}Q.\label{eq:intro-P1}
\end{equation}
We remark that $P$ and $P_{1}$ \emph{do not involve any cut-offs}
as mentioned in Section~\ref{subsec:Discussions}. The profiles $P$
and $P_{1}$ do not belong to $L^{2}$ due to their spatial tails;
these profiles will be used only when decomposing solutions in
spaces allowing these tails (e.g., $\dot{H}^{2}$ and $\dot{H}^{1}$,
respectively). In the rigorous analysis, we will in fact use $P_{1}=-(ib+\eta)\frac{y}{2}Q+(i\nu+\mu)\frac{1}{2}Q$
with an additional parameter $\mu$. The reason is technical and we
will not discuss it here.

In order to close the formal dynamical system, we need to identify
the dynamics of $b,\eta,\nu$. We substitute our $P$- and $P_{1}$-ansatzes
\eqref{eq:intro-P}--\eqref{eq:intro-P1} into the $w_{1}$-equation
\eqref{eq:w1-equ}, and keep $b_{s},\eta_{s},\nu_{s}$ and the quadratic
terms in $b,\eta,\nu$. Here, we take advantage of the nonlinear adapted
derivative to simplify the computation. As a result, we obtain 
\[
\Big(b_{s}+\frac{3}{2}b^{2}+\frac{1}{2}\eta^{2}\Big)\frac{i}{2}yQ-(\eta_{s}+b\eta)\frac{1}{2}yQ+(\nu_{s}+b\nu)\frac{i}{2}Q=(\text{cubic and higher}).
\]
This nonlinear computation is surprising because \emph{the quadratic
terms only involve three special directions} $yQ,iyQ,iQ$, which can
be easily handled by requiring 
\begin{equation}
b_{s}+\frac{3}{2}b^{2}+\frac{1}{2}\eta^{2}=0,\quad\eta_{s}+b\eta=0,\quad\nu_{s}+b\nu=0.\label{eq:1. modulation first equ}
\end{equation}
The formal modulation laws are now given by \eqref{eq:adiabatic-ansatz}
and \eqref{eq:1. modulation first equ}.

We derive our blow-up dynamics using the formal modulation laws \eqref{eq:adiabatic-ansatz}
and \eqref{eq:1. modulation first equ}. First, we observe the conserved
quantities 
\begin{equation}
\frac{b^{2}+\eta^{2}}{\lmb^{3}},\quad\frac{\eta}{\lmb},\quad\frac{\nu}{\lmb}.\label{eq:intro-formal-conserved}
\end{equation}
The conservation of the first two quantities implies that \emph{blow-up
($\lmb\to0$) is possible only if $\eta/\lmb\equiv0$}. Assuming $\eta\equiv0$,
the dynamics of $\lmb$ is reduced to $\lmb\lmb_{t}+b=0$ with $b^{2}/\lmb^{3}$
conserved, which implies that the blow-up is accompanied with the
rate 
\[
\lmb(t)\sim(T-t)^{2}.
\]
See the last paragraph of Section~\ref{subsec:Formal-modulation-laws}
for the integration of \eqref{eq:adiabatic-ansatz} and \eqref{eq:1. modulation first equ}.

\vspace{5bp}

\emph{4. Trapped solutions.} In the previous steps, we have studied
the linearized dynamics and derived the formal modulation laws that lead
to blow-up. From now on, we consider the problem of constructing nonlinear
blow-up solutions.

For the nonlinear solution $v(t)$ and its nonlinear adapted derivative
$v_{1}=\bfD_{v}v$, we will use the decompositions \eqref{eq:intro-v-decom-1}--\eqref{eq:intro-v1-decom}
introduced in Step~2. We impose six orthogonality conditions on $\eps$
so that the mapping $\eps\mapsto A_{Q}L_{Q}\eps$ (in Step 1) is one-to-one.
In more technical terms, we choose local orthogonality conditions
to have coercivity estimates for $A_{Q}L_{Q}$.

Ideally, we hope to work only with the decompositions \eqref{eq:intro-v-decom-2}
and \eqref{eq:intro-v1-decom} because they provide more refined information
on the structure of solutions. However, since we did not truncate the
modified profiles (to avoid serious cut-off problems in the
dynamical analysis), these decompositions make sense only in local
$L^{2}$ spaces (or together with more regularity, say in $\dot{H}^{2}$
and $\dot{H}^{1}$, respectively). Then, the fact that the nonlinear
analysis cannot be closed in $\dot{H}^{2}$ forces us to use the (global)
$L^{2}$ space decomposition. For this $L^2$ space, we use a simpler and
rougher decomposition, \eqref{eq:intro-v-decom-1}. In other words, we decompose the
solutions \emph{depending on topologies}.

We will construct blow-up solutions by \emph{bootstrapping} some controls
on the modulation parameters and the radiation, together with a \emph{topological
argument}. The first set of controls is 
\begin{equation}
b(s)>0,\qquad\frac{b^{2}(s)}{\lmb^{3}(s)}\sim1,\qquad|\eta(s)|\ll b(s).\label{eq:intro-control-1}
\end{equation}
The positivity $b(s)>0$ is required to have $\lmb(s)$ decreasing,
the smallness $|\eta(s)|\ll b(s)$ is a condition that forces the
blow-up scenario in Step~3, and $b^{2}(s)/\lmb^{3}(s)$ will be almost
conserved. Here, the control $|\eta(s)|\ll b(s)$ cannot be bootstrapped
forward in time due to the formal conservation of $\eta(s)/\lmb(s)$
and $b(s)\sim\lmb^{3/2}(s)$; it will be achieved for nongeneric
initial data (hence adding codimension one) via a topological argument.
Another control we assume in the construction is 
\begin{equation}
|\nu(s)|\ll b^{1-}(s),\label{eq:intro-control-2}
\end{equation}
which also cannot be bootstrapped forward in time. At this moment,
our blow-up construction seems only codimension two stable, but applying
Galilean boosts to the blow-up solutions formally recovers the codimension
one stability of the blow-up dynamics. The reason for assuming \eqref{eq:intro-control-2}
is to simplify the blow-up analysis at the technical level. Our final
set of controls is imposed on the radiation, which we roughly state
here as 
\begin{equation}
\|\wt{\eps}(s)\|_{L^{2}}\leq\delta\ll1\quad\text{and}\quad\|\eps(s)\|_{\dot{H}^{2}}\lesssim\lmb^{2}(s).\label{eq:intro-control-3}
\end{equation}
Note that we use $\wt{\eps}(s)$ for the $L^{2}$ bound because we
decompose $v(t)$ depending on topologies. The decay estimate $\|\eps(s)\|_{\dot{H}^{2}}\lesssim\lmb^{2}(s)$
is equivalent to saying that the radiation part of $v(t)$ (i.e.,
the last term of \eqref{eq:intro-v-decom-2}) is bounded in $\dot{H}^{2}$
(as $t\to T$). As the final data of the radiation of $v(t)$ (i.e.,
$v(T)$ in some sense) is nontrivial in general, the decay estimate
$\|\eps(s)\|_{\dot{H}^{2}}\lesssim\lmb^{2}(s)$ is optimal, and one
cannot ask for stronger bounds.

We call a nonlinear solution a \emph{trapped solution} if it satisfies
the controls \eqref{eq:intro-control-1}--\eqref{eq:intro-control-3}
for all time. Once we develop the tools for controlling the dynamics
of parameters, it is not difficult to show that trapped solutions
blow up in finite time with sharp descriptions provided in the main
theorems. Henceforth, we explain how we control the dynamics and how
we construct trapped solutions.

When we control the modulation parameters, we introduce \emph{refined
modulation parameters} $\td b,\td{\eta},\td{\nu}$, which are carefully
designed corrections of $b,\eta,\nu$ motivated from $N_{g}(\calL_{Q}i)^{\perp}$
(at the level of the nonlinear adapted derivative $w_{1}$; see \eqref{eq:RefinedModulationDefinition}).
Compared to the very rough estimate $|b_{s}+\frac{3}{2}b^{2}|\lesssim\|\eps\|_{\dot{H}^{2}}\lesssim\lmb^{2}\sim b^{4/3}$,
these refined parameters enjoy much better time derivative estimates,
e.g., 
\[
\Big|\lmb^{3/2}\rd_{s}\Big(\frac{\td b}{\lmb^{3/2}}\Big)\Big|\lesssim b^{2+}.
\]
This fact allows us to close the bootstrap by working only in the
$H^{2}$-topology.

In order to close the bootstrapping, we need to propagate the smallness of $\eps$. For this, we utilize higher-order energy conservation laws. Using this argument in \cite{JeongKim2024arXiv}, the authors constructued quantized blow-up rates for \eqref{CMdnls-gauged} with even data. The conservation laws directly yield the following energy estimates \emph{without bootstrapping}:
\begin{align*}
    \|w_1(s)\|_{L^2}\lesssim_{v_0}\lambda(s),\quad \|\td\bfD_{w(s)}w_{1}(s)\|_{L^2}\lesssim_{v_0} \lambda^2(s) \quad\textnormal{for all}\quad s\geq s(0).
\end{align*}
In particular, since $w\approx Q$, the operator $\td\bfD_{w(s)}$ can be approximated by $\td\bfD_{Q}$. Combining this observation with the repulsive properties \eqref{eq:DQ BQ decomp} and \eqref{eq:BQBQstar equal I}, we obtain
\begin{align*}
    \|A_Qw_{1}\|_{L^2}=\|B_Q^*A_Qw_{1}\|_{L^2}=\|\td\bfD_Qw_{1}\|_{L^2}\lesssim \lambda^{2}.
\end{align*}
Finally, noting that $A_Qw_1\approx A_QL_Q\eps $ and applying the coercivity of $A_QL_Q$, we deduce the desired $\dot H^2$-level estimate for $\eps(s)$.

Finally, the two unstable conditions, $|\eta(s)|\ll b(s)$ and $|\nu(s)|\ll b^{1-}(s)$
for all time, are achieved by a topological argument in a similar
spirit of \cite{CoteMartelMerle2011TopologicalExample}. In more detail,
for any parameters $(\lmb_{0},\gmm_{0},b_{0},x_{0},\eps_{0})$ satisfying
the above controls at initial time, we prepare a two-dimensional set
of initial data that are almost tangent to the directions of $\eta$-
and $\nu$-variations, i.e., $(1+y^{2})Q$ and $iyQ$. We then use
the Brouwer fixed point theorem to prove the existence of special
initial data whose evolution satisfies the unstable conditions for
all time and hence become trapped. This constructs a codimension two
set of initial data leading to the desired blow-up. Applying the Galilean
invariance then improves codimension two to codimension one. This
gives Theorem~\ref{TheoremCodimension2Blowup}.

\vspace{5bp}

\emph{5. Chiral initial data}. To construct chiral finite-time blow-up
solutions, we need one more ingredient. The dynamical analysis (modulation
estimates and energy estimates) is the same, but we need to construct
\emph{chiral} modified profiles $\Rc_{b,\eta,\nu}\in H_{+}^{\infty}(\R)$
that well approximate $Q+P$ in the gauge transformed side. This is
achieved by introducing suitable truncations of the function $y\mapsto y$
whose multiplication preserves chirality. Moreover, by simply smoothing
out the low frequency parts of the profiles, one can find initial
data in $\mathcal{S}_{+}(\bbR)$. This gives Theorem~\ref{TheoremChiralBlowup}.

\noindent \vspace{5bp}

\noindent %
\mbox{%
\textbf{Organization of the paper.}~%
}In Section~\ref{sec:notation and preliminaries}, we collect the
notation and preliminary facts for the analysis. In Section~\ref{Section 3 Linearization},
we develop a linear theory of the linearized operator of \eqref{CMdnls-gauged}
around the soliton. In Section~\ref{SectionFormalBlowupLaw}, we
derive the formal modulation laws and define the modified profiles
$P$ and $P_{1}$. In Section~\ref{sec:Trapped-solution}, we set
up the bootstrap controls for the modulation parameters and the radiation
(leading to the definition of trapped solutions) and prove the main
theorems. The nonlinear estimates used in Section~\ref{sec:Trapped-solution}
are finally proved in Section~\ref{Section Proof of nonlinear estimates}.

\noindent \vspace{5bp}

\noindent %
\mbox{\textbf{Acknowledgements.~}}K. Kim's research was supported by Huawei Young Talents Programme at IHES, the New Faculty Startup Fund from Seoul National University, the POSCO Science Fellowship of POSCO TJ Park Foundation, and the National Research Foundation of Korea (NRF) grant funded by the Korea government (MSIT) RS-2025-00523523. T.~Kim and S.~Kwon are partially supported by the National Research Foundation of Korea, NRF-2019R1A5A1028324 and NRF-2022R1A2C1091499.

\section{Notation and preliminaries}

\label{sec:notation and preliminaries}

In this section, we record some preliminaries. We provide frequently
used notation, properties of the Hilbert transform, and also introduce
relevant Sobolev and adapted function spaces.

\subsection{\label{subsec:Notation}Notation and definitions}

For quantities $A\in\bbC$ and $B\geq0$, we write $A\lesssim B$
if $|A|\leq CB$ holds for some implicit constant $C$. For $A,B\geq0$,
we write $A\sim B$ if $A\lesssim B$ and $B\lesssim A$. If $C$
is allowed to depend on some parameters, then we write them as subscripts
of $\lesssim,\sim,\gtrsim$ to indicate the dependence.

We write $\langle x\rangle\coloneqq(1+x^{2})^{\frac{1}{2}}$. We define
the smooth even cut-off $\chi_{R}$ by $\chi_{R}(x)=\chi(R^{-1}x)$,
$\chi(x)=1$ for $|x|\leq1$, and $\chi(x)=0$ for $|x|\geq2$. We
also denote the sharp cut-off on a set $A$ by $\mathbf{1}_{A}$.
Denote by $\delta_{jk}$ the Kronecker-Delta symbol, i.e., $\delta_{jk}=1$
if $j=k$ and $\delta_{jk}=0$ if $j\neq k$.

Denote by $L^{p}(\mathbb{R})$ and $H^{s}(\bbR)$ the standard $L^{p}$
and Sobolev spaces on $\bbR$. As we work on $\R$, we often omit
$\R$ when there is no confusion.

The Fourier transform (on $\bbR$) is denoted by 
\begin{align*}
\mathcal{F}(f)(\xi)=\widehat{f}(\xi)\coloneqq\int_{\mathbb{R}}f(x)e^{-ix\xi}dx.
\end{align*}
The inverse Fourier transform is given by $\mathcal{F}^{-1}(f)(x)\coloneqq\tfrac{1}{2\pi}\int_{\R}\wt f(\xi)e^{ix\xi}d\xi$.
We denote $|D|$ by the Fourier multiplier with symbol $|\xi|$, that
is, $|D|\coloneqq\mathcal{F}^{-1}|\xi|\mathcal{F}$. We denote by
$\mathcal{H}$ the Hilbert transform: 
\begin{align}
\mathcal{H}f\coloneqq\left(\frac{1}{\pi}\text{p.v.}\frac{1}{x}\right)*f=\mathcal{F}^{-1}(-i\text{sgn}(\xi))\mathcal{F}f.\label{eq:HilbertFourierFormula}
\end{align}
Another equivalent formula for the Hilbert transform is 
\begin{align}
\mathcal{H}f(x)=\frac{1}{\pi}\int_{0+}^{\infty}\frac{f(x-y)-f(x+y)}{y}dy.\label{eq:HilbertIntegralFormula}
\end{align}
In view of \eqref{eq:HilbertFourierFormula}, we have $\partial_{x}\mathcal{H}f=|D|f=\mathcal{H}\partial_{x}f$
for $f\in H^{1}$. 
We denote by $\Pi_{+}$ the Cauchy--Szeg\H{o} projection from $L^{2}(\mathbb{R})$
onto the Hardy space $L_{+}^{2}(\mathbb{R})=\{f\in L^{2}(\bbR):\mathrm{supp}\wt f\subset[0,\infty)\}$:
\begin{align*}
\Pi_{+}f=\mathcal{F}^{-1}\textbf{1}_{\xi>0}\mathcal{F}f.
\end{align*}
Then, we have $\Pi_{+}=\frac{1}{2}(1+i\mathcal{H}).$ We denote $D=-i\rd_{x}$,
$D_{+}=D\Pi_{+}$, and $|D|=\calH\rd_{x}$. We denote the Hardy-Sobolev
space by $H_{+}^{s}\coloneqq H^{s}\cap L_{+}^{2}$, where $H^{s}$
is the usual Sobolev space based on $L^{2}$. We also denote the Schwartz
space $\mathcal{S}$ and $\mathcal{S}_{+}\coloneqq\mathcal{S}\cap L_{+}^{2}$.

We mostly use the \emph{real inner product} defined by 
\begin{align*}
(f,g)_{r}\coloneqq\Re\int_{\R}f\overline{g}dx.
\end{align*}
We also use the modulated form of a function $f$ by parameters $\lmb$,
$\gmm$, and $x$: 
\begin{equation}
[f]_{\la,\ga,x}\coloneqq\frac{e^{i\ga}}{\la^{1/2}}f\lr{\frac{\cdot-x}{\la}}.\tag{\ref{eq:notation renornmalized}}
\end{equation}
We denote by $\Lambda_{s}$ the $\dot{H}^{s}$-scaling generator in
$\R$ as 
\begin{align*}
\Lambda_{s}f & \coloneqq\frac{d}{d\lambda}\bigg|_{\lambda=1}\lambda^{\frac{1}{2}-s}f(\lambda\cdot)=(\frac{1}{2}-s+x\partial_{x})f,\\
\Lmb & \coloneqq\Lmb_{0}.
\end{align*}

Given a function $f$ on $\bbR$, $f_{o}$ denotes the odd part and
$f_{e}$ denotes the even part: 
\begin{align*}
f_{o}(x)=\frac{f(x)-f(-x)}{2},\quad\text{and}\quad f_{e}(x)=\frac{f(x)+f(-x)}{2}.
\end{align*}
Hence $f$ is odd if $f=f_{o}$ and $f$ is even if $f=f_{e}$.

We denote $\langle f\rangle_{k}$ by 
\begin{align*}
\langle f\rangle_{-k}(x)\coloneqq\sup_{0\leq j\leq k}\langle x\rangle^{-j}|\partial_{x}^{k-j}f(x)|.
\end{align*}
In addition to the standard Sobolev spaces $H^{s}$ and $\dot{H}^{s}$,
we use \emph{the adapted Sobolev spaces $\dot{\mathcal{H}}^{1}$ and
$\dot{\mathcal{H}}^{2}$} whose norms are defined by 
\begin{align}
\|f\|_{\dot{\mathcal{H}}^{1}}^{2} & \coloneqq\|\partial_{x}f\|_{L^{2}}^{2}+\left\Vert \langle x\rangle^{-1}f\right\Vert _{L^{2}}^{2},\tag{\ref{eq:DefinitionAdaptedH1norm}}\\
\|f\|_{\dot{\mathcal{H}}^{2}}^{2} & \coloneqq\|\partial_{xx}f\|_{L^{2}}^{2}+\left\Vert \langle x\rangle^{-1}\langle f\rangle_{-1}\right\Vert _{L^{2}}^{2}.\tag{\ref{eq:DefinitionAdaptedH2norm}}
\end{align}
Note that $\dot{\mathcal{H}}^{1}\hookrightarrow\dot{H}^{1}$ and $\dot{\mathcal{H}}^{2}\hookrightarrow\dot{H}^{2}$
but $\dot{\mathcal{H}}^{1}\cap L^{2}=H^{1},\quad\dot{\mathcal{H}}^{2}\cap L^{2}=H^{2}$.
We also note that $\dot{\mathcal{H}}^{1}$ and $\dot{\mathcal{H}}^{2}$
are compactly embedded in $L_{\textnormal{loc}}^{2}$ and $H_{\textnormal{loc}}^{1}$
respectively.

In the following, we record formulas for profiles and linear operators
(and their adjoints) appearing in Sections~\ref{Section 3 Linearization}--\ref{SectionFormalBlowupLaw}
for the reader's convenience. We will use profile functions $Q,P,\widehat{P}_{1},P_{1}$
: 
\begin{align}
Q=Q(x) & \coloneqq\frac{\sqrt{2}}{\sqrt{1+x^{2}}},\nonumber \\
P=P(x;b,\eta,\nu) & \coloneqq-ib\frac{x^{2}}{4}Q-\eta\frac{1+x^{2}}{4}Q+i\nu\frac{x}{2}Q,\tag{\ref{eq:DefinitionProfileP}}\\
\widehat{P}_{1}=\widehat{P}_{1}(x;b,\eta,\nu) & \coloneqq-(ib+\eta)\frac{x}{2}Q+i\nu\frac{1}{2}Q,\tag{\ref{eq:DefinitionProfileP1hat}}\\
P_{1}=P_{1}(x;b,\eta,\nu,\mu) & \coloneqq-(ib+\eta)\frac{x}{2}Q+(i\nu+\mu)\frac{1}{2}Q.\nonumber 
\end{align}

For linear operators (at $v$ or $Q$), we denote 
\begin{align}
\mathbf{D}_{v}f & \coloneqq\partial_{x}f+\frac{1}{2}\mathcal{H}(|v|^{2})f,\tag{\ref{eq:DefinitionDv}}\\
\widetilde{\mathbf{D}}_{v}f & \coloneqq\partial_{x}f+\frac{1}{2}v\mathcal{H}(\overline{v}f),\tag{\ref{eq:DefinitionDvtilde}}\\
L_{v}f & \coloneqq\partial_{x}f+\frac{1}{2}\mathcal{H}(|v|^{2})f+v\mathcal{H}(\text{Re}(\overline{v}f)),\tag{\ref{eq:DefinitionLv}}\\
L_{v}^{*}f & \coloneqq-\partial_{x}f+\frac{1}{2}\mathcal{H}(|v|^{2})f-v\mathcal{H}(\text{Re}(\overline{v}f)),\tag{\ref{eq:DefinitionLvstar}}
\end{align}
and 
\begin{align*}
\widetilde{L}_{Q}f\coloneqq\mathbf{D}_{Q}f+Q^{-1}\mathcal{H}\Re(Q^{3}f).
\end{align*}
Some expression of nonlinear parts are denoted by 
\begin{align*}
N_{v}(f)=f\mathcal{H}(\Re(\overline{v}f))+\frac{1}{2}(v+f)\mathcal{H}(|f|^{2}),
\end{align*}
and
\begin{align*}
	\td N_Q(f)
	=f\mathcal{H}(\Re(Qf))+\frac{1}{2}f\mathcal{H}(|f|^{2})+\frac{1}{4}yQ\mathcal{H}(yQ^2|f|^{2})
	+\frac{1}{4}Q\mathcal{H}(Q^2|f|^{2}). 
\end{align*}
A second-order operator $H_{v}$ is given by 
\begin{align}
H_{v}f\coloneqq -\partial_{xx}f+\frac{1}{4}|v|^{4}f-v|D|\overline{v}f.\tag{\ref{eq:DefinitionHv}}
\end{align}
Here, we use the convention that
\begin{align*}
    |D|\ol vf= |D|(\ol v f),
\end{align*}
so as to avoid excessive parentheses in similar expressions involving other operators such as $\mathcal{H}$.
Another important linear operator $A_{Q}$ is defined by 
\begin{align}
A_{Q}f\coloneqq\partial_{x}(x-\mathcal{H})\langle x\rangle^{-1}f.\tag{\ref{eq:DefinitionAQ}}
\end{align}
Then, its adjoint operator $A_{Q}^{*}$ with respect to $(\cdot,\cdot)_{r}$
is given by 
\begin{align*}
A_{Q}^{*}f=-\langle x\rangle^{-1}(x+\mathcal{H})\partial_{x}f.
\end{align*}
We also write $A_{Q}=\partial_{x}B_{Q}$ and $A_{Q}^{*}=-B_{Q}^{*}\partial_{x}$where
\begin{align}
B_{Q}f\coloneqq(x-\mathcal{H})\langle x\rangle^{-1}f,\quad B_{Q}^{*}f\coloneqq\langle x\rangle^{-1}(x+\mathcal{H})f.\tag{\ref{eq:DefinitionBQ}}
\end{align}

\subsection{Useful identities for $\mathcal{H}$}

We will use some commutator identities related to the Hilbert transform
$\mathcal{H}$. Here, we only record the formulas and postpone the
proofs to Appendix~\ref{AppendixPreliminaryProof}. 
\begin{lem}
\label{LemmaHilbertUsefulEquation} Let $f,g\in H^{\frac{1}{2}+}$.
Then, 
\begin{align}
fg=\mathcal{H}f\cdot\mathcal{H}g-\mathcal{H}(f\cdot\mathcal{H}g+\mathcal{H}f\cdot g)\label{eq:HilbertProductRule}
\end{align}
in a pointwise sense. 
\end{lem}

\begin{lem}
\label{LemmaCommuteHilbert} For $f\in\langle x\rangle^{-1}L^{2}$,
we have 
\begin{align}
[x,\mathcal{H}]f(x)=\frac{1}{\pi}\int_{\mathbb{R}}f(y)dy.\label{eq:CommuteHilbert}
\end{align}
If $f\in H^{1}$ and $\partial_{x}f\in\langle x\rangle^{-1}L^{2}$,
then we have 
\begin{align}
\partial_{x}[x,\mathcal{H}]f=[x,\mathcal{H}]\partial_{x}f=0,\quad\text{i.e.,}\quad[\mathcal{H}\partial_{x},x]f=\mathcal{H}f.\label{eq:CommuteHilbertDerivative}
\end{align}
The above equalities are satisfied in almost everywhere sense. 
\end{lem}

\eqref{eq:CommuteHilbert} implies that if $f\in\langle x\rangle^{-1}L^{2}$,
then $\mathcal{H}f\sim\tfrac{1}{x}$ at $|x|\to\infty$. This means
that even if $f$ is $C_{c}^{\infty}$, the decay of $\mathcal{H}f$
may not be good. The above identities can be used to exchange $\mathcal{H}$
and decay as follows. In fact, for $f\in L^{2}$, using $\langle x\rangle^{-2}(1+x^{2})=1$,
we have 
\begin{align}
\mathcal{H}(f)=\mathcal{H}\left(\langle x\rangle^{-2}f\right)+x\mathcal{H}\left(x\langle x\rangle^{-2}f\right)-\tfrac{1}{\pi}{\textstyle \int_{\R}x\langle x\rangle^{-2}fdx}.\label{eq:HilbertDecayExchange1}
\end{align}
Applying again \eqref{eq:CommuteHilbert} to \eqref{eq:HilbertDecayExchange1},
we deduce 
\begin{align}
\mathcal{H}(f)=\langle x\rangle^{2}\mathcal{H}\left(\langle x\rangle^{-2}f\right)-x\cdot\tfrac{1}{\pi}{\textstyle \int_{\R}\langle x\rangle^{-2}fdx}-\tfrac{1}{\pi}{\textstyle \int_{\R}x\langle x\rangle^{-2}fdx}.\label{eq:HilbertDecayExchange2}
\end{align}
In these two equalities \eqref{eq:HilbertDecayExchange1} and \eqref{eq:HilbertDecayExchange2},
we move the decay to $\mathcal{H}$. Reversing this calculation,
we can also take the decay out of $\mathcal{H}$.

\section{Linearization of \eqref{CMdnls-gauged} around $Q$}

\label{Section 3 Linearization}

In this section, we analyze the linearized operators around solitons.
We will examine the self-dual nature of the linearized operator $i\mathcal{L}_{Q}=iL_{Q}^{*}L_{Q}$,
and we outline the spectral properties of $\mathcal{L}_{Q}$ and $L_{Q}$
in Section~\ref{SubsectionLQ}. Following this, our aim is to identify
a \textit{conjugation identity} similar to the one discussed in previous work (\cite{KimKwon2020blowup}, \cite{KimKwonOh2020blowup}), emerging
from the non-commutativity of $[L_{Q},i]$ and insights gained from
the self-dual Chern--Simons--Schrödinger equation (CSS). In Section~\ref{SubsectionLaxStructure},
we will discuss efforts to connect the Lax structure introduced in
\cite{GerardLenzmann2022} with the nonlinear conjugation identity
in (CSS), and the challenges encountered. Section~\ref{SubsectionConjugationAQ}
introduces an alternative operator $A_{Q}$ and its conjugation identity
to address the issues highlighted in Section~\ref{SubsectionLaxStructure},
along with a modified zero-order operator $B_{Q}$ to further address
the limitations of $A_{Q}$. We will also discuss the coercivity for
$L_{Q}$ and $A_{Q}$ at the end of Section~\ref{SubsectionConjugationAQ}.

In our analysis throughout the paper, we handle the equation \eqref{CMdnls-gauged}
instead of the original \eqref{CMdnls}. We recall from the introduction
the complete square form of energy 
\begin{align*}
E(v)=\frac{1}{2}\int_{\mathbb{R}}\Big|\partial_{x}v+\frac{1}{2}\mathcal{H}(|v|^{2})v\Big|^{2}dx=\frac{1}{2}\int_{\bbR}|\bfD_{v}v|^{2}dx
\end{align*}
and the nonlinear Bogomol'nyi operator 
\begin{align*}
v\mapsto\bfD_{v}v\coloneqq\partial_{x}v+\frac{1}{2}\mathcal{H}(|v|^{2})v.
\end{align*}
We separately denote the operator $\mathbf{D}_{v}$ by 
\begin{align}
\mathbf{D}_{v}\coloneqq\partial_{x}+\tfrac{1}{2}\mathcal{H}(|v|^{2}).\label{eq:DefinitionDv}
\end{align}
Now, we linearize the Bogomol'nyi operator $v\mapsto\mathbf{D}_{v}v$.
We write 
\begin{align}
\mathbf{D}_{v+\eps}(v+\eps)=\mathbf{D}_{v}v+L_{v}\eps+N_{v}(\eps),\label{eq:DwLinearlize Lecompose}
\end{align}
where the linearized operator $L_{v}$ and the nonlinear part $N_{v}(\eps)$
are given by 
\begin{align}
L_{v}\eps & \coloneqq\partial_{x}\eps+\tfrac{1}{2}\mathcal{H}(|v|^{2})\eps+v\mathcal{H}(\text{Re}(\overline{v}\eps)),\label{eq:DefinitionLv}\\
N_{v}(\eps) & \coloneqq\eps\mathcal{H}(\Re(\overline{v}\eps))+\tfrac{1}{2}\eps\mathcal{H}(|\eps|^{2}).
\end{align}
The adjoint operator $L_{v}^{*}$ of $L_{v}$ with respect to the
real $L^{2}$-inner product is given by 
\begin{align}
L_{v}^{*}\eps & \coloneqq-\partial_{x}\eps+\tfrac{1}{2}\mathcal{H}(|v|^{2})\eps-v\mathcal{H}(\text{Re}(\overline{v}\eps)).\label{eq:DefinitionLvstar}
\end{align}
Equipped with $L_{v}$ and $L_{v}^{*}$, we derive the self-dual form
of \eqref{CMdnls-gauged}. 
\begin{lem}[Self-dual form of gauge transformed CM-DNLS]
\label{LemmaSelfdualform} \eqref{CMdnls-gauged} is equivalent to
\begin{align}
\partial_{t}v+iL_{v}^{*}\mathbf{D}_{v}v=0.\label{eq:CMdnlsSelfdualform}
\end{align}
\end{lem}

\begin{proof}
By the Hamiltonian form of \eqref{CMdnls-gauged}, we have $i\partial_{t}v=\nabla E$
where $\nabla$ is a functional derivative with respect to $(\cdot,\cdot)_{r}$.
By \eqref{eq:DwLinearlize Lecompose}, we have 
\begin{align*}
\text{Re}\int\frac{\delta E}{\delta v}\overline{\phi}dx=\frac{d}{d\eps}\bigg|_{\eps=0}E(v+\eps\phi) & =\frac{1}{2}\frac{d}{d\eps}\bigg|_{\eps=0}\int_{\R}|\mathbf{D}_{v}v+L_{v}\eps\phi+N_{v}(\eps\phi)|^{2}dx\\
 & =\text{Re}\int_{\R}\mathbf{D}_{v}v\overline{L_{v}\phi}dx\\
 & =\text{Re}\int_{\R}L_{v}^{*}\mathbf{D}_{v}v\overline{\phi}dx.
\end{align*}
Thus, we have \eqref{eq:CMdnlsSelfdualform} in a weak sense. In fact,
we can obtain this using \eqref{eq:HilbertProductRule}. From a
direct computation, we have 
\begin{align}
L_{v}^{*}\mathbf{D}_{v}v=-\partial_{xx}v-|D|(|v|^{2})v+\left(\tfrac{1}{4}|\mathcal{H}(|v|^{2})|^{2}v-\tfrac{1}{2}\mathcal{H}(\mathcal{H}(|v|^{2})|v|^{2})v\right).\label{eq:CMDNLS gauged Selfdualform}
\end{align}
Thanks to \eqref{eq:HilbertProductRule} with $f=g=|v|^{2}$, the
last term of \eqref{eq:CMDNLS gauged Selfdualform} becomes 
\begin{align*}
\left(\tfrac{1}{4}|\mathcal{H}(|v|^{2})|^{2}-\tfrac{1}{2}\mathcal{H}(\mathcal{H}(|v|^{2})|v|^{2})\right)v=\tfrac{1}{4}|v|^{4}v.
\end{align*}
\end{proof}
A self-dual form \eqref{eq:CMdnlsSelfdualform} at the nonlinear PDE
level was first introduced in the self-dual Chern--Simons--Schrödinger
equation \cite{KimKwon2019}.

Next, we linearize \eqref{CMdnls-gauged} at the static solution $Q$.
Denote the linearized operator by $i\calL_{Q}$ and then the linearized
equation is 
\[
\rd_{t}\eps+i\calL_{Q}\eps=0.
\]
From \eqref{CMdnls-gauged} we find a second-order operator formula
\begin{align*}
\calL_{Q}f=\partial_{xx}f+|D|(Q^{2})f+2Q|D|\Re(Qf)-\tfrac{1}{4}Q^{4}f-\Re(Q^{4}f).
\end{align*}
 On the other hand, using the self-dual form \eqref{eq:CMdnlsSelfdualform}
and $\bfD_{Q}Q=0$, we also arrive at the self-dual factorization
formula 
\begin{align}
i\mathcal{L}_{Q}=iL_{Q}^{*}L_{Q},\qquad & L_{Q}=\partial_{x}+\tfrac{1}{2}\mathcal{H}(Q{}^{2})+Q\calH(\text{Re}(Q\,\cdot)).\label{eq:LinearizedOperatorSelfdualform}
\end{align}
Note that $\mathcal{L}_{Q}$ and $L_{Q}$ are $\R$-linear. i.e. $[i,L_{Q}]\ne0\ne[i,\mathcal{L}_{Q}]$. 

\subsection{Generalized kernel of linearized operator $i\mathcal{L}_{Q}$}

\label{SubsectionLQ} In this subsection, we discuss the algebraic
structure of $i\mathcal{L}_{Q}$. More precisely, we identify the
(formal) generalized kernel of $i\mathcal{L}_{Q}$ for the preparation
of the modulation analysis.

First of all, from various symmetries, one can derive identities of
linearized operator $i\mathcal{L}_{Q}$ as follows; 
\begin{align}
i\mathcal{L}_{Q}(iQ) & =0, & \text{(phase rotation)}\nonumber \\
i\mathcal{L}_{Q}(\Lambda Q) & =0, & \text{(scaling)}\nonumber \\
i\mathcal{L}_{Q}(Q_{x}) & =0, & \text{(space translation)}\label{eq:generalizedKernel}\\
i\mathcal{L}_{Q}(ixQ) & =2Q_{x}, & \text{(Galilean boost)}\nonumber \\
i\mathcal{L}_{Q}(ix^{2}Q) & =4\Lambda Q. & \text{(pseudoconformal symmetry)}\nonumber 
\end{align}
For the sake of the reader's convenience, we explain how one can derive
the identities \eqref{eq:generalizedKernel} from continuous symmetries.
Each continuous symmetry is applied to $Q$ and then differentiates
it with respect to the symmetry parameters. Assume that we have a continuous
family of solutions $v^{(a)}(t,x)$ to \eqref{CMdnls-gauged} with
$v^{(0)}(t,x)=Q(x)$. For example, if we want to obtain an equality
related to Galilean symmetry, $v^{(a)}(t,x)$ would be expressed as
\begin{align*}
v^{(a)}(t,x)=e^{iax-ita^{2}}Q(x-2at)
\end{align*}
Substituting $v^{(a)}$ into the equation for \eqref{eq:CMdnlsSelfdualform},
we have 
\begin{align*}
\partial_{t}v^{(a)}+iL_{v^{(a)}}^{*}\mathbf{D}_{v^{(a)}}v^{(a)}=0,
\end{align*}
and differentiating with respect to $a$ at $a=0$, we obtain 
\begin{align*}
\partial_{t}(\partial_{a}v^{(a)})|_{a=0}+i\mathcal{L}_{Q}(\partial_{a}v^{(a)})|_{a=0}=0.
\end{align*}
Other identities in \eqref{eq:generalizedKernel} are obtained by
a process similar to that of each symmetry. (To derive the identity \eqref{eq:generalizedKernel}
for the pseudo-conformal symmetry, one may use the continuous version
\eqref{eq:pseudo-conf-conti}.)

We note that $iQ$, $\Lambda Q$, $Q_{x}$ belong to the $L^{2}$-kernel.
However, $ixQ$, $ix^{2}Q$ lack decay and are not found in $L^{2}$.
We understand from \eqref{eq:generalizedKernel} that they are merely
formal generalized kernel elements. Similarly, on the way, we will
continue to encounter technical issues due to the slow decay of $Q$.
A more detailed explanation of these issues will be provided at the
end of this section.

In order to rigorously state the kernel of $\mathcal{L}_{Q}$ and
$L_{Q}$, let us first discuss a function space where we will obtain
the kernel element. We will work with the \textit{adapted function
spaces} $\dot{\mathcal{H}}^{1}$ and $\dot{\mathcal{H}}^{2}$ with
the following norms; 
\begin{align}
\|f\|_{\dot{\mathcal{H}}^{1}}^{2} & \coloneqq\|\partial_{x}f\|_{L^{2}}^{2}+\left\Vert \langle x\rangle^{-1}f\right\Vert _{L^{2}}^{2},\label{eq:DefinitionAdaptedH1norm}\\
\|f\|_{\dot{\mathcal{H}}^{2}}^{2} & \coloneqq\|\partial_{xx}f\|_{L^{2}}^{2}+\left\Vert \langle x\rangle^{-1}\langle f\rangle_{-1}\right\Vert _{L^{2}}^{2}.\label{eq:DefinitionAdaptedH2norm}
\end{align}
$\dot{\mathcal{H}}^{1}$ and $\dot{\mathcal{H}}^{2}$ are defined
as closures under each norm of the Schwartz class. By their definitions,
one observes $\dot{\mathcal{H}}^{1}\subset\dot{H}^{1}$,$\dot{\mathcal{H}}^{2}\subset\dot{H}^{2}$
and $\dot{\mathcal{H}}^{1}\cap L^{2}=H^{1}$ and $\dot{\mathcal{H}}^{2}\cap L^{2}=H^{2}$.
The main motivation for these adapted function spaces is related to
the coercivity property of linear operators such as $L_{Q}$. In the
modulation analysis of this paper, we will perform a higher-order
energy estimate of the radiation part. We will use nonlinear higher-order variables that are close to the so-called \textit{adapted
derivative}s, $L_{Q}$ and $A_{Q}$ of the radiation. For this, we need to
redesign the function spaces at the $\dot{H}^{1}$ or $\dot{H}^{2}$ levels to
restore the coercivity properties of the adapted derivatives. See Section~\ref{SubsectionConjugationAQ}
for more details.

Due to self-duality, $\mathcal{L}_{Q}=L_{Q}^{*}L_{Q}$, we
observe that $\mathcal{L}_{Q}f=0$ if and only if $L_{Q}f=0$ in a
weak sense. Thus, it suffices to find the kernel of $L_{Q}$ instead
of $\mathcal{L}_{Q}$. We show that in the adapted function space
$\dot{\mathcal{H}}^{1}$, the kernel of $L_{Q}$ is determined from
the symmetry identities in \eqref{eq:generalizedKernel}. 
\begin{prop}[Kernel of $L_{Q}$]
\label{PropKernel mathcalLQ} We have $\textnormal{ker}L_{Q}=\textnormal{ker}\mathcal{L}_{Q}=\textnormal{span}_{\mathbb{R}}\{iQ,\Lambda Q,Q_{x}\}$
on $\dot{\mathcal{H}}^{1}$. 
\end{prop}

\begin{proof}
We will show $\ker L_{Q}\subset\textnormal{span}_{\mathbb{R}}\{iQ,\Lambda Q,Q_{x}\}$.
We recall that 
\begin{align}
iQ=i\tfrac{\sqrt{2}}{(1+x^{2})^{1/2}},\quad\Lambda Q=\tfrac{1-x^{2}}{\sqrt{2}(1+x^{2})^{3/2}},\quad Q_{x}=-\tfrac{\sqrt{2}x}{(1+x^{2})^{3/2}}.
\end{align}
We want to decompose the null space of $L_{Q}$ into the pure real
space and the pure imaginary space. Let $v=v_{1}+iv_{2}$ with real-valued functions $v_{1}$ and $v_{2}$. Then, $L_{Q}v=0$ becomes
\begin{align*}
L_{1}v_{1}\coloneqq\partial_{x}v_{1}+\tfrac{1}{2}\mathcal{H}(Q^{2})v_{1}+Q\mathcal{H}(Qv_{1})=0,\quad L_{2}v_{2}\coloneqq\partial_{x}v_{2}+\tfrac{1}{2}\mathcal{H}(Q^{2})v_{2}=0.
\end{align*}
We note that $L_{2}=\mathbf{D}_{Q}$ and $L_{1}=L_{Q}$ in real-valued function spaces. Thus, we reduce to showing 
\begin{itemize}
\item[(i)] $\ker L_{1}=\text{span}_{\mathbb{R}}\{\Lambda Q,Q_{x}\}$, 
\item[(ii)] $\ker L_{2}=\text{span}_{\mathbb{R}}\{Q\}$. 
\end{itemize}
For (ii), suppose $f\in\ker L_{2}$. Recalling the formula $\frac{1}{2}\mathcal{H}(Q^{2})=\frac{x}{1+x^{2}}$,
$f$ satisfies 
\begin{align*}
L_{2}f=f_{x}+\tfrac{x}{1+x^{2}}f=0.
\end{align*}
By an elementary computation, we deduce that $f=\frac{c}{\sqrt{1+x^{2}}}$
for $c\in\mathbb{R}$, which implies (ii).

To show (i), we first need to simplify the equation $L_{1}f=0$. We
rewrite $L_{1}f$ as 
\begin{align}
L_{1}f=\partial_{x}f+\tfrac{x}{1+x^{2}}f+\tfrac{2}{\sqrt{1+x^{2}}}\mathcal{H}\left(\tfrac{f}{\sqrt{1+x^{2}}}\right)=0.\label{eq:L11Equation}
\end{align}
Set $\frac{f}{\sqrt{1+x^{2}}}=g$. Then, \eqref{eq:L11Equation} becomes
\begin{align}
0=(1+x^{2})g_{x}+2xg+2\mathcal{H}(g)=\partial_{x}((1+x^{2})g)+2\mathcal{H}(g).\label{eq:L11Equation-2}
\end{align}
Then, we have 
\begin{align}
-\xi\partial_{\xi\xi}\widehat{g}+(\xi-2\text{sgn}(\xi))\widehat{g}=0.\label{eq:L11Equation-Fourier}
\end{align}
In addition, set $g=(1+x^{2})h$. By \eqref{eq:L11Equation-2}, $\widehat{h}$
satisfies 
\begin{align}
(-\xi\partial_{\xi\xi}+(\xi-2\text{sgn}(\xi)))(1-\partial_{\xi\xi})\widehat{h}=0\label{eq:L11Equation-FourierModify}
\end{align}
in a distribution sense. In addition, since $f\in\dot{\mathcal{H}}^{1}$,
we have 
\begin{align}
\|\widehat{h}\|_{H^{2}}<\infty.\label{eq:L11Equation-FourierModifyNormbound}
\end{align}
We note that the weak solutions $\widehat{h}$ are also classical
solutions. We separate the region $\xi>0$ and $\xi<0$, and solve
the fourth order ODE \eqref{eq:L11Equation-FourierModify} for each
region, and collect solutions $\widehat{h}(\xi)$ such that $\widehat{h}(\xi)\to0$
as $|\xi|\to\infty$. We first discuss when $\xi>0$. In view of $g=(1+x^{2})h$,
we first solve \eqref{eq:L11Equation-Fourier} to obtain two independent solutions
$\widehat{g}_{1}$ and $\widehat{g}_{2}$. Then, by solving $(1-\partial_{\xi\xi})\widehat{h}=\widehat{g}_{1}$
and $(1-\partial_{\xi\xi})\widehat{h}=\widehat{g}_{2}$, we can determine
all the solutions for \eqref{eq:L11Equation-FourierModify}. Since
$\xi>0$, we rewrite \eqref{eq:L11Equation-Fourier} as 
\begin{align}
-\xi\partial_{\xi\xi}\widehat{g}+(\xi-2)\widehat{g}=-(\partial_{\xi}-1)(\xi\partial_{\xi}-(1-\xi))\widehat{g}=0.\label{eq:L11Equation-Fourier decompose}
\end{align}
Solving \eqref{eq:L11Equation-Fourier decompose}, we have two solutions
$\widehat{g}_{1}=\xi e^{-\xi}$ and $\widehat{g}_{2}$ that satisfies
\begin{align*}
(\xi\partial_{\xi}-(1-\xi))\widehat{g}_{2}=e^{\xi}.
\end{align*}
Taking $\widetilde{\widehat{g}}$ by $\widetilde{\widehat{g}}=\widehat{g}_{2}(\xi)\cdot\xi^{-1}e^{\xi}$,
$\widetilde{\widehat{g}}$ satisfies 
\begin{align*}
(\partial_{\xi}\widetilde{\widehat{g}})e^{-\xi}\xi^{2}=e^{\xi}.
\end{align*}
Therefore, we conclude $\widehat{g}_{2}(\xi)\cdot\xi^{-1}e^{\xi}=\int_{1}^{\xi}(\xi^{\prime})^{-2}e^{2\xi^{\prime}}d\xi^{\prime}+\widehat{g}_{2}(1)\cdot e$.
Therefore, we conclude that $\widehat{g}_{2}(\xi)\sim e^{\xi}$ as
$\xi\to\infty$. This implies that there are no solutions $\widehat{h}$
such that $(1-\partial_{\xi\xi})\widehat{h}=\widehat{g}_{2}$ and
$\widehat{h}$ satisfies \eqref{eq:L11Equation-FourierModifyNormbound}
in the region $\xi>0$. Thus, it suffices to find the solution $\widehat{h}$
such that 
\begin{align*}
(1-\partial_{\xi\xi})\widehat{h}=C_{1}\widehat{g}_{1}(\xi)=C_{1}\xi e^{-\xi}.
\end{align*}
By a computation, we obtain 
\begin{align*}
\widehat{h}=C_{2}e^{\xi}+C_{3}e^{-\xi}+\tfrac{1}{4}C_{1}e^{-\xi}(\xi+\xi^{2})\quad\text{for }\xi>0,
\end{align*}
and $C_{2}=0$ due to \eqref{eq:L11Equation-FourierModifyNormbound}.
Similarly, we can apply the same argument for the region $\xi<0$,
and then we have 
\begin{align*}
\widehat{h}=C_{5}e^{\xi}+\tfrac{1}{4}C_{4}e^{\xi}(\xi-\xi^{2})\quad\text{for }\xi<0,
\end{align*}
From \eqref{eq:L11Equation-FourierModifyNormbound} with the Sobolev embedding,
we have $\widehat{h}\in C^{1}$. That is, we have $\widehat{h}(0+)=\widehat{h}(0-)$
and $\partial_{\xi}\widehat{h}(0+)=\partial_{\xi}\widehat{h}(0-)$.
This implies $C_{3}=C_{5}$ and $C_{4}=C_{1}-8C_{3}$. Thus, the solutions
to \eqref{eq:L11Equation-FourierModify} with \eqref{eq:L11Equation-FourierModifyNormbound}
become 
\begin{align*}
\widehat{h}(\xi)=\begin{cases}
C_{3}e^{-\xi}+\tfrac{1}{4}C_{1}e^{-\xi}(\xi+\xi^{2}), & \xi>0,\\
C_{3}e^{\xi}+\tfrac{1}{4}(C_{1}-8C_{3})e^{\xi}(\xi-\xi^{2}), & \xi<0.
\end{cases}
\end{align*}
However, since $L_{1}$ has a real-valued function as its domain
and $C_{1}$ and $C_{3}$ may be complex numbers, it is still not
enough to claim that the kernel dimension is $2$. Since $h$ is real
valued, we have $\widehat{h}(-\xi)=\overline{\widehat{h}(\xi)}$.
That is, we have $\overline{C_{3}}=C_{3}$ and $-\overline{C_{1}}=C_{1}-8C_{3}$.
Taking $C_{3}=c_{1}$ and $C_{1}=c_{3}+ic_{2}$ with real values $c_{1},c_{2},c_{3}$,
we deduce $c_{3}=4c_{1}$ and 
\begin{align}
\widehat{h}(\xi)=\begin{cases}
c_{1}e^{-\xi}+\tfrac{1}{4}(4c_{1}+ic_{2})e^{-\xi}(\xi+\xi^{2}), & \xi>0,\\
c_{1}e^{\xi}+\tfrac{1}{4}(-4c_{1}+ic_{2})e^{\xi}(\xi-\xi^{2}), & \xi<0,
\end{cases}\label{eq:LQKernelProof h hat}
\end{align}
with real numbers $c_{1}$ and $c_{2}$. Therefore, we conclude that the
kernel of $L_{1}$ is 2-dimensional and $\ker L_{1}=\text{span}\{\Lambda Q,Q_{x}\}$.
We note that the Fourier transforms of $(1+x^{2})^{-3/2}\Lambda Q$
and $(1+x^{2})^{-3/2}Q_{x}$ are 
\begin{align}
\mathcal{F}(\langle x\rangle^{-3}\Lambda Q)(\xi)=\tfrac{\pi}{4\sqrt{2}}e^{-|\xi|}(\xi^{2}+|\xi|+1),\quad\mathcal{F}(\langle x\rangle^{-3}Q_{x})(\xi)=i\tfrac{\pi}{4\sqrt{2}}\xi e^{-|\xi|}(1+|\xi|),\label{eq:Lambda Q Qsubx Fourier explicit formula}
\end{align}
and this coincides with \eqref{eq:LQKernelProof h hat}.

Thus, we conclude claims (i) and (ii), and we finish the proof. 
\end{proof}
Finally, we would like to check \emph{formally} the generalized kernel
of $i\mathcal{L}_{Q}$. Let $\rho=Q^{-1}$. From \eqref{eq:DefinitionLv},
we formally have 
\begin{align}
L_{Q}\rho=xQ+Q\mathcal{H}(1),\label{eq:LQrho formal}
\end{align}
and using the convention $\mathcal{H}(1)=0$, we have $L_{Q}\rho=xQ$.
In addition, we have $L_{Q}^{*}xQ=Q$, and this implies that 
\begin{align*}
i\mathcal{L}_{Q}Q^{-1}=i\mathcal{L}_{Q}\tfrac{1}{2}(1+x^{2})Q=iQ.
\end{align*}
Therefore, we have the generalized kernel relations of $i\mathcal{L}_{Q}$
as 
\begin{align}
 & \begin{split}i\mathcal{L}_{Q}Q_{x} & =0,\\
i\mathcal{L}_{Q}ixQ & =2Q_{x},
\end{split}
 & \begin{split}i\mathcal{L}_{Q}\Lambda Q & =0,\\
i\mathcal{L}_{Q}ix^{2}Q & =4\Lambda Q,
\end{split}
 & \begin{split}i\mathcal{L}_{Q}iQ & =0,\\
i\mathcal{L}_{Q}(1+x^{2})Q & =2iQ.
\end{split}
\label{eq:GeneralizedKernelRelation}
\end{align}
We also write down some algebraic identities for $L_{Q}$ for the
reader's convenience: 
\begin{align}
L_{Q}(ixQ)=\textbf{D}_{Q}(ixQ)=iQ,\quad L_{Q}(ix^{2}Q)=\textbf{D}_{Q}(ix^{2}Q)=2ixQ.\label{eq:LQ ixQ ix2Q}
\end{align}

\begin{rem}
\label{Remark H1 Warning remark} We note that $\mathcal{H}(1)$ in
\eqref{eq:LQrho formal} is not well defined as a single-valued function.
We can presume $\mathcal{H}(1)=0$ formally when following the computations
in the sequel. However, to be more precise, we will slightly avoid
these issues and rigorously justify the computations in practical
analysis in later sections. For related discussions, see Remark~\ref{RemarkLQ Well Definedness} and Section~\ref{sec:proof of lemma 5.17}
\end{rem}

\begin{rem}
\label{Remark Invariant genKernel and Transversality} We \emph{formally}
define $N_{g}(i\mathcal{L}_{Q})$ and $N_{g}(\mathcal{L}_{Q}i)^{\perp}$
by 
\begin{align*}
N_{g}(i\mathcal{L}_{Q}) & \coloneqq\textnormal{span}_{\mathbb{R}}\{iQ,\Lambda Q,Q_{x},ixQ,ix^{2}Q,(1+x^{2})Q\},\\
N_{g}(\mathcal{L}_{Q}i)^{\perp} & \coloneqq\{i(1+x^{2})Q,x^{2}Q,xQ,iQ_{x},i\Lambda Q,Q\}^{\perp}.
\end{align*}
We recognize that these two spaces remain invariant under linear
flow $\partial_{t}+i\mathcal{L}_{Q}$. Furthermore, the $6\times6$
matrix constructed from the $L^{2}$ inner products (after suitable
cut-offs) between $iQ,\Lambda Q,Q_{x},ixQ,ix^{2}Q,(1+x^{2})Q$ and
$i(1+x^{2})Q,x^{2}Q,xQ,iQ_{x},i\Lambda Q,Q$ is found to have a nonzero
determinant. In essence, $N_{g}(i\mathcal{L}_{Q})$ and $N_{g}(\mathcal{L}_{Q}i)$
are transversal. See Lemma~\ref{LemmaZkTransversality} for a precise statement. Therefore, the direct sum 
\begin{align*}
N_{g}(i\mathcal{L}_{Q})\oplus N_{g}(\mathcal{L}_{Q}i)^{\perp}
\end{align*}
formally spans the whole function space. Consequently, linearized
evolution is divided into parts belonging to $N_{g}(i\mathcal{L}_{Q})$
and $N_{g}(\mathcal{L}_{Q}i)^{\perp}$. Motivated by this, we consider
a decomposition of the nonlinear solution in the following form, 
\begin{align*}
v(\cdot)=\frac{e^{i\gamma}}{\lambda^{1/2}}[Q+P(b,\eta,\nu)+\eps]\left(\frac{\cdot-x}{\lambda}\right).
\end{align*}
In view of the generalized kernel relations, we choose $P$ to satisfy
$\partial_{b}P=-i\frac{y^{2}}{4}Q$, $\partial_{\eta}P=-\frac{1+y^{2}}{4}Q$,
and $\partial_{\nu}P=i\frac{y}{2}Q$. Concretely, we will define $P(b,\eta,\nu)$
by the linear part in the modified profile, i.e., 
\begin{align}
P(b,\eta,\nu)\coloneqq-ib\tfrac{y^{2}}{4}Q-\eta\tfrac{1+y^{2}}{4}Q+i\nu\tfrac{y}{2}Q.\label{eq:DefinitionProfileP}
\end{align}
\end{rem}

\subsection{Conjugation identity and Lax pair structure}

\label{SubsectionLaxStructure} In this subsection, we present some
algebraic features of \eqref{CMdnls-gauged}. From the experiences
from earlier work, for example, the Chern--Simons--Schr\"odinger equation, we
attempt to induce a repulsive dynamics of $\varepsilon$ by finding
a conjugation identity. In \cite{KimKwon2020blowup,KimKwonOh2020blowup}
the authors use \emph{method of nonlinear conjugation} to obtain
a new factorization of $L_{Q}iL_{Q}^{*}$ that enables them to uncover
a repulsive structure in the conjugated equation. On the other hand,
it is known that \eqref{CMdnls} enjoys a Lax pair structure \eqref{eq:LaxPair}
on chiral solutions in $L_{+}^{2}$ and if one takes the gauge transform
$\mathcal{G}$ on it, we have a Lax pair for \eqref{CMdnls-gauged}.
The goal here is to present a relation between the Lax pair and the
nonlinear conjugation identity. From this we will verify that the
Lax pair is extended to nonchiral solutions in $H^{1}$. 

Recall the linearized equation at $Q$, of \eqref{CMdnls-gauged}.
\begin{align}
\partial_{t}\eps+i\mathcal{L}_{Q}\eps=\partial_{t}\eps+iL_{Q}^{*}L_{Q}\eps\approx0.\label{eq:LinearizedEquation}
\end{align}
When we do modulation analysis, a crucial part of analysis is to propagate
the smallness of the $\varepsilon$-part. On the one hand, we obtain coercivity
estimates by modulating the generalized kernel directions of the linearized
operator and imposing orthogonality conditions on $\varepsilon$.
On the other hand, we take the adapted derivatives of $\varepsilon$ and
find a repulsive dynamics of the higher derivative of $\varepsilon.$.
Following the idea of \cite{RodnianskiSterbenz2010} in the self-dual
case, we take an adapted derivative, $L_{Q}$, and we obtain a linear
equation for $\eps_{1}=L_{Q}\eps$ 
\begin{align*}
(\partial_{t}+L_{Q}iL_{Q}^{*})\eps_{1}\approx0,\quad\eps_{1}\coloneqq L_{Q}\eps. & {\color{black}}
\end{align*}
At first glance, this equation does not look like a Hamiltonian form
due to the non-commutativity $[L_{Q},i]\ne0$. However, a direct computation
shows that $L_{Q}iL_{Q}^{\ast}=iH_{Q}$ for some \emph{$\bbC$-linear}
operator $H_{Q}$ and hence it \emph{is} of Hamiltonian form. A similar
algebra was first observed in the context of the Chern--Simons--Schrödinger
equation \cite{KimKwon2020blowup,KimKwonOh2020blowup} and the authors
naturally derive this fact by the method of \emph{nonlinear conjugation}.
In an attempt to find a conjugation identity, we investigate
the nonlinear equation for $\bfD_{v}v$. Inspired by the Lax pair
for \eqref{CMdnls} in Hardy space $L_{+}^{2}$ by Gérard and
Lenzmann \cite{GerardLenzmann2022}, we define two operators $\widetilde{\mathbf{D}}_{v}$
and $H_{v}$ by 
\begin{align}
\widetilde{\mathbf{D}}_{v}f & \coloneqq\partial_{x}f+\tfrac{1}{2}v\mathcal{H}(\overline{v}f),\label{eq:DefinitionDvtilde}\\
H_{v}f & \coloneqq -\partial_{xx}f+\tfrac{1}{4}|v|^{4}f-v|D|\overline{v}f.\label{eq:DefinitionHv}
\end{align}
Note that when applying $v$ in place of $f$, $v\mapsto\td{\bfD}_{v}v=\bfD_{v}v$
is also the Bogomol'nyi operator. One can check through the gauge
transform that $\mathcal{L}_{\textnormal{Lax}}$ and $\mathcal{P}_{\textnormal{Lax}}$
in \eqref{eq:LaxPair} correspond to $-i\widetilde{\mathbf{D}}_{v}$
and $-iH_{v}$ respectively. See Proposition~\ref{Appendix Proposition Laxpair upto gauge}.
Using \eqref{eq:DefinitionHv}, \eqref{CMdnls-gauged} can be expressed
as 
\begin{align}
\partial_{t}v+iH_{v}v=0.\label{eq:v-equ using Hv}
\end{align}
These motivate us to utilize the Lax pair structure to derive the
equation for $\bfD_{v}v(=\td{\bfD}_{v}v)$. We commute $\widetilde{\mathbf{D}}_{v}$
and \eqref{eq:CMDNLS gauged Selfdualform} to obtain a commutator
formula. From this we can obtain a Lax pair for \eqref{CMdnls-gauged}.
\begin{prop}[Unconditional Lax pair for \eqref{CMdnls-gauged}]
\label{PropositionUnconditionalLax} For $v\in C([0,T];H^{s})$ which
solves \eqref{CMdnls-gauged} with $s\geq0$ sufficiently large, we
have 
\begin{align}
\partial_{t}(-i\widetilde{\mathbf{D}}_{v})=[-iH_{v},-i\widetilde{\mathbf{D}}_{v}].\label{eq:LaxEqu Unconditional}
\end{align}
\end{prop}

An important observation is that \eqref{eq:LaxEqu Unconditional}
will verify a Lax pair structure for \eqref{CMdnls} \emph{without
the chiral condition}. The proof of Proposition~\ref{PropositionUnconditionalLax}
is a direct calculation, so we pass it to Appendix~\ref{AppendixUnconditionalLaxProof}.
As a standard consequence of Lax pairs, we have a hierarchy of conservation laws.  
\begin{cor}[Hierarchy of conservation laws] \label{cor:hierarchy}
	For $j\in \bbN$ with $0\leq j\leq 2s$, the quantities
	\begin{align}
		I_j(v)\coloneqq(\td \bfD_v^j v,v)_r \label{eq:hierarchy}
	\end{align}
	are conserved.
\end{cor}
\begin{proof}
	Applying \eqref{eq:LaxEqu Unconditional} to \eqref{eq:v-equ using Hv}, for any $j\in \bbN$, we derive
    \begin{align}
        \partial_{t}(\td \bfD_v^jv)+iH_{v}\td \bfD_v^jv=0. \label{eq:v1equationTilde}
    \end{align}
    We also have
    \begin{align*}
       \partial_t I_j(v)=&-(iH_{v}\td \bfD_v^jv,v)_r-(\td \bfD_v^jv,iH_{v}v)_r
       \\
       =&-(iH_{v}\td \bfD_v^jv,v)_r+(iH_{v}\td \bfD_v^jv,v)_r=0.
    \end{align*}
    That is, we have $\partial_t I_j(v)=0$, and this completes the proof.
\end{proof}

As in \cite{KimKwon2020blowup,KimKwonOh2020blowup}, we take a conjugation
by $\tilde{\D}_{v}$ on \eqref{eq:v-equ using Hv} to obtain a nonlinear
conjugation identity, as well as a linear conjugation identity by
linearization.
\begin{prop}[Conjugation identity]
\label{PropositionConjugationIndentity} Let $v$ solve \eqref{CMdnls-gauged}.
Then, we have 
\begin{equation}
    \begin{aligned}
        0 & =\partial_{t}(\widetilde{\mathbf{D}}_{v}v)+iH_{v}\widetilde{\mathbf{D}}_{v}v \\
 & =\partial_{t}(\widetilde{\mathbf{D}}_{v}v)+i\widetilde{\mathbf{D}}_{v}^{*}\widetilde{\mathbf{D}}_{v}\widetilde{\mathbf{D}}_{v}v
 +\tfrac{i}{2}\widetilde{\mathbf{D}}_{v}v\mathcal{H}(\overline{v}\widetilde{\mathbf{D}}_{v}v)
 -\tfrac{i}{2}v\mathcal{H}(|\widetilde{\mathbf{D}}_{v}v|^{2}).
    \end{aligned} \label{eq:NonlinearConjugationIdentity}
\end{equation}
At the linearized level, we have a linear conjugation identity 
\begin{align}
L_{Q}iL_{Q}^{*}=iH_{Q}=i\widetilde{\mathbf{D}}_{Q}^{*}\widetilde{\mathbf{D}}_{Q}=-i\widetilde{\mathbf{D}}_{Q}\widetilde{\mathbf{D}}_{Q}=i(-\partial_{xx}+\tfrac{1}{4}Q{}^{4}-Q|D|Q).\label{eq:ConjugationIdentityBad}
\end{align}
\end{prop}

An analogous argument as in \cite{KimKwon2020blowup,KimKwonOh2020blowup}
gives a natural linear conjugation identity $L_{Q}iL_{Q}^{\ast}=i\td{\bfD}_{Q}^{\ast}\td{\bfD}_{Q}$.
However, as opposed to the case of (CSS) this conjugation identity
does not meet our need. In (CSS), the linear conjugation identity
naturally derived from the method of nonlinear conjugation gave a
factorization $L_{Q}iL_{Q}^{\ast}=iA_{Q}^{\ast}A_{Q}$ \emph{with
an additional property}: the supersymmetric conjugate $\td H_{Q}=A_{Q}A_{Q}^{\ast}$
in (CSS) becomes \emph{repulsive} Schrödinger operator. In contrast,
we do \emph{not} have this additional property for \eqref{CMdnls-gauged}.
We have $\td{\bfD}_{Q}^{\ast}=-\td{\bfD}_{Q}$ so the supersymmetric
conjugate still gives rise to the same operator, which has (formal)
kernel elements $Q$ and $xQ$. These are obvious obstructions to monotonicity. The linearized identity \eqref{eq:ConjugationIdentityBad} is motivated from 
\begin{align*}
L_{Q}iL_{Q}^{*}L_{Q}=iH_{Q}L_{Q}=i\widetilde{\mathbf{D}}_{Q}^{*}\widetilde{\mathbf{D}}_{Q}L_{Q}.
\end{align*}
The second identity follows directly from the linearization of \eqref{eq:NonlinearConjugationIdentity}, whereas the first identity is obtained by comparing the linearization of the first line of \eqref{eq:NonlinearConjugationIdentity} with the result of applying $L_Q$ to \eqref{eq:LinearizedEquation}.

We close this subsection with the proof of Lax structure \eqref{eq:NonlinearConjugationIdentity}
and \eqref{eq:ConjugationIdentityBad}.
\begin{proof}[Proof of Proposition~\ref{PropositionConjugationIndentity}]
We have 
\begin{align}
-\widetilde{\mathbf{D}}_{v}^{*}\widetilde{\mathbf{D}}_{v}h= & \partial_{xx}h+\tfrac{1}{2}v_{x}\mathcal{H}(\overline{v}h)+\tfrac{1}{2}v|D|(\overline{v}h)+\tfrac{1}{2}v\mathcal{H}(\overline{v}\partial_{x}h)+\tfrac{1}{4}v\mathcal{H}[|v|^{2}\mathcal{H}(\overline{v}h)]\nonumber \\
= & \partial_{xx}h+\tfrac{1}{2}v_{x}\mathcal{H}(\overline{v}h)+v|D|(\overline{v}h)-\tfrac{1}{2}v\mathcal{H}(\overline{\partial_{x}v}h)+\tfrac{1}{4}v\mathcal{H}[|v|^{2}\mathcal{H}(\overline{v}h)].\label{eq:NonlinearConjugationIdentity1}
\end{align}
From \eqref{eq:HilbertProductRule} with $f=|v|^{2}$ and $g=\overline{v}h$,
we have 
\begin{align}
v\mathcal{H}[|v|^{2}\mathcal{H}(\overline{v}h)]=-|v|^{4}h+v\mathcal{H}(|v|^{2})\mathcal{H}(\overline{v}h)-\mathcal{H}(\mathcal{H}(|v|^{2})\overline{v}h).\label{eq:NonlinearConjugationIdentity2}
\end{align}
Inserting \eqref{eq:NonlinearConjugationIdentity2} into \eqref{eq:NonlinearConjugationIdentity1},
we have 
\begin{align}
\eqref{eq:NonlinearConjugationIdentity1}=-H_{v}h+\tfrac{1}{2}\widetilde{\mathbf{D}}_{v}v\mathcal{H}(\overline{v}h)-\tfrac{1}{2}v\mathcal{H}(\overline{\widetilde{\mathbf{D}}_{v}v}h).\label{eq:NonlinearConjugationIdentityGoal}
\end{align}
We now prove \eqref{eq:ConjugationIdentityBad} by a direct calculation. We note that the
linearization of the map $v\mapsto\widetilde{\mathbf{D}}_{v}v$ is
$L_{v}$ since $\widetilde{\mathbf{D}}_{v}v=\mathbf{D}_{v}v$. Inserting $Q$ into \eqref{eq:NonlinearConjugationIdentityGoal}
directly, we have 
\begin{align}
iH_{Q}=i\widetilde{\mathbf{D}}_{Q}^{*}\widetilde{\mathbf{D}}_{Q}.\label{eq:LinearConjugationIdentityBadPart1}
\end{align}
By the definition of $\widetilde{\mathbf{D}}_{Q}$, we observe that $\widetilde{\mathbf{D}}_{Q}^*=-\widetilde{\mathbf{D}}_{Q}$, which yields 
\begin{align*}
    i\widetilde{\mathbf{D}}_{Q}^{*}\widetilde{\mathbf{D}}_{Q}=-i\widetilde{\mathbf{D}}_{Q}\widetilde{\mathbf{D}}_{Q}.
\end{align*} 
For $L_{Q}iL_{Q}^{*}=iH_{Q}$, we first rewrite $L_{Q}$ and $L_{Q}^{*}$
as the matrix form 
\begin{align}
L_{Q}f=\begin{pmatrix}L_{11} & 0\\
0 & L_{22}
\end{pmatrix}\begin{pmatrix}\Re f\\
\Im f
\end{pmatrix},\quad L_{Q}^{*}f=\begin{pmatrix}L_{11}^{*} & 0\\
0 & L_{22}^{*}
\end{pmatrix}\begin{pmatrix}\Re f\\
\Im f
\end{pmatrix},\label{eq:LQ MatrixForm}
\end{align}
where $L_{11}=\mathbf{D}_{Q}+Q\mathcal{H}(Q\ \cdot\ )$ and $L_{22}=\mathbf{D}_{Q}$.
We also note that 
\begin{align}
i=\begin{pmatrix}0 & -1\\
1 & 0
\end{pmatrix}.\label{eq:i MatrixForm}
\end{align}
Thus, by \eqref{eq:LQ MatrixForm} and \eqref{eq:i MatrixForm}, we
have 
\begin{align*}
L_{\mathcal{Q}}iL_{\mathcal{Q}}^{*}=i\begin{pmatrix}0 & 1\\
-1 & 0
\end{pmatrix}\begin{pmatrix}L_{11} & 0\\
0 & L_{22}
\end{pmatrix}\begin{pmatrix}0 & -1\\
1 & 0
\end{pmatrix}\begin{pmatrix}L_{11}^{*} & 0\\
0 & L_{22}^{*}
\end{pmatrix}=i\begin{pmatrix}L_{22}L_{11}^{*} & 0\\
0 & L_{11}L_{22}^{*}
\end{pmatrix}.
\end{align*}
We have 
\begin{align*}
L_{22}L_{11}^{*}f=\mathbf{D}_{Q}(\mathbf{D}_{Q}^{*}f-Q\mathcal{H}(Qf)) & =\mathbf{D}_{Q}\mathbf{D}_{Q}^{*}f-(\mathbf{D}_{Q}Q)\mathcal{H}(Qf)-Q|D|Qf\\
 & =\mathbf{D}_{Q}\mathbf{D}_{Q}^{*}f-Q|D|Qf
\end{align*}
since $\mathbf{D}_{Q}Q=0$. Similarly, we also have 
\begin{align*}
L_{11}L_{22}^{*}f=\mathbf{D}_{Q}\mathbf{D}_{Q}^{*}f+Q\mathcal{H}(Q\mathbf{D}_{Q}^{*}f) & =\mathbf{D}_{Q}\mathbf{D}_{Q}^{*}f+Q\mathcal{H}((\mathbf{D}_{Q}Q)f)-Q|D|Qf\\
 & =\mathbf{D}_{Q}\mathbf{D}_{Q}^{*}f-Q|D|Qf.
\end{align*}
We note that we used $Q\partial_{x}f=\partial_{x}(Qf)-(\partial_{x}Q)f$.
Therefore, we deduce that 
\begin{align}
L_{Q}iL_{Q}^{*}f=i\mathbf{D}_{Q}\mathbf{D}_{Q}^{*}f-Q|D|Qf.\label{eq:LinearConjugationIdentityBadPart2-1}
\end{align}
We have 
\begin{align}
\mathbf{D}_{Q}\mathbf{D}_{Q}^{*}=-\partial_{xx}+\tfrac{1}{4}\left(\mathcal{H}(Q^{2})\right)^{2}+\tfrac{1}{2}\mathcal{H}(\partial_{x}(Q^{2})).\label{eq:LinearConjugationIdentityBadPart2-2}
\end{align}
From $\mathcal{H}(Q^{2})=xQ^{2}=\frac{2x}{1+x^{2}}$, we have $\partial_{x}(Q^{2})=-Q^{2}\mathcal{H}(Q^{2})$.
By this and \eqref{eq:HilbertProductRule} with $f=g=Q^{2}$, we have
\begin{align*}
\eqref{eq:LinearConjugationIdentityBadPart2-2}=-\partial_{xx}+\tfrac{1}{4}Q^{4},
\end{align*}
and this implies that 
\begin{align}
\eqref{eq:LinearConjugationIdentityBadPart2-1}=i(-\partial_{xx}f+\tfrac{1}{4}Q^{4}f-Q|D|Qf)=iH_{Q}f.\label{eq:LinearConjugationIdentityBadPart2}
\end{align}
By \eqref{eq:LinearConjugationIdentityBadPart1} and \eqref{eq:LinearConjugationIdentityBadPart2},
we conclude \eqref{eq:ConjugationIdentityBad}. 
\end{proof}

\subsection{Conjugation identity with repulsivity and higher-order conservation laws\label{SubsectionConjugationAQ}}
Following the previous discussion, we aim to search for a conjugation
identity that introduces repulsivity. Afterward, we will explain how the higher-order conservation laws are used to observe that the higher-order nonlinear variables enjoy the repulsivity. To this end, we propose an operator
$A_{Q}$, distinct from $\widetilde{\mathbf{D}}_{v}$, which adheres
to the conjugation identity through equation $L_{Q}iL_{Q}^{*}=iA_{Q}^{*}A_{Q}$.
In addition, we will discuss another operator, $B_{Q}$, linked to
$A_{Q}$. $A_{Q}$ is a (nonlocal) first order differential operator,
but $B_{Q}$ is a zeroth-order operator that will resolve coercivity
issues that arise when using $A_{Q}$. A surprising property we find
here is that 
\begin{align*} 
A_{Q}A_{Q}^{\ast}=-i\rd_{xx}.
\end{align*}
We will explain in more
detail the conjugation identity after introducing it. 
\begin{prop}[Conjugation identity with repulsivity]
\label{prop:AQ-BQ definition}We define $A_{Q}$ by 
\begin{align}
A_{Q}f=\partial_{x}(x-\mathcal{H})(\langle x\rangle^{-1}f).\label{eq:DefinitionAQ}
\end{align}
We also denote the adjoint operator of $A_{Q}$ as $A_{Q}^{*}$, and
$A_{Q}^{*}$ is given by 
\begin{align*}
A_{Q}^{*}f=-\langle x\rangle^{-1}(x+\mathcal{H})\partial_{x}f.
\end{align*}
Then, the following hold true: 
\begin{enumerate}
\item (Conjugation identity) We have 
\begin{align}
L_{Q}iL_{Q}^{*}=iH_{Q}=iA_{Q}^{*}A_{Q}\label{eq:Conjugation Identity}
\end{align}
for sufficiently good $f$, i.e., $f\in H^{2}$. 
\item (Repulsivity) For $\partial_{x}f\in H^{1}$, 
\begin{align}
A_{Q}A_{Q}^{*}f=-\partial_{xx}f,\label{eq:RepulsivityAQAQstarequality}
\end{align}
\item (Kernel for $A_{Q}$ on $\dot{\mathcal{H}}^{1}$) Let $v\in\dot{\mathcal{H}}^{1}$
with $A_{Q}v=0$. Then, 
\begin{align}
v\in\textnormal{span}_{\mathbb{C}}\{Q,xQ\}.\label{eq:kernel AQ}
\end{align}
\item We have an algebraic equality, 
\begin{align}
A_{Q}[Q\mathcal{H}(f)]=A_{Q}[\langle x\rangle^{2}Q\mathcal{H}(\langle x\rangle^{-2}f)].\label{eq:LemmaHilbertDecayExchangeInAQ}
\end{align}
\end{enumerate}
We further define $B_{Q}$ and its adjoint $B_{Q}^{*}$ by 
\begin{align}
B_{Q}f=(x-\mathcal{H})(\langle x\rangle^{-1}f),\quad B_{Q}^{*}f=\langle x\rangle^{-1}(x+\mathcal{H})f.\label{eq:DefinitionBQ}
\end{align}
For $B_{Q}$, the following hold true: 
\begin{enumerate}
\item[(5)] \setcounter{enumi}{5} $A_{Q}=\partial_{x}B_{Q}$, $A_{Q}^{*}=-B_{Q}^{*}\partial_{x}$, and $A_{Q}^{*}A_{Q}=-B_{Q}^{*}\partial_{xx}B_{Q}$.
In addition, we have
\begin{align}
	\widetilde{\D}_Q=B_Q^*\partial_xB_Q.  \label{eq:DQ BQ decomp}
\end{align}
\item We have 
\begin{align}
B_{Q}^{*}B_{Q}f=f-\frac{1}{2\pi}Q\int_{\R}Qfdx,\quad B_{Q}B_{Q}^{*}=I.\label{eq:BQBQstar equal I}
\end{align}
\item (Kernel for $B_{Q}$ on $\dot{\mathcal{H}}^{1}$) Let $v\in\dot{\mathcal{H}}^{1}$
with $B_{Q}v=0$. Then, 
\begin{align*}
v\in\textnormal{span}_{\mathbb{C}}\{Q\}.
\end{align*}
\end{enumerate}
\end{prop}

The proof of this proposition is rather straightforward and will
be given later. The key point for introducing $A_{Q}$ is that it
gives a new factorization of $L_{Q}iL_{Q}^{*}$, \eqref{eq:Conjugation Identity}
such that its supersymmetric conjugation becomes repulsive as in \eqref{eq:RepulsivityAQAQstarequality}.
We also remark that $A_{Q}$ is $\C$-linear, that is, $[A_{Q},i]=0$. The derivation of $A_{Q}$ was
motivated by examining the decomposition on the Fourier side, with
an emphasis on the fact that differentiations, the Hilbert transform,
and multiplication by $x$ behave as local operators in the Fourier
side. Note that $A_{Q}$ is further decomposed into a usual derivative
$\partial_{x}$ and $B_{Q}$ and the repulsivity $A_{Q}A_{Q}^{\ast}=-\rd_{xx}$
is related to $B_{Q}B_{Q}^{*}=I$.
\begin{rem}
The search for a new factorization \eqref{eq:Conjugation Identity}
is motivated by a similar conjugation identity in the context of
Chern--Simons--Schrödinger equations (CSS) \cite{KimKwon2020blowup,KimKwonOh2020blowup}.
As explained in Section~\ref{SubsectionLaxStructure}, the favorable
$A_{Q}$ operator in (CSS) arose naturally via nonlinear conjugation,
but this is not the case for \eqref{CMdnls-gauged}. This is why we
believe that the discovery of $A_{Q}$ is highly nontrivial compared
to the case of (CSS).
\end{rem}

Applying the conjugation to $A_{Q}$ and using $[A_{Q},i]=0$ along
with the conjugation identity \eqref{eq:Conjugation Identity} and
\eqref{eq:RepulsivityAQAQstarequality}, we transform \eqref{eq:LinearizedEquation}
into 
\begin{align}
(\partial_{t}+iA_{Q}A_{Q}^{*})\eps_{2}=(\partial_{t}-i\partial_{xx})\eps_{2}\approx0,\quad\eps_{2}=A_{Q}\eps_{1}.\label{eq:eps2 equ}
\end{align}
However, we will face a coercivity issue in this setting. If we were
to perform the energy-Morawetz argument in the blow-up analysis in
this setting, we would need to control the energy of \eqref{eq:eps2 equ},
$\|\partial_{x}\eps_{2}\|_{L^{2}}^{2}=\|A_{Q}^{*}\eps_{2}\|_{L^{2}}^{2}$,
which requires a coercivity estimate for the operator $A_{Q}^{*}A_{Q}L_{Q}$.
However, there are ($\dot{H}^{3}$-)kernel elements of $A_{Q}^{*}A_{Q}L_{Q}$,
$x(1+x^{2})Q$ and $ix^{3}Q$, which do not belong to $N_{g}(i\calL_{Q})$.
To avoid this coercivity issue, we set aside $\partial_{x}$ from
$A_{Q}$. This is our main motivation for introducing $B_{Q}$, a
zeroth-order operator with $A_{Q}=\partial_{x}B_{Q}$.

Similarly to \eqref{eq:eps2 equ}, now we apply the conjugation to
$B_{Q}$ along with \eqref{eq:LinearizedEquation} to derive 
\begin{align}
(\partial_{t}-i\partial_{xx})B_{Q}\eps_{1}\approx0.\label{eq:BQ eps1 equ}
\end{align}
Then, the associated energy for \eqref{eq:BQ eps1 equ} is $\|\partial_{x}B_{Q}\eps_{1}\|_{L^{2}}^{2}=\|A_{Q}\eps_{1}\|_{L^{2}}^{2}$,
and it suffices to establish the coercivity for $A_{Q}L_{Q}$. Moreover,
we will demonstrate in Proposition~\ref{PropKernelAQLQ} that the
kernel of $A_{Q}L_{Q}$ aligns with $N_{g}(i\calL_{Q})$, as indicated
by \eqref{eq:GeneralizedKernelRelation}.

It is widely acknowledged that this one-dimensional free Schr\"odinger flow lacks monotonicity because of the zero resonance. If we were to restrict \emph{odd} solutions, these free 1D Schr\"odinger waves with odd initial conditions can be interpreted as radially symmetric free 3D Schrödinger waves through the substitution $v \mapsto rv $. This substitution permits the use of the Morawetz estimate in three dimensions. In contrast, for general one-dimensional solutions, an alternative approach is employed. In fact, we apply the higher-order conservation laws \eqref{eq:hierarchy} for $j=2,4$. More specifically, we can restore the linearized equation \eqref{eq:BQ eps1 equ} from the $\bfD_vv$ equation. We note that, starting from \eqref{eq:v1equationTilde} with $j=1$, the $\bfD_vv$ equation is represented as
\begin{align*}
    (\partial_t +iH_v)\D_vv= 0.
\end{align*}
Since $\bfD_QQ=0$, linearizing $\bfD_vv$ around $Q$ is the same as considering $\bfD_vv$ itself, thus we linearize this equation around $Q$ as follows.
\begin{align*}
    (\partial_t +iH_Q)\D_vv\approx 0.
\end{align*}
Next, by taking $B_Q$ and applying \eqref{eq:Conjugation Identity}, \eqref{eq:RepulsivityAQAQstarequality}, and \eqref{eq:BQBQstar equal I}, it is quite notable that we deduce 
\begin{equation}
    (\partial_t -i\partial_{xx})B_Q\D_vv\approx 0
\end{equation}
this, which exactly matches \eqref{eq:BQ eps1 equ}. This implies that $B_Q\D_vv$ approximately solves the 1D free Schr\"odinger flow and motivates us to view the $B_Q\D_vv$-variable. On the other hand, the Hamiltonian of this flow is $\|\partial_xB_Q\D_vv\|_{L^2}$. Thanks to \eqref{eq:BQBQstar equal I}, the operator $B_Q^*$ is an isometry on $L^2$, so we expect that $\|B_Q^*\partial_xB_Q\D_vv\|_{L^2}$ is almost conserved. Then by \eqref{eq:DQ BQ decomp}, we connect to the higher-order conservation laws in \eqref{eq:hierarchy} as  
\begin{align*}
    \|B_Q^*\partial_xB_Q\D_vv\|_{L^2}^2
    =\|\td{\D}_Q\D_vv\|_{L^2}^2
    \approx
    \|\td{\D}_v\D_vv\|_{L^2}^2=I_4(v),
\end{align*}
which gives a energy control in $B_Q\D_vv$-variable. Similarly, we also observe that 
\begin{align*}
    \|B_Q\D_vv\|_{L^2}^2\lesssim \|\bfD_vv\|^2_{L^2}=-I_2(v).
\end{align*}
Indeed, one can systematically use the higher-order conservation laws $I_j(v)$ for higher-order energy control. In a forthcoming paper, this argument is used in \cite{JeongKim2024arXiv} to construct blow-up solutions with quantized blow-up rates.

\begin{proof}[Proof of Proposition~\ref{prop:AQ-BQ definition} ]
$(1)$: By \eqref{eq:ConjugationIdentityBad}, we have 
\begin{align*}
-iL_{Q}iL_{Q}^{*}f=-\partial_{xx}f+\tfrac{1}{4}Q^{4}f-Q|D|Qf=H_{Q}f.
\end{align*}
By Lemma~\ref{LemmaCommuteHilbert}, we have $[x,\mathcal{H}]\partial_{xx}\langle x\rangle^{-1}f\equiv0$
for $f\in H^{2}$. Using this, we conclude 
\begin{align*}
A_{Q}^{*}A_{Q}f & =-\langle x\rangle^{-1}(x\partial_{x}+|D|)(x\partial_{x}+1-|D|)(\langle x\rangle^{-1}f)\\
 & =-\langle x\rangle^{-1}((1+x^{2})\partial_{xx}+2x\partial_{x}+2|D|-[x,\mathcal{H}]\partial_{xx})(\langle x\rangle^{-1}f)\\
 & =-\langle x\rangle^{-1}(\langle x\rangle\partial_{xx}\langle x\rangle-\langle x\rangle^{-2}+2|D|)(\langle x\rangle^{-1}f)\\
 & =-\partial_{xx}f+\tfrac{1}{4}Q^{4}f-Q|D|Qf=H_{Q}f.
\end{align*}

$(2)$, $(6)$: By Lemma~\ref{LemmaHilbertUsefulEquation}, we can
obtain \eqref{eq:RepulsivityAQAQstarequality}. In fact, we have 
\begin{align}
B_{Q}B_{Q}^{*}= & (x-\mathcal{H})\langle x\rangle^{-2}(x+\mathcal{H})h\nonumber \\
 & =(x^{2}\langle x\rangle^{-2}h+x\langle x\rangle^{-2}\mathcal{H}h)-[\mathcal{H}\left(x\langle x\rangle^{-2}h\right)+\mathcal{H}\left(\langle x\rangle^{-2}\mathcal{H}h\right)].\label{eq:BQBQstarIproof1}
\end{align}
Using $\mathcal{H}(\langle x\rangle^{-2})=x\langle x\rangle^{-2}$
and \eqref{eq:HilbertProductRule} with $f=x\langle x\rangle^{-2}$
and $g=h$ to latter two terms of \eqref{eq:BQBQstarIproof1}, we
deduce 
\begin{align*}
\eqref{eq:BQBQstarIproof1}=(x^{2}\langle x\rangle^{-2}h+x\langle x\rangle^{-2}\mathcal{H}h)+[\langle x\rangle^{-2}h-\mathcal{H}\left(\langle x\rangle^{-2}\right)\mathcal{H}h]=h.
\end{align*}
This leads to $B_{Q}B_{Q}^{*}=I$ and \eqref{eq:RepulsivityAQAQstarequality}.

$(3)$: Due to $[A_{Q},i]=0$, the kernel of $A_{Q}$ is a span over
$\mathbb{C}$. The proof is similar to that of Proposition~\ref{PropKernel mathcalLQ}.
Let $\frac{v}{\sqrt{1+x^{2}}}=(1+x^{2})h$, then we have 
\begin{align*}
(\xi\partial_{\xi}+|\xi|)(1-\partial_{\xi\xi})\widehat{h}=0
\end{align*}
with $\widehat{h}\in H^{2}$. That is, $\widehat{h}\in C^{1}$. Therefore,
we conclude 
\begin{align*}
\widehat{h}(\xi)=\begin{cases}
C_{1}e^{-\xi}+(C_{2}+C_{1})\xi e^{-\xi} & \xi\geq0\\
C_{1}e^{\xi}+(C_{2}-C_{1})\xi e^{\xi} & \xi<0.
\end{cases}
\end{align*}
Thus, the null space of $A_{Q}$ on $\dot{\mathcal{H}}^{1}$ is two-dimensional. Since we know $Q$ and $xQ$ are kernel elements of $A_{Q}$
in $\dot{\mathcal{H}}^{1}$, we finish the proof. We also note that
$\mathcal{F}^{-1}[(1+|\xi|)e^{-|\xi|}]=cQ^{4}$ and $\mathcal{F}^{-1}[\xi e^{-|\xi|}]=c^{\prime}xQ^{4}$
for some constants $c$ and $c^{\prime}$.

$(4)$: By \eqref{eq:HilbertDecayExchange2} and $[x,\mathcal{H}]=c\int_{\R}$,
we have 
\begin{align*}
A_{Q}[Q\mathcal{H}(f)]= & A_{Q}[\langle x\rangle^{2}Q\mathcal{H}(\langle x\rangle^{-2}f)]\\
 & -A_{Q}[Q\cdot c{\textstyle \int_{\R}\left(x\langle x\rangle^{-2}f\right)dx]-A_{Q}[Q\cdot xc{\textstyle \int_{\R}\left(\langle x\rangle^{-2}f\right)dx].}}
\end{align*}
We remark that $\int_{\bbR}\left(x\langle x\rangle^{-2}f\right)dx$
and $\int_{\bbR}\left(\langle x\rangle^{-2}f\right)dx$ are constants
in $x$. Since $A_{Q}Q=A_{Q}xQ=0$, we deduce \eqref{eq:LemmaHilbertDecayExchangeInAQ}.

$(5)$: Except for \eqref{eq:DQ BQ decomp}, it comes directly from the definition of $A_{Q}$ and $B_{Q}$. For \eqref{eq:DQ BQ decomp}, a direct computation yields
\begin{align}
	B_Q^*\partial_xB_Qf=&
	\tfrac{x}{\langle x \rangle}\partial_x(\tfrac{x}{\langle x \rangle}f)
	+\tfrac{1}{\langle x \rangle}\calH \partial_x(-\calH)(\tfrac{1}{\langle x \rangle}f) \label{eq:DQ BQ decomp pf 1}
	\\
	&+\tfrac{1}{\langle x \rangle}\calH\partial_x(\tfrac{x}{\langle x \rangle}f)
	-\tfrac{x}{\langle x \rangle}\partial_x\calH(\tfrac{1}{\langle x \rangle}f). \label{eq:DQ BQ decomp pf 2}
\end{align}
Using $\calH \partial_x(-\calH)=-\calH^2 \partial_x=\partial_x$, we have
\begin{align*}
	\eqref{eq:DQ BQ decomp pf 1}=\partial_x f.
\end{align*}
To handle \eqref{eq:DQ BQ decomp pf 2}, we use \eqref{eq:CommuteHilbertDerivative} to derive
\begin{align*}
	\eqref{eq:DQ BQ decomp pf 2}=\tfrac{1}{\langle x \rangle}\partial_x[\calH,x](\tfrac{1}{\langle x \rangle}f)
	+
	\tfrac{1}{\langle x \rangle}\calH(\tfrac{1}{\langle x \rangle}f)
	=\tfrac{1}{\langle x \rangle}\calH(\tfrac{1}{\langle x \rangle}f)=\tfrac{1}{2}Q\calH(Qf). 
\end{align*}
Therefore, we conclude \eqref{eq:DQ BQ decomp}.

$(7)$: We omit the proof as it is similar to the proof of $(3)$. 
\end{proof}
\begin{rem}
As stated previously, we can find \eqref{eq:Conjugation Identity}
by observing the equation on the Fourier side. Set $h=\frac{f}{\sqrt{1+x^{2}}}$.
Then, we have 
\begin{align*}
H_{Q}f & =-\partial_{xx}f+\tfrac{1}{4}Q^{4}f-Q|D|Qf\\
 & =-h_{xx}\sqrt{1+x^{2}}-\tfrac{2x}{\sqrt{1+x^{2}}}h_{x}-\tfrac{2}{\sqrt{1+x^{2}}}|D|h.
\end{align*}
By multiplying $\sqrt{1+x^{2}}$, we obtain 
\begin{align*}
-h_{xx}(1+x^{2})-2xh_{x}-2|D|h=-\partial_{x}((1+x^{2})h_{x})-2|D|h.
\end{align*}
Taking Fourier transform, we have 
\begin{align*}
-i\xi(1-\partial_{\xi}^{2})(i\xi\widehat{h})-2|\xi|\widehat{h} & =-\xi^{2}\widehat{h}_{\xi\xi}-2\xi\widehat{h}_{\xi}+(\xi^{2}-2|\xi|)\widehat{h}\\
 & =-(-\partial_{\xi}\xi+|\xi|)(-\xi\partial_{\xi}-|\xi|)\widehat{h}.
\end{align*}
Thus, we formally deduce that 
\begin{align*}
-\partial_{x}((1+x^{2})h_{x})-2|D|h=-(x\partial_{x}+|D|)(\partial_{x}x-|D|)h=\langle x\rangle A_{Q}^{*}A_{Q}f.
\end{align*}
\end{rem}

It remains to show the coercivity of $A_{Q}L_{Q}$. For technical
reasons, we will also use a modification of $L_{Q}$, denoted by $\widetilde{L}_{Q}$,
to rigorously describe the coercivity estimate at the $H^{2}$ level.
Define $\widetilde{L}_{Q}$ by 
\begin{align}
\widetilde{L}_{Q}f\coloneqq\mathbf{D}_{Q}f+Q^{-1}\mathcal{H}\Re(Q^{3}f).\label{eq:tildeLQ definition}
\end{align}

\begin{rem}[Introduction of $\widetilde{L}_{Q}$]
\label{RemarkLQ Well Definedness} As mentioned above, we will perform
an energy estimate at $\dot{H}^{2}$-level and hence we expect $L_{Q}\eps$
to belong to $\dot{\mathcal{H}}^{1}$. However, $L_{Q}f$ may not
be defined for general $\dot{\calH}^{2}$ functions due to the presence
of $\mathcal{H}(L^{\infty})$, as mentioned before in Remark~\ref{Remark H1 Warning remark}.
For example, a part of $P$, $(1+x^{2})Q$ belongs to $\dot{\mathcal{H}}^{2}$,
but $L_{Q}(1+x^{2})Q=2xQ+2Q\mathcal{H}(1)$. We modify $L_{Q}$
to $\widetilde{L}_{Q}$ to avoid this issue. Now $\widetilde{L}_{Q}$
is well defined on $\dot{\mathcal{H}}^{2}$ and $A_{Q}L_{Q}=A_{Q}\widetilde{L}_{Q}$. 
\end{rem}

By \eqref{eq:HilbertDecayExchange2}, we compute 
\begin{align}
(\widetilde{L}_{Q}-L_{Q})f=\tfrac{1}{\pi}xQ{\textstyle \int_{\R}\langle x\rangle^{-2}\Re(Qf)dx}+\tfrac{1}{\pi}Q{\textstyle \int_{\R}x\langle x\rangle^{-2}\Re(Qf)dx}.\label{eq:LQ LQtilde Difference}
\end{align}
The expression in the integral of \eqref{eq:LQ LQtilde Difference}
is valid for $\dot{\mathcal{H}}^{1}$-functions and using \eqref{eq:kernel AQ}
we obtain the following lemma. 
\begin{lem}
\label{LemmaAQLQ equal AQtildeLQ} For any $f\in\dot{\mathcal{H}^{1}}\cap\dot{H}^{2}$,
We have 
\begin{align*}
A_{Q}L_{Q}f=A_{Q}\widetilde{L}_{Q}f & .
\end{align*}
\end{lem}

We will substitute the second-order adapted derivative $A_{Q}L_{Q}$
with $A_{Q}\widetilde{L}_{Q}$ when the acting function $\eps$ belongs
to $\dot{\mathcal{H}}^{2}$ . To elaborate this process, it is important
to understand the strategy to which we use various decompositions of nonlinear
solutions, depending on the spaces where the radiation parts $\eps$
or $\widehat{\eps}$ belong. This strategy is one of our novelty
in the analysis of controlling the radiation part. As mentioned in
Remark~\ref{Remark Invariant genKernel and Transversality}, the
nonlinear solution is decomposed as $\lambda^{\frac{1}{2}}e^{-i\gamma}v(\lambda y+x)=(Q+P+\eps)(y)$.
However, this decomposition is not applicable when the solution is
measured in $H^{1}$ since $P$ does not belong to $H^{1}$. When
we control the solution in $H^{1}$, for instance in modulation estimates
of $\lambda(t),\gamma(t),x(t)$ or in the energy estimates, we use
a rougher decomposition $\lambda^{\frac{1}{2}}e^{-i\gamma}v(\lambda y+x)=Q+\widehat{\eps}$
where $\widehat{\eps}=P+\eps\in H^{2}$. On the other hand, in the
$\dot{\mathcal{H}^{2}}$ topology, we can use the finer decomposition
$Q+P+\eps$ as $P,\eps\in\dot{\mathcal{H}}^{2}$. For further details, refer to Section~\ref{Subsection Decomposition}.

From the decomposition in the form of $Q+\widehat{\eps}$, we are
allowed to apply $A_{Q}L_{Q}$ to $\widehat{\eps}$, whereas $A_{Q}L_{Q}$
is not applicable to $\eps$ since $P,\eps\notin H^{1}$. From Lemma~\ref{LemmaAQLQ equal AQtildeLQ},
we have $A_{Q}L_{Q}\widehat{\eps}=A_{Q}\widetilde{L}_{Q}\widehat{\eps}$.
Applying the decomposition for the $\dot{\mathcal{H}^{2}}$-topology,
we can write $A_{Q}\widetilde{L}_{Q}\widehat{\eps}=A_{Q}\widetilde{L}_{Q}P+A_{Q}\widetilde{L}_{Q}\eps$,
as $P,\eps\in\dot{\mathcal{H}}^{2}$. By the definition of $P$ and
\eqref{eq:GeneralizedKernelRelation}, we have $A_{Q}\widetilde{L}_{Q}P=0$
and therefore arrive at 
\begin{equation}
A_{Q}L_{Q}\widehat{\eps}=A_{Q}\widetilde{L}_{Q}\eps.\label{eq:AQLQ epsilon hat AQtildeLQ epsilon}
\end{equation}
From this we take the seamless transition from $A_{Q}L_{Q}\widehat{\eps}$
to $A_{Q}\widetilde{L}_{Q}\eps$. This kind of approach will be frequently
used in Section~\ref{Section Proof of nonlinear estimates}. We will
often find some cases where $Ta$ and $Tb$ are not defined by
themselves, but $T(a+b)$ is well defined. For example, $T=A_{Q}L_{Q}$,
$T=[\partial_{x},\mathcal{H}]$, etc. See also Remark~\ref{Remark commuting Hilbert derivative}.

Now, we find the kernel of $A_{Q}\widetilde{L}_{Q}$ in $\dot{\calH}^{2}$.
In view of \eqref{eq:GeneralizedKernelRelation}, we have 
\begin{align*}
\text{span}_{\mathbb{R}}\{iQ,\Lambda Q,Q_{x},ixQ,ix^{2}Q,(1+x^{2})Q\}\subset\ker A_{Q}\widetilde{L}_{Q}.
\end{align*}
We show that the above indeed forms the kernel of $A_{Q}\widetilde{L}_{Q}$. 
\begin{prop}[Kernel of $A_{Q}\widetilde{L}_{Q}$ in $\dot{\mathcal{H}}^{2}$]
\label{PropKernelAQLQ} Let $v\in\dot{\mathcal{H}}^{2}$ be with $A_{Q}\widetilde{L}_{Q}v=0$.
Then, 
\begin{align*}
v\in\textnormal{span}_{\mathbb{R}}\{iQ,\Lambda Q,Q_{x},ixQ,ix^{2}Q,(1+x^{2})Q\}.
\end{align*}
\end{prop}

\begin{proof}
We check the real part and the imaginary part separately. Assume $v$ is real-valued.
Then, $v$ satisfies 
\begin{align*}
A_{Q}(\textbf{D}_{Q}v+Q^{-1}\mathcal{H}(Q^{3}v))=0.
\end{align*}
Set $h=\langle x\rangle^{-3}v$, then $h$ satisfies 
\begin{align*}
\partial_{x}(x-\mathcal{H})\langle x\rangle^{-1}\widetilde{L}_{Q}v=\partial_{x}(x-\mathcal{H})[\partial_{x}((1+x^{2})h)+2xh+2\mathcal{H}h]
\end{align*}
Since $v\in\dot{\mathcal{H}}^{2}$, $h\in H^{2}$ and $xh\in L^{2}$.
Taking the Fourier transform, we have 
\begin{align*}
(\xi\partial_{\xi}+|\xi|)(\xi(1-\partial_{\xi\xi})+2\partial_{\xi}-2\text{sgn}(\xi))\widehat{h}(\xi)=0
\end{align*}
weakly. For $\xi>0$ and $\xi<0$, respectively, classical solutions
are given by 
\begin{align*}
\widehat{h}(\xi)=\begin{cases}
C_{1}e^{-\xi}+C_{2}e^{\xi}+C_{3}e^{-\xi}(1+\xi)\xi, & \xi>0,\\
C_{4}e^{-\xi}+C_{5}e^{\xi}+C_{6}e^{\xi}(1-\xi)\xi, & \xi<0.
\end{cases}
\end{align*}
Since $\widehat{h}\in H^{1}\subset C_{0}$, we have $C_{2}=C_{4}=0$
and $C_{1}=C_{5}$. Moreover, since $h$ is real-valued, we have $\widehat{h}(-\xi)=\overline{\widehat{h}(\xi)}$.
From this, we deduce that $C_{1}$ is real-valued, $-\overline{C_{3}}=C_{6}$.
Thus, 
\begin{align}
\widehat{h}(\xi)=\begin{cases}
c_{1}e^{-\xi}+(c_{2}+ic_{3})e^{-\xi}(1+\xi)\xi, & \xi>0,\\
c_{1}e^{\xi}+(-c_{2}+ic_{3})e^{\xi}(1-\xi)\xi, & \xi<0,
\end{cases}\label{eq:AQLQKernelProof h hat}
\end{align}
for real $c_{1},c_{2},c_{3}$. This implies that the real-valued kernel
for $A_{Q}\widetilde{L}_{Q}$ is at most $3$-dimensional. Since $\Lambda Q,Q_{x},(1+x^{2})Q$
are the kernel elements, we obtain the kernel for real-valued functions.
We note that the Fourier transforms of $\langle x\rangle^{-3}\Lambda Q$
and $\langle x\rangle^{-3}Q_{x}$ are given by \eqref{eq:Lambda Q Qsubx Fourier explicit formula},
and the Fourier transform of $\langle x\rangle^{-3}(1+x^{2})Q=\sqrt{2}\langle x\rangle^{-2}$
is given by 
\begin{align*}
\mathcal{F}(\sqrt{2}\langle x\rangle^{-2})(\xi)=\sqrt{\pi}e^{-|\xi|}.
\end{align*}
These are identical to \eqref{eq:AQLQKernelProof h hat}.

Now, we assume that $v$ is pure imaginary. Since the kernel of $A_{Q}$
is $\text{span}_{\mathbb{C}}\{yQ,Q\}$ and $\textbf{D}_{Q}$ commutes
with $i$, it suffices to find $v$ such that $\textbf{D}_{Q}v=0$,
$\textbf{D}_{Q}v=iyQ$, and $\textbf{D}_{Q}v=iQ$. We know that $\ker\textbf{D}_{Q}=\text{span}_{\mathbb{C}}\{Q\}$
on $L^{2}$. By direct calculation, we have $\textbf{D}_{Q}(xQ)=Q$
and $\textbf{D}_{Q}(x^{2}Q)=2xQ$. This implies that the pure imaginary
part of the kernel of $A_{Q}$ is given by $\text{span}_{\mathbb{R}}\{iQ,ixQ,ix^{2}Q\}$.
This finishes the proof.
\end{proof}

Finally, we state the coercivity of $L_{Q}$ and $A_{Q}\widetilde{L}_{Q}$.
This step is the main motivation for the adapted function spaces $\dot{\mathcal{H}}^{1}$
and $\dot{\mathcal{H}}^{2}$. We want to find a norm comparable to
$\|L_{Q}f\|_{L^{2}}$ at $\dot{H}^{1}$-level. Recall $L_{Q}=\partial_{x}+\tfrac{x}{1+x^{2}}+Q\mathcal{H}(\Re(Q\cdot))$.
Note that $L_{Q}\approx\rd_{x}$ when $x$ is small and $L_{Q}\approx\rd_{x}+\frac{1}{x}$
when $x$ is large. In view of Hardy's inequality for the operators
$\rd_{x}+\frac{k}{x}$ ($k\in\bbR$), $\|L_{Q}f\|_{L^{2}}$ can control
$\|f\|_{\dot{\calH}^{1}}$ modulo a compact term; one has \emph{subcoercivity
estimate} $\|L_{Q}f\|_{L^{2}}+\|\chf_{|x|\aleq1}f\|\sim\|f\|_{\dot{\calH}^{1}}$
(Lemma~\ref{LemmaAppendix LQ subcoer}). The norm of $\dot{\calH}^{1}$
was built in this way. The compactly supported term can be deleted
(hence giving \emph{coercivity estimate}) by restricting $\dot{\calH}^{1}$
to a subspace that is transversal to the kernel of $L_{Q}$. A similar
argument works for $A_{Q}\tilde L_{Q}f$ at the $\dot{H}^{2}$-level
to build the $\dot{\calH}^{2}$ norm via the subcoercivity estimate
(Lemma~\ref{lem:subcoercivity AQLQ}). The associated coercivity
estimate is obtained in a similar way by restricting to a subspace.

For notational convenience, denote the elements of $N_{g}(i\mathcal{L}_{Q})$
(or the kernel of $A_{Q}\widetilde{L}_{Q}$) by 
\begin{align}
 & \begin{split}\mathcal{K}_{1} & =\Lmb Q,\\
\mathcal{K}_{4} & =ix^{2}Q,
\end{split}
 & \begin{split}\mathcal{K}_{2} & =iQ,\\
\mathcal{K}_{5} & =(1+x^{2})Q,
\end{split}
 & \begin{split}\mathcal{K}_{3} & =Q_{x},\\
\mathcal{K}_{6} & =ixQ.
\end{split}
\label{eq:3.3 def ker elem calK}
\end{align}
We state the coercivity estimates related to the linear operators
$L_{Q}$ and $A_{Q}\td L_{Q}$: 
\begin{prop}[Coercivity for $L_{Q}$ at $\dot{\mathcal{H}}^{1}$ level]
\label{PropCoercivityLQ} Let $\psi_{1},\psi_{2},\psi_{3}$ be elements
of the dual space $(\dot{\mathcal{H}}^{1})^{*}$. If the $3\times3$
matrix $(a_{ij})$ defined by $a_{ij}=(\psi_{i},\mathcal{K}_{j})_{r}$
has a nonzero determinant, then we have a coercivity estimate 
\begin{align*}
\|v\|_{\dot{\mathcal{H}}^{1}}\lesssim_{\psi_{1},\psi_{2},\psi_{3}}\|L_{Q}v\|_{L^{2}}\lesssim\|v\|_{\dot{\mathcal{H}}^{1}},\quad\forall v\in\dot{\mathcal{H}}^{1}\cap\{\psi_{1},\psi_{2},\psi_{3}\}^{\perp}.
\end{align*}
\end{prop}

\begin{prop}[Coercivity for $A_{Q}\widetilde{L}_{Q}$ at $\dot{\mathcal{H}}^{2}$
level]
\label{PropCoercivityAQLQ} Let $\psi_{1},\psi_{2},\psi_{3},\psi_{4},\psi_{5},\psi_{6}$
be elements of the dual space $(\dot{\mathcal{H}}^{2})^{*}$. If the
$6\times6$ matrix $(a_{ij})$ defined by $a_{ij}=(\psi_{i},\mathcal{K}_{j})_{r}$
has a nonzero determinant, then we have a coercivity estimate 
\begin{align*}
\|v\|_{\dot{\mathcal{H}}^{2}}\lesssim_{\psi_{1},\psi_{2},\psi_{3},\psi_{4},\psi_{5},\psi_{6}}\|A_{Q}\widetilde{L}_{Q}v\|_{L^{2}}\lesssim\|v\|_{\dot{\mathcal{H}}^{2}},\ \forall v\in\dot{\mathcal{H}}^{2}\cap\{\psi_{i}:1\leq i\leq6\}^{\perp}.
\end{align*}
\end{prop}

The above propositions are standard consequences of the associated
subcoercivity estimates (see Lemmas~\ref{LemmaAppendix LQ subcoer}
and \ref{lem:subcoercivity AQLQ} in Appendix~\ref{AppendixSubcoercivity})
and the characterization of the kernels (given in Propositions~\ref{PropKernel mathcalLQ}
and \ref{PropKernelAQLQ}), and hence are omitted. We may refer to
the proof in \cite[Lemma A.15]{KimKwon2020blowup}.

\section{\label{SectionFormalBlowupLaw}Formal modulation law and modified
profile}

In this section, we discuss the blow-up profile along with the formal
dynamics of the modulation parameters. As mentioned in Remark~\ref{Remark Invariant genKernel and Transversality},
we will decompose the nonlinear solution $v(t)$ into the blow-up
profile and the radiation part; 
\begin{align}
v(t,x)=\frac{e^{i\gamma(t)}}{\lambda(t)^{1/2}}[Q+P(b(t),\eta(t),\nu(t);\cdot)+\eps(t,\cdot)]\left(\frac{x-x(t)}{\lambda(t)}\right)\label{eq:v Decomposition form}
\end{align}
where $Q+P$ is the modified blow-up profile such that $P(b,\eta,\nu)\to0$
as $t\to T$. After all the analysis, we hope that the radiation part
(i.e., the last term of \eqref{eq:v Decomposition form}) converges
to some function $v^{\ast}$ called the \emph{asymptotic profile}.
In this section, we ignore $\eps$ (i.e., $\eps\equiv0$) and conduct
a formal computation for modulation laws of $\lambda(t)$, $\gmm(t)$,
$x(t)$, $b(t)$, $\eta(t)$, $\nu(t)$ together with the derivation
of the right modified profile. This procedure often called a \emph{tail
computation} \cite{RaphaelRodnianski2012}.

\subsection{\label{subsec:Formal-modulation-laws}Formal modulation laws}

We work on renormalized coordinates $(s,y)$ instead of the original
ones $(t,x)$, sending the blow-up time to $+\infty$ and zooming
in on the soliton scale $\la(t)$. Define $(s,y)$ and the renormalized
soliton $w$ by 
\begin{align*}
\frac{ds}{dt}=\frac{1}{\lambda^{2}},\quad y=\frac{x-x(t)}{\lambda},\quad w(s,y)=\lambda^{\frac{1}{2}}e^{-i\gamma}v(t,\lambda y+x(t))|_{t=t(s)}=[v]_{\la^{-1},-\ga,-x},
\end{align*}
where we denoted 
\begin{equation}
[f]_{\la,\ga,x}=\frac{e^{i\ga}}{\la^{1/2}}f\lr{\frac{\cdot-x}{\la}}.\label{eq:notation renornmalized}
\end{equation}
Modulation parameters $\lambda(t)$, $\gmm(t)$, $x(t)$, $b(t)$,
$\eta(t)$, $\nu(t)$ are also identified with functions of $s$;
we abuse notation and simply write $\la(s)=\la(t(s))$, $\ga(s)=\ga(t(s))$,
and so on.

The (renormalized) \emph{nonlinear adapted derivative $w_{1}$ }is
defined by 
\begin{align*}
w_{1}=\textbf{D}_{w}w=\lambda^{\frac32} e^{-i\gamma}(\textbf{D}_{v}v)(t,\lambda y+x(t))|_{t=t(s)}=\lmb[\mathbf{D}_{v}v]_{\la^{-1},-\ga,-x}.
\end{align*}
Renormalizing \eqref{eq:CMdnlsSelfdualform} and \eqref{eq:v1equationTilde} with $j=1$,
we derive the equations of $w$ and $w_{1}$ as 
\begin{align}
(\partial_{s}-\frac{\lambda_{s}}{\lambda}\Lambda+\gamma_{s}i-\frac{x_{s}}{\lambda}\partial_{y})w+iL_{w}^{*}w_{1} & =0,\label{eq:w-equ}\\
(\partial_{s}-\frac{\lambda_{s}}{\lambda}\Lambda_{-1}+\gamma_{s}i-\frac{x_{s}}{\lambda}\partial_{y})w_{1}+iH_{w}w_{1} & =0,\label{eq:w1-equ}
\end{align}
where $H_{w}=-\partial_{yy}+\frac{1}{4}|w|^{4}-w|D|\overline{w}$.
Since $H_{w}w=L_{w}^{*}\mathbf{D}_{w}w$, we can rewrite \eqref{eq:w-equ}
as 
\begin{align*}
(\partial_{s}-\frac{\lambda_{s}}{\lambda}\Lambda+\gamma_{s}i-\frac{x_{s}}{\lambda}\partial_{y})w+iH_{w}w=0.
\end{align*}
One caveat is that $H_{w}w_{1}\ne L_{w}^{*}\mathbf{D}_{w}w_{1}$.
The nonlinear adapted derivative $w_{1}$ is a nonlinear version of the
adapted derivative $L_{Q}w$. This was first introduced in the context
of (CSS) \cite{KimKwonOh2020blowup}. Since it is a more canonical
derivative of $w$, it simplifies many parts of analysis, such as
Lemma~\ref{LemmaNonlinearEstimate1}, etc.

Next, we rewrite \eqref{eq:w1-equ} by decomposing it into the linear
and the nonlinear parts; 
\begin{align}
(\partial_{s}-\frac{\lambda_{s}}{\lambda}\Lambda_{-1}+\gamma_{s}i-\frac{x_{s}}{\lambda}\partial_{y})w_{1}+iH_{Q}w_{1}=i\textnormal{NL}_{1},\label{eq:w1-equ-modify}
\end{align}
with 
\begin{align}
\textnormal{NL}_{1}=(H_{Q}-H_{w})w_{1}.\label{eq:DefinitionNL1}
\end{align}

Now, we construct the modified profiles. We will work for the solution
pair $(w,w_{1})$ and find the profiles for each of them. From the
linear analysis in Section~\ref{SubsectionLQ}, the following form
of the modified profile $Q+P$ is suggested: 
\[
Q+P=Q-ib\tfrac{y^{2}}{4}Q-\eta\tfrac{1+y^{2}}{4}Q+i\nu\tfrac{y}{2}Q+\text{h.o.t}.
\]
It will turn out that the linear expansion of the modified profile
$P(b,\eta,\nu)$ suffices for our analysis: 
\begin{equation}
P\coloneqq-ib\frac{y^{2}}{4}Q-\eta\frac{1+y^{2}}{4}Q+i\nu\frac{y}{2}Q\label{eq:definition P}
\end{equation}
Next, motivated by $w_{1}=\bfD_{w}w$, we hope to define $P_{1}$
such that $P_{1}\approx\bfD_{(Q+P)}(Q+P)$. Linearizing, we want to
define $\wt P_{1}=L_{Q}P$. However, since $L_{Q}(\tfrac{1+y^{2}}{4}Q)=\tfrac{1}{2}yQ+\tfrac{1}{2}Q\mathcal{H}(1)$
by computation, $L_{Q}P$ is not rigorously defined. Hence, we discard
the term $\mathcal{H}(1)$ (or use the formal convention $\calH(1)=0$)
and define the profile $\widehat{P}_{1}$ for $w_{1}$ as 
\begin{align}
\widehat{P}_{1} & \coloneqq-(ib+\eta)\frac{y}{2}Q+i\nu\frac{1}{2}Q.\label{eq:DefinitionProfileP1hat}
\end{align}
Note that $\wt P_{1}$ \emph{degenerates} in order 1 in $b,\eta,\nu$.

Now, we search for the formal dynamical laws of $\lambda,\gamma$,
and $x$. Plugging $w=Q+P$ in \eqref{eq:w-equ}, and collecting the
first order terms in $b,\eta,\nu$ (and $\la_{s},\ga_{s},x_{s}$),
we have 
\begin{align*}
-\frac{\lambda_{s}}{\lambda}\Lambda w & \to-\frac{\lambda_{s}}{\lambda}\Lambda Q,\quad\gamma_{s}iw\to\gamma_{s}iQ,\quad-\frac{x_{s}}{\lambda}\partial_{y}w\to-\frac{x_{s}}{\lambda}\partial_{y}Q,\\
L_{w}^{*}w_{1} & \to ib\Lambda Q-\frac{\eta}{2}Q-i\nu Q_{y}.
\end{align*}
Then, we obtain 
\begin{align*}
\eqref{eq:w-equ}\approx-(\frac{\lambda_{s}}{\lambda}+b)\Lambda Q+(\gamma_{s}-\frac{\eta}{2})iQ-(\frac{x_{s}}{\lambda}-\nu)Q_{y}+\text{h.o.t}.
\end{align*}
Thus, we arrive at the first modulation laws for $\la,\ga,x$; 
\begin{align}
\frac{\lambda_{s}}{\lambda}+b=0,\quad\gamma_{s}-\frac{\eta}{2}=0,\quad\frac{x_{s}}{\lambda}-\nu=0.\label{eq:FirstModulationEqu}
\end{align}

Now, we find the formal modulation laws for $b,\eta$ and $\nu$.
We use the $w_{1}$-equation \eqref{eq:w1-equ} and the ansatz \eqref{eq:FirstModulationEqu}.
Collecting the quadratic terms in $b,\eta,\nu$ and $b_{s},\eta_{s},\nu_{s}$,
we have 
\begin{align}
\begin{split}\partial_{s}w_{1} & \to-(ib_{s}+\eta_{s})\tfrac{y}{2}Q+i\tfrac{1}{2}\nu_{s}Q,\\
b\Lambda_{-1}w_{1} & \to-(ib^{2}+b\eta)\tfrac{1}{2}\Lambda_{-1}(yQ)+ib\nu\tfrac{1}{2}\Lambda_{-1}Q,\\
\gamma_{s}iw_{1} & \to(b\eta-i\eta^{2})\tfrac{1}{4}yQ-\eta\nu\tfrac{1}{4}Q,\\
-\nu\partial_{y}w_{1} & \to(ib\nu+\eta\nu)\tfrac{1}{2}\partial_{y}(yQ)-i\nu^{2}\tfrac{1}{2}Q_{y},\\
iH_{w}w_{1} & \to iH(Q,P)\widehat{P}_{1},
\end{split}
\label{eq:QuadraticModulation calculation}
\end{align}
where 
\begin{align*}
H(Q,P)f\coloneqq |Q|^{2}\Re(\overline{Q}P)f-Q|D|\overline{P}f-P|D|\overline{Q}f.
\end{align*}
Now we compute the last term of \eqref{eq:QuadraticModulation calculation}. 
\begin{lem}
\label{lem:H(P,Q)}Formally, we have 
\begin{equation}
\begin{aligned}iH(Q,P)\widehat{P}_{1}= & -(b^{2}+\eta^{2})\tfrac{i}{4}Q\mathcal{H}(1)+ib^{2}\tfrac{1}{2}\Lambda_{\frac{1}{2}}yQ+b\eta\tfrac{1}{2}\Lambda_{\frac{1}{2}}yQ\\
 & -ib\nu\left(\tfrac{1}{2}\Lambda Q+\tfrac{1}{2}\partial_{y}(yQ)\right)+\eta\nu\left(\tfrac{1}{4}Q-\tfrac{1}{2}\partial_{y}(yQ)\right)+i\nu^{2}\tfrac{1}{2}Q_{y}.
\end{aligned}
\label{eq:QuadraticModulation calculation nonlinearpart}
\end{equation}
\end{lem}

We will use the convention $Q\mathcal{H}(1)=0$. We drop the $\mathcal{H}(1)$
term when we derive the formal modulation laws. This issue will
not appear in the rigorous analysis in nonlinear estimates. See Section~\ref{sec:proof of lemma 5.17} for more details. 
\begin{proof}[Proof of Lemma~\ref{lem:H(P,Q)}]
The proof is a direct computation using $\mathcal{H}\partial_{y}=|D|$.
One use various formulas in Appendix~\ref{AppendixAlgebraicIdentity},
such as $\mathcal{H}\partial_{y}(y^{3}Q^{2}),\mathcal{H}\partial_{y}(y^{2}Q^{2})$,
etc. First, we calculate 
\begin{align*}
iQ^{2}\Re(QP)\widehat{P}_{1}=(-b\eta+i\eta^{2})\tfrac{1}{2}\tfrac{1}{1+y^{2}}yQ+\eta\nu\tfrac{1}{2}\tfrac{1}{1+y^{2}}Q.
\end{align*}
We deduce 
\begin{align*}
-iQ|D|\overline{P}\widehat{P}_{1}= & -iQ|D|\left(\overline{-ib\tfrac{y^{2}}{4}Q-\eta\tfrac{1+y^{2}}{4}Q+i\tfrac{\nu}{2}yQ}\right)\left(-(ib+\eta)\tfrac{y}{2}Q+i\tfrac{\nu}{2}Q\right)\\
= & -b^{2}i\tfrac{1}{8}Q|D|(y^{3}Q^{2})-b\eta\tfrac{1}{8}Q|D|(y^{3}Q^{2}-2y)-\eta^{2}i\tfrac{1}{4}Q|D|y\\
 & +b\nu i\tfrac{3}{8}Q|D|(y^{2}Q^{2})+\eta\nu\tfrac{1}{4}Q|D|(y^{2}Q^{2}-1)-\nu^{2}i\tfrac{1}{4}Q|D|yQ^{2}.
\end{align*}
Thus, we have 
\begin{equation}
	\begin{aligned}
		-iQ|D|\overline{P}\widehat{P}_{1}= & -(b^{2}+\eta^{2})\tfrac{i}{4}Q\mathcal{H}(1)+(i\tfrac{b^{2}}{2}+\tfrac{b\eta}{2})\tfrac{1}{(1+y^{2})^{2}}yQ\\
		& -(i\tfrac{3b\nu}{4}+\tfrac{\eta\nu}{2})\tfrac{1-y^{2}}{(1+y^{2})^{2}}Q-i\nu^{2}\tfrac{y}{(1+y^{2})^{2}}Q.
	\end{aligned} \label{eq:QDPP1}
\end{equation}
Similarly, using $|D|(yQ^{2})=yQ^{4}$, we obtain 
\begin{align*}
-iP|D|Q\widehat{P}_{1}= & i\tfrac{b^{2}}{2}\tfrac{y^{2}}{(1+y^{2})^{2}}yQ+\tfrac{b\eta}{2}\tfrac{1+2y^{2}}{(1+y^{2})^{2}}yQ-i\tfrac{\eta^{2}}{2}\tfrac{1}{1+y^{2}}yQ\\
 & -i\tfrac{b\nu}{4}\tfrac{(5-y^{2})y^{2}}{(1+y^{2})^{2}}Q-\tfrac{\eta\nu}{4}\tfrac{1+4y^{2}-y^{4}}{(1+y^{2})^{2}}Q+i\tfrac{\nu^{2}}{2}\tfrac{(1-y^{2})y}{(1+y^{2})^{2}}Q.
\end{align*}
By collecting $b^{2}$-terms, we have 
\begin{align*}
 & -b^{2}\tfrac{i}{4}Q\mathcal{H}(1)+i\tfrac{b^{2}}{2}\tfrac{1}{(1+y^{2})^{2}}yQ+i\tfrac{b^{2}}{2}\tfrac{y^{2}}{(1+y^{2})^{2}}yQ\\
 & =-b^{2}\tfrac{i}{4}Q\mathcal{H}(1)+i\tfrac{b^{2}}{2}\Lambda_{\frac{1}{2}}yQ.
\end{align*}
Collecting $b\eta$-terms, 
\begin{align*}
 & -b\eta\tfrac{1}{2}\tfrac{1}{1+y^{2}}yQ+\tfrac{b\eta}{2}\tfrac{1}{(1+y^{2})^{2}}yQ+\tfrac{b\eta}{2}\tfrac{1+2y^{2}}{(1+y^{2})^{2}}yQ\\
 & =b\eta\tfrac{1}{2}\Lambda_{\frac{1}{2}}yQ.
\end{align*}
Collecting $\eta^{2}$-terms, 
\begin{align*}
 & -\eta^{2}\tfrac{i}{4}Q\mathcal{H}(1)+i\tfrac{\eta^{2}}{2}\tfrac{1}{1+y^{2}}yQ-i\tfrac{\eta^{2}}{2}\tfrac{1}{1+y^{2}}yQ\\
= & -\eta^{2}\tfrac{i}{4}Q\mathcal{H}(1).
\end{align*}
Collecting $b\nu$-terms, 
\begin{align*}
 & -i\tfrac{3b\nu}{4}\tfrac{1-y^{2}}{(1+y^{2})^{2}}Q-i\tfrac{b\nu}{4}\tfrac{(5-y^{2})y^{2}}{(1+y^{2})^{2}}Q\\
 & =-ib\nu\left(\tfrac{1}{2}\Lambda Q+\tfrac{1}{2}\partial_{y}(yQ)\right).
\end{align*}
Collecting $\eta\nu$-terms, 
\begin{align*}
 & \eta\nu\tfrac{1}{2}\tfrac{1}{1+y^{2}}Q-\tfrac{\eta\nu}{2}\tfrac{1-y^{2}}{(1+y^{2})^{2}}Q-\tfrac{\eta\nu}{4}\tfrac{1+4y^{2}-y^{4}}{(1+y^{2})^{2}}Q\\
 & =\eta\nu\left(\tfrac{1}{4}Q-\tfrac{1}{2}\partial_{y}(yQ)\right).
\end{align*}
Finally, we collect $\nu^{2}$-terms, 
\begin{align*}
 & -i\nu^{2}\tfrac{y}{(1+y^{2})^{2}}Q+i\tfrac{\nu^{2}}{2}\tfrac{(1-y^{2})y}{(1+y^{2})^{2}}Q\\
 & =i\nu^{2}\tfrac{1}{2}Q_{y}.
\end{align*}
This completes the proof. 
\end{proof}
By \eqref{eq:QuadraticModulation calculation} and \eqref{eq:QuadraticModulation calculation nonlinearpart}
with the convention $\mathcal{H}(1)=0$, we deduce 
\begin{align}
\eqref{eq:w1-equ}\approx-\left(b_{s}+\frac{3b^{2}}{2}+\frac{\eta^{2}}{2}\right)\frac{i}{2}yQ-(\eta_{s}+b\eta)\frac{1}{2}yQ+(\nu_{s}+b\nu)\frac{i}{2}Q.\label{eq:FormalBlowupLawApproxEqu}
\end{align}
Therefore, we have derive the formal modulation laws as 
\begin{align}
 & \begin{split}\frac{\lambda_{s}}{\lambda}+b=0,\\
b_{s}+\frac{3}{2}b^{2}+\frac{1}{2}\eta^{2}=0,
\end{split}
 & \begin{split}\gamma_{s}-\frac{\eta}{2}=0,\\
\eta_{s}+b\eta=0,
\end{split}
 & \begin{split}\frac{x_{s}}{\lambda}+\nu=0,\\
\nu_{s}+b\nu=0,
\end{split}
\label{eq:FormalBlowupLawConclusion}
\end{align}
so that \eqref{eq:FormalBlowupLawApproxEqu} almost vanishes. A noteworthy
point here is that we are able to obtain the second modulation laws
\emph{without tail computations}. It was surprising to us that \eqref{eq:FormalBlowupLawApproxEqu}
is expressed as a linear combination of three special directions $iyQ$,
$yQ$, and $iQ$, belonging to the kernel of $A_{Q}$. When we conduct
the energy estimates later, we apply $A_{Q}$ to these directions
and they will even disappear. Typically, for example in (WM) or (CSS)
\cite{RaphaelRodnianski2012,KimKwonOh2020blowup}, the (nonlinear)
quadratic terms do not necessarily belong to the kernel of $A_{Q}$,
so these contribute as quadratic errors $O(b^{2},\eta^{2},\nu^{2})$
in the energy estimate. To avoid these large errors, one adds quadratic
terms (and possibly more) in the modified profile such as $b^{2}T_{2,0,0}$,
$b\eta T_{1,1,0}$, etc., and modify the modulation laws incorporating
effects from these corrections. In our case, we do not need to add
higher order corrector profiles such as $T_{2,0,0}$, $T_{1,1,0}$,
and so on. In other words, although the modified profiles $Q+P$ and
$\wt P_{1}$ are linear in $b,\eta,\nu$, they solve the evolution
equation even at the quadratic order. This is one of the remarkable
points of \eqref{CMdnls-gauged} at the technical level.

Before closing this section, we integrate \eqref{eq:FormalBlowupLawConclusion}
to obtain the dynamics of modulation parameters. First, we have the
following conserved quantities: 
\begin{equation}
\frac{b^{2}+\eta^{2}}{\la^{3}}\equiv\ell\in[0,+\infty),\qquad\frac{\eta}{\la}\equiv\eta_{0}\in\bbR,\qquad\frac{\nu}{\lmb}\equiv\nu_{0}\in\bbR.\label{eq:conserved quantity of modulation}
\end{equation}
With these conservation laws, we can rewrite the remaining parts of
the system as 
\[
\Big(\frac{b}{\lmb}\Big)_{t}+\frac{\ell}{2}=0,\qquad\lmb=\frac{1}{\ell}\Big\{\Big(\frac{b}{\lmb}\Big)^{2}+(\eta_{0})^{2}\Big\},\qquad\gamma_{t}-\frac{\eta_{0}}{2\lmb}=0,\qquad x_{t}+\nu_{0}=0.
\]
In particular, $x(t)=x(t_{0})-\nu_{0}(t-t_{0})$ for any $t,t_{0}\in\bbR$.
If $\ell^{2}=0$, then $b\equiv\eta\equiv0$ implies that $\lmb\equiv\lmb^{\ast}$
and $\gmm\equiv\gmm^{\ast}$ are constant in time, so we obtain the
traveling wave dynamics. If $\ell\neq0$, then we have $\frac{b(t)}{\lmb(t)}=\frac{\ell}{2}(T-t)$
for some $T\in\bbR$ and $\lmb(t)$ and $\gmm(t)$ are given by 
\begin{equation}
\begin{aligned}\lmb(t) & =\frac{\ell}{4}(t-T)^{2}+\frac{(\eta_{0})^{2}}{\ell},\\
\gmm(t) & =\gmm^{\ast}+\begin{cases}
0 & \text{if }\eta_{0}=0,\\
{\displaystyle \sgn(\eta_{0})\Big\{\arctan\Big(\frac{\ell}{2|\eta_{0}|}(t-T)\Big)+\frac{\pi}{2}\Big\}} & \text{if }\eta_{0}\neq0,
\end{cases}
\end{aligned}
\label{eq:lmb-gmm-formal-law}
\end{equation}
for some $\gmm^{\ast}\in\bbR$ and $\sgn(\eta_{0})$ denotes the sign
of $\eta_{0}$. Now we see that the blow-up ($\lmb(t)\to0$) scenario
is codimension one stable and its instability mechanism is given by
rotational instability as illustrated in the comment (rotation instability)
in Section~\ref{subsec:Discussions}.

\subsection{Modification of the profile for $w_{1}$}

\label{Subsection P1 Further modify}It turns out that the choice
$\wt P_{1}$ is not sufficient for our analysis. The reason is rather
technical, so we defer the explanation to Section~\ref{SubsectionInterpolationEstimates}.
Here, we only define the modified profile $P_{1}$ with an additional
parameter $\mu$, and compute the quadratic terms in the evolution
equation, similar to the derivation of \eqref{eq:FormalBlowupLawApproxEqu}.
Define $P_{1}$ by 
\begin{align*}
P_{1}=\wt P_{1}+\mu\tfrac{1}{2}Q=-(ib+\eta)\tfrac{y}{2}Q+(i\nu+\mu)\tfrac{1}{2}Q,
\end{align*}
where $\mu$ will be determined by \eqref{eq:definition of mu} in
Section~\ref{SubsectionInterpolationEstimates}. We note that $P_{1}$
still belongs to the kernel of $A_{Q}$. As we will see in Section~\ref{SubsectionInterpolationEstimates},
this correction is to compensate for a drawback of nonlinear adapted derivative.

We compute the modulation equation for this $P_{1}$ including $\mu$.
It is essential to check if there is a notable cancellation as before, i.e. to check if the quadratic terms belong to the kernel
of $A_{Q}$. More precisely, we claim that 
\begin{align}
\begin{split}\eqref{eq:w1-equ}\approx & -\left(b_{s}+\frac{3b^{2}}{2}+\frac{\eta^{2}}{2}\right)\frac{i}{2}yQ-(\eta_{s}+b\eta)\frac{1}{2}yQ\\
 & +(\nu_{s}+b\nu)\frac{i}{2}Q+(\mu_{s}+b\mu)\frac{1}{2}Q.
\end{split}
\label{eq:FormalBlowupLawApproxEquModified}
\end{align}
Without $\mu$, we saw this in Section~\ref{subsec:Formal-modulation-laws}.
Thus, it suffices to check the quadratic terms involving $\mu$.

In the linear part for $w_{1}$, we have 
\begin{align*}
 & \quad(\partial_{s}-\frac{\lambda_{s}}{\lambda}\Lambda_{-1}+\gamma_{s}i-\frac{x_{s}}{\lambda}\partial_{y})(\tfrac{\mu}{2}Q)\\
 & \approx\tfrac{\mu_{s}}{2}Q+b\mu\tfrac{1}{2}\Lambda_{-1}Q+\tfrac{\eta\mu}{4}iQ-\nu\mu\tfrac{1}{2}Q_{y}.
\end{align*}
In $i\frac{1}{4}|w|^{4}w_{1}$, we have 
\begin{align*}
(-i\tfrac{\eta}{2}Q^{2})\cdot(\tfrac{\mu}{2}Q)=-i\tfrac{\eta\mu}{4}Q^{3}
\end{align*}
In $-iw|D|\overline{w}w_{1}$, we compute 
\begin{align*}
 & \quad-iQ|D|(\overline{i\tfrac{\nu}{2}yQ})(\tfrac{\mu}{2}Q)-i\cdot(i\tfrac{\nu}{2}yQ)|D|Q(\tfrac{\mu}{2}Q)\\
 & \quad-iQ|D|(\overline{-ib\tfrac{y^{2}}{4}Q-\eta\tfrac{1+y^{2}}{4}Q})(\tfrac{\mu}{2}Q)-i\cdot(-ib\tfrac{y^{2}}{4}Q-\eta\tfrac{1+y^{2}}{4}Q)|D|Q(\tfrac{\mu}{2}Q)\\
 & =\tfrac{\nu\mu}{2}Q^{2}Q_{y}+\tfrac{\nu\mu}{4}yQ^{3}+\tfrac{\nu\mu}{2}y^{2}Q^{2}Q_{y}-\tfrac{b\mu}{4}\tfrac{1-y^{2}}{(1+y^{2})^{2}}Q-\tfrac{b\mu}{4}\tfrac{(1-y^{2})y^{2}}{(1+y^{2})^{2}}Q+i\tfrac{\eta\mu}{4}\tfrac{1-y^{2}}{1+y^{2}}Q.
\end{align*}
Now, we collect the quadratic terms. 
\begin{itemize}
\item $b\mu$ : Using $\Lambda_{-1}Q=\tfrac{1-y^{2}}{2(1+y^{2})}Q+Q$, we
compute 
\begin{align*}
b\mu\tfrac{1}{2}\Lambda_{-1}Q-\tfrac{b\mu}{4}\tfrac{1-y^{2}}{(1+y^{2})^{2}}Q-\tfrac{b\mu}{4}\tfrac{(1-y^{2})y^{2}}{(1+y^{2})^{2}}Q=b\mu\tfrac{1}{2}Q.
\end{align*}
\item $\eta\mu:$ we have by a direct computation 
\begin{align*}
\tfrac{\eta\mu}{4}iQ-i\tfrac{\eta\mu}{4}Q^{3}+i\tfrac{\eta\mu}{4}\tfrac{1-y^{2}}{1+y^{2}}Q=0.
\end{align*}
\item $\nu\mu:$using $yQ^{3}=-2Q_{y}$, we compute 
\begin{align*}
-\nu\mu\tfrac{1}{2}Q_{y}+\tfrac{\nu\mu}{2}Q^{2}Q_{y}+\tfrac{\nu\mu}{4}yQ^{3}+\tfrac{\nu\mu}{2}y^{2}Q^{2}Q_{y}=0.
\end{align*}
\end{itemize}
In conclusion, the quadratic terms involving $\mu$ are summarized
as 
\begin{align}
+(\mu_{s}+b\mu)\tfrac{1}{2}Q.\label{eq:FormalBlowupLawApproxEquModify-mu}
\end{align}
Combining this with \eqref{eq:FormalBlowupLawApproxEqu}, we obtain
\eqref{eq:FormalBlowupLawApproxEquModified}.

With those modified profiles $Q+P$ and $P_{1}$, we will fix the
decomposition of $w$ and $w_{1}$ in Section~\ref{Subsection Decomposition}
and Section~\ref{SubsectionInterpolationEstimates}.

\section{Trapped solution\label{sec:Trapped-solution}}

From the formal modulation laws \eqref{eq:FormalBlowupLawConclusion}
in Section~\ref{SectionFormalBlowupLaw}, we anticipate the existence
of finite-time blow-up solutions to \eqref{CMdnls-gauged} with the
blow-up profile $Q$. This section aims to construct these blow-up
solutions by bootstrapping. We will decompose the solution $v(t)$
in the form \eqref{eq:v Decomposition form}, $v=[Q+P(b,\eta,\nu)+\eps]_{\lambda,\gamma,x}$
with the modulation parameters $\la,\ga,x,b,\eta,\nu$. Inspired from
the invariant subspace decomposition, Remark~\ref{Remark Invariant genKernel and Transversality},
we will fix the decomposition through orthogonality conditions on
$\eps$. By bootstrapping, we will justify the modulation laws \eqref{eq:FormalBlowupLawConclusion}. This will be conducted by
exploiting a conservation law at the $\dot{H}^{2}$-level, which ensures the propagation of the smallness of $\eps$. We follow the
road map developed in other dispersive wave models, such as wave maps,
Schrödinger maps, Chern--Simons--Schrödinger equations, etc. \cite{RodnianskiSterbenz2010,RaphaelRodnianski2012,MerleRaphaelRodnianski2013Invention,MerleRaphaelRodnianski2015CambJMath,KimKwon2020blowup,KimKwonOh2020blowup,Kim2022CSSrigidityArxiv}

As we saw in Section~\ref{SectionFormalBlowupLaw}, we will view
\eqref{CMdnls-gauged} as a system of the (renormalized) solution
$w$ and its nonlinear adapted derivative $w_{1}$. This strategy
was first introduced in \cite{KimKwonOh2020blowup}. We will decompose
$(w,w_{1})$ using the modified profiles ($Q+P,P_{1})$. One of the
main novelties of this work is to exploit different decompositions
for $(w,w_{1})$, depending on the norms we measure. This is because
$P$ and $P_{1}$ have slow decays (or even growth) and do not belong
to $L^{2}$; see Section~\ref{Subsection Decomposition} and Section~\ref{SubsectionInterpolationEstimates}
for more details. This strategy allows us to avoid introducing any
cut-offs in the modified profiles. Due to the slow decay of $P$ and
$P_{1}$, if there were cut-offs, they would create large cut-off
errors. On the other hand, in the energy estimates, we will exploit the higher-order conservation laws associated with nonlinear variables $w$ and $w_1$. This takes advantage of the integrability and the nonlinear adapted derivatives. This strategy extremely simplifies the higher-order energy argument, and this corresponds to exploiting a repulsive dynamics of the $\eps$-part, \eqref{eq:BQ eps1 equ}, as explained in Section~\ref{SubsectionConjugationAQ}.

Section~\ref{sec:Trapped-solution} is organized as follows. In Section~\ref{Subsection Decomposition},
we fix the decomposition of solutions by imposing orthogonality conditions.
We also discuss the strategy to use two different decompositions depending
on the topology. In Section~\ref{Subsection ExistenceTrappedSol},
we reduce our main Theorem~\ref{TheoremCodimension2Blowup} to bootstrap
propositions and topological lemmas, Propositions~\ref{PropositionMainBootstrap},
\ref{PropSharpDescription} and Lemmas~\ref{LemmaTopologicalPhiContinuous},
\ref{LemmaPsiNonExistenceBrouwer}. In Section~\ref{SubsectionInterpolationEstimates}--\ref{SubsectionBootstrapDescription},
we are under the bootstrap assumption \eqref{eq:BootstapAssumption}.
In Section~\ref{SubsectionInterpolationEstimates}, we explain further
modification of the profile $P_{1}$. Also, we obtain various interpolation
estimates originated from the bootstrap bounds. In Section~\ref{SubsectionModulationEstimate},
we prove modulation estimates that justify formal modulation laws
\eqref{eq:FormalBlowupLawConclusion}. Additionally, we improve modulation
estimates for $b_{s},\eta_{s},\nu_{s}$, by introducing refined modulation
parameters $\widetilde{b},\widetilde{\eta},\widetilde{\nu}$. These
estimates will be proved assuming the nonlinear estimate, Lemma~\ref{LemmaNonlinearEstimate1}.
We set aside nonlinear estimates used in Section~\ref{sec:Trapped-solution}
to Section~\ref{Section Proof of nonlinear estimates}. Next, in Section~\ref{SubsectionBootstrapDescription},
we prove the bootstrap proposition, Proposition~\ref{PropositionMainBootstrap}
and a sharp description of solutions, Proposition~\ref{PropSharpDescription}.
In addition, we also prove topological lemmas, Lemma~\ref{LemmaTopologicalPhiContinuous}
and Lemma~\ref{LemmaPsiNonExistenceBrouwer}. Up to this point, we
construct codimension one finite time blow-up solutions (Theorem~\ref{TheoremCodimension2Blowup}).
Finally, in Section~\ref{SubsectionChiralBlowup}, we prove Theorem~\ref{TheoremChiralBlowup};
we construct finite time blow-up solutions in Hardy space, i.e., the
\emph{chiral blow-up solutions}. This will be done with a suitable
modification of the earlier argument.

\subsection{Decomposition\label{Subsection Decomposition}}

In this subsection, our aim is to decompose the nonlinear solution $v(t,x)$
as mentioned above. We hope to decompose $v(t,x)$ in the form:
\begin{align}
v(t,x)=\frac{e^{i\gamma(t)}}{\lambda(t)^{1/2}}[Q+P(b(t),\eta(t),\nu(t);\cdot)+\eps(t,\cdot)]\left(\frac{x-x(t)}{\lambda(t)}\right).\tag{\ref{eq:v Decomposition form}}
\end{align}
Recall that we did not introduce any cut-offs in the definition of
$P$, in order to avoid (potentially dangerous) cut-off errors in
the analysis. However, the decomposition \eqref{eq:v Decomposition form}
does not make sense in $L^{2}$ due to the spatial growth (as $|x|$)
of $P$. In view of nonlinear estimates, we need a $L^{2}$ decomposition
also, in which case we will use a rough decomposition in $H^{1}$
\[
v(t,x)=[Q+\widehat{\eps}]_{\lambda(t),\gamma(t),x(t)}(x).
\]
In higher topology $\dot{\mathcal{H}}^{2}$, as $P\in\dot{\calH}^{2}$,
we are able to use the refined decomposition 
\[
v(t,x)=[Q+P(\cdot;b(t),\eta(t),\nu(t))+\eps]_{\lambda(t),\gamma(t),x(t)}(x).
\]
In view of the formal invariant subspace decomposition in Remark~\ref{Remark Invariant genKernel and Transversality},
we formally put $\widehat{\eps}\in\textnormal{ker}(\mathcal{L}_{Q}i)^{\perp}=\{iQ_{y},i\Lambda Q,Q\}^{\perp}$
and $\eps\in N_{g}(\mathcal{L}_{Q}i)^{\perp}=\{y^{2}Q,i(1+y^{2})Q,yQ,iQ_{y},i\Lambda Q,Q\}^{\perp}$
by orthogonality conditions. However, we will impose the truncated
ones because of the slow decay of generalized kernel elements. Recall
the generalized kernel elements $\mathcal{K}_{j}$ given by \eqref{eq:3.3 def ker elem calK}. 
\begin{lem}[Transversality]
\label{LemmaZkTransversality} There exist profiles $\mathcal{Z}_{k}\in C_{c}^{\infty}$
that are transversal to $N_{g}(i\mathcal{L}_{Q})$, i.e. $\mathcal{Z}_{k}'s$
satisfy the following transversality; the $6\times6$ matrix $\mathcal{A}_{jk}\coloneqq(\mathcal{K}_{j},\mathcal{Z}_{k})_{r}$
becomes a nonsingular diagonal matrix. 
\end{lem}

\begin{proof}
We fix $R_{0}$ to be a universal constant, say $R_{0}=10000$. We
define $\mathcal{Z}_{k}$ by 
\begin{align*}
\mathcal{Z}_{1} & =y^{2}Q\chi_{R_{0}}-\tfrac{((1+y^{2})Q,y^{2}Q\chi_{R_{0}})_{r}}{2(yQ,yQ\chi_{R_{0}})_{r}}\mathbf{D}_{Q}^{*}(yQ\chi_{R_{0}}),\\
\mathcal{Z}_{2} & =i(1+y^{2})Q\chi_{R_{0}}-\tfrac{(y^{2}Q,(1+y^{2})Q\chi_{R_{0}})_{r}}{2(yQ,yQ\chi_{R_{0}})_{r}}\mathbf{D}_{Q}^{*}(iyQ\chi_{R_{0}}),\\
\mathcal{Z}_{3} & =yQ\chi_{R_{0}},\\
\mathcal{Z}_{4} & =\mathbf{D}_{Q}^{*}(iyQ\chi_{R_{0}}),\\
\mathcal{Z}_{5} & =\mathbf{D}_{Q}^{*}(yQ\chi_{R_{0}})-\tfrac{(\Lambda Q,\mathbf{D}_{Q}^{*}(yQ\chi_{R_{0}}))_{r}}{(\Lambda Q,\chi_{R_{0}})_{r}}Q\chi_{R_{0}},\\
\mathcal{Z}_{6} & =\mathbf{D}_{Q}^{*}(iQ\chi_{R_{0}}).
\end{align*}
Then, by direct computation, we have 
\begin{align}
(\mathcal{K}_{j},\mathcal{Z}_{k})=c_{j}\delta_{jk},\label{eq:transversality}
\end{align}
where $\delta_{jk}$ is the Kronecker delta, and 
\[
c_{1}=(\Lambda Q,y^{2}Q\chi_{R_{0}})_{r}-\tfrac{((1+y^{2})Q,y^{2}Q\chi_{R_{0}})_{r}}{2(yQ,yQ\chi_{R_{0}})_{r}}(\D_{Q}(\Lambda Q),yQ\chi_{R_{0}})=O(R_{0})-O(R_{0}^{2})\ne0.
\]
Similarly, we compute 
\begin{align*}
c_{2}=(1,\chi_{R_{0}})_{r}=O(R_{0}), & \quad c_{3}=(Q_{y},yQ\chi_{R_{0}})_{r}=O(1),\\
c_{4}=2(yQ,yQ\chi_{R_{0}})_{r}=O(R_{0}), & \quad c_{6}=(Q,Q\chi_{R_{0}})_{r}=O(1),
\end{align*}
\[
c_{5}=(2yQ,yQ\chi_{R_{0}})_{r}-\tfrac{(Q^{2}Q_{y},yQ\chi_{R_{0}})_{r}}{(\Lambda Q,Q\chi_{R_{0}})_{r}}2(1,\chi_{R_{0}})_{r}=O(R_{0})-O(R_{0}^{2})\neq0.
\]
This finishes the proof. 
\end{proof}
Denote by $\widehat{\eps}\in\widehat{\mathcal{Z}}^{\perp}$ and $\eps\in\mathcal{Z}^{\perp}$
when $\wt{\eps}$ and $\eps$ satisfy the orthogonality conditions:
\begin{align*}
\begin{gathered}\widehat{\mathcal{Z}}^{\perp}\coloneqq\{\widehat{\eps}\in L^{2}:(\widehat{\eps},\mathcal{Z}_{k})_{r}=0,k=1,2,3\},\\
\mathcal{Z}^{\perp}\coloneqq\{\eps\in\langle y\rangle^{2}L^{2}:(\eps,\mathcal{Z}_{k})_{r}=0,k=1,2,\cdots,6\}.
\end{gathered}
\end{align*}
Now we define the sets on which decompositions will be performed.
First, we denote the soliton tube $\mathcal{O}_{dec}$ for a small
parameter $\delta_{dec}>0$ to be determined later. 
\begin{align*}
\mathcal{O}_{dec}\coloneqq\{[Q+\widehat{\eps}_{0}]_{\lambda_{0},\gamma_{0},x_{0}}(x)\in L^{2}:(\lambda_{0},\gamma_{0},x_{0})\in\R_{+}\times\R/2\pi\Z\times\R,\,\|\widehat{\eps}_{0}\|_{L^{2}}<\delta_{dec}\}.
\end{align*}
To define the set of initial data, we first denote the set of coordinates
\[
\mathcal{U}_{init}\coloneqq\{(\lambda_{0},\gamma_{0},x_{0},b_{0},\eta_{0},\nu_{0},\tilde{\eps}_{0})\in\mathbb{R}_{+}\times\mathbb{R}/2\pi\mathbb{Z}\times\mathbb{R}^{4}\times\mathcal{Z}^{\perp}:\text{satisfying \eqref{eq:InitialConditionBootstrap}}\},
\]
where $b^{*}>0$ will be chosen small depending on $\delta_{dec}$,
a universal small constant $\kappa>0$ ($\kappa=\frac{1}{100}$ will
suffice), and the assumptions for $\mathcal{U}_{init}$ are given
by 
\begin{align}
0<b_{0}<\frac{1}{2}b^{*},\quad|\eta_{0}|<\frac{1}{2}b_{0}^{1+\kappa},\quad|\nu_{0}|<\frac{1}{2}b_{0}^{1-\kappa},\quad\frac{1}{1.1}<\frac{b_{0}^{2}}{\lambda_{0}^{3}}<1.1,\quad\|\widetilde{\eps}_{0}\|_{H^{2}}<b_{0}^{2}.\label{eq:InitialConditionBootstrap}
\end{align}
The initial data set $\mathcal{O}_{init}\subset\mathcal{O}_{dec}$
is then defined by 
\begin{align*}
\mathcal{O}_{init}\coloneqq\{[Q+P(b_{0},\eta_{0},\nu_{0})\chi_{\delta_{dec}\la_{0}^{-1}}+\widetilde{\eps}_{0}]_{\lambda_{0},\gamma_{0},x_{0}}(x):(\lambda_{0},\gamma_{0},x_{0},b_{0},\eta_{0},\nu_{0},\widetilde{\eps}_{0})\in\mathcal{U}_{init}\}.
\end{align*}
We will also express codimension two coordinates of initial data.
Define 
\begin{align*}
\widetilde{\mathcal{U}}_{init}\coloneqq & \Big\{(\lmb_{0},\gmm_{0},x_{0},b_{0},\td{\eps}_{0})\in\mathbb{R}_{+}\times\mathbb{R}/2\pi\mathbb{Z}\times\mathbb{R}^{2}\times\mathcal{Z}^{\perp}:\\
 & \qquad\qquad0<2b_{0}<b^{*},\quad\|\widetilde{\eps}_{0}\|_{H^{2}}<b_{0}^{2},\quad\frac{1}{1.1}<\frac{b_{0}^{2}}{\lambda_{0}^{3}}<1.1\Big\}.
\end{align*}

In the above definition of $\calO_{init}$, the cut-off $\chi_{\delta_{dec}\la_{0}^{-1}}$
is introduced to ensure $\calO_{init}\subset H^{2}$. However, in
the later dynamical analysis, we will switch to the decomposition
without cut-offs to avoid any dangerous cut-off errors. 
\begin{lem}[Decomposition]
\label{LemmaDecomposition} There exist $\delta',\delta'',\text{ and }\delta_{dec}>0$
such that the following hold. 
\begin{enumerate}
\item[(1)] (The decomposition on $L^{2}$) There exists a map $(\mathbf{G},\wt{\eps}):\mathcal{O}_{dec}\to\mathbb{R}_{+}\times\mathbb{R}/2\pi\mathbb{Z}\times\mathbb{R}\times\wt{\mathcal{Z}}^{\perp}$
satisfying 
\begin{align*}
v=[Q+\widehat{\eps}]_{\lambda,\gamma,x},\quad(\widehat{\eps},\mathcal{Z}_{k})_{r}=0\quad\text{for }k=1,2,3,
\end{align*}
with $\textnormal{\textbf{G}}=(\lambda,\gamma,x)$, and $\|\widehat{\eps}\|_{L^{2}}<\delta^{\prime}$. 
\item[(2)] (The decomposition on $\langle y\rangle^{2}L^{2}$) Let $(\mathbf{G},\wt{\eps})$
be the map obtained in (1). Then, there is a further decomposition
map $(\mathbf{G},\mathbf{H},\eps):\mathcal{O}_{dec}\to\mathbb{R}_{+}\times\mathbb{R}/2\pi\mathbb{Z}\times\mathbb{R}\times\R^{3}\times\mathcal{Z}^{\perp}$
satisfying 
\begin{align*}
v=[Q+P(b,\eta,\nu)+\eps]_{\lambda,\gamma,x},\quad P(b,\eta,\nu)+\eps=\wt{\eps}
\end{align*}
with $\mathbf{G}=(\lambda,\gamma,x)$, $\mathbf{H}=(b,\eta,\nu)$,
$(\eps,\mathcal{Z}_{k})_{r}=0$ for $k=1,\cdots,6$, and $\|\langle y\rangle^{-2}\eps\|_{L^{2}}<\delta^{\prime\prime}$. 
\item[(3)] ($C^{1}$-regularity) The map $v\mapsto(\lambda,\gamma,x,b,\eta,\nu)$
for each decomposition is $C^{1}$. 
\item[(4)] (Initial data set) We have $\mathcal{O}_{init}\subseteq\mathcal{O}_{dec}$,
and the statements of (1), (2), and (3) also hold when we replace
$\mathcal{O}_{dec}$ by $\mathcal{O}_{init}$. Moreover, there exists
a homeomorphism in $H^{2}$ from $\overline{\mathcal{U}_{init}}$
to $\overline{\mathcal{O}_{init}}$. 
\end{enumerate}
\end{lem}

The proof of Lemma~\ref{LemmaDecomposition} is rather standard and
relegated to Appendix~\ref{AppendixDecomposition}. The proof of
(1) is a standard application of Implicit Function Theorem. Once we
obtain a decomposition $(\mathbf{G},\wt{\eps})$, then we can impose
orthogonality conditions $(\eps,\mathcal{Z}_{k})=0$ for $k=4,5,6$
to fix decomposition $\wt{\eps}=P(b,\eta,\nu)+\eps$. This process
is explicit since $P$ is linear in $b,\eta,\nu$ \eqref{eq:definition P}.
Then, by \eqref{eq:transversality}, $(\eps,\mathcal{Z}_{k})=0$ for
$k=1,2,3$ follow. The second decomposition is done in a weighted
space $\langle y\rangle^{2}L^{2}$ since $P(b,\eta,\nu)$ lacks spatial
decay and so does $\eps$. This is also the main reason why we need
two different decompositions for different topologies.
\begin{rem}
\label{RemarkDecomposeModulationParamterEqual} As mentioned above,
the main reason for introducing a cut-off in the definition of $\mathcal{O}_{init}$
is to ensure $\calO_{init}\subset H^{2}$. Moreover, one can choose
compactly supported initial data (despite being non-chiral). In
the later dynamical analysis, we will use the decomposition in Lemma~\ref{LemmaDecomposition},
which does not involve any cut-offs. After this decomposition, $\eps\notin L^{2}$
in general, due to the spatial tail of $P$.

We also remark that \emph{the modulation parameters do not change}
when decomposing $Q+P\chi_{\delta_{dec}\la_{0}^{-1}}+\tilde{\eps}_{0}$
(with parameters in $\mathcal{U}_{init}$) into the form $Q+P+\eps_{0}$
as in Lemma~\ref{LemmaDecomposition}. Indeed, we have $\eps_{0}=\widetilde{\eps}_{0}-P(b_{0},\eta_{0},\nu_{0})(1-\chi_{\delta_{dec}\la_{0}^{-1}})$
and this $\eps_{0}$ still satisfies the same orthogonality conditions
$(\eps_{0},\mathcal{Z}_{k})=0$ for all $k=1,\dots,6$ because $\mathcal{Z}_{k}$s
are compactly supported (say $\text{supp}\mathcal{Z}_{k}\subset[-2R_0,2R_0]$)
and $\delta_{dec}\la_{0}^{-1}\gtrsim\delta_{dec}b_{0}^{-\frac{2}{3}}\gg R_0$
for $0<b_{0}<b^{*}\ll\delta_{dec}$. 
\end{rem}

\subsection{Trapped solutions and precise statement of Theorem~\ref{TheoremCodimension2Blowup}\label{Subsection ExistenceTrappedSol}}

In this subsection, we state the precise version (Theorem~\ref{thm:precise statement codimensionone blowup})
of Theorem~\ref{TheoremCodimension2Blowup} and reduce its proof
to Propositions~\ref{PropositionMainBootstrap}, \ref{PropSharpDescription}
and Lemmas~\ref{LemmaTopologicalPhiContinuous}, \ref{LemmaPsiNonExistenceBrouwer}.
The core of this process is the main bootstrap Proposition~\ref{PropositionMainBootstrap},
from which we can confirm that the smallness (time decay) for $\eps$
and the controls for the stable modulation parameters $\la$ and $b$
propagate. In addition, it also justifies the modulation laws \eqref{eq:FormalBlowupLawApproxEqu}.
Finally, we use a topological argument to show the existence of blow-up
solutions.

Here, we use the conservation laws from the integrability, Corollary \ref{cor:hierarchy}. Denoting $v_1:=\D_vv$, we have
\begin{align*}
	\|v_{1}(0)\|_{L^2}=\|v_1(t)\|_{L^2},\quad \|\td \bfD_{v_0} v_{1}(0)\|_{L^2}=\|\td \bfD_{v(t)} v_1(t)\|_{L^2}
\end{align*}
for any time $t$ in the maximal interval of existence. Renormalizing them, we derive
\begin{align*}
	\|v_{1}(0)\|_{L^2}=\frac{\|w_1(s)\|_{L^2}}{\lambda(s)},\quad 
	\|\td \bfD_{v_0} v_{1}(0)\|_{L^2}=\frac{\|\td \bfD_{w(s)} w_1(s)\|_{L^2}}{\lambda^2(s)},
\end{align*}
which yield 
\begin{align}
	\|w_1(s)\|_{L^2}\lesssim_{v_0} \lambda(s),\quad \|\td \bfD_{w(s)} w_1(s)\|_{L^2}\lesssim_{v_0} \lambda^2(s) \label{eq:energy scale bounds}
\end{align}
for any $s\in [s(0),\infty)$. Note that \eqref{eq:energy scale bounds} hold true \textit{without any bootstrap argument}. We will show the behaviors of the modulation parameters by bootstrapping.

We state \emph{bootstrap assumptions} as follows:
\begin{align}
\begin{split}  0<b<&b^{*},\quad|\eta|<\lambda^{\frac{3}{2}(1+\kappa)},\quad|\nu|<\lambda^{\frac{3}{2}(1-\kappa)},\\
 &\frac{1}{2}\leq\frac{b^{2}}{\lambda^{3}}\leq2,\quad  \|\widehat{\eps}\|_{L^{2}}\leq\frac{1}{10}\delta_{dec}.
\end{split}
\label{eq:BootstapAssumption}
\end{align}
Recall that $0<\delta_{dec}\ll1$ is a small universal constant in Lemma~\ref{LemmaDecomposition}, $\kappa=\frac{1}{100}$ is as in \eqref{eq:InitialConditionBootstrap}, and $R_0$ is given in the proof of Lemma~\ref{LemmaZkTransversality}.
We will choose a small constant $b^{*}$ depending on $\delta_{dec}$ and
$R_0$.

We call the solution $v$ a \emph{trapped solution} if it satisfies
\eqref{eq:BootstapAssumption} for $0\le t<T$, where $[0,T)$ is
its maximal forward lifespan. We are now ready to provide a rigorous
statement of Theorem~\ref{TheoremCodimension2Blowup}. 
\begin{thm}[Codimension one blow-up construction]
\label{thm:precise statement codimensionone blowup} There exists a constant
$b^{*}>0$ satisfying the following properties. Let $(\la_{0},\ga_{0},x_{0},b_{0},\tilde{\eps}_{0})\in\tilde{\mathcal{U}}_{init}$.
Then, there exists $(\eta_{0},\nu_{0})\in(-\tfrac{1}{2}b^{1+\kappa},\tfrac{1}{2}b^{1+\kappa})\times(-\tfrac{1}{2}b^{1-\kappa},\tfrac{1}{2}b^{1-\kappa})$
such that the solution $v(t)$ to \eqref{CMdnls-gauged} starting
from the initial data 
\begin{equation}
v_{0}=[Q+P(b_{0},\eta_{0},\nu_{0})\chi_{\delta_{dec}\la_{0}^{-1}}+\widetilde{\eps}_{0}]_{\lambda_{0},\gamma_{0},x_{0}}\in\mathcal{O}_{init}\label{eq:initial data form v_0}
\end{equation}
satisfies: 
\begin{enumerate}
\item (Existence of trapped solutions) $v(t)$ is a trapped solution and
blows up in finite time $T=T(v_{0})\in(0,\infty)$. 
\item (Sharp description of the blow-up) There exist $\ell=\ell(v_{0})\in(0,\infty)$,
$\gamma^{*}=\gamma^{*}(v_{0})\in\mathbb{R}$, $x^{*}=x^{*}(v_{0})\in\mathbb{R}$,
and $v^{*}=v^{*}(v_{0})\in H^{1}$ such that 
\[
v(t,x)-\frac{e^{i\gamma^{*}}}{\sqrt{\ell(T-t)^{2}}}Q\left(\frac{x-x^{*}}{\ell(T-t)^{2}}\right)\to v^{*}\quad\text{in }L^{2}\quad\text{as }t\to T.
\]
\end{enumerate}
Moreover, by applying the Galilean boost to $v(t)$ with any parameter
$c\in\bbR$, 
\[
v(t,x)\mapsto v^{c}(t,x)=e^{icx-itc^{2}}v(t,x-2ct),\quad(c\in\mathbb{R}),
\]
we further obtain an one-parameter family of blow-up solutions. 
\end{thm}

\begin{rem}[Codimension one stability]
\label{rem:codim-one}In the first part of Theorem~\ref{thm:precise statement codimensionone blowup},
we obtain a codimension two set of initial data leading to the blow-up
dynamics as described in the theorem. Although it is not explicitly
mentioned, the spatial center $x(t)$ of the soliton part of any trapped
solution has \emph{vanishing velocity} as $t\to T$ in view of $|\nu|\lesssim b^{1-\kappa}\ll\la$;
see also \eqref{eq:gmm-x-deriv-in-time}. Hence, the Galilean boosted
solutions $v^{c}(t)$ have the asymptotic velocity $2c$ and are
genuinely different from $v(t)$. This is why Theorem~\ref{thm:precise statement codimensionone blowup}
indeed demonstrates the `codimension one' stability of the blow-up.

One can further study the regularity issues of this codimension one
set as in \cite{Collot2018MemAmer,KimKwon2020blowup}, but we do not
pursue this direction here. 
\end{rem}

\begin{rem}[Parameter dependence]
We can track in what order we choose small parameters $\delta_{dec}$,
$\kappa$, and $b^{*}$. We will presume the size dependence of those
parameters: 
\begin{align}
0<b^{\ast} \ll\delta_{dec}^{\frac{100}{\kappa}\cdot\frac{1}{4}}\ll1.\label{eq:5.2 rem para dependecy}
\end{align}
This dependency is mainly used for proving the main bootstrap Proposition~\ref{PropositionMainBootstrap}.
We note that the universal constant $\kappa=\frac{1}{100}$ was
fixed. 
\end{rem}

In the sequel, we will prove Theorem~\ref{thm:precise statement codimensionone blowup}
assuming Propositions~\ref{PropositionMainBootstrap}, \ref{PropSharpDescription}
and Lemmas~\ref{LemmaTopologicalPhiContinuous}, \ref{LemmaPsiNonExistenceBrouwer}.
We prove these propositions and lemmas in Section~\ref{SubsectionBootstrapDescription}.

Before starting the proof, we write the \emph{initial bound condition}
for the initial modulation parameters. 
\begin{align}
\begin{split}  0<&b_{0}<\frac{1}{2}b^{*},\quad|\eta_{0}|<\frac{1}{2}b_{0}^{1+\kappa},\quad|\nu_{0}|<\frac{1}{2}b_{0}^{1-\kappa},\\
 & \frac{1}{1.1}<\frac{b_{0}^{2}}{\lambda_{0}^{3}}<1.1,\quad \|\widehat{\eps}_{0}\|_{L^{2}}\leq\frac{1}{40}\delta_{dec}.
\end{split}
\label{eq:5.2 ini boots assump}
\end{align}

We note that if the initial data after being decomposed according
to Lemma~\ref{LemmaDecomposition} satisfies \eqref{eq:5.2 ini boots assump},
then it also verifies the bootstrap assumption \eqref{eq:BootstapAssumption}
at initial time. We now check that the initial data in $\mathcal{O}_{init}$
satisfies the initial bound condition \eqref{eq:5.2 ini boots assump}
(after being decomposed according to Lemma~\ref{LemmaDecomposition}). 
\begin{lem}[Initial bound condition]
\label{lem:initial bound}Any $v_{0}\in\mathcal{O}_{init}$ (see
\eqref{eq:initial data form v_0}) satisfies the initial bound condition
\eqref{eq:5.2 ini boots assump} when it is decomposed according to
Lemma~\ref{LemmaDecomposition}.
\end{lem}
\begin{proof}
	We begin with the initial data of the form 
	\begin{align*}
		v_{0}=[Q+P(b_{0},\eta_{0},\nu_{0})\chi_{\delta_{dec}\la_{0}^{-1}}+\widetilde{\eps}_{0}]_{\lambda_{0},\gamma_{0},x_{0}}.
	\end{align*}
	Using the decomposition Lemma~\ref{LemmaDecomposition}, we write
	$v_{0}=[Q+P+\eps]_{\la,\ga,x}.$ Observing the truncation in $v_{0}$,
	the parameters $\la,\ga,x,b,\eta,\nu$ in the decomposition remain
	unchanged (see Remark~\ref{RemarkDecomposeModulationParamterEqual}).
	So, from \eqref{eq:InitialConditionBootstrap}, the modulation parameters automatically satisfy \eqref{eq:5.2 ini boots assump}. We will check
the second line of \ref{eq:5.2 ini boots assump}. For $\wt{\eps}_{0}$,
we have 
	\begin{align*}
		\|\wt{\eps}_{0}\|_{L^{2}}=\|P(b_{0},\eta_{0},\nu_{0})\chi_{\delta_{dec}\la_{0}^{-1}}+\widetilde{\eps}_{0}\|_{L^{2}}=O(\delta_{dec}^{\frac{3}{2}}+b_{0}^{2})\leq\tfrac{1}{40}\delta_{dec}.
	\end{align*}
	Hence, we conclude that the initial data $v_{0}$ with
the parameters in $\mathcal{U}_{init}$ satisfy \eqref{eq:5.2 ini boots assump}. 
\end{proof}

We now prove our main blow-up construction theorem for \eqref{CMdnls-gauged}.

\begin{proof}[Proof of Theorem~\ref{thm:precise statement codimensionone blowup}
assuming Propositions~\ref{PropositionMainBootstrap}, \ref{PropSharpDescription},
and Lemmas~\ref{LemmaTopologicalPhiContinuous}, \ref{LemmaPsiNonExistenceBrouwer}]

Let $(\lambda_{0},\gamma_{0},x_{0},b_{0},\tilde{\eps}_{0})\in\tilde{\mathcal{U}}_{init}$.
For any $(\eta_{0},\nu_{0})$ varying in the range $|\eta_{0}|\leq\frac{1}{2}b_{0}^{1+\kappa}$
and $|\nu_{0}|\leq\frac{1}{2}b_{0}^{1-\kappa}$, define the initial
data $v_{0}\in\mathcal{O}_{init}$ using the formula \eqref{eq:initial data form v_0}.
Let $v(t)$ be the (forward) maximal solution to \eqref{CMdnls-gauged}
with the initial data $v_{0}$ and lifespan $[0,T)$. Our goal is
to show that there exists $(\eta_{0},\nu_{0})$ such that $v(t)$
becomes a trapped solution. Define the exit time of the trapped regime:
\begin{align}
T_{exit}\coloneqq\sup\{\tau\in[0,T):v(\tau')\in\calO_{dec}\text{ and }\eqref{eq:BootstapAssumption}\text{ holds for all }\tau^{\prime}\in[0,\tau]\}\in(0,T].\label{eq:DefinitionTexit}
\end{align}

Note also that 
\begin{align*}
(\lambda_{0},\gamma_{0},x_{0},b_{0},\eta_{0},\nu_{0})=(\lambda(0),\gamma(0),x(0),b(0),\eta(0),\nu(0))
\end{align*}
by Remark~\ref{RemarkDecomposeModulationParamterEqual}. As \eqref{eq:BootstapAssumption}
is satisfied at $t=0$ by Lemma~\ref{lem:initial bound}, we have
$T_{exit}\in(0,T]$. Our goal (the existence of trapped solution)
is then equivalent to showing that $T=T_{exit}$ for some $(\eta_{0},\nu_{0})$.

We argue by contradiction. Suppose $T_{exit}<T$ for all such $(\eta_{0},\nu_{0})$.
Now we state the main bootstrap proposition, which will be proved
in Section~\ref{SubsectionBootstrapDescription}. 
\begin{prop}[Main bootstrap]
\label{PropositionMainBootstrap} Let the parameter $b^{*}$ satisfy the parameter dependence \eqref{eq:5.2 rem para dependecy}.
Let $v(t)$ be a solution to \eqref{CMdnls-gauged} that admits the
decomposition as in Lemma~\ref{LemmaDecomposition} with the parameters
$(\lambda,\gamma,x,b,\eta,\nu,\eps)(t)$. Suppose that the initial
bound condition \eqref{eq:5.2 ini boots assump} holds at $t=0$.
If the bootstrap assumption \eqref{eq:BootstapAssumption} holds for
$\forall$$t\in[0,t^{*}]$ for some $t^{*}>0$, then the following holds
for $\forall t\in[0,t^{*}]$: 
\begin{align}
\begin{split} & 0<b<1.2b_{0}<0.6b^{*},\quad\frac{1}{1.2}\leq\frac{b^{2}}{\lambda^{3}}\leq1.2, \quad  \|\widehat{\eps}\|_{L^{2}}\leq\frac{1}{20}\delta_{dec}.
\end{split}
\label{eq:5.2 Bootstrap prop goal}
\end{align}
\end{prop}

We already showed in Lemma~\ref{lem:initial bound} that the decomposition
$(\lambda,\gamma,x,b,\eta,\nu,\eps)(0)$ for $v_{0}\in\mathcal{O}_{init}$
satisfies \eqref{eq:5.2 ini boots assump}, so we can apply Proposition~\ref{PropositionMainBootstrap}.
We note that mass conservation, the smallness of parameters from \eqref{eq:BootstapAssumption},
and $\|\widehat{\eps}_{0}\|_{L^{2}}\leq\frac{1}{40}\delta_{dec}$
also imply that $v(T_{exit})\in\calO_{dec}$. As a result, there are
only three possible scenarios at the exit time $t=T_{exit}$: 
\begin{align*}
b=0,\quad|\eta|=\lambda^{\frac{3}{2}(1+\kappa)},\quad\text{or}\quad|\nu|=\lambda^{\frac{3}{2}(1-\kappa)}.
\end{align*}
If $b=0$ at the exit time, then we also have $\eta=\nu=0$ so $v(T_{exit})$
is a rescaled $Q$, implying that $v(t)=[Q]_{\la_{0},\ga_{0},x_{0}}$
is a static solution. In particular, we have $v_{0}\notin\mathcal{O}_{init}$
so we have a contradiction. Hence, we have either 
\begin{align}
|\eta|=\lambda^{\frac{3}{2}(1+\kappa)}\quad\text{or}\quad|\nu|=\lambda^{\frac{3}{2}(1-\kappa)}\quad\text{at }t=T_{exit}.\label{eq:TexitBootstrapfail}
\end{align}

To bring out the contradiction, we will use a topological argument
in two dimensions based on the Brouwer fixed point theorem. This argument
is a higher-dimensional version of the continuity argument and was first
used in \cite{CoteMartelMerle2011TopologicalExample}. Let 
\begin{align}
\widetilde{T}_{exit}\coloneqq\inf\{0<t<T_{exit}:\left|\tfrac{10\widetilde{\eta}}{\lambda^{\frac{3}{2}(1+\kappa)}}\right|>0.9\text{ or }\left|\tfrac{10\widetilde{\nu}}{\lambda^{\frac{3}{2}(1-\kappa)}}\right|>0.9\},\label{eq:DefinitionTexitTilde}
\end{align}
where $\widetilde{\eta}$ and $\widetilde{\nu}$ are the refined modulation
parameters defined by \eqref{eq:RefinedModulationDefinition} in Section~\ref{SubsectionModulationEstimate}.
These parameters are close to the original modulation parameters $\eta$
and $\nu$ but they enjoy better dynamical controls compared to $\eta$
and $\nu$. By \eqref{eq:TexitBootstrapfail} and the proximity estimate
for the refined modulation parameters \eqref{eq:RefinedModulationDiffer},
we see that $0\leq\widetilde{T}_{exit}<T_{exit}$ exists. Therefore,
for each $(\eta_{0},\nu_{0})$ with $|\eta_{0}|\leq0.1\lambda_{0}^{\frac{3}{2}(1+\kappa)}$
and $|\nu_{0}|\leq0.1\lambda_{0}^{\frac{3}{2}(1-\kappa)}$, since \eqref{eq:BootstapAssumption} is also satisfied, we can
define a map $\Psi:[-1,1]^{2}\to[-1.1,1.1]^{2}\setminus[-0.9,0.9]^{2}$
by 
\begin{align}
\Psi\Big(\frac{10\eta_{0}}{\lmb_{0}^{\frac{3}{2}(1+\kappa)}},\frac{10\nu_{0}}{\lmb_{0}^{\frac{3}{2}(1-\kappa)}}\Big)=\bigg(\Big(\frac{10\td{\eta}}{\lmb^{\frac{3}{2}(1+\kappa)}}\Big)(\widetilde{T}_{exit}),\Big(\frac{10\td{\nu}}{\lmb^{\frac{3}{2}(1-\kappa)}}\Big)(\widetilde{T}_{exit})\bigg).\label{eq:PsiDefinitionTopology}
\end{align}
In summary, when $\lambda_{0},\gamma_{0},x_{0},b_{0},\tilde{\eps}_{0}$
are fixed, the maps $(\eta_{0},\nu_{0})\mapsto v_{0}\mapsto v(\widetilde{T}_{exit})\mapsto(\widetilde{\eta}(\widetilde{T}_{exit}),\widetilde{\nu}(\widetilde{T}_{exit}))$
are well defined.

We will show that $\Psi$ is continuous and is almost an identity
on the boundary $\rd[-1,1]^{2}$. Then, we will be able to apply the
Brouwer fixed point theorem to reach a contradiction. In the following
lemmas, for their use in Section~\ref{SubsectionChiralBlowup} also,
we prove slightly more general statements. Denote by $\mathcal{I}\subset\mathbb{R}_{+}\times\mathbb{R}/2\pi\mathbb{Z}\times\mathbb{R}^{4}$
the set of $(\lambda_{0},\gamma_{0},x_{0},b_{0},\eta_{0},\nu_{0},\eps_{0})$
satisfying 
\begin{align}
0<b_{0}<\tfrac{1}{2}b^{*},\quad|\eta_{0}|<0.11\lambda_{0}^{\frac{3}{2}(1+\kappa)},\quad|\nu_{0}|<0.11\lambda_{0}^{\frac{3}{2}(1-\kappa)},\quad\tfrac{1}{1.1}<\tfrac{b_{0}^{2}}{\lambda_{0}^{3}}<1.1.\label{eq:5.2. def set I}
\end{align}
Compared to the initial data condition \eqref{eq:InitialConditionBootstrap},
note that the condition \eqref{eq:5.2. def set I} requires stronger
bounds for $\eta_{0}$ and $\nu_{0}$. Proceeding as in the previous
paragraph, one can define another map $\Phi:\mathcal{I}\to[-1.2,1.2]^{2}\setminus[-0.9,0.9]^{2}$
by\footnote{The $(\eta_{0},\nu_{0})$-component of the domain of $\Phi$ is slightly
larger than that of $\Psi$ for a later purpose in Section \ref{SubsectionChiralBlowup}.} 
\begin{align}
\Phi(\lambda_{0},\gamma_{0},x_{0},b_{0},\eta_{0},\nu_{0},\eps_{0})\coloneqq((\tfrac{10\widetilde{\eta}}{\lambda^{\frac{3}{2}(1+\kappa)}})(\widetilde{T}_{exit}),(\tfrac{10\widetilde{\nu}}{\lambda^{\frac{3}{2}(1-\kappa)}})(\widetilde{T}_{exit})).\label{eq:5.2 def Phi}
\end{align}
We will use $\Phi$ defined on $\mathcal{I}$ in Section~\ref{SubsectionChiralBlowup}. However, we here use $\Phi$ restricted to a smaller set $\td{\mathcal{I}}\subset \mathcal{I}$, which we define as the set of tuples 
$(\lambda_{0},\gamma_{0},x_{0},b_{0},\eta_{0},\nu_{0},\eps_{0})$ 
satisfying
\begin{align*}
    0<b_{0}<\tfrac{1}{2}b^{*},\quad|\eta_{0}|<0.1\lambda_{0}^{\frac{3}{2}(1+\kappa)},\quad|\nu_{0}|<0.1\lambda_{0}^{\frac{3}{2}(1-\kappa)},\quad\tfrac{1}{1.1}<\tfrac{b_{0}^{2}}{\lambda_{0}^{3}}<1.1
\end{align*}
With this definition, one can check for $(\lmb_{0},\gmm_{0},x_{0},b_{0},\eta_{0},\nu_{0},\eps_{0})\in\td{\mathcal{I}}$
that 
\begin{align}
\Psi\Big(\frac{10\eta_{0}}{\lmb_{0}^{\frac{3}{2}(1+\kappa)}},\frac{10\nu_{0}}{\lmb_{0}^{\frac{3}{2}(1-\kappa)}}\Big)=\Phi(\lambda_{0},\gamma_{0},x_{0},b_{0},\eta_{0},\nu_{0},\eps_{0}).\label{eq:5.2 Psi Phi relation}
\end{align}
\begin{lem}
\label{LemmaTopologicalPhiContinuous} The map $\Phi:\mathcal{I}\to[-1.2,1.2]^{2}\setminus[-0.9,0.9]^{2}$
is continuous. In particular, $\Psi:[-1,1]^{2}\to[-1.1,1.1]^{2}\setminus[-0.9,0.9]^{2}$
is continuous. 
\end{lem}

Lemma~\ref{LemmaTopologicalPhiContinuous} will be a consequence
of the refined modulation estimates \eqref{eq:RefinedModulationBlowupLaw},
which are discussed in Section~\ref{SubsectionModulationEstimate}.
Next, we will use a variant of the Brouwer fixed point theorem. 
\begin{lem}
\label{LemmaPsiNonExistenceBrouwer} There exists no function $\Psi$
that satisfies the following: 
\begin{enumerate}
\item $\Psi:[-1,1]^{2}\to[-1.2,1.2]^{2}\setminus[-0.9,0.9]^{2}$. \label{StateLemmaPsiNonExistenceBrouwer1} 
\item Let $S$ be a square, $S=\{(x,y)\in\mathbb{R}^{2}:|x+y|+|x-y|=2\}$.
For any $z\in S$, $\Psi(z)\in B_{0.2}(z)$. \label{StateLemmaPsiNonExistenceBrouwer2} 
\item $\Psi$ is continuous. \label{StateLemmaPsiNonExistenceBrouwer3} 
\end{enumerate}
\end{lem}

We postpone the proof of these two lemmas to Section~\ref{SubsectionBootstrapDescription}.

By the definition \eqref{eq:PsiDefinitionTopology} of $\Psi$ and
\eqref{eq:RefinedModulationDiffer}, our $\Psi$ satisfies (\ref{StateLemmaPsiNonExistenceBrouwer1})
and (\ref{StateLemmaPsiNonExistenceBrouwer2}) in Lemma~\ref{LemmaPsiNonExistenceBrouwer}
for $b^{*}\ll1$. Lemma~\ref{LemmaTopologicalPhiContinuous} says
that $\Psi$ also satisfies (\ref{StateLemmaPsiNonExistenceBrouwer3}).
Therefore, the assumptions of Lemma~\ref{LemmaPsiNonExistenceBrouwer}
are satisfied with our $\Psi$ and we get a contradiction. As we began
our contradiction argument from $T_{exit}<T$, we conclude that $T_{exit}=T$
for some $(\eta_{0},\nu_{0})$ and we obtain a trapped solution $v$
with this $(\eta_{0},\nu_{0})$.

The rest of statement of Theorem~\ref{thm:precise statement codimensionone blowup}
is summarized in the following proposition, which will be proved in
Section~\ref{SubsectionBootstrapDescription} after we complete the
proof of the bootstrap proposition.
\begin{prop}[Sharp description of blow-up solutions]
\label{PropSharpDescription} Let $v$ be a trapped solution. Then,
it blows up in finite time as described in Theorem~\ref{TheoremCodimension2Blowup}. 
\end{prop}

So far, we have shown that for each $(\lambda_{0},\gamma_{0},x_{0},b_{0},\widetilde{\eps}_{0})\in\widetilde{\mathcal{U}}_{init}\subset\mathbb{R}_{+}\times\mathbb{R}/2\pi\mathbb{Z}\times\mathbb{R}^{2}\times\mathcal{Z}^{\perp}$,
there exists $(\eta_{0},\nu_{0})\in(-\frac{1}{2}b_{0}^{1+\kappa},\frac{1}{2}b_{0}^{1+\kappa})\times(-\frac{1}{2}b_{0}^{1-\kappa},\frac{1}{2}b_{0}^{1-\kappa})$
so that the solution $v(t)$ starting from the initial data $v_{0}$
is a finite-time blow-up solution as described in Theorem~\ref{TheoremCodimension2Blowup}.
In view of this, we have constructed a codimension two set of initial
data that leads finite-time blow-up. For each such $v(t)$, we can
add an one parameter family of blow-up solutions $\{v^{c}(t)\}_{c\in\R}$
using Galilean invariance. This finishes the proof of Theorem~\ref{thm:precise statement codimensionone blowup}.
\end{proof}
\begin{quote}
\emph{In the rest of the paper except Section}~\emph{\ref{SubsectionChiralBlowup},
we assume the bootstrap assumption \eqref{eq:BootstapAssumption}.}
\end{quote}

\subsection{\label{SubsectionInterpolationEstimates}Modification of the profile
for $w_{1}$ and consequences of coercivity}

In this subsection, we handle several technical preparations for modulation
analysis. Under the bootstrap assumption \eqref{eq:BootstapAssumption},
before entering the modulation decomposition of nonlinear solution
$v(t,x)$, we will take further modification of the profile for technical reasons. Moreover, we will find interpolation-type bounds
of solutions that mostly follow from the coercivity estimates in Section~\ref{SubsectionConjugationAQ}
and the bootstrap assumption \eqref{eq:BootstapAssumption}.

First, we explain why the profile $\widehat{P}_{1}$ is not sufficient
and how we further modify it. Originally, we have decomposed $w_{1}$
as 
\begin{align*}
w_{1}=\widehat{P}_{1}+\widehat{\eps}_{1}=-(ib+\eta)\tfrac{y}{2}Q+i\tfrac{\nu}{2}Q+\widehat{\eps}_{1}.
\end{align*}
Note that $iyQ,yQ,iQ$ belong to the kernel of $A_{Q}$. However,
as $A_{Q}$ is $\bbC$-linear, another kernel direction $Q$ is missing
in $\wt P_{1}$. This fact generates an issue related to coercivity
estimates for $A_{Q}\wt{\eps}_{1}$ and indeed prevents us to obtain
the scaling invariant bound $\|\widehat{\eps}_{1}\|_{\dot{\mathcal{H}}^{1}}\lesssim\lambda^{2}$
for $\widehat{\eps}_{1}$ in the energy estimate. To resolve this
issue, we will identify the main $Q$-directional component of $\widehat{\eps}_{1}$
and put it into $P_{1}$ as a further correction profile.

To identify the main $Q$-component of $\wt{\eps}_{1}$, we begin
with $w_{1}=\bfD_{w}w$, which implies $\wt P_{1}+\wt{\eps}_{1}=L_{Q}\widehat{\eps}+N_{Q}(\widehat{\eps})$.
As $Q$ is even, we take the even part to observe 
\begin{align*}
\widehat{\eps}_{1,e} & =L_{Q}(\eps_{o})+(N_{Q}(\widehat{\eps}))_{e}\\
 & =\mathbf{D}_{Q}(\eps_{o})+Q\mathcal{H}\Re(Q\eps_{o})+(N_{Q}(\widehat{\eps}))_{e}.
\end{align*}
We also remark that $Q\eps_{o}\in L^{2}$ since $QP_{o}=\tfrac{1}{2}i\nu yQ^{2}\in L^{2}$,
and then $\mathcal{H}\Re(Q\eps_{o})$ is well defined. Similarly to
the proof of \eqref{eq:LemmaHilbertDecayExchangeInAQ}, using $(1+y^{2})Q^{2}=2$
and \eqref{eq:CommuteHilbert}, we have 
\begin{align*}
Q\mathcal{H}\Re(Q\eps_{o}) & =\tfrac{1}{2}Q\mathcal{H}(1+y^{2})\Re(Q^{3}\eps_{o})\\
 & =\tfrac{1}{2}Q\mathcal{H}\Re(Q^{3}\eps_{o})+\tfrac{1}{2}yQ\mathcal{H}y\Re(Q^{3}\eps_{o})-\tfrac{1}{2\pi}Q{\textstyle \int_{\R}\Re(yQ^{3}\eps_{o})dy}\\
 & =Q^{-1}\mathcal{H}\Re(Q^{3}\eps_{o})-\tfrac{1}{2\pi}Q{\textstyle \int_{\R}\Re(yQ^{3}\eps_{o})dy.}
\end{align*}
We note that, in the second equality, we used \eqref{eq:CommuteHilbert}
with $\int Q^{3}\eps_{o}dy=0$ since \eqref{eq:CommuteHilbert} is
well defined for $f=Q^{3}\eps_{o}$ and the integration of odd functions
is $0$. In addition, a similar procedure applies to the nonlinear part $(N_{Q}(\widehat{\eps}))_{e}$, obtaining
\begin{equation}
	\begin{aligned}
		N_{Q}(\wt\eps)=&\wt\eps\mathcal{H}(\Re(Q\wt\eps))+\tfrac{1}{2}(Q+\wt\eps)\mathcal{H}(|\wt\eps|^{2})
		\\
		=&\td N_Q(\wt\eps)-\tfrac{1}{4\pi}Q {\textstyle\int_\bbR} yQ^2|\wt\eps|^2dy.
	\end{aligned} \label{eq:NQ td NQ relation}
\end{equation}
where
\begin{align}
	\td N_Q(\wt\eps)
	=\wt\eps\mathcal{H}(\Re(Q\wt\eps))+\tfrac{1}{2}\wt\eps\mathcal{H}(|\wt\eps|^{2})+\tfrac{1}{4}yQ\mathcal{H}(yQ^2|\wt\eps|^{2})
	+\tfrac{1}{4}Q\mathcal{H}(Q^2|\wt\eps|^{2}). \label{eq:def td NQ}
\end{align}
Therefore, we have 
\begin{align}
\widehat{\eps}_{1,e}= & \mathbf{D}_{Q}(\eps_{o})+Q^{-1}\mathcal{H}\Re(Q^{3}\eps_{o})-\tfrac{1}{2\pi}Q{\textstyle \int_{\R}\Re(yQ^{3}\eps_{o})dy+(N_{Q}(\widehat{\eps}))_{e},\nonumber}\\
= & \widetilde{L}_{Q}\eps_{o}
-\tfrac{1}{2\pi}Q{\textstyle \int_{\R}\Re(yQ^{3}\eps_{o})dy
	-\tfrac{1}{4\pi}Q {\textstyle\int_\bbR} yQ^2|\wt\eps|^2dy
	+(\td N_{Q}(\widehat{\eps}))_{e}\nonumber}\\
= & \widetilde{L}_{Q}\eps_{o}+\tfrac{1}{2}\mu Q+(\td N_{Q}(\widehat{\eps}))_{e},\label{eq:epsilon tilde L expression}
\end{align}
where the parameter $\mu=\mu(\eps)$ is given by

\begin{equation}
\mu\coloneqq-\frac{1}{\pi}\int_{\R}\Re(yQ^{3}\eps_{o})dy
-\frac{1}{2\pi} \int_\bbR yQ^2|\wt\eps|^2dy
.\label{eq:definition of mu}
\end{equation}
As the kernel elements of $A_{Q}\tilde L_{Q}$ are already completely
modulated (using six parameters), the $Q$-direction cannot arise
from $\tilde L_{Q}f$; we do not need to care $\tilde L_{Q}\eps_{o}$.
Therefore, we arrive at the final decomposition $w_{1}=P_{1}+\eps_{1}$
with 
\begin{align}
\eps_{1}\coloneqq\wt{\eps}{}_{1}-\tfrac{1}{2}\mu Q,\quad P_{1}=\wt P_{1}+\tfrac{1}{2}\mu Q=-(ib+\eta)\tfrac{y}{2}Q+(i\nu+\mu)\tfrac{1}{2}Q.\label{eq:DefinitionEpsilon1afterModify}
\end{align}

\begin{rem}
In fact, $\eps_{1}$ cannot be rigorously expressed using $L_{Q}\eps$,
as explained in Remark~\ref{RemarkLQ Well Definedness}. But it can
be represented using $\widetilde{L}_{Q}\eps$: 
\begin{align}
\eps_{1}=\widetilde{L}_{Q}\eps+\td N_{Q}(\widehat{\eps})-yQ\cdot\tfrac{1}{2\pi}{\textstyle \int_{\R}\Re(Q^{3}\eps_{e})dy}.\label{eq:5.3 Remark goal}
\end{align}
In this expression, we still have an odd profile in the $yQ$-direction.
Unlike $Q$, this direction is already included in the generalized
kernel of $L_{Q}$ and does not create an issue in our modulation
setting. Indeed, we can verify on the way the contribution of this
term is negligible as 
\begin{align*}
\left|\tfrac{1}{2\pi}{\textstyle \int_{\R}\Re(Q^{3}\eps_{e})dy}\right|\lesssim\lambda^{2}\sim b^{\frac{4}{3}}\ll b.
\end{align*}
In contrast, the $Q$-direction could be critical in size since $|\mu|\lesssim \lambda$
(see \eqref{eq:mu estimate by b}) and we need the correction from
$\wt{\eps}_{1}$to $\eps_{1}$.

In the following, we verify \eqref{eq:5.3 Remark goal}. We note that
$\D_{Q}P_{e}=\widehat{P}_{1,o}$. From the decomposition \eqref{eq:DefinitionEpsilon1afterModify}
and \eqref{eq:epsilon tilde L expression}, we have 
\begin{align*}
\eps_{1} & =(L_{Q}\widehat{\eps}_{e}-\widehat{P}_{1,o})+\widetilde{L}_{Q}\eps_{o}+\td N_{Q}(\widehat{\eps})\\
 & =\mathbf{D}_{Q}\eps+Q\mathcal{H}\Re(Q\widehat{\eps}_{e})+Q^{-1}\mathcal{H}\Re(Q^{3}\eps_{o})+\td N_{Q}(\widehat{\eps}).
\end{align*}
For $Q\mathcal{H}\Re(Q\widehat{\eps}_{e})$, using $(1+y^{2})Q^{2}=2$
and \eqref{eq:CommuteHilbert}, we have 
\begin{align}
2Q\mathcal{H}\Re(Q\widehat{\eps}_{e})= & Q\mathcal{H}\Re(Q^{3}\widehat{\eps}_{e})+yQ\mathcal{H}\Re(yQ^{3}\widehat{\eps}_{e})\nonumber \\
= & Q\mathcal{H}\Re(Q^{3}\eps_{e})+yQ\mathcal{H}\Re(yQ^{3}\eps_{e})+Q\mathcal{H}\Re(Q^{3}P_{e})+yQ\mathcal{H}\Re(yQ^{3}P_{e}).\label{eq:5.3 Remark 1}
\end{align}
For the profile part, we have 
\begin{align}
Q\mathcal{H}\Re(Q^{3}P_{e})+yQ\mathcal{H}\Re(yQ^{3}P_{e})=-\tfrac{\eta}{2}(Q\mathcal{H}Q^{2}+yQ\mathcal{H}yQ^{2})=0.\label{eq:5.3 Remark 2}
\end{align}
By \eqref{eq:5.3 Remark 1}, \eqref{eq:5.3 Remark 2}, and again \eqref{eq:CommuteHilbert},
we have 
\begin{align*}
Q\mathcal{H}\Re(Q\widehat{\eps}_{e})=Q^{-1}\mathcal{H}\Re(Q^{3}\eps_{e})-\tfrac{1}{2\pi}yQ{\textstyle \int_{\R}\Re(Q^{3}\eps_{e})dy.}
\end{align*}
Hence, we arrive at \eqref{eq:5.3 Remark goal}. 
\end{rem}

\begin{rem}[Kernel of $A_{Q}L_{Q}$]
The technical issue mentioned above arises from the use of nonlinear adapted
derivative $\D_{w}w=w_{1}$ instead of linear adapted derivatives
such as $L_{Q}\eps$ or $A_{Q}L_{Q}\eps$. The problematic (real)
$Q$-direction belongs to the kernel of $A_{Q}$ but $Q$ cannot be
achieved by taking $L_{Q}$ to the generalized kernels \eqref{eq:GeneralizedKernelRelation},
which we completely modulated by introducing six modulation parameters.
The presence of a nontrivial $Q$-direction of $\wt{\eps}_{1}$ is
an artifact of using nonlinear adapted derivative $w_{1}$. If we
were to use the linear one $L_{Q}\eps$ (or $\td L_{Q}\eps$), then
we would not face this issue. However, the nonlinear adapted derivative
$\D_{w}w$ has a lot of benefits; it solves a natural higher-order
equation \eqref{eq:w1-equ} and so essentially simplifies the nonlinear
analysis in the sequel.
\end{rem}

We are ready to fix the decompositions for $w$ and $w_{1}$. As mentioned
above, we will use two different decompositions for $(w,w_{1})$. 
\begin{align}
\begin{split}w & =Q+\widehat{\eps}\\
 & =Q+P+\eps,\\
w_{1} & =w_{1}\text{ itself}\\
 & =P_{1}+\eps_{1}.
\end{split}
\begin{split}(\text{on }H^{1})\\
(\text{on }\dot{\mathcal{H}}^{2})\\
(\text{on }L^{2})\\
(\text{on }\dot{\mathcal{H}}^{1})
\end{split}
\label{eq:DecomposeDependTop}
\end{align}
where $\widehat{\eps}\in L^{2}$, $\eps\in\dot{\mathcal{H}}^{2}$,
and $\eps_{1}\in\dot{\mathcal{H}}^{1}$. The decomposition that we use
depends on the norms we measure. For example, since $P,P_{1}\notin L^{2}$,
we cannot use $P+\eps$ or $P_{1}+\eps_{1}$ decomposition when measured
in $L^{2}$-norm. Since $P\in\dot{\mathcal{H}}^{2}$ and $P_{1}\in\dot{\mathcal{H}}^{1}$,
we use suitable decompositions when the solutions are posed in $\dot{\mathcal{H}^{1}}$
or $\dot{\mathcal{H}}^{2}$. Their refined decompositions are valid
only in these spaces allowing spatial growths. Along with this strategy,
we are able to \emph{avoid cut-offs} on the profile. Because the profiles
$P$ and $P_{1}$ lack decay and even grow in space, if we were to
use cut-offs on the profiles, it will create serious cut-off errors
in the dynamical analysis of $\eps$ and $\eps_{1}$. In addition,
we can use a number of explicit identities related to $\mathcal{H}$ and
$Q$ because we did not introduce cut-offs. We believe that this strategy
is robust and applicable to many related cases when the blow-up profile
has bad polynomial decay.

To perform the nonlinear analysis, we measure solutions in various
norms. In the following, we prove some relevant estimates, which basically originate from \eqref{eq:energy scale bounds} 
and the coercivity estimates (Proposition~\ref{PropCoercivityLQ}
and Proposition~\ref{PropCoercivityAQLQ}) by the bootstrap assumption \eqref{eq:BootstapAssumption} and interpolation arguments. 

\begin{lem}[Coercivity estimates]
\label{LemmaNonlinearCoercivity} The following estimates hold: 
\begin{enumerate}
\item We have 
\begin{align}
\|\widehat{\eps}\|_{\dot{\mathcal{H}}^{1}}\sim\|w_{1}\|_{L^{2}}\lesssim\lambda,\label{eq:CoercivityEstimateH1}
\end{align}
from which we deduce 
\begin{align}
\|\widehat{\eps}\|_{L^{\infty}}\lesssim\lambda^{\frac{1}{2}}.\label{eq:InterpolationEstimateW,inf}
\end{align}
\item We have 
\begin{align}
	\|A_Qw_1\|_{L^2}\lesssim \lambda^2, \label{eq: AQw1 esti}
\end{align}
and
\begin{align}
\|\eps_{1}\|_{\dot{\mathcal{H}}^{1}}+\|\eps\|_{\dot{\mathcal{H}}^{2}} & \lesssim\lambda^{2}.\label{eq:CoercivityEstimateH2epsilon}
\end{align}
\item We have $\eps_{1,e}\in L^{2}$ and $\eps_{o}\in\dot{H}^{1}$ with
\begin{align}
\|\eps_{o}\|_{\dot\calH^1}+\|\eps_{1,e}\|_{L^{2}}\lesssim\lambda.\label{eq:5.3 interpolation additional decay}
\end{align}
\item We also have the interpolation type estimates, 
\begin{align}
\begin{split} & \|Q\widehat{\eps}_{e}\|_{L^{\infty}}+\|\partial_{y}\widehat{\eps}_{e}\|_{L^{\infty}}\lesssim\lambda^{\frac{3}{2}},\\
 & \|Q\widehat{\eps}_{o}\|_{L^{\infty}}+\|\partial_{y}\widehat{\eps}_{o}\|_{L^{\infty}}\lesssim|\nu|\lesssim\lambda^{\frac{3}{2}(1-\kappa)}.
\end{split}
\label{eq:InterpolationEstimateW1,inf}
\end{align}
\item In addition, we have interpolation estimates which do not depend on
parity, 
\begin{align}
\|\eps_{1}\|_{L^{\infty}}+\|Q\eps\|_{L^{\infty}}+\|\partial_{y}\eps\|_{L^{\infty}}\lesssim\lambda^{\frac{3}{2}},\label{eq:InterpolationEstimateW1,inf +}
\end{align}
\item We estimate $\mu$ defined in \eqref{eq:definition of mu} 
\begin{equation}
|\mu|\lesssim \lambda.\label{eq:mu estimate by b}
\end{equation}
\end{enumerate}
\end{lem}

\begin{proof}
	We first remark that from \eqref{eq:energy scale bounds}, we have
	\begin{align}
		\|w_1\|_{L^2}\lesssim \lambda,\quad \|\td \bfD_{w} w_1\|_{L^2}\lesssim \lambda^2. \label{eq:energy bound coerc}
	\end{align}
We use the coercivity $\|L_{Q}\widehat{\eps}\|_{L^{2}}\sim\|\widehat{\eps}\|_{\dot{\mathcal{H}}^{1}}$ (Proposition~\ref{PropCoercivityLQ})
and $\|\eps\|_{\dot{\mathcal{H}}^{2}}\sim\|A_{Q}\widetilde{L}_{Q}\eps\|_{L^{2}}$ (Proposition~\ref{PropCoercivityAQLQ})
with \eqref{eq:AQLQ epsilon hat AQtildeLQ epsilon} and $A_{Q}L_{Q}\widehat{\eps}=A_{Q}\widetilde{L}_{Q}\eps$
without mentioning. However, we will not use the coercivity for $\eps_{1}$
and $A_{Q}$. 

Recall that 
\begin{align*}
w_{1}=\mathbf{D}_{w}w & =\mathbf{D}_{Q}Q+L_{Q}\widehat{\eps}+N_{Q}(\widehat{\eps})\\
 & =L_{Q}\widehat{\eps}+N_{Q}(\widehat{\eps})
\end{align*}
where 
\begin{align*}
N_{Q}(\widehat{\eps})=\widehat{\eps}\mathcal{H}(\Re(Q\widehat{\eps}))+\tfrac{1}{2}(Q+\widehat{\eps})\mathcal{H}(|\widehat{\eps}|^{2}).
\end{align*}

\textbf{Proof of \eqref{eq:CoercivityEstimateH1} and \eqref{eq:InterpolationEstimateW,inf}.}
We estimate 
\begin{align}
    \|w_{1}-L_{Q}\widehat{\eps}\|_{L^{2}}=\|N_{Q}(\widehat{\eps})\|_{L^{2}} & \lesssim
    \|\widehat{\eps}\|_{L^{\infty}}^{2}\|\widehat{\eps}\|_{L^{2}}
    +\|Q\calH(|\wt\eps|^2)\|_{L^2}\nonumber \\
     & \lesssim
     \|\widehat{\eps}\|_{\dot{\mathcal{H}}^{1}}\|\widehat{\eps}\|_{L^{2}}^{2}
     +\|Q\calH(|\wt\eps|^2)\|_{L^2}.\label{eq:N_Q wt epsilon}
\end{align}
For the last term of \eqref{eq:N_Q wt epsilon}, we apply \eqref{eq:CommuteHilbert} with $(1+y^2)\langle y\rangle^{-2}=1$. This yields
\begin{align*}
    Q\calH(|\wt\eps|^2)
    =Q\calH(\tfrac{1+y^2}{\langle y\rangle^{2}}|\wt\eps|^2)
    =&Q\calH(\tfrac{1}{\langle y\rangle^{2}}|\wt\eps|^2)+Q\calH(\tfrac{y^2}{\langle y\rangle^{2}}|\wt\eps|^2)
    \\
    =&Q\calH(\tfrac{1}{\langle y\rangle^{2}}|\wt\eps|^2)+yQ\calH(\tfrac{y}{\langle y\rangle^{2}}|\wt\eps|^2)-\tfrac{1}{\pi}Q {\textstyle\int}\tfrac{y}{\langle y\rangle^{2}}|\wt\eps|^2.
\end{align*}
Since $\|\tfrac{1}{\langle y\rangle^{2}}\wt\eps\|_{L^2} \leq \|\wt\eps\|_{L^2}$, we obtain
\begin{align*}
    \|Q\calH(|\wt\eps|^2)\|_{L^2}
    \lesssim& \|\tfrac{1}{\langle y\rangle^{2}}\wt\eps\|_{L^2}^2
    +\|\widehat{\eps}\|_{L^{\infty}}\|\tfrac{1}{\langle y\rangle^{2}}\wt\eps\|_{L^2}
    +\|\tfrac{1}{\langle y\rangle^{2}}\wt\eps\|_{L^2}\|\wt\eps\|_{L^2}
    \\
    \lesssim& \|\widehat{\eps}\|_{\dot{\mathcal{H}}^{1}}\|\widehat{\eps}\|_{L^{2}}
    +\|\wt\eps\|_{\dot H^1}^{\frac12}\|\widehat{\eps}\|_{L^2}^{\frac12}\|\tfrac{1}{\langle y\rangle^{2}}\wt\eps\|_{L^2}
    +\|\widehat{\eps}\|_{\dot{\mathcal{H}}^{1}}\|\widehat{\eps}\|_{L^{2}}
    \\
    \lesssim& \|\widehat{\eps}\|_{\dot{\mathcal{H}}^{1}}\|\widehat{\eps}\|_{L^{2}}.
\end{align*}
Consequently,
\begin{align*}
    \|w_{1}-L_{Q}\widehat{\eps}\|_{L^{2}}
    \lesssim \|\widehat{\eps}\|_{\dot{\mathcal{H}}^{1}}(\|\widehat{\eps}\|_{L^{2}}^{2}+\|\widehat{\eps}\|_{L^{2}}),
\end{align*}
which implies
\begin{align*}
    (1-\|\widehat{\eps}\|_{L^{2}})\|\widehat{\eps}\|_{\dot{\mathcal{H}}^{1}}
    \lesssim\|w_{1}\|_{L^{2}}
    \lesssim(1+\|\widehat{\eps}\|_{L^{2}})\|\widehat{\eps}\|_{\dot{\mathcal{H}}^{1}}.
\end{align*}
Since $\|\widehat{\eps}\|_{L^2}\ll1$ by \eqref{eq:BootstapAssumption}, together with $\|L_{Q}\widehat{\eps}\|_{L^{2}}\sim\|\widehat{\eps}\|_{\dot{\mathcal{H}}^{1}}$ and \eqref{eq:energy bound coerc}, we obtain \eqref{eq:CoercivityEstimateH1}. Interpolating this
and $\|\widehat{\eps}\|_{L^{2}}\ll1$, we deduce \eqref{eq:InterpolationEstimateW,inf}.

\textbf{Proof of \eqref{eq: AQw1 esti}.}
Next, we show \eqref{eq: AQw1 esti}. Thanks to \eqref{eq:DQ BQ decomp} and \eqref{eq:BQBQstar equal I}, we have
\begin{align}
	A_Qw_1=B_QB_Q^*A_Qw_1=B_Q \td \bfD_{Q} w_1. \label{eq:AQ bfDQ relation}
\end{align}
We claim
\begin{align}
	\|B_Q(\td \bfD_{w} w_1-\td \bfD_{Q} w_1)\|_{L^2}\lesssim \lambda^2. \label{eq:coerc H2 claim}
\end{align}
We note that, by the definition of $B_Q$, we derive
\begin{align}
	\|B_Qf\|_{L^2}\lesssim \|f\|_{L^2}. \label{eq: BQ L2 estimate}
\end{align}
Thanks to the definition of $\td \bfD_w$, we have
\begin{align*}
	2(\td \bfD_{w} w_1-\td \bfD_{Q} w_1)
	=&w\calH(\ol w w_1)-Q\calH(Qw_1)
	\\
	=&Q\calH(\ol{\wt\eps} w_1)+\wt\eps\calH(Qw_1)+\wt\eps\calH(\ol{\wt\eps} w_1).
\end{align*}
We first estimate the $\wt\eps\calH(\ol{\wt\eps} w_1)$ part. By \eqref{eq: BQ L2 estimate}, the isometry property of $\mathcal{H}$, \eqref{eq:CoercivityEstimateH1}, and \eqref{eq:InterpolationEstimateW,inf}, we have
\begin{align*}
	\|B_Q[\wt\eps\calH(\ol{\wt\eps} w_1)]\|_{L^2}
	\lesssim \|\wt\eps\calH(\ol{\wt\eps} w_1)\|_{L^2}
	\leq \|\wt\eps\|_{L^\infty}^2\|w_1\|_{L^2}\lesssim \lambda^2.
\end{align*}
For the second term $\wt\eps\calH(Qw_1)$, using \eqref{eq:CommuteHilbert} with $(1+y^2)\langle y\rangle^{-2}=1$, we derive
\begin{align*}
	\wt\eps\calH(Qw_1)=\tfrac{1}{\langle y\rangle^{2}}\wt\eps\calH(Qw_1)+\tfrac{y}{\langle y\rangle^{2}}\wt\eps\calH(yQw_1)+\tfrac{y}{\langle y\rangle^{2}}\wt\eps\cdot \tfrac{1}{\pi}\textstyle{\int_{\bbR}} Qw_1 dy.
\end{align*}
Thus, we have
\begin{align*}
	\|B_Q[\wt\eps\calH(Qw_1)]\|_{L^2}
	&\lesssim \|\wt\eps\calH(Qw_1)\|_{L^2}
	\\
	&\lesssim \|Q\wt\eps\|_{L^\infty}(\|\calH(Qw_1)\|_{L^2}+\|\calH(yQw_1)\|_{L^2})+\|\wt\eps\|_{\dot\calH^1}|\textstyle{\int_{\bbR}} Qw_1 dy|
	\\
	&\lesssim \|\wt\eps\|_{\dot\calH^1}\|w_1\|_{L^2}
	\\
	&\lesssim \lambda^2.
\end{align*}
We note that $\|Q\wt\eps\|_{L^\infty}\lesssim \|Q\wt\eps\|_{L^2}^{\frac12}\|Q\wt\eps\|_{\dot H^1}^{\frac12}\lesssim \|\wt\eps\|_{\dot\calH^1}$.
We control the first term $Q\calH(\ol{\wt\eps} w_1)$. We use $B_QQ=0$ and \eqref{eq:CommuteHilbert} with $(1+y^2)\langle y\rangle^{-2}=1$. We deduce
\begin{align*}
	B_Q[Q\calH(\ol{\wt\eps} w_1)]
    =&B_Q[Q\calH(\tfrac{1}{\langle y\rangle^2}\ol{\wt\eps} w_1)]
	+B_Q[Q\calH(\tfrac{y^2}{\langle y\rangle^2}\ol{\wt\eps} w_1)]
    \\
    =&B_Q[Q\calH(\tfrac{1}{\langle y\rangle^2}\ol{\wt\eps} w_1)]
	+B_Q[yQ\calH(\tfrac{y}{\langle y\rangle^2}\ol{\wt\eps} w_1)
    -\tfrac{1}{\pi}Q {\textstyle \int}\tfrac{y}{\langle y\rangle^2}\ol{\wt\eps} w_1]
    \\
    =&B_Q[Q\calH(\tfrac{1}{\langle y\rangle^2}\ol{\wt\eps} w_1)]
	+B_Q[yQ\calH(\tfrac{y}{\langle y\rangle^2}\ol{\wt\eps} w_1)].
\end{align*}
Thus, we have
\begin{align*}
	\|B_Q[Q\calH(\ol{\wt\eps} w_1)]\|_{L^2}
	\lesssim \|Q\wt\eps\|_{L^\infty}\|w_1\|_{L^2}
	\lesssim \|\wt\eps\|_{\dot\calH^1}\|w_1\|_{L^2}
	\lesssim \lambda^2.
\end{align*}
Therefore, we conclude our claim \eqref{eq:coerc H2 claim}. By \eqref{eq:energy bound coerc}, together with \eqref{eq:AQ bfDQ relation}--\eqref{eq: BQ L2 estimate}, we deduce \eqref{eq: AQw1 esti}:
\begin{align*}
	\|A_Qw_1\|_{L^2}\leq \|B_Q(\td \bfD_{w} w_1-\td \bfD_{Q} w_1)\|_{L^2}+\|B_Q\td \bfD_{w} w_1\|_{L^2}\lesssim \lambda^2.
\end{align*} 

\textbf{Proof of \eqref{eq:CoercivityEstimateH2epsilon}.}
Now, we prove \eqref{eq:CoercivityEstimateH2epsilon}. Thanks to \eqref{eq:5.3 Remark goal} and \eqref{eq:kernel AQ}, we have
\begin{align*}
A_{Q}\eps_{1}=A_{Q}\widetilde{L}_{Q}\eps+A_{Q}\td N_{Q}(\widehat{\eps}),
\end{align*}
where $\widetilde{L}_{Q}$ and $\td N_Q(\wt\eps)$ are given by \eqref{eq:tildeLQ definition} and \eqref{eq:def td NQ} respectively. This yields
\begin{align*}
\|A_{Q}\eps_{1}-A_{Q}\widetilde{L}_{Q}\eps\|_{L^{2}}=\|A_{Q} N_{Q}(\widehat{\eps})\|_{L^{2}}
=\|A_{Q} \td N_{Q}(\widehat{\eps})\|_{L^{2}},
\end{align*}
We decompose $\td N_Q$ into two parts,
\begin{align*}
	\td N_Q(\wt\eps)=\td N_1+ \td N_2,
\end{align*}
where
\begin{align*}
	\td N_1&=\wt\eps\mathcal{H}(\Re(Q\wt\eps))+\tfrac{1}{2}\wt\eps\mathcal{H}(|\wt\eps|^{2}),
	\\
	\td N_2&=\tfrac{1}{4}yQ\mathcal{H}(yQ^2|\wt\eps|^{2})
	+\tfrac{1}{4}Q\mathcal{H}(Q^2|\wt\eps|^{2}).
\end{align*}
For $\td N_1$, thanks to Lemma~\ref{LemmaAppendix AQ Subcoercivity}, we have 
\begin{equation}
    \begin{aligned}
        \|A_{Q}\td N_1\|_{L^{2}}\lesssim\|\td N_1\|_{\dot{\mathcal{H}}^{1}}
    \lesssim & \|\partial_{y}(\widehat{\eps}\mathcal{H}(\Re(Q\widehat{\eps})))\|_{L^{2}}+\|\partial_{y}(\widehat{\eps}\mathcal{H}(|\widehat{\eps}|^{2}))\|_{L^{2}}\\
     & +\|Q\widehat{\eps}\mathcal{H}(\Re(Q\widehat{\eps}))\|_{L^{2}}+\|Q\widehat{\eps}\mathcal{H}(|\widehat{\eps}|^{2})\|_{L^{2}}.
    \end{aligned}\label{eq:td N1 aux}
\end{equation}
Using $\partial_{y}\Re(QP)=0$ and Sobolev inequality $\|f\|_{L^{\infty}}\lesssim\|f\|_{L^{2}}^{\frac{1}{2}}\|\partial_{y}f\|_{L^{2}}^{\frac{1}{2}}$,
we estimate 
\begin{align*}
\|\partial_{y}(\widehat{\eps}\mathcal{H}(\Re(Q\widehat{\eps})))\|_{L^{2}} & \leq\|\partial_{y}\widehat{\eps}\|_{L^{2}}\|\mathcal{H}(Q\widehat{\eps})\|_{L^{\infty}}+\|\widehat{\eps}\|_{L^{\infty}}\|\partial_{y}\Re(Q\widehat{\eps})\|_{L^{2}}\\
 & \lesssim\|\partial_{y}\widehat{\eps}\|_{L^{2}}\|Q\widehat{\eps}\|_{L^{2}}^{\frac{1}{2}}\|\partial_{y}(Q\widehat{\eps})\|_{L^{2}}^{\frac{1}{2}}+\|\widehat{\eps}\|_{L^{\infty}}\|\partial_{y}(Q\eps)\|_{L^{2}}\\
 & \lesssim\|\widehat{\eps}\|_{\dot{\mathcal{H}}^{1}}^{2}+\|\widehat{\eps}\|_{L^{\infty}}\|\eps\|_{\dot{\mathcal{H}}^{2}}\\
 & \lesssim\lambda^{2}+\lambda^{\frac{1}{2}}\|A_{Q}\widetilde{L}_{Q}\eps\|_{L^{2}}.
\end{align*}
In the third line of the preceding inequalities, we used $\|\partial_{y}(Q\widehat{\eps})\|_{L^{2}}\lesssim\|Q\|_{L^{\infty}}\|\partial_{y}\widehat{\eps}\|_{L^{2}}$.
By a similar argument, we have $\|Q\widehat{\eps}\mathcal{H}(\Re(Q\widehat{\eps}))\|_{L^{2}}\lesssim\|\widehat{\eps}\|_{\dot{\mathcal{H}}^{1}}^{2}\lesssim\lambda^{2}$.
Moreover, 
\begin{align*}
\|\partial_{y}(\widehat{\eps}\mathcal{H}(|\widehat{\eps}|^{2}))\|_{L^{2}}\leq\|\partial_{y}\widehat{\eps}\|_{L^{2}}\||\widehat{\eps}|^{2}\|_{L^{\infty}}+\|\widehat{\eps}\|_{L^{\infty}}\|\partial_{y}(|\widehat{\eps}|^{2})\|_{L^{2}}\lesssim\lambda^{2}.
\end{align*}
and we also have $\|Q\widehat{\eps}\mathcal{H}(|\widehat{\eps}|^{2})\|_{L^{2}}\lesssim\lambda^{2}$
following a similar argument. 
Thus, we arrive at
\begin{align}
	\|A_{Q}\td N_1\|_{L^{2}}\lesssim \lambda^{2}+\lambda^{\frac{1}{2}}\|A_{Q}\widetilde{L}_{Q}\eps\|_{L^{2}}. \label{eq:td N1 goal}
\end{align}

For $\td N_2$, using \eqref{eq:CommuteHilbert} with $A_Q(yQ)=0$, we derive
\begin{align*}
	A_Q\td N_2=\tfrac{1}{4}A_Q[y^2Q\mathcal{H}(Q^2|\wt\eps|^{2})]
	+\tfrac{1}{4}A_Q[Q\mathcal{H}(Q^2|\wt\eps|^{2})].
\end{align*}
Since the second term, $\tfrac{1}{4}A_Q[Q\mathcal{H}(Q^2|\wt\eps|^{2})]$, can be handled in a similar way, we focus on the first term, $\tfrac{1}{4}A_Q[y^2Q\mathcal{H}(Q^2|\wt\eps|^{2})]$. By the definition of $A_Q$, we have
\begin{align*}
	A_Q[y^2Q\mathcal{H}(Q^2|\wt\eps|^{2})]
	=
	\partial_y[\tfrac{y^3}{\langle y \rangle}Q\mathcal{H}(Q^2|\wt\eps|^{2})]-|D|[\tfrac{y^2}{\langle y \rangle}Q\mathcal{H}(Q^2|\wt\eps|^{2})].
\end{align*}
For the latter term $|D|[\frac{y^2}{\langle y \rangle}Q\mathcal{H}(Q^2|\wt\eps|^{2})]$, we have
\begin{equation}
	\begin{aligned}
		\||D|[\tfrac{y^2}{\langle y \rangle}Q\mathcal{H}(Q^2|\wt\eps|^{2})]\|_{L^2}=
		&\|\partial_y[\tfrac{y^2}{\langle y \rangle}Q\mathcal{H}(Q^2|\wt\eps|^{2})]\|_{L^2}
		\\
		\lesssim &\|Q\wt\eps\|_{L^\infty}\|\wt\eps\|_{\dot \calH^1}
		\lesssim \|\wt\eps\|_{\dot \calH^1}^2\lesssim \lambda^2.
	\end{aligned} \label{eq:td N2 1}
\end{equation}
For the term $ \partial_y[\frac{y^3}{\langle y \rangle}Q\mathcal{H}(Q^2|\wt\eps|^{2})]$, if the derivative falls on $ \frac{y^3}{\langle y \rangle}Q$, then we observe that $ \partial_y(\frac{y^3}{\langle y \rangle}Q) \in L^\infty$, so the estimate proceeds similarly as before. However, if the derivative falls on $ \mathcal{H}(Q^2|\wt\eps|^{2})$, then since $ \frac{y^3}{\langle y \rangle}Q \notin L^\infty$, an additional step is required. Using \eqref{eq:CommuteHilbertDerivative}, we have
\begin{align*}
	\tfrac{y^3}{\langle y \rangle}Q\partial_y\calH(Q^2|\wt\eps|^{2})
	=\tfrac{y^2}{\langle y \rangle}Q\calH[y\partial_y(Q^2|\wt\eps|^{2})].
\end{align*}
Thus, we obtain 
\begin{align}
	\|\partial_y[\tfrac{y^3}{\langle y \rangle}Q\mathcal{H}(Q^2|\wt\eps|^{2})]\|_{L^2}
	\lesssim \|Q\wt\eps\|_{L^\infty}\|\wt\eps\|_{\dot \calH^1}
	\lesssim \|\wt\eps\|_{\dot \calH^1}^2\lesssim \lambda^2. \label{eq:td N2 2}
\end{align}
From \eqref{eq:td N2 1} and \eqref{eq:td N2 2}, we have
\begin{align*}
	\|A_Q[y^2Q\mathcal{H}(Q^2|\wt\eps|^{2})]\|_{L^2}\lesssim \lambda^2.
\end{align*}
Since the second term $\tfrac{1}{4}A_Q[Q\mathcal{H}(Q^2|\wt\eps|^{2})]$ can be estimated in a similar manner, we obtain the bound
\begin{align}
	\|A_Q\td N_2\|_{L^2}\lesssim \lambda^2. \label{eq:td N2 goal}
\end{align}
Therefore, collecting \eqref{eq:td N1 goal} and \eqref{eq:td N2 goal}, we deduce 
\begin{align*}
\|A_{Q}\td N_{Q}(\widehat{\eps})\|_{L^{2}}\lesssim\lambda^{2}+\lambda^{\frac{1}{2}}\|A_{Q}\widetilde{L}_{Q}\eps\|_{L^{2}},
\end{align*}
and so 
\begin{align}
\|A_{Q}\eps_{1}-A_{Q}\widetilde{L}_{Q}\eps\|_{L^{2}}\lesssim\lambda^{2}+\lambda^{\frac{1}{2}}\|A_{Q}\widetilde{L}_{Q}\eps\|_{L^{2}}.\label{eq:5.3 proof coercivity 1}
\end{align}
By \eqref{eq: AQw1 esti} and $A_{Q}P_{1}=0$, we have $\|A_{Q}\eps_{1}\|_{L^{2}}\lesssim\lambda^{2}$.
Thus, from \eqref{eq:5.3 proof coercivity 1}, we have 
\begin{align}
(1-C\lambda^{\frac{1}{2}})\|A_{Q}\widetilde{L}_{Q}\eps\|_{L^{2}}\lesssim\lambda^{2}+\|A_{Q}\eps_{1}\|_{L^{2}}\lesssim\lambda^{2}.\label{eq:5.3 proof coercivity 2}
\end{align}
From $\|\eps\|_{\dot{\mathcal{H}}^{2}}\sim\|A_{Q}\widetilde{L}_{Q}\eps\|_{L^{2}}$,
we conclude $\|\eps\|_{\dot{\mathcal{H}}^{2}}\lesssim\lambda^{2}$.

Next, we prove $\|\eps_{1}\|_{\dot{\mathcal{H}}^{1}}\lesssim\lambda^{2}$.
It suffices to show $\|Q\eps_{1}\|_{L^{2}}+\|\partial_{y}\eps_{1}\|_{L^{2}}\lesssim\lambda^{2}.$
We recall \eqref{eq:5.3 Remark goal}, 
\begin{align*}
\eps_{1} & =\widetilde{L}_{Q}\eps+\td N_{Q}(\widehat{\eps})-yQ\cdot\tfrac{1}{2\pi}{\textstyle \int_{\R}\Re(Q^{3}\eps_{e})dy}\\
 & =\D_{Q}\eps-yQ\cdot\tfrac{1}{2\pi}{\textstyle \int_{\R}\Re(Q^{3}\eps_{e})dy+\td N_{Q}(\widehat{\eps})+Q^{-1}\mathcal{H}\Re(Q^{3}\eps).}
\end{align*}
By $\|\eps\|_{\dot{\mathcal{H}}^{2}}\lesssim\lambda^{2}$, we obtain
\begin{align*}
\|\mathbf{D}_{Q}\eps\|_{\dot{\mathcal{H}}^{1}}\lesssim\lambda^{2}.
\end{align*}
For the profile term, we have 
\begin{align*}
{\textstyle \left\Vert yQ\cdot\frac{1}{2\pi}\int_{\R}\Re(Q^{3}\eps_{e})dy\right\Vert _{\dot{\mathcal{H}}^{1}}} & \lesssim\|yQ\|_{\dot{\mathcal{H}}^{1}}{\textstyle \int_{\R}|Q^{3}\eps|dy}\\
 & \lesssim\|yQ\|_{\dot{\mathcal{H}}^{1}}\|Q\|_{L^{2}}\|Q^{2}\eps\|_{L^{2}}\lesssim\lambda^{2}.
\end{align*}
For the nonlinear term $\td N_{Q}(\widehat{\eps})=\td N_1+\td N_2$, we follow the proof of \eqref{eq:td N1 goal} and \eqref{eq:td N2 goal}. Indeed, by \eqref{eq:td N1 aux}, the proof of \eqref{eq:td N1 goal}, and \eqref{eq:5.3 proof coercivity 2}, we have
\begin{align*}
    \|\td N_1\|_{\dot{\mathcal{H}}^{1}}
    \lesssim \lambda^2+\lambda^{\frac{1}{2}}\|A_{Q}\widetilde{L}_{Q}\eps\|_{L^{2}}
    \lesssim \lambda^2. 
\end{align*}
For $\td N_2$, thanks to \eqref{eq:CommuteHilbert}, we have
\begin{align*}
    \td N_2
    =
    \tfrac{1}{4}[y^2Q\mathcal{H}(Q^2|\wt\eps|^{2})+Q\mathcal{H}(Q^2|\wt\eps|^{2})]
    -\tfrac{1}{4\pi}yQ{\textstyle \int}Q^2|\wt\eps|^{2}.
\end{align*}
By the proof of \eqref{eq:td N2 goal}, we have
\begin{align*}
    \|\tfrac{1}{4}[y^2Q\mathcal{H}(Q^2|\wt\eps|^{2})+Q\mathcal{H}(Q^2|\wt\eps|^{2})]\|_{\dot{\mathcal{H}}^{1}}
    \lesssim \lambda^2.
\end{align*}
For $-\tfrac{1}{4\pi}yQ{\textstyle \int}Q^2|\wt\eps|^{2}$, we obtain
\begin{align*}
    \|\tfrac{1}{4\pi}yQ{\textstyle \int}Q^2|\wt\eps|^{2}\|_{\dot\calH^1}
    \lesssim \|yQ\|_{\dot\calH^1}\|Q\wt\eps\|_{L^2}
    \lesssim \|\wt\eps\|_{\dot\calH^1}
    \lesssim \lambda^2,
\end{align*}
which implies
\begin{align*}
    \|\td N_2\|_{\dot{\mathcal{H}}^{1}}\lesssim \lambda^2.
\end{align*}
Thus, we arrive at
\begin{align*}
    \|\td N_{Q}(\widehat{\eps})\|_{\dot{\mathcal{H}}^{1}}\lesssim \lambda^2.
\end{align*}
For the last term $Q^{-1}\mathcal{H}\Re(Q^{3}\eps_{o})$, we have
\begin{align*}
\|Q^{-1}\mathcal{H}\Re(Q^{3}\eps)\|_{\dot{\mathcal{H}}^{1}} & \lesssim\|\mathcal{H}\Re(Q^{3}\eps)\|_{L^{2}}+\|Q^{-1}\partial_{y}\mathcal{H}\Re(Q^{3}\eps)\|_{L^{2}}\\
 & \lesssim\lambda^{2}+\|Q^{-1}\partial_{y}\mathcal{H}\Re(Q^{3}\eps)\|_{L^{2}}.
\end{align*}
By \eqref{eq:CommuteHilbertDerivative} with $Q^{-1}=\frac{1}{2}(1+y^{2})Q$,
we have 
\begin{align}
Q^{-1}\partial_{y}\mathcal{H}\Re(Q^{3}\eps)=\tfrac{1}{2}Q\mathcal{H}\partial_{y}\Re(Q^{3}\eps)+\tfrac{1}{2}yQ\mathcal{H}[y\partial_{y}\Re(Q^{3}\eps)],\label{eq:5.3 Lemma proof 1}
\end{align}
with $f=\Re(Q^{3}\eps)$. By \eqref{eq:5.3 Lemma proof 1}, we justify
$Q^{-1}\partial_{y}\mathcal{H}\Re(Q^{3}\eps)\in L^{2}$, and we have
\begin{align*}
\|Q^{-1}\partial_{y}\mathcal{H}\Re(Q^{3}\eps)\|_{L^{2}}\lesssim\|\eps\|_{\dot{\mathcal{H}}^{2}}\lesssim\lambda^{2}.
\end{align*}
This proves $\|\eps_{1}\|_{\dot{\mathcal{H}}^{1}}\lesssim\lambda^{2}$,
and we conclude \eqref{eq:CoercivityEstimateH2epsilon}.

\textbf{Proof of \eqref{eq:5.3 interpolation additional decay} and \eqref{eq:mu estimate by b}.}
For \eqref{eq:5.3 interpolation additional decay}, we have 
\begin{align*}
\partial_{y}\eps_{o} & =\partial_{y}\widehat{\eps}_{o}-\partial_{y}P_{o}=\partial_{y}\widehat{\eps}_{o}-\tfrac{i\nu}{2}\partial_{y}(yQ)\in L^{2},\\
Q\eps_{o} & =Q\widehat{\eps}_{o}-QP_{o}=Q\widehat{\eps}_{o}-\tfrac{i\nu}{2}yQ^2 \in L^{2},\\
\eps_{1,e} & =w_{1,e}-P_{1,e}=w_{1,e}-\tfrac{i\nu+\mu}{2}Q\in L^{2}.
\end{align*}
Thus, we have 
\begin{align}
\|\eps_{o}\|_{\dot\calH^1} & \lesssim\|\widehat{\eps}_{o}\|_{\dot\calH^1}+|\nu|\lesssim\lambda, \label{eq:eps odd H1}
\\
\|\eps_{1,e}\|_{L^{2}} & \lesssim\|w_{1,e}\|_{L^{2}}+|\nu|+|\mu|\lesssim\lambda+|\mu|.
\end{align}
To conclude \eqref{eq:5.3 interpolation additional decay} we have to estimate
the size of $\mu$, \eqref{eq:mu estimate by b}. We recall the definition of $\mu$, \eqref{eq:definition of mu},
\begin{align*}
	\mu=-\tfrac{1}{\pi}{\textstyle \int_{\bbR}}\Re(yQ^{3}\eps_{o})dy
	-\tfrac{1}{2\pi} {\textstyle \int_{\bbR}} yQ^2|\wt\eps|^2dy.
\end{align*}
By interpolating \eqref{eq:eps odd H1} and $\|\eps\|_{\dot{\mathcal{H}}^{2}}\lesssim\lambda^{2}$, we have
\begin{align*}
	|{\textstyle \int_{\bbR}}\Re(yQ^{3}\eps_{o})dy|
	\lesssim& \|Q^{\frac{3}{2}-}\eps_o\|_{L^2}\lesssim \lambda^{\frac{3}{2}-}\ll \lambda,
	\\
	{\textstyle \int_{\bbR}} yQ^2|\wt\eps|^2dy
	\lesssim& \|Q\wt\eps\|_{L^2}\|\wt\eps\|_{L^2}\lesssim \lambda,
\end{align*}
which prove \eqref{eq:mu estimate by b}.

\textbf{Proof of \eqref{eq:InterpolationEstimateW1,inf}.}
We now show \eqref{eq:InterpolationEstimateW1,inf}. We decompose $\partial_{y}\widehat{\eps}$ as 
\begin{align*}
\partial_{y}\widehat{\eps}=(\partial_{y}P)\chi_{\lambda^{-1}}+\mathring{\eps},\quad\mathring{\eps}=(\partial_{y}P)(1-\chi_{\lambda^{-1}})+\partial_{y}\eps.
\end{align*}
Thus, using $\mathring{\eps}=\partial_{y}\widehat{\eps}-(\partial_{y}P)\chi_{\lambda^{-1}}$, we have
\begin{align*}
\|\mathring{\eps}\|_{L^{2}}\leq\|\partial_{y}\widehat{\eps}\|_{L^{2}}+\|(\partial_{y}P)\chi_{\lambda^{-1}}\|_{L^{2}}\lesssim\lambda+b\lambda^{-\frac{1}{2}}+|\nu|\lesssim\lambda.
\end{align*}
In addition, 
\begin{align*}
\|\partial_{y}\mathring{\eps}\|_{L^{2}}\leq\|\partial_{yy}\eps\|_{L^{2}}+\|\partial_{y}[(\partial_{y}P)(1-\chi_{\lambda^{-1}})]\|_{L^{2}}\lesssim\lambda^{2}+b\lambda^{\frac{1}{2}}+|\nu|\lambda^{1-}\lesssim\lambda^{2}.
\end{align*}
Therefore, we deduce 
\begin{align*}
\|\mathring{\eps}\|_{L^{\infty}}\lesssim\|\mathring{\eps}\|_{L^{2}}^{\frac{1}{2}}\|\partial_{y}\mathring{\eps}\|_{L^{2}}^{\frac{1}{2}}\lesssim\lambda^{\frac{3}{2}}.
\end{align*}
So, we have 
\begin{align}
\begin{split}\|\partial_{y}\widehat{\eps}_{e}\|_{L^{\infty}} & \leq\|(\partial_{y}P_{e})\chi_{\lambda^{-1}}\|_{L^{\infty}}+\|\mathring{\eps}\|_{L^{\infty}}\lesssim\lambda^{\frac{3}{2}},\\
\|\partial_{y}\widehat{\eps}_{o}\|_{L^{\infty}} & \leq\|(\partial_{y}P_{o})\chi_{\lambda^{-1}}\|_{L^{\infty}}+\|\mathring{\eps}\|_{L^{\infty}}\lesssim|\nu|\lesssim\lambda^{\frac{3}{2}(1-\kappa)}.
\end{split}
\label{eq:5.3 Lemma proof L inf goal 3}
\end{align}
By a similar argument, we also deduce 
\begin{align}
 & \|Q\widehat{\eps}_{e}\|_{L^{\infty}}\lesssim\lambda^{\frac{3}{2}},\quad\|Q\widehat{\eps}_{o}\|_{L^{\infty}}\lesssim|\nu|\lesssim\lambda^{\frac{3}{2}(1-\kappa)}.\label{eq:5.3 Lemma proof L inf goal 4}
\end{align}
From \eqref{eq:5.3 Lemma proof L inf goal 3} and \eqref{eq:5.3 Lemma proof L inf goal 4},
we conclude \eqref{eq:InterpolationEstimateW1,inf}.

\textbf{Proof of \eqref{eq:InterpolationEstimateW1,inf +}.}
Now, we prove \eqref{eq:InterpolationEstimateW1,inf +}. For the $\eps_1$ part, we decompose $w_{1}$ as 
\begin{align*}
	w_{1} & =\chi_{\lambda^{-1}}P_{1}+\mathring{\eps}_{1},\quad\mathring{\eps}_{1}=(1-\chi_{\lambda^{-1}})P_{1}+\eps_{1}.
\end{align*}
Then, 
\begin{align}
	\|\mathring{\eps}_{1}\|_{L^{2}}\leq\|w_{1}\|_{L^{2}}+\|\chi_{\lambda^{-1}}P_{1}\|_{L^{2}}\lesssim\lambda+|\nu|+|\mu|\lesssim \lambda,\label{eq:5.3 Lemma proof e ring 1}
\end{align}
and 
\begin{align}
	\|\mathring{\eps}_{1}\|_{\dot{H}^{1}}\leq\|\eps_{1}\|_{\dot{H}^{1}}+\|(1-\chi_{\lambda^{-1}})P_{1}\|_{\dot{H}^{1}}\lesssim\lambda^{2}+(|\nu|+|\mu|)\lambda^{\frac{3}{2}}\sim\la^{2}.\label{eq:5.3 Lemma proof e ring 2}
\end{align}
Thus, interpolating \eqref{eq:5.3 Lemma proof e ring 1} and \eqref{eq:5.3 Lemma proof e ring 2},
we deduce $\|\mathring{\eps}_{1}\|_{L^{\infty}}\lesssim\lambda^{\frac{3}{2}}$,
and then we have 
\begin{align}
	\begin{split}
		\|\eps_{1,o}\|_{L^{\infty}}\leq&\|(\chi_{\lambda^{-1}}-1)P_{1,o}\|_{L^{\infty}}+\|\mathring{\eps}_{1}\|_{L^{\infty}}\lesssim b+\lambda^{\frac{3}{2}} \sim\lambda^{\frac{3}{2}},\\
		\|\eps_{1,e}\|_{L^{\infty}}\leq&\|(\chi_{\lambda^{-1}}-1)P_{1,e}\|_{L^{\infty}}+\|\mathring{\eps}_{1}\|_{L^{\infty}}\lesssim (|\nu|+|\mu|)\lambda + \lambda^{\frac{3}{2}}\lesssim\lambda^{\frac{3}{2}}.
	\end{split}
	\label{eq:proof L inf no parity 1}
\end{align}
For the $\eps$ part of \eqref{eq:InterpolationEstimateW1,inf +}, we have 
\begin{align*}
	\|\partial_{y}\eps_{e}\|_{L^{\infty}}\leq\|\partial_{y}\widehat{\eps}_{e}\|_{L^{\infty}}+\|\partial_{y}P_{e}\|_{L^{\infty}}\lesssim\lambda^{\frac{3}{2}}.
\end{align*}
Since $\partial_{y}P_{o}\in L^{2}$, we have 
\begin{align*}
	\|\partial_{y}\eps_{o}\|_{L^{2}}\leq\|\partial_{y}P_{o}\|_{L^{2}}+\|\partial_{y}\widehat{\eps}_{o}\|_{L^{2}}\lesssim\lambda.
\end{align*}
Interpolating $\|\partial_{y}\eps_{o}\|_{L^{2}}\lesssim\lambda$ and
$\|\partial_{yy}\eps_{o}\|_{L^{2}}\lesssim\lambda^{2}$, we deduce
\begin{align}
	\|\partial_{y}\eps_{o}\|_{L^{\infty}}\lesssim\lambda^{\frac{3}{2}}, \label{eq:proof L inf no parity 2}
\end{align}
which proves $\|\partial_{y}\eps\|_{L^{\infty}}\lesssim\lambda^{\frac{3}{2}}$.
By a similar argument, we show 
\begin{align}
	\|Q\eps\|_{L^{\infty}}\lesssim\lambda^{\frac{3}{2}}. \label{eq:proof L inf no parity 3}
\end{align}
Collecting \eqref{eq:proof L inf no parity 1}--\eqref{eq:proof L inf no parity 3}, we conclude \eqref{eq:InterpolationEstimateW1,inf +}.

Hence, we finish the proof. 
\end{proof}

\subsection{Modulation estimates}

\label{SubsectionModulationEstimate} In this subsection, we prove
the modulation estimates that justify the formal modulation laws \eqref{eq:FormalBlowupLawConclusion}.
First, we obtain the modulation estimates for $\la,\ga,x$ from the
equation for $w$ \eqref{eq:w-equ}. For this, we use decomposition
$w=Q+\wt{\eps}$, write the equation for $\wt{\eps}$ and use the
orthogonality for $\widehat{\eps}$, i.e., $(\widehat{\eps},\mathcal{Z}_{k})$
with $k=1,2,3$. To obtain satisfactory modulation estimates for $b,\eta,\nu$,
we need to introduce refined modulation parameters $\tilde b,\tilde{\eta},\tilde{\nu}$,
which are comparable to the original ones. This is \emph{different}
from the usual approach, where one sticks to estimate $b_{s},\eta_{s},\nu_{s}$
through differentiating $(\eps,\mathcal{Z}_{k})$ in $s$. Typically,
the main term becomes $(iH_{Q}\eps,\mathcal{Z}_{k})$, whose best
possible bound is $\lmb^{2}$ in our bootstrap scheme due to the cut-off
error of $\mathcal{Z}_{k}$. However, this bound for the modulation
laws (e.g., $|b_{s}+\tfrac{3}{2}b^{2}+\tfrac{1}{2}\eta^{2}|\aleq\lmb^{2}$)
is \emph{not} sufficient to justify $b_{s}+\tfrac{3}{2}b^{2}+\tfrac{1}{2}\eta^{2}\approx0$
since $\la^{3}\sim b^{2}$. Equipped with the refined modulation parameters,
we will measure a rescaled time derivative $(\partial_{s}-c\tfrac{\la_{s}}{\la})$. 
\begin{lem}[First modulation estimate]
\label{LemmaModulationEstimate}We have 
\begin{align}
\left|\frac{\lambda_{s}}{\lambda}+b\right|+\left|\gamma_{s}-\frac{\eta}{2}\right|+\left|\frac{x_{s}}{\lambda}-\nu\right|\lesssim\lambda^{2},\label{eq:modulationEstimate1}
\end{align}
In particular, we have 
\begin{align}
\left|\frac{\lambda_{s}}{\lambda}\right|\lesssim b,\quad|\gamma_{s}|\lesssim\eta,\quad\left|\frac{x_{s}}{\lambda}\right|\lesssim|\nu|.\label{eq:modulationEstimate1+}
\end{align}
\end{lem}

\begin{proof}
Write the equation for $\wt{\eps}$ using the formula of $P_{1}$
and $L_{Q}^{*}P_{1}$; 
\begin{align*}
(\partial_{s}-\tfrac{\lambda_{s}}{\lambda}\Lambda+\gamma_{s}i-\tfrac{x_{s}}{\lambda}\partial_{y})\widehat{\eps} & +(iL_{w}^{*}w_{1}-iL_{Q}^{*}P_{1})\\
 & =\left(\tfrac{\lambda_{s}}{\lambda}+b\right)\Lambda Q-\left(\gamma_{s}-\tfrac{\eta}{2}\right)iQ+\left(\tfrac{x_{s}}{\lambda}-\nu\right)Q_{y}.
\end{align*}
Differentiating the orthogonality conditions, $\partial_{s}(\widehat{\eps},\mathcal{Z}_{k})_{r}=0$
for $k=1,2,3$, gives 
\begin{align}
\left(\tfrac{\lambda_{s}}{\lambda}+b\right)(\Lambda Q, & \mathcal{Z}_{k})_{r}-\left(\gamma_{s}-\tfrac{\eta}{2}\right)(iQ,\mathcal{Z}_{k})_{r}+\left(\tfrac{x_{s}}{\lambda}-\nu\right)(Q_{y},\mathcal{Z}_{k})_{r}\nonumber \\
= & \tfrac{\lambda_{s}}{\lambda}(\widehat{\eps},\Lambda\mathcal{Z}_{k})_{r}-\gamma_{s}(\widehat{\eps},i\mathcal{Z}_{k})_{r}+\tfrac{x_{s}}{\lambda}(\widehat{\eps},\partial_{y}\mathcal{Z}_{k})_{r}+((iL_{w}^{*}w_{1}-iL_{Q}^{*}P_{1}),\mathcal{Z}_{k})_{r}.\label{eq:5.4.1}
\end{align}
We have 
\begin{align*}
|(\widehat{\eps},\Lambda\mathcal{Z}_{k})_{r}|=|(Q\widehat{\eps},Q^{-1}\Lambda\mathcal{Z}_{k})_{r}|\leq\|Q\widehat{\eps}\|_{L^{\infty}}\|Q^{-1}\Lambda\mathcal{Z}_{k}\|_{L^{1}}\lesssim b.
\end{align*}
Similarly, we have 
\begin{align*}
|(\widehat{\eps},i\mathcal{Z}_{k})_{r}|\lesssim b.
\end{align*}
Hence, the first three terms in RHS of \eqref{eq:5.4.1} are absorbed
in LHS. To finish the proof of \eqref{eq:modulationEstimate1}, we
reduce to showing 
\begin{equation}
|((iL_{w}^{*}w_{1}-iL_{Q}^{*}P_{1}),\mathcal{Z}_{k})_{r}|\lesssim\la^{2}.\label{eq:claim 5.2.1}
\end{equation}
Expand $L_{w}^{*}w_{1}-L_{Q}^{*}P_{1}$ as 
\begin{align*}
L_{w}^{*}w_{1}-L_{Q}^{*}P_{1}= & -\partial_{y}\eps_{1}+\tfrac{1}{2}\mathcal{H}(|w|^{2}-|Q|^{2})P_{1}+\tfrac{1}{2}\mathcal{H}(|w|^{2})\eps_{1}\\
 & -\widehat{\eps}\mathcal{H}\Re(\overline{w}w_{1})-Q\mathcal{H}\Re(\overline{w}w_{1}-QP_{1}).
\end{align*}
Now, we estimate terms in the form $(i(\text{RHS)},\mathcal{Z}_{k})$.
First, for $\partial_{y}\eps_{1}$, we estimate 
\begin{align*}
|(i\partial_{y}\eps_{1},\mathcal{Z}_{k})_{r}|\lesssim\|\partial_{y}\eps_{1}\|_{L^{2}}\lesssim\lambda^{2}.
\end{align*}
For $\mathcal{H}(|w|^{2}-|Q|^{2})P_{1}$, using \eqref{eq:mu estimate by b},
\begin{align*}
|(i\mathcal{H}(|w|^{2}-|Q|^{2})P_{1},\mathcal{Z}_{k})_{r}| & \lesssim\|P_{1}\|_{L^{\infty}}\|\mathcal{H}(|w|^{2}-|Q|^{2})\|_{L^{\infty}}\|\mathcal{Z}_{k}\|_{L^{1}}\\
 & \lesssim (b+|\eta|+|\nu|+|\mu|)(\|\mathcal{H}(Q\widehat{\eps})\|_{L^{\infty}}+\|\mathcal{H}(|\widehat{\eps}|^{2})\|_{L^{\infty}}).
\end{align*}
By an elementary inequality $\|f\|_{L^\infty}\lesssim \|f\|_{L^2}^{\frac{1}{2}}\|f\|_{\dot H^1}^{\frac{1}{2}}$, we have
\begin{align*}
    \|\mathcal{H}(Q\widehat{\eps})\|_{L^{\infty}}&\lesssim \|\wt\eps\|_{\dot \calH^1}\lesssim \lambda,
    \\
    \|\mathcal{H}(|\widehat{\eps}|^2)\|_{L^{\infty}}&\lesssim \|\wt\eps\|_{L^\infty}\|\wt\eps\|_{L^2}^{\frac{1}{2}}\|\wt\eps\|_{\dot H^1}^{\frac{1}{2}} \lesssim \lambda.
\end{align*}
Thus, we have
\begin{align*}
    |(i\mathcal{H}(|w|^{2}-|Q|^{2})P_{1},\mathcal{Z}_{k})_{r}|\lesssim \lambda^2.
\end{align*}
For $\mathcal{H}(|w|^{2})\eps_{1}$, 
\begin{align*}
|(i\mathcal{H}(|w|^{2})\eps_{1},\mathcal{Z}_{k})_{r}|\lesssim & |(i\mathcal{H}(Q^{2})\eps_{1},\mathcal{Z}_{k})_{r}|+|(i\mathcal{H}(\Re(Q\widehat{\eps}))\eps_{1},\mathcal{Z}_{k})_{r}|\\
 & +|(i\mathcal{H}(|\widehat{\eps}|^{2})\eps_{1},\mathcal{Z}_{k})_{r}|\\
\lesssim & (\|Q\eps_{1}\|_{L^{2}}+\|\eps_{1}\|_{L^{\infty}}\|Q\widehat{\eps}\|_{L^{2}}+\|\eps_{1}\|_{L^{\infty}}\||\widehat{\eps}|^{2}\|_{L^{2}})\|\mathcal{Z}_{k}\|_{L^{2}}\\
\lesssim & \lambda^{2}.
\end{align*}
For $\widehat{\eps}\mathcal{H}\Re(\overline{w}w_{1})$, we estimate
\begin{align*}
|(i\widehat{\eps}\mathcal{H}\Re(\overline{w}w_{1}),\mathcal{Z}_{k})_{r}| & =|(i(Q\widehat{\eps})\mathcal{H}\Re(\overline{w}w_{1}),Q^{-1}\mathcal{Z}_{k})_{r}|\\
 & \lesssim\|Q\widehat{\eps}\|_{L^{\infty}}\|\mathcal{H}\Re(\overline{w}w_{1})\|_{L^{2}}\|Q^{-1}\mathcal{Z}_{k}\|_{L^{2}}\\
 & \lesssim b^{1-}\|w\|_{L^{2}}\|w_{1}\|_{L^{\infty}}\lesssim b^{2-}\lesssim\lambda^{2}.
\end{align*}
Here we applied \eqref{eq:InterpolationEstimateW1,inf} to obtain $|Q\widehat{\eps}|_{L^{\infty}} \lesssim b^{1-}$.

For $Q\mathcal{H}\Re(\overline{w}w_{1}-QP_{1})$, using $\overline{w}w_{1}-QP_{1}=\overline{w}\eps_{1}+\overline{\widehat{\eps}}P_{1},$we
have 
\begin{align*}
|(iQ\mathcal{H}\Re(\overline{w}w_{1}-QP_{1}),\mathcal{Z}_{k})_{r}|\lesssim\|\mathcal{H}\Re(\overline{w}\eps_{1})\|_{L^{\infty}}+\|\mathcal{H}\Re(\overline{\widehat{\eps}}P_{1})\|_{L^{\infty}}.
\end{align*}
Then, we further decompose as 
\begin{align*}
\|\mathcal{H}\Re(\overline{w}\eps_{1})\|_{L^{\infty}} & \leq\|\mathcal{H}\Re(Q\eps_{1})\|_{L^{\infty}}+\|\mathcal{H}\Re(\overline{\widehat{\eps}}\eps_{1})\|_{L^{\infty}}\\
 & \lesssim\|Q\eps_{1}\|_{L^{2}}^{\frac{1}{2}}\|Q\eps_{1}\|_{\dot{H}^{1}}^{\frac{1}{2}}+\|\widehat{\eps}\eps_{1}\|_{L^{2}}^{\frac{1}{2}}\|\widehat{\eps}\eps_{1}\|_{\dot{H}^{1}}^{\frac{1}{2}}.
\end{align*}
We easily check that
\begin{align*}
    \|Q\eps_{1}\|_{\dot{H}^{1}}
    \leq&
    \|Q_y\eps_{1}\|_{L^2}+\|Q \partial_y\eps_{1}\|_{L^2}
    \\
    \leq& \|Q^{-1}Q_y\|_{L^\infty} \|Q\eps_{1}\|+\|Q\|_{L^\infty}\|\partial_y\eps_{1}\|_{L^2}
    \lesssim \|\eps_{1}\|_{\dot\calH^1},
\end{align*}
which implies
\begin{align*}
    \|Q\eps_{1}\|_{L^{2}}^{\frac{1}{2}}\|Q\eps_{1}\|_{\dot{H}^{1}}^{\frac{1}{2}}\lesssim\|\eps_{1}\|_{\dot{\mathcal{H}}^{1}}\lesssim\lambda^{2}.
\end{align*}
By \eqref{eq:CoercivityEstimateH1}--\eqref{eq:CoercivityEstimateH2epsilon} and \eqref{eq:InterpolationEstimateW1,inf +}, we have 
\begin{align*}
\|\widehat{\eps}\eps_{1}\|_{L^{2}}  \lesssim\|\widehat{\eps}\|_{L^{2}}\|\eps_{1}\|_{L^{\infty}}
\lesssim\lambda^{\frac{3}{2}},
\end{align*}
and 
\begin{align*}
\|\partial_{y}(\widehat{\eps}\eps_{1})\|_{L^{2}} & \leq\|(\partial_{y}\widehat{\eps})\eps_{1}\|_{L^{2}}+\|\widehat{\eps}(\partial_{y}\eps_{1})\|_{L^{2}}\\
 & \lesssim\|\partial_{y}\widehat{\eps}\|_{L^{2}}\|\eps_{1}\|_{L^{\infty}}
 +\|\widehat{\eps}\|_{L^{\infty}}\|\partial_{y}\eps_{1}\|_{L^{2}}\\
 & \lesssim\lambda^{\frac{5}{2}}.
\end{align*}
From the above two, we obtain $\|\widehat{\eps}\eps_{1}\|_{L^{2}}^{\frac{1}{2}}\|\widehat{\eps}\eps_{1}\|_{\dot{H}^{1}}^{\frac{1}{2}}\lesssim\lambda^{2}$.
Hence, $\|\mathcal{H}\Re(\overline{w}\eps_{1})\|_{L^{\infty}}\lesssim\la^{2}$.
To estimate $\|\mathcal{H}\Re(\overline{\widehat{\eps}}P_{1})\|_{L^{\infty}}$,
using the Sobolev inequality, 
\[
\|\mathcal{H}\Re(\overline{\widehat{\eps}}P_{1})\|_{L^{\infty}}\lesssim\|\overline{\widehat{\eps}}P_{1}\|_{L^{2}}^{\frac{1}{2}}\|\Re(\overline{\widehat{\eps}}P_{1})\|_{\dot{H}^{1}}^{\frac{1}{2}},
\]
and for $\overline{\widehat{\eps}}P_{1}$, we can argue similarly
to $\widehat{\eps}\eps_{1}$, but using \eqref{eq:mu estimate by b}. This yields:
\begin{equation}
	\begin{aligned}
		\|\widehat{\eps}P_{1}\|_{L^{2}} & \lesssim\|\widehat{\eps}\|_{L^{2}}\|P_{1,o}\|_{L^{\infty}}+(|\nu|+|\mu|)\|Q\widehat{\eps}\|_{L^{2}}\\
		& \lesssim b+\lambda^2\lesssim\lambda^{\frac{3}{2}}.
	\end{aligned}
	\label{eq:first mod nonlin 0}
\end{equation}
We estimate by \eqref{eq:CoercivityEstimateH1}
\begin{equation}
	\begin{aligned}
		\|\partial_{y}\Re(\ol{\widehat{\eps}}P_{1})\|_{L^{2}} & \lesssim\|\partial_{y}\widehat{\eps}\|_{L^{2}}\|P_{1,o}\|_{L^{\infty}}
		+\|Q\widehat{\eps}\|_{L^{2}}\|Q^{-1}\partial_yP_{1,o}\|_{L^{\infty}}
		+
		\|\partial_y\Re(\ol{\widehat{\eps}}P_{1,e})\|_{L^{2}}\\
		& \lesssim\la b +
		\|\partial_y\Re(\ol{\widehat{\eps}}P_{1,e})\|_{L^{2}}
		\\
		&\lesssim\lambda^{\frac{5}{2}} +
		\|\partial_y\Re(\ol{\widehat{\eps}}P_{1,e})\|_{L^{2}}.
	\end{aligned} \label{eq:first mod nonlin 1}
\end{equation}
For $\|\partial_y\Re(\ol{\widehat{\eps}}P_{1,e})\|_{L^{2}}$, we observe that
\begin{align*}
	\Re (\ol{P_{o}}P_{1,e})=yQ^2\cdot O(\nu^2), 
\end{align*}
and there is no $\mu$ part. Thus, we have
\begin{equation}
	\begin{aligned}
		\|\partial_y\Re(\ol{\widehat{\eps}}P_{1,e})\|_{L^{2}}
		\leq&
		\|\partial_y\Re(\ol{\widehat{\eps}}P_{1,e}-\ol{P_{o}}P_{1,e})\|_{L^{2}}+\|\partial_y\Re (\ol{P_{o}}P_{1,e})\|_{L^2}
		\\
		\lesssim& \lambda^{3-}+\|\partial_y\Re(\ol{(\widehat{\eps}-P_{o})}P_{1,e})\|_{L^{2}}.
	\end{aligned}
	\label{eq:first mod nonlin 2}
\end{equation}
For the last term, by \eqref{eq:CoercivityEstimateH2epsilon}, \eqref{eq:mu estimate by b}, and the decomposition $\wt\eps=P+\eps$, we deduce
\begin{equation}
	\begin{aligned}
		\|\partial_y\Re(\ol{(\widehat{\eps}-P_{o})}P_{1,e})\|_{L^{2}}
		\lesssim (|\nu|+|\mu|)\|\wt\eps-P_o\|_{\dot \calH^2}
		\lesssim& \lambda(\|\eps\|_{\dot \calH^2}+\|P_e\|_{\dot \calH^2})
		\\
		\lesssim & \lambda^3+\lambda b\lesssim \lambda^{\frac52}.
	\end{aligned} \label{eq:first mod nonlin 3}
\end{equation}
Hence, by combining \eqref{eq:first mod nonlin 0}--\eqref{eq:first mod nonlin 3}, we have 
\[\|\overline{\widehat{\eps}}P_{1}\|_{L^{2}}^{\frac{1}{2}}\|\Re(\overline{\widehat{\eps}}P_{1})\|_{\dot{H}^{1}}^{\frac{1}{2}}\lesssim\lambda^{2}
\]
and so 
\[
|(iQ\mathcal{H}\Re(\overline{w}w_{1}-QP_{1}),\mathcal{Z}_{k})_{r}|\lesssim\la^{2}.
\]
This completes the proof of \eqref{eq:claim 5.2.1}. 
\end{proof}
Now, we introduce the refined modulation parameters $\widetilde{b},\widetilde{\eta},\text{ and }\widetilde{\nu}$.
As explained above, modulation estimates, obtained by differentiating
orthogonality conditions (that is, $(\eps,\mathcal{Z}_{k})_{r}$ with
$k=4,5,6$) are insufficient to justify formal laws \eqref{eq:FormalBlowupLawConclusion}.
Moreover, an attempt to add a correctional term to $b$ still creates
unwanted terms $b_{s}$ and technical issues. They originate mainly
from the lack of decay of $P$ and the technical issues from the Hilbert
transform in nonlinear terms. Here, we search for quantities written
only in terms of $w$ or $w_{1}$, but comparable to the original
modulation parameters.

Let the normalizing constant $A$ be defined by 
\begin{align*}
A\coloneqq\lim_{R\to\infty}\frac{\|\frac{1}{2}yQ\sqrt{\chi_{R}}\|_{L^{2}}^{2}}{R}.
\end{align*}
We define the refined modulation parameters: 
\begin{align}
\widetilde{b}\coloneqq\frac{(w_{1},i\frac{1}{2}yQ\chi_{R_{1}})_{r}}{AR_{1}},\quad\widetilde{\eta}\coloneqq\frac{(w_{1},\frac{1}{2}yQ\chi_{R_{1}})_{r}}{AR_{1}},\quad\widetilde{\nu}\coloneqq\frac{(w_{1},i\frac{1}{2}Q\chi_{R_{1}})_{r}}{\|\frac{1}{2}Q\|_{L^{2}}^{2}}\label{eq:RefinedModulationDefinition}
\end{align}
with the cut-off radius $R_{1}=\lambda^{-\frac{3}{4}}\sim b^{-\frac{1}{2}}$.
This radius corresponds to the self-similar scale in the $(t,x)$
coordinates. One can observe that $\frac{|(\eps_{1},i\frac{1}{2}yQ\chi_{R_{1}})_{r}|}{AR_{1}}\ll b$
since $\eps_{1}$ is a radiation and using the $P_{1}$-formula we
expect 
\[
\tilde b\approx\frac{(P_{1},i\frac{1}{2}yQ\chi_{R_{1}})_{r}}{AR_{1}}=b\frac{(i\tfrac{1}{2}yQ,i\frac{1}{2}yQ\chi_{R_{1}})_{r}}{AR_{1}}\approx b,
\]
where we used the definition of $A$. Similar arguments apply to $\tilde{\eta}$
and $\tilde{\nu}$. With these corrections, we have satisfactory modulation
estimates. 
\begin{lem}[Refined modulation estimate]
\label{LemmaModulationEstimateRefined} Let $\widetilde{b}$, $\widetilde{\eta}$,
and $\td{\nu}$ be the refined modulation parameters defined by \eqref{eq:RefinedModulationDefinition}.
Then, we have the proximity estimates 
\begin{align}
|\widetilde{b}-b|\lesssim b^{1+\frac{1}{12}},\quad|\widetilde{\eta}-\eta|\lesssim b^{1+\frac{1}{12}},\quad|\widetilde{\nu}-\nu|\lesssim b^{1+\frac{1}{12}},\label{eq:RefinedModulationDiffer}
\end{align}
and the modulation estimates 
\begin{align}
\left|\left(\partial_{s}-\frac{3}{2}\frac{\lambda_{s}}{\lambda}\right)\widetilde{b}+\frac{1}{2}\eta^{2}\right|\lesssim b^{2+\frac{1}{12}-},\quad\left|\left(\partial_{s}-\frac{\lambda_{s}}{\lambda}\right)\widetilde{\eta}\right|\lesssim b^{2+\frac{1}{12}-},\quad\left|\left(\partial_{s}-\frac{\lambda_{s}}{\lambda}\right)\widetilde{\nu}\right|\lesssim b^{2}.\label{eq:RefinedModulationBlowupLaw}
\end{align}
\end{lem}

To prove this lemma, we need to estimate the nonlinear term in $w_{1}$
equation (see \eqref{eq:w1-equ-modify} and \eqref{eq:DefinitionNL1})
as the refined modulation parameters are defined with $w_{1}$. We
state the first nonlinear estimate for $\text{NL}_1$ in \eqref{eq:DefinitionNL1}. We decompose it into even and odd parts and extract their main contributions.
\begin{lem}
\label{LemmaNonlinearEstimate1} We have a estimate for the odd part;
\begin{align}
	\|\textnormal{NL}_{1,o}+\tfrac{1}{2}(\nu^2-i\nu\mu)Q_y
	\|_{L^{2}}\lesssim b^{2}.\label{eq:NL1Estimate}
\end{align}
Moreover, for the even part, we have 
\begin{align}
	\|\textnormal{NL}_{1,e}+[i(b\Lambda_{-1}+\tfrac{\eta}{2}i-\nu\partial_{y})P_{1}+b\nu\tfrac{1}{2}Q-b\mu\tfrac{i}{2}Q]_{e}\|_{L^{2}}\lesssim b^{2}.\label{eq:NL1eWithoutProfileEstimate}
\end{align}
\end{lem}

We postpone the proof of Lemma~\ref{LemmaNonlinearEstimate1} to
Section~\ref{Section Proof of nonlinear estimates}. If we collect
the quadratic terms in $b,\eta,\nu$ from $\textnormal{NL}_{1}$,
then it should be 
\begin{align}
-i(b\Lambda_{-1}+\tfrac{\eta}{2}i-\nu\partial_{y})P_{1}-i(\tfrac{3b^{2}}{2}+\tfrac{\eta^{2}}{2})\tfrac{i}{2}yQ-ib\eta\tfrac{1}{2}yQ+ib\nu\tfrac{i}{2}Q+ib\mu\tfrac{1}{2}Q.\label{eq:NL1EstimateProfileTerm}
\end{align}
Note that the profiles in \eqref{eq:NL1Estimate} and \eqref{eq:NL1eWithoutProfileEstimate} are the odd part and the even part
of \eqref{eq:NL1EstimateProfileTerm}, respectively. These terms are only of size $b^{2-}$ due to a slightly weaker bound $|\nu|\le b^{1-\kappa}$ and \eqref{eq:mu estimate by b}.

\begin{proof}[Proof of Lemma~\ref{LemmaModulationEstimateRefined} assuming Lemma~\ref{LemmaNonlinearEstimate1}]
We first show the proximity estimates \eqref{eq:RefinedModulationDiffer}.
We can rewrite 
\begin{align*}
\widetilde{b} & =b\frac{\|\frac{1}{2}yQ\sqrt{\chi_{R_{1}}}\|_{L^{2}}^{2}}{AR_{1}}+\frac{(\eps_{1},i\frac{1}{2}yQ\chi_{R_{1}})_{r}}{AR_{1}},\\
\widetilde{\eta} & =\eta\frac{\|\frac{1}{2}yQ\sqrt{\chi_{R_{1}}}\|_{L^{2}}^{2}}{AR_{1}}+\frac{(\eps_{1},\frac{1}{2}yQ\chi_{R_{1}})_{r}}{AR_{1}},\\
\widetilde{\nu} & =\nu\frac{\|\frac{1}{2}Q\sqrt{\chi_{R_{1}}}\|_{L^{2}}^{2}}{\|\frac{1}{2}Q\|_{L^{2}}^{2}}+\frac{(\eps_{1},i\frac{1}{2}Q\chi_{R_{1}})_{r}}{\|\frac{1}{2}Q\|_{L^{2}}^{2}}.
\end{align*}
We claim that 
\begin{align}
\Big|1-\frac{\|\frac{1}{2}yQ\sqrt{\chi_{R_{1}}}\|_{L^{2}}^{2}}{AR_{1}}\Big|+\Big|1-\frac{\|\frac{1}{2}Q\sqrt{\chi_{R_{1}}}\|_{L^{2}}^{2}}{\|\frac{1}{2}Q\|_{L^{2}}^{2}}\Big| & =O(\lambda^{\frac{3}{4}}),\label{eq:refine modulation proof claim 1}\\
\frac{|(\eps_{1},i\frac{1}{2}Q\chi_{R_{1}})_{r}|}{\|\frac{1}{2}Q\|_{L^{2}}^{2}}+\frac{|(\eps_{1},i\frac{1}{2}yQ\chi_{R_{1}})_{r}|}{AR_{1}}+\frac{|(\eps_{1},\frac{1}{2}yQ\chi_{R_{1}})_{r}|}{AR_{1}} & \lesssim b^{1+\frac{1}{12}}.\label{eq:refined modulation proof claim 2}
\end{align}
For \eqref{eq:refined modulation proof claim 2}, we estimate 
\begin{align*}
\frac{|(\eps_{1},i\frac{1}{2}Q\chi_{R_{1}})_{r}|}{\|\frac{1}{2}Q\|_{L^{2}}^{2}}\leq\frac{\|Q\eps_{1}\|_{L^{2}}\|i\frac{1}{2}\chi_{R_{1}}\|_{L^{2}}}{\|\frac{1}{2}Q\|_{L^{2}}^{2}}\lesssim\la^{2}R_{1}^{\tfrac{1}{2}}\sim b^{1+\tfrac{1}{12}}.
\end{align*}
And other two terms are estimated similarly. To prove \eqref{eq:refine modulation proof claim 1},
we have 
\begin{align*}
\partial_{R}\left(\frac{\|yQ\sqrt{\chi_{R}}\|_{L^{2}}^{2}}{R}\right) & =-\frac{1}{R^{2}}\int_{\R}y^{2}Q^{2}\chi_{R}dy-\frac{1}{R^{2}}\int_{\R}y^{3}Q^{2}\frac{1}{R}\partial_{y}\chi\left(\frac{y}{R}\right)dy\\
 & =\frac{1}{2R^{2}}\int_{\R}y^{2}Q^{4}\chi_{R}=O(R^{-2}),
\end{align*}
and this implies that 
\begin{align*}
1-\frac{\|\frac{1}{2}yQ\sqrt{\chi_{R_{1}}}\|_{L^{2}}^{2}}{AR_{1}}=\int_{R_{1}}^{\infty}\partial_{R}\left(\frac{\|\tfrac{1}{2}yQ\sqrt{\chi_{R}}\|_{L^{2}}^{2}}{AR}\right)dR=O(R_{1}^{-1})=O(\lambda^{\frac{3}{4}}).
\end{align*}
$\Big|1-\frac{\|\frac{1}{2}Q\sqrt{\chi_{R_{1}}}\|_{L^{2}}^{2}}{\|\frac{1}{2}Q\|_{L^{2}}^{2}}\Big|=O(\lambda^{\frac{3}{4}})$
is a straightforward computation. This proves \eqref{eq:refine modulation proof claim 1},\eqref{eq:refined modulation proof claim 2},
and so we conclude \eqref{eq:RefinedModulationDiffer}.

Now, we move to \eqref{eq:RefinedModulationBlowupLaw}. Note that
these are motivated by formal modulation laws.

1. $\widetilde{b}$: We compute using equation \eqref{eq:w1-equ-modify}
for $w_{1}$ to have 
\begin{align}
\left(\partial_{s}-\tfrac{3}{2}\tfrac{\lambda_{s}}{\lambda}\right)\widetilde{b}=\tfrac{1}{AR}\big[ & ((\partial_{s}-\tfrac{3}{2}\tfrac{\lambda_{s}}{\lambda})w_{1},i\tfrac{1}{2}yQ\chi_{R})_{r}\nonumber \\
 & -\tfrac{R_{s}}{R}((w_{1},i\tfrac{1}{2}y^{2}Q(\partial_{y}\chi_{R}))_{r}+(w_{1},i\tfrac{1}{2}yQ\chi_{R})_{r})\big].\label{eq: 5.4.2 Refined b 1}
\end{align}
We reorganize the first term of the second line of the RHS of \eqref{eq: 5.4.2 Refined b 1};
\begin{align*}
(w_{1},i\tfrac{1}{2}y^{2}Q(\partial_{y}\chi_{R}))_{r} & =(w_{1},\partial_{y}(i\tfrac{1}{2}y^{2}Q\chi_{R}))_{r}-(w_{1},\partial_{y}(i\tfrac{1}{2}y^{2}Q)\chi_{R})_{r}\\
 & =-(y\partial_{y}w_{1},i\tfrac{1}{2}yQ\chi_{R})_{r}-(w_{1},iyQ\chi_{R})_{r}-(w_{1},i\tfrac{1}{2}y^{2}Q_{y}\chi_{R})_{r}.
\end{align*}
Thus, we write the second line of \eqref{eq: 5.4.2 Refined b 1} 
\begin{align}
(w_{1},i\tfrac{1}{2}y^{2}Q(\partial_{y}\chi_{R}) & )_{r}+(w_{1},i\tfrac{1}{2}yQ\chi_{R})_{r}\nonumber \\
 & =-(y\partial_{y}w_{1},i\tfrac{1}{2}yQ\chi_{R})_{r}-(w_{1},i\tfrac{1}{2}yQ\chi_{R})_{r}-(w_{1},i\tfrac{1}{2}y^{2}Q_{y}\chi_{R})_{r}\nonumber \\
 & =-(y\partial_{y}w_{1},i\tfrac{1}{2}yQ\chi_{R})_{r}-(w_{1},i\tfrac{1}{4}yQ^{3}\chi_{R})_{r}.\label{eq: 5.4.2 Refined b 2}
\end{align}
The first line of \eqref{eq: 5.4.2 Refined b 1} is written by $w_{1}$-equation;
\begin{align}
 & ((\partial_{s}-\tfrac{3}{2}\tfrac{\lambda_{s}}{\lambda})w_{1},i\tfrac{1}{2}yQ\chi_{R})_{r}\nonumber \\
 & =\big((\tfrac{\lambda_{s}}{\lambda}y\partial_{y}-\gamma_{s}i+\tfrac{x_{s}}{\lambda}\partial_{y})w_{1}-iH_{Q}w_{1}+i\textnormal{NL}_{1},i\tfrac{1}{2}yQ\chi_{R}\big)_{r}.\label{eq: 5.4.2 Refined b 3}
\end{align}
Therefore, by \eqref{eq: 5.4.2 Refined b 2} and \eqref{eq: 5.4.2 Refined b 3},
we have 
\begin{align*}
AR\left(\partial_{s}-\tfrac{3}{2}\tfrac{\lambda_{s}}{\lambda}\right)\widetilde{b}= & \left(\tfrac{\lambda_{s}}{\lambda}+\tfrac{R_{s}}{R}\right)(y\partial_{y}w_{1},i\tfrac{1}{2}yQ\chi_{R})_{r}+\tfrac{R_{s}}{R}(w_{1},i\tfrac{1}{4}yQ^{3}\chi_{R})_{r}\\
 & +((-\gamma_{s}i+\tfrac{x_{s}}{\lambda}\partial_{y})w_{1}-iH_{Q}w_{1}+i\textnormal{NL}_{1},i\tfrac{1}{2}yQ\chi_{R})_{r}.
\end{align*}
Now, we take $R=R_{1}=\lambda^{-\frac{3}{4}}$. Then we have 
\begin{align*}
\left|\tfrac{\lambda_{s}}{\lambda}+\tfrac{R_{s}}{R}\right|=\left|\tfrac{(\lambda R)_{s}}{\lambda R}\right|\sim b,\quad\left|\tfrac{R_{s}}{R}\right|\sim b.
\end{align*}
For $(y\partial_{y}w_{1},i\frac{1}{2}yQ\chi_{R})_{r}$, 
\begin{align*}
(y\partial_{y}w_{1},i\tfrac{1}{2}yQ\chi_{R})_{r} & =(y\partial_{y}P_{1},i\tfrac{1}{2}yQ\chi_{R})_{r}+(y\partial_{y}\eps_{1},i\tfrac{1}{2}yQ\chi_{R})_{r}\\
 & =\tfrac{b}{4}(y\partial_{y}(yQ),yQ\chi_{R})_{r}+O(\lambda^{2}R^{\frac{3}{2}}).
\end{align*}
Since $\partial_{y}(yQ)=\tfrac{Q^{3}}{2}$, we have 
\begin{align*}
|(y\partial_{y}w_{1},i\tfrac{1}{2}yQ\chi_{R})_{r}|\lesssim b+\lambda^{2}R^{\frac{3}{2}}.
\end{align*}
Using similar arguments, we estimate 
\begin{align*}
 & |(w_{1},i\tfrac{1}{4}yQ^{3}\chi_{R})_{r}|\lesssim b,\\
 & ((-\gamma_{s}i+\tfrac{x_{s}}{\lambda}\partial_{y})w_{1},i\tfrac{1}{2}yQ\chi_{R})_{r}=-\tfrac{1}{2}|\eta|^{2}AR+O(|\eta|\lambda^{2}R^{\frac{3}{2}})+O(|\nu|^{2}+|\nu|\lambda^{2}R^{\frac{1}{2}}).
\end{align*}
On the other hands, we have 
\begin{align*}
(H_{Q}w_{1},yQ\chi_{R})_{r}=(H_{Q}\eps_{1},yQ\chi_{R})_{r}=-(\eps_{1},H_{Q}(yQ(1-\chi_{R})))_{r},
\end{align*}
 using $H_{Q}(Q)=H_{Q}(yQ)=0$. We also have 
\begin{align*}
|(\eps_{1},H_{Q}(yQ(1-\chi_{R})))_{r}| & =|(\partial_{y}\eps_{1},\partial_{y}(yQ(1-\chi_{R})))_{r}+(Q\partial_{y}\eps_{1},Q^{3}(yQ(1-\chi_{R})))_{r}\\
 & \quad+(Q\eps_{1},|D|Q(yQ(1-\chi_{R})))_{r}|\\
 & \lesssim\lambda^{2}R^{-\frac{1}{2}}.
\end{align*}
In addition, by \eqref{eq:NL1Estimate}, we have 
\begin{align*}
|(i\textnormal{NL}_{1},i\tfrac{1}{2}yQ\chi_{R})_{r}|\leq&
|(\textnormal{NL}_{1,o}+\tfrac{1}{2}(\nu^2-i\nu\mu)Q_y,\tfrac{1}{2}yQ\chi_{R})_{r}|
\\
&+
|(\tfrac{1}{2}(\nu^2-i\nu\mu)Q_y,\tfrac{1}{2}yQ\chi_{R})_{r}|
\\
\lesssim& b^{2}R^{\frac{1}{2}}+|\nu|^2+|\nu\mu|
\lesssim b^{2}R^{\frac{1}{2}}+b^{1-}\lambda.
\end{align*}
Therefore, choosing $R=R_{1}=\lambda^{-\frac{3}{4}}\sim b^{-\frac{1}{2}}$ and
dividing by $AR$, we have 
\begin{align*}
\left|\left(\partial_{s}-\frac{3}{2}\frac{\lambda_{s}}{\lambda}\right)\widetilde{b}+\frac{1}{2}\eta^{2}\right|\lesssim b^{2+\frac{1}{12}-}.
\end{align*}

2. $\widetilde{\eta}$: We have 
\begin{align}
\left(\partial_{s}-\tfrac{\lambda_{s}}{\lambda}\right)\widetilde{\eta}=\tfrac{1}{AR}\big[ & ((\partial_{s}-\tfrac{\lambda_{s}}{\lambda})w_{1},\tfrac{1}{2}yQ\chi_{R})_{r}\nonumber \\
 & -\tfrac{R_{s}}{R}((w_{1},\tfrac{1}{2}y^{2}Q(\partial_{y}\chi_{R}))_{r}+(w_{1},\tfrac{1}{2}yQ\chi_{R})_{r})\big].\label{eq: 5.4.2 Refined eta 1}
\end{align}
Using an argument similar to $\widetilde{b}$,
the second line of \eqref{eq: 5.4.2 Refined eta 1} is estimated as
\begin{align*}
(w_{1},\tfrac{1}{2}y^{2}Q(\partial_{y}\chi_{R}))_{r}+(w_{1},\tfrac{1}{2}yQ\chi_{R})_{r}=-(y\partial_{y}w_{1},\tfrac{1}{2}yQ\chi_{R})_{r}-(w_{1},\tfrac{1}{4}yQ^{3}\chi_{R})_{r},
\end{align*}
and using $\partial_{y}(yQ)=\tfrac{Q^{3}}{2}$, 
\begin{align*}
 & |(y\partial_{y}w_{1},\tfrac{1}{2}yQ\chi_{R})_{r}|\lesssim|\eta|+\lambda^{2}R^{\frac{3}{2}},\quad|(w_{1},\tfrac{1}{4}yQ^{3}\chi_{R})_{r}|\lesssim|\eta|.
\end{align*}
Thus, we have
\begin{align}
    |\tfrac{R_{s}}{R}((w_{1},\tfrac{1}{2}y^{2}Q(\partial_{y}\chi_{R}))_{r}+(w_{1},\tfrac{1}{2}yQ\chi_{R})_{r})|
    \lesssim b(|\eta|+\lambda^{2}R^{\frac{3}{2}}). \label{eq:second eta goal 1}
\end{align}
The first line of the RHS of \eqref{eq: 5.4.2 Refined eta 1} is written
as 
\begin{align}
 & ((\partial_{s}-\tfrac{\lambda_{s}}{\lambda})w_{1},\tfrac{1}{2}yQ\chi_{R})_{r}\nonumber \\
 & =((\tfrac{\lambda_{s}}{\lambda}(\Lambda_{-1}-1)-\gamma_{s}i+\tfrac{x_{s}}{\lambda}\partial_{y})w_{1}-iH_{Q}w_{1}+i\textnormal{NL}_{1},\tfrac{1}{2}yQ\chi_{R})_{r}.\label{eq: 5.4.2 Refined eta 2}
\end{align}
We note that $\Lambda_{-1}-1=\Lambda$. 
For the first term of \eqref{eq: 5.4.2 Refined eta 2}, we have 
\begin{align*}
((\tfrac{\lambda_{s}}{\lambda}\Lambda-\gamma_{s}i+\tfrac{x_{s}}{\lambda}\partial_{y})w_{1},\tfrac{1}{2}yQ\chi_{R})_{r}= & ((\tfrac{\lambda_{s}}{\lambda}\Lambda-\gamma_{s}i+\tfrac{x_{s}}{\lambda}\partial_{y})\eps_{1},\tfrac{1}{2}yQ\chi_{R})_{r}\\
 & +((\tfrac{\lambda_{s}}{\lambda}\Lambda-\gamma_{s}i+\tfrac{x_{s}}{\lambda}\partial_{y})P_{1},\tfrac{1}{2}yQ\chi_{R})_{r},
\end{align*}
and 
\begin{align*}
|((\tfrac{\lambda_{s}}{\lambda}\Lambda-\gamma_{s}i+\tfrac{x_{s}}{\lambda}\partial_{y})\eps_{1},\tfrac{1}{2}yQ\chi_{R})_{r}|\lesssim b\lambda^{2}R^{\frac{3}{2}}+|\nu|\lambda^{2}R^{\frac{1}{2}}.
\end{align*}
We have 
\begin{align*}
((\tfrac{\lambda_{s}}{\lambda} & \Lambda-\gamma_{s}i+\tfrac{x_{s}}{\lambda}\partial_{y})P_{1},\tfrac{1}{2}yQ\chi_{R})_{r}\\
= & (((\tfrac{\lambda_{s}}{\lambda}+b)\Lambda-(\gamma_{s}-\tfrac{\eta}{2})i+(\tfrac{x_{s}}{\lambda}-\nu)\partial_{y})P_{1},\tfrac{1}{2}yQ\chi_{R})_{r}\\
 & +((-b\Lambda-\tfrac{\eta}{2}i+\nu\partial_{y})P_{1},\tfrac{1}{2}yQ\chi_{R})_{r},
\end{align*}
and 
\begin{align*}
|([(\tfrac{\lambda_{s}}{\lambda}+b)\Lambda-(\gamma_{s}-\tfrac{\eta}{2})i+(\tfrac{x_{s}}{\lambda}-\nu)\partial_{y}]P_{1},\tfrac{1}{2}yQ\chi_{R})_{r}|\lesssim\lambda^{2}(bR+|\mu|),
\end{align*}
\begin{align*}
|([-b\Lambda-\tfrac{\eta}{2}i+\nu\partial_{y}]P_{1},\tfrac{1}{2}yQ\chi_{R})_{r}| & =|(-b\eta y\partial_{y}(\tfrac{1}{2}yQ)+\nu\mu\partial_{y}(\tfrac{1}{2}Q),\tfrac{1}{2}yQ\chi_{R})_{r}|\\
 & \sim|(-b\eta yQ^{3}+\nu\mu\tfrac{1}{2}yQ^{3},\tfrac{1}{2}yQ\chi_{R})_{r}|\\
 & \lesssim b|\eta|+|\nu\mu|.
\end{align*}
For the second and third terms of \eqref{eq: 5.4.2 Refined eta 2},
again using an argument similar to $\widetilde{b}$, using \eqref{eq:NL1Estimate},
\begin{align*}
|(-iH_{Q}w_{1}+i\textnormal{NL}_{1},\tfrac{1}{2}yQ\chi_{R})_{r}|
\lesssim\lambda^{2}R^{-\frac{1}{2}}+b^{2}R^{\frac{1}{2}}+b^{1-\kappa}\lambda.
\end{align*}
Thus, we have
\begin{align}
    |((\partial_{s}-\tfrac{\lambda_{s}}{\lambda})w_{1},\tfrac{1}{2}yQ\chi_{R})_{r}|
    \lesssim 
    b\lambda^{2}R^{\frac{3}{2}}+|\nu|\lambda^{2}R^{\frac{1}{2}}
    +\lambda^{2}R^{-\frac{1}{2}}+b^{2}R^{\frac{1}{2}}+b^{1-\kappa}\lambda. \label{eq:second eta goal 2}
\end{align}
Therefore, from \eqref{eq: 5.4.2 Refined eta 1}, \eqref{eq:second eta goal 1}, and \eqref{eq:second eta goal 2}, choosing $R=R_{1}=\lambda^{-\frac{3}{4}}\sim b^{-\frac{1}{2}}$ and
dividing by $AR$,
we conclude 
\begin{align*}
\left|\left(\partial_{s}-\frac{\lambda_{s}}{\lambda}\right)\widetilde{\eta}\right|\lesssim b^{2+\frac{1}{12}-}.
\end{align*}

3. $\widetilde{\nu}$: We write 
\begin{align}
\|\tfrac{1}{2}Q\|_{L^{2}}^{2}\left(\partial_{s}-\tfrac{\lambda_{s}}{\lambda}\right)\widetilde{\nu}=\big((\partial_{s}-\tfrac{\lambda_{s}}{\lambda})w_{1},i\tfrac{1}{2}Q\chi_{R}\big)_{r}-\tfrac{R_{s}}{R}(w_{1},i\tfrac{1}{2}yQ(\partial_{y}\chi_{R}))_{r}.\label{eq: 5.4.2 Refined bu 1}
\end{align}
For the second term of \eqref{eq: 5.4.2 Refined bu 1}, we estimate
\begin{align*}
|(w_{1},i\tfrac{1}{2}yQ(\partial_{y}\chi_{R}))_{r}| & =|(\eps_{1},i\tfrac{1}{2}yQ(\partial_{y}\chi_{R}))_{r}+(P_{1,e},i\tfrac{1}{2}yQ(\partial_{y}\chi_{R}))_{r}|\\
 & \lesssim\|\eps_{1}\|_{\dot{\mathcal{H}^{1}}}R^{\frac{1}{2}}+|\nu|\tfrac{|\log R|}{R}.
\end{align*}
We note that the $\mu$ part of $P_{1,e}$ is real-valued, so it vanishes in the inner product $(P_{1,e},i\tfrac{1}{2}yQ(\partial_{y}\chi_{R}))_{r}$.

The first term of \eqref{eq: 5.4.2 Refined bu 1} is expanded by $w_{1}$-equation,
\begin{align}
 & ((\partial_{s}-\tfrac{\lambda_{s}}{\lambda})w_{1},i\tfrac{1}{2}Q\chi_{R})_{r}\nonumber \\
 & =((\tfrac{\lambda_{s}}{\lambda}(\Lambda_{-1}-1)-\gamma_{s}i+\tfrac{x_{s}}{\lambda}\partial_{y})w_{1}-iH_{Q}w_{1}+i\textnormal{NL}_{1},i\tfrac{1}{2}Q\chi_{R})_{r}.\label{eq:5.4 Lemma second mouldation nu 1}
\end{align}
Furthermore, using $\Lambda=\Lambda_{-1}-1$, we reorganize 
\begin{align}
\eqref{eq:5.4 Lemma second mouldation nu 1}= & -(\eps_{1},H_{Q}\tfrac{1}{2}Q\chi_{R})_{r}\label{eq:5.4 Lemma second mouldation nu 1-1}\\
 & +((\tfrac{\lambda_{s}}{\lambda}\Lambda-\gamma_{s}i+\tfrac{x_{s}}{\lambda}\partial_{y})w_{1}-i[i(b\Lambda+\tfrac{\eta}{2}i-\nu\partial_{y})P_{1}]_{e},i\tfrac{1}{2}Q\chi_{R})_{r}\label{eq:5.4 Lemma second mouldation nu 1-2}\\
 & +(\textnormal{NL}_{1}+[i(b\Lambda_{-1}+\tfrac{\eta}{2}i-\nu\partial_{y})P_{1}+b\nu\tfrac{1}{2}Q-b\mu\tfrac{i}{2}Q]_{e},\tfrac{1}{2}Q\chi_{R})_{r}.\label{eq:5.4 Lemma second mouldation nu 1-3}
\end{align}
Here, we have used $[ib\Lambda P_{1}]_{e}=[ib\Lambda_{-1}P_{1}+b\nu\tfrac{1}{2}Q-b\mu\tfrac{i}{2}Q]_{e}$.
For \eqref{eq:5.4 Lemma second mouldation nu 1-1}, we have 
\begin{align*}
|(\eps_{1},H_{Q}\tfrac{1}{2}Q\chi_{R})_{r}|=|(\eps_{1},-H_{Q}\tfrac{1}{2}Q(1-\chi_{R}))_{r}|\lesssim R^{-1}\|\eps_{1}\|_{\dot{\mathcal{H}}^{1}}.
\end{align*}
We now control \eqref{eq:5.4 Lemma second mouldation nu 1-2}. First,
since the odd part is canceled, 
\begin{align*}
 & ((\tfrac{\lambda_{s}}{\lambda}\Lambda-\gamma_{s}i+\tfrac{x_{s}}{\lambda}\partial_{y})w_{1}+[(b\Lambda+\tfrac{\eta}{2}i-\nu\partial_{y})P_{1}]_{e},i\tfrac{1}{2}Q\chi_{R})_{r}\\
 & =([(\tfrac{\lambda_{s}}{\lambda}\Lambda-\gamma_{s}i+\tfrac{x_{s}}{\lambda}\partial_{y})w_{1}]_{e}+[(b\Lambda+\tfrac{\eta}{2}i-\nu\partial_{y})P_{1}]_{e},i\tfrac{1}{2}Q\chi_{R})_{r}.
\end{align*}
We have 
\begin{align*}
 & [(\tfrac{\lambda_{s}}{\lambda}\Lambda-\gamma_{s}i+\tfrac{x_{s}}{\lambda}\partial_{y})w_{1}]_{e}+[(b\Lambda+\tfrac{\eta}{2}i-\nu\partial_{y})P_{1}]_{e}\\
 & =[(\tfrac{\lambda_{s}}{\lambda}+b)\Lambda-(\gamma_{s}-\tfrac{\eta}{2})i]P_{1,e}+(\tfrac{x_{s}}{\lambda}-\nu)\partial_{y}P_{1,o}+[(\tfrac{\lambda_{s}}{\lambda}\Lambda-\gamma_{s}i+\tfrac{x_{s}}{\lambda}\partial_{y})\eps_{1}]_{e}.
\end{align*}
We estimate 
\begin{align*}
\left|\left([(\tfrac{\lambda_{s}}{\lambda}+b)\Lambda-(\gamma_{s}-\tfrac{\eta}{2})i]P_{1,e}+(\tfrac{x_{s}}{\lambda}-\nu)\partial_{y}P_{1,o},i\tfrac{1}{2}Q\chi_{R}\right)_{r}\right|\lesssim \lambda^3\sim b^2.
\end{align*}
Moreover, 
\begin{align*}
|((\tfrac{\lambda_{s}}{\lambda}\Lambda-\gamma_{s}i+\tfrac{x_{s}}{\lambda}\partial_{y})\eps_{1},i\tfrac{1}{2}Q\chi_{R})_{r}|\lesssim b\lambda^{2}R^{\frac{1}{2}}+|\nu|\lambda^{2}.
\end{align*}
Thus, we have 
\begin{align*}
|\eqref{eq:5.4 Lemma second mouldation nu 1-2}|\lesssim b^{2}\lambda^{\frac{1}{3}-}+b\lambda^{2}R^{\frac{1}{2}}.
\end{align*}
For the last term \eqref{eq:5.4 Lemma second mouldation nu 1-3},
we use \eqref{eq:NL1eWithoutProfileEstimate} to obtain 
\begin{align}
(\textnormal{NL}_{1}+[i(b\Lambda_{-1}+\tfrac{\eta}{2}i-\nu\partial_{y})P_{1}+b\nu\tfrac{1}{2}Q-b\mu\tfrac{i}{2}Q]_{e},\tfrac{1}{2}Q\chi_{R})_{r}\lesssim b^{2}.\label{eq:5.4 Lemma second mouldation nu 2}
\end{align}
Taking $R=R_{1}=\lambda^{-\frac{3}{4}}\sim b^{-\frac{1}{2}}$, we
conclude 
\begin{align*}
\left|\left(\partial_{s}-\frac{\lambda_{s}}{\lambda}\right)\widetilde{\nu}\right|\lesssim b^{2+\frac{1}{12}-}+b^{2}\lesssim b^2.
\end{align*}
This finishes the proof of \eqref{eq:RefinedModulationBlowupLaw}.
We note that the worst bound $b^{2}$ comes from \eqref{eq:5.4 Lemma second mouldation nu 2}. 
\end{proof}

\subsection{Proof of main bootstrap proposition and topological lemmas.}

\label{SubsectionBootstrapDescription} In this subsection, we provide
the proofs of Proposition~\ref{PropositionMainBootstrap}, Lemma~\ref{LemmaTopologicalPhiContinuous},
Lemma~\ref{LemmaPsiNonExistenceBrouwer}, and Proposition~\ref{PropSharpDescription}.
The rest
of arguments are rather standard and obtained basically by integrating
the modulation ODEs in time.

We begin with the almost conserved quantity. 
\begin{lem}
\label{Lemmab-lambdaRelation} We have 
\begin{align}
\frac{b^{2}(t)}{\lambda^{3}(t)}=(1+O(b_{0}^{2\kappa}))\frac{b_{0}^{2}}{\lambda_{0}^{3}}.\label{eq:almost conserved quantity}
\end{align}
\end{lem}

\begin{proof}
We first note that by proximity \eqref{eq:RefinedModulationDiffer},
it suffices to show it for $\tilde b$, 
\begin{align*}
\frac{\widetilde{b}^{2}(t)}{\lambda^{3}(t)}=(1+O(\widetilde{b}_{0}^{2\kappa}))\frac{\widetilde{b}_{0}^{2}}{\lambda_{0}^{3}}.
\end{align*}
Thanks to \eqref{eq:RefinedModulationBlowupLaw}, we have 
\begin{align*}
\partial_{s}\log\frac{\lambda^{3}}{\widetilde{b}^{2}}=-3\left(\frac{\lambda_{s}}{\lambda}+\widetilde{b}\right)+\frac{4}{\widetilde{b}}\left(\frac{\widetilde{b}_{s}}{2}+\frac{3\widetilde{b}^{2}}{4}\right)=O(b^{1+2\kappa}).
\end{align*}
Using \eqref{eq:RefinedModulationBlowupLaw} in the form 
\begin{align*}
\frac{b^{1+2\kappa}}{\lambda^{2}}(1+O(b^{2\kappa}))=-\frac{2\widetilde{b}_{t}}{3\widetilde{b}^{1-2\kappa}},
\end{align*}
and integrating this in $t$, we have 
\begin{align*}
|\frac{\lambda(t)^{3}}{\widetilde{b}(t)^{2}}\frac{\widetilde{b}_{0}^{2}}{\lambda_{0}^{3}}-1|\lesssim\int_{0}^{t}\frac{b^{1+2\kappa}}{\lambda^{2}}d\tau\lesssim_{\kappa}b_{0}^{2\kappa},
\end{align*}
which finishes the proof of \eqref{eq:almost conserved quantity}. 
\end{proof}
\begin{proof}[Proof of Proposition~\ref{PropositionMainBootstrap}]
We will show the bootstrap conclusion \eqref{eq:5.2 Bootstrap prop goal}. We first verify $0<b(t)<2b_{0}$. By \eqref{eq:RefinedModulationBlowupLaw}
and \eqref{eq:modulationEstimate1}, we have $\widetilde{b}_{s}=-b^{2}(\frac{3}{2}+O(b^{0+}))$.
Thus, $\widetilde{b}(t)\leq\widetilde{b}_{0}$. By \eqref{eq:RefinedModulationDiffer}
and the initial bound \eqref{eq:5.2 ini boots assump}, we have 
\begin{align*}
b(t)(1-C(b^{*})^{\frac{1}{12}})\leq\widetilde{b}(t)\leq\widetilde{b}_{0}\leq b_{0}(1+C(b^{*})^{\frac{1}{12}}).
\end{align*}
Since $b_{*}\ll1$, we deduce 
\begin{align*}
0<b(t)\leq\Big(\frac{1+C(b^{*})^{\frac{1}{12}}}{1-C(b^{*})^{\frac{1}{12}}}\Big)b_{0}\leq1.1b_{0}.
\end{align*}

Next, we close $\frac{1}{1.2}\leq\frac{b^{2}}{\lambda^{3}}\leq1.2$.
Indeed, by Lemma~\ref{Lemmab-lambdaRelation} with the initial bound
\eqref{eq:5.2 ini boots assump}, we have 
\begin{align*}
\frac{1}{1.2}\leq\frac{b^{2}(t)}{\lambda^{3}(t)}=(1+O(b_{0}^{2\kappa}))\frac{b_{0}^{2}}{\lambda_{0}^{3}}\leq1.2.
\end{align*}

We show $\|\widehat{\eps}\|_{L^{2}}<\frac{1}{20}\delta_{dec}$. By
\eqref{eq:InterpolationEstimateW,inf}, we have 
\begin{align*}
\|w\|_{L^{2}}^{2}=\|Q+\widehat{\eps}\|_{L^{2}}^{2} & =\|Q\|_{L^{2}}^{2}+\|\widehat{\eps}\|_{L^{2}}^{2}+2(Q,\widehat{\eps})\\
 & =\|Q\|_{L^{2}}^{2}+\|\widehat{\eps}\|_{L^{2}}^{2}+O(\lambda^{\frac{1}{2}-}).
\end{align*}
Thanks to the mass conservation law and $b(t)\leq2b_{0}$, we have
\begin{align*}
\|\widehat{\eps}(t)\|_{L^{2}}^{2}\leq\|\widehat{\eps}_{0}\|_{L^{2}}^{2}+O(b_{0}^{\frac{1}{3}-}).
\end{align*}
By the initial bound \eqref{eq:5.2 ini boots assump} and $b^{*}\ll\delta_{dec}^{6+}$,
we obtain $\|\widehat{\eps}\|_{L^{2}}<\frac{1}{20}\delta_{dec}$.
This finishes the proof. 
\end{proof}

Now, we prove topological Lemmas~\ref{LemmaTopologicalPhiContinuous}
and \ref{LemmaPsiNonExistenceBrouwer}. 
\begin{proof}[Proof of Lemma~\ref{LemmaTopologicalPhiContinuous}]
It suffice to show the continuity of $\Phi$, as $\Psi$ is given
by \eqref{eq:5.2 Psi Phi relation}. For the convenience of notation,
we denote $\mathcal{V}_{0}=(\lambda_{0},\gamma_{0},x_{0},b_{0},\eta_{0},\nu_{0},\eps_{0})$.
We also denote $\widetilde{T}_{exit}^{\mathcal{V}_{0}}$ by the exit
time $\widetilde{T}_{exit}$ with the initial data $\mathcal{V}_{0}$.
By the definition of the refined modulation parameters \eqref{eq:RefinedModulationDefinition}
and the continuous dependence, we have the continuity of $\mathcal{V}_{0}\mapsto(\frac{10\widetilde{\eta}}{\lambda^{\frac{3}{2}(1+\kappa)}},\frac{10\widetilde{\nu}}{\lambda^{\frac{3}{2}(1-\kappa)}})^{\mathcal{V}_{0}}(t^{\prime})$
for each $0\leq t^{\prime}<T$. It suffices to show the continuity
of $(\frac{10\widetilde{\eta}}{\lambda^{\frac{3}{2}(1+\kappa)}},\frac{10\widetilde{\nu}}{\lambda^{\frac{3}{2}(1-\kappa)}})^{\mathcal{V}_{0}}(0)\mapsto\tilde T_{exit}^{\mathcal{V}_{0}}$.
When $(\frac{10\widetilde{\eta}}{\lambda^{\frac{3}{2}(1+\kappa)}},\frac{10\widetilde{\nu}}{\lambda^{\frac{3}{2}(1-\kappa)}})^{\mathcal{V}_{0}}(0)\in(-0.9,0.9)^{2}$,
we have $\widetilde{T}_{exit}^{\mathcal{V}_{0}}>0$. By the continuous
dependence of solutions and the continuity of decomposition map, we
can show the continuity of $\mathcal{V}_{0}\mapsto\widetilde{T}_{exit}^{\mathcal{V}_{0}}$
when $(\frac{10\widetilde{\eta}}{\lambda^{\frac{3}{2}(1+\kappa)}},\frac{10\widetilde{\nu}}{\lambda^{\frac{3}{2}(1-\kappa)}})^{\mathcal{V}_{0}}(0)\in(-0.9,0.9)^{2}$.
Hence, $\Phi$ is continuous in this case, interior region. On the
other hand, on exterior region, $(\frac{10\widetilde{\eta}}{\lambda^{\frac{3}{2}(1+\kappa)}},\frac{10\widetilde{\nu}}{\lambda^{\frac{3}{2}(1-\kappa)}})^{\mathcal{V}_{0}}(0)\in[-1,1]^{2}\setminus[-0.9,0.9]^{2}$,
we have $\widetilde{T}_{exit}^{\mathcal{V}_{0}}\equiv0$ by its definition
and so $\Phi$ is obviously continuous on it. Hence, we are left to
show the continuity of $\Phi$ at the boundary region, $(\frac{10\widetilde{\eta}}{\lambda^{\frac{3}{2}(1+\kappa)}},\frac{10\widetilde{\nu}}{\lambda^{\frac{3}{2}(1-\kappa)}})^{\mathcal{V}_{0}}(0)\in[-0.9,0.9]^{2}\setminus(-0.9,0.9)^{2}$.
For this case, we have to use the modulation estimates \eqref{eq:RefinedModulationBlowupLaw}
in the sense that the vector field $(\partial_{s}(\tfrac{\widetilde{\eta}}{\lambda^{\frac{3}{2}+}}),\partial_{s}(\tfrac{\widetilde{\nu}}{\lambda^{\frac{3}{2}-}}))^{\mathcal{V}_{0}}$
is pointing outward direction. If it were to happen that $(\frac{10\widetilde{\eta}}{\lambda^{\frac{3}{2}(1+\kappa)}},\frac{10\widetilde{\nu}}{\lambda^{\frac{3}{2}(1-\kappa)}})(t)$
does not immediately exit the boundary, and exit at the other side
by crossing through the interior of $(-0.9,0.9)^{2}$, then the continuity
of $\Phi$ would break down. We show that this indeed does not happen
by using the modulation estimates. We compute $\partial_{s}(\frac{\widetilde{\eta}}{\lambda^{\frac{3}{2}+}}),\partial_{s}(\frac{\widetilde{\nu}}{\lambda^{\frac{3}{2}-}})$
by 
\begin{align}
\partial_{s}(\tfrac{\widetilde{\eta}}{\lambda^{\frac{3}{2}+}}) & =\tfrac{\widetilde{\eta}}{\lambda^{\frac{3}{2}+}}\cdot b(\tfrac{1}{2}+\kappa+O(b^{\frac{1}{12}-})),\label{eq:TopologicalArgumentTimeDerivativeEta}\\
\partial_{s}(\tfrac{\widetilde{\nu}}{\lambda^{\frac{3}{2}-}}) & =\tfrac{\widetilde{\nu}}{\lambda^{\frac{3}{2}+}}\cdot b(\tfrac{1}{2}-\kappa+O(b^{\frac{1}{12}-})).\label{eq:TopologicalArgumentTimeDerivativeNu}
\end{align}
We note that this implies that $\frac{\widetilde{\eta}}{\lambda^{\frac{3}{2}+}}$
or $\frac{\widetilde{\nu}}{\lambda^{\frac{3}{2}-}}$ increases(decreases)
when $\frac{\widetilde{\eta}}{\lambda^{\frac{3}{2}+}}>0(<0)$ or $\frac{\widetilde{\nu}}{\lambda^{\frac{3}{2}-}}>0(<0)$,
respectively since $0<b_{0}\ll1.$ In other words, the vector field
point uniformly outgoing near the boundary and so it exit immediately,
i.e., $\tilde T_{exit}^{\mathcal{V}_{0}}=0.$ If $(\frac{10\widetilde{\eta}_{0}}{\lambda_{0}^{\frac{3}{2}+}},\frac{10\widetilde{\nu}_{0}}{\lambda_{0}^{\frac{3}{2}-}})^{\mathcal{V}_{0}}$
is on the boundary, $[-0.9,0.9]^{2}\setminus(-0.9,0.9)^{2}$, it is
immediate we have one of the following four situations: 
\begin{align*}
\tfrac{10\widetilde{\eta}_{0}}{\lambda_{0}^{\frac{3}{2}+}}=\pm0.9,\quad\tfrac{10\widetilde{\nu}_{0}}{\lambda_{0}^{\frac{3}{2}-}}=\pm0.9.
\end{align*}
We suppose that $\frac{10\widetilde{\eta}_{0}}{\lambda_{0}^{\frac{3}{2}+}}=0.9$.
Let $\mathcal{V}_{0}^{\prime}\in B_{\epsilon_{0}}(\mathcal{V}_{0})$
and $((\frac{10\widetilde{\eta}_{0}}{\lambda_{0}^{\frac{3}{2}+}})^{\mathcal{V}_{0}^{\prime}},(\frac{10\widetilde{\nu}_{0}}{\lambda_{0}^{\frac{3}{2}-}})^{\mathcal{V}_{0}^{\prime}})\in(-0.9,0.9)^{2}$
for a sufficiently small $\epsilon_{0}>0$. From \eqref{eq:TopologicalArgumentTimeDerivativeEta}
and \eqref{eq:TopologicalArgumentTimeDerivativeNu}, we have 
\begin{align}
2\left\{ \partial_{s}\left(\tfrac{\widetilde{\eta}}{\lambda^{\frac{3}{2}+}}\right)^{\mathcal{V}_{0}^{\prime}}\right\} (t)\geq\left|\left\{ \partial_{s}\left(\tfrac{\widetilde{\nu}}{\lambda^{\frac{3}{2}-}}\right)^{\mathcal{V}_{0}^{\prime}}\right\} (t)\right|.\label{eq:TopologicalArgumentTimeDerivative}
\end{align}
By integrating \eqref{eq:TopologicalArgumentTimeDerivative} from
$0$ to $\widetilde{T}_{exit}^{\mathcal{V}_{0}^{\prime}}$, we deduce
\begin{align}
2\left(\left(\tfrac{\widetilde{\eta}}{\lambda^{\frac{3}{2}+}}\right)^{\mathcal{V}_{0}^{\prime}}(\widetilde{T}_{exit}^{\mathcal{V}_{0}^{\prime}})-\tfrac{\widetilde{\eta}_{0}^{\prime}}{(\lambda_{0}^{\prime})^{\frac{3}{2}+}}\right)\geq\left|\left(\tfrac{\widetilde{\nu}}{\lambda^{\frac{3}{2}-}}\right)^{\mathcal{V}_{0}^{\prime}}(\widetilde{T}_{exit}^{\mathcal{V}_{0}^{\prime}})-\tfrac{\widetilde{\nu}_{0}^{\prime}}{(\lambda_{0}^{\prime})^{\frac{3}{2}-}}\right|.\label{eq:TopologicalArgumentTimeDerivativeEstimate}
\end{align}
For $(\frac{10\widetilde{\eta}_{0}^{\prime}}{\lambda_{0}^{\frac{3}{2}+}},\frac{10\widetilde{\nu}_{0}^{\prime}}{\lambda_{0}^{\frac{3}{2}-}})\in B_{\epsilon_{0}}(\frac{10\widetilde{\eta}_{0}}{\lambda_{0}^{\frac{3}{2}+}},\frac{10\widetilde{\nu}_{0}}{\lambda_{0}^{\frac{3}{2}-}})$,
$\left(\frac{\widetilde{\eta}}{\lambda^{\frac{3}{2}+}}\right)^{\mathcal{V}_{0}^{\prime}}$
increases in time by \eqref{eq:TopologicalArgumentTimeDerivativeEta}.
Therefore, we have 
\begin{align}
10\times\text{LHS of }\eqref{eq:TopologicalArgumentTimeDerivativeEstimate}\leq2\Big(0.9-\tfrac{10\widetilde{\eta}_{0}^{\prime}}{(\lambda_{0}^{\prime})^{\frac{3}{2}+}}\Big)<2\epsilon_{0}.\label{eq:TopologicalArgumentTimeDerivativeEstimate2}
\end{align}
Thus, by \eqref{eq:TopologicalArgumentTimeDerivativeEstimate} and
\eqref{eq:TopologicalArgumentTimeDerivativeEstimate2}, we conclude
\begin{align*}
\bigg(\left(\tfrac{\widetilde{\eta}}{\lambda^{\frac{3}{2}+}}\right)^{\mathcal{V}_{0}^{\prime}}(\widetilde{T}_{exit}^{\mathcal{V}_{0}^{\prime}}),\left(\tfrac{\widetilde{\nu}}{\lambda^{\frac{3}{2}-}}\right)^{\mathcal{V}_{0}^{\prime}}(\widetilde{T}_{exit}^{\mathcal{V}_{0}^{\prime}})\bigg)\in B_{3\eps}\left(\tfrac{\widetilde{\eta}_{0}}{\lambda_{0}^{\frac{3}{2}+}},\tfrac{\widetilde{\nu}_{0}}{\lambda_{0}^{\frac{3}{2}-}}\right).
\end{align*}
Since the map $\mathcal{V}_{0}\mapsto(\widetilde{\eta}_{0},\widetilde{\nu}_{0})$
is continuous, we deduce that $\Phi$ is continuous at $\mathcal{V}_{0}=(\lambda_{0},\gamma_{0},x_{0},b_{0},\eta_{0},\nu_{0},\eps_{0})$
which satisfies $\frac{10\widetilde{\eta}_{0}}{\lambda_{0}^{\frac{3}{2}+}}=0.9$.
Similarly, we can deduce the continuity when $\mathcal{V}_{0}$ satisfies
$\frac{10\widetilde{\eta}_{0}}{\lambda_{0}^{\frac{3}{2}+}}=-0.9$
or $\frac{10\widetilde{\nu}_{0}}{\lambda_{0}^{\frac{3}{2}-}}=\pm0.9$.
Therefore, $\Phi$ is continuous on the boundary of $(-0.9,0.9)^{2}$,
and we conclude that $\Phi$ is continuous on $\mathcal{I}$. 
\end{proof}
Now, we prove the nonexistence of the above $\Psi$. 
\begin{proof}[Proof of Lemma~\ref{LemmaPsiNonExistenceBrouwer}]
We note that $S$ is a square, $S=\{(x,y):|x+y|+|x-y|=2\}$. Define
$P:\mathbb{R}^{2}\setminus\{0\}\to S$ as a projection function by
$P(x,y)\coloneqq$the intersection of the square $S$ and the ray
that comes from the origin passing through $(x,y)$. Then, we deduce $P\circ\Psi:[-1,1]^{2}\to S$
is a continuous function since both $\Psi$ and $P$ are continuous.
In addition, according to the definition of $\Psi$ and \eqref{eq:RefinedModulationDiffer},
for any $(x,y)\in S$, $\Psi(x,y)\in B_{0.2}(x,y)$. Thus, we have
$P\circ\Psi(x,y)\in S\cap B_{0.2}(x,y)$. Now we apply the Brouwer
fixed point theorem to $-P\circ\Psi$. Since $-P\circ\Psi:[-1,1]^{2}\to S\subset[-1,1]^{2}$
and $-P\circ\Psi$ are continuous. By the Brouwer fixed point theorem,
we can find a fixed point $(x_{0},y_{0})\in[-1,1]^{2}$. Since $-P\circ\Psi:[-1,1]^{2}\to S$,
$(x_{0},y_{0})\in S$. However, we have $-P\circ\Psi(x_{0},y_{0})\in B_{0.2}(-x_{0},-y_{0})$,
and this implies that $-P\circ\Psi(x_{0},y_{0})\neq(x_{0},y_{0}).$
Hence, this is a contradiction, and we conclude such $\Psi$ does
not exist. 
\end{proof}
Next, we prove the sharp description of blow-up solutions. 
\begin{proof}[Proof of Proposition~\ref{PropSharpDescription}]
Due to the bootstrap conclusion, the trapped solution $v(t)$ satisfies
\eqref{eq:5.2 Bootstrap prop goal} for $0\le t<T$. We first observe
that 
\begin{align*}
\left(\frac{\lambda(t)}{\lambda_{0}}\right)^{\frac{8}{5}}\leq\frac{\widetilde{b}(t)}{\widetilde{b}_{0}},
\end{align*}
because using \eqref{eq:modulationEstimate1} and Lemma~\ref{LemmaModulationEstimateRefined}
we have 
\begin{align*}
\partial_{s}\log\frac{\widetilde{b}}{\lambda^{\frac{8}{5}}}=\frac{8}{5}\widetilde{b}+\frac{\widetilde{b}_{s}}{\widetilde{b}}-\frac{8}{5}\left(\frac{\lambda_{s}}{\lambda}+\widetilde{b}\right)=\left(\frac{1}{10}+O(b^{0+})\right)b\geq0.
\end{align*}
Then, we estimate 
\begin{align*}
\partial_{t}\lambda^{\frac{2}{5}}=-\frac{2}{5}\frac{b}{\lambda^{\frac{8}{5}}}+\frac{2}{5}\frac{1}{\lambda^{\frac{8}{5}}}\left(\frac{\lambda_{s}}{\lambda}+b\right)=-\frac{2}{5}\frac{\widetilde{b}}{\lambda^{\frac{8}{5}}}\left(1+O((b^{*})^{\frac{1}{12}})\right)\leq-\frac{1}{5}\frac{\widetilde{b}_{0}}{\lambda_{0}^{\frac{8}{5}}},
\end{align*}
and verify that the solution $v(t,x)$ blows up in finite time, i.e.,
$T<+\infty$. By a standard blow-up criterion, we have $\lim_{t\to T}\|v(t)\|_{\dot{H}^{1}}=\infty$
and $\lambda(T)\coloneqq\lim_{t\to T}\lambda(t)=0$. By Lemma~\ref{Lemmab-lambdaRelation}
($b^{2}\sim\lambda^{3}$), $|\eta|<b^{1+}$, and $|\nu|<b^{1-}$,
we have $b(T)\coloneqq\lim_{t\to T}b(t)=0$, $\eta(T)\coloneqq\lim_{t\to T}\eta(t)=0$,
and $\nu(T)\coloneqq\lim_{t\to T}\nu(t)=0$. From the proof of Lemma~\ref{Lemmab-lambdaRelation},
we have 
\begin{align*}
\frac{b^{2}(t)}{\lambda^{3}(t)}=\ell(1+O(b^{0+})),\quad\text{with}\quad\ell\coloneqq\lim_{t\to T}\frac{b^{2}(t)}{\lambda^{3}(t)}\in(0,\infty).
\end{align*}
We claim the asymptotics of the parameters $\lambda$ and $b$: 
\begin{align}
\begin{split}\lambda(t)=\ell\cdot(T-t)^{2}(1+o_{t\to T}(1)),\\
b(t)=\ell^{2}\cdot(T-t)^{3}(1+o_{t\to T}(1)).
\end{split}
\label{eq:asymptoticModulationParam}
\end{align}
From \eqref{eq:RefinedModulationBlowupLaw}, we have 
\begin{align*}
\partial_{t}\widetilde{b}=\frac{\widetilde{b}_{s}+\frac{3}{2}\widetilde{b}^{2}}{\lambda^{2}}-\frac{3}{2}\frac{\widetilde{b}^{2}}{\lambda^{2}}=-\left(\frac{3}{2}+O(b^{0+})\right)\frac{\widetilde{b}^{2}}{\lambda^{2}}=-\ell^{\frac{2}{3}}\left(\frac{3}{2}+O(b^{0+})\right)b^{\frac{2}{3}},
\end{align*}
and we deduce 
\begin{align*}
b(t)=\ell^{2}\cdot(T-t)^{3}(C+o_{t\to T}(1))\quad\text{for some }C>0.
\end{align*}
Using this and \eqref{eq:modulationEstimate1}, 
\begin{align*}
\partial_{t}(\lambda^{\frac{1}{2}})=\frac{\lambda_{t}}{\lambda^{\frac{1}{2}}}=-\frac{b}{\lambda^{\frac{3}{2}}}(1+O(b^{\frac{1}{3}}))=-\sqrt{l}(1+O(b^{\frac{1}{3}})),
\end{align*}
and this implies 
\begin{align*}
\lambda(t)=\ell\cdot(T-t)^{2}(1+o_{t\to T}(1)).
\end{align*}
By the definition of $\ell$, we have $C=1$ and we conclude
\eqref{eq:asymptoticModulationParam}.

By \eqref{eq:modulationEstimate1+}, we have 
\begin{align*}
|\gamma_{s}|\lesssim b^{1+},\qquad|x_{s}|\lesssim b^{\frac{5}{3}-}.
\end{align*}
Applying \eqref{eq:asymptoticModulationParam}, we obtain 
\begin{equation}
|\gmm_{t}|\aleq_{\ell}(T-t)^{-1+},\qquad|x_{t}|\aleq_{\ell}(T-t)^{1-}.\label{eq:gmm-x-deriv-in-time}
\end{equation}
Hence $\gamma(t)$ and $x(t)$ converge to some $\gamma^{*}$ and
$x^{*}$ as $t\to T$.

Finally, we prove that the solution $v$ decomposes as in Theorem~\ref{thm:precise statement codimensionone blowup}
. We first claim the outer $L^{2}$-convergence. Fix large $R>0.$
We have 
\begin{align*}
\|(1-\chi_{R})L_{v}^{*}\mathbf{D}_{v}v\|_{L^{2}}= & \lambda^{-2}(t)\|(1-\chi_{R\lambda^{-1}(t)})L_{w}^{*}w_{1}\|_{L^{2}}\\
\leq & \lambda^{-2}(t)(\|(1-\chi_{R\lambda^{-1}(t)})(L_{w}^{*}w_{1}-L_{Q}^{*}P_{1})\|_{L^{2}}\\
 & +\|(1-\chi_{R\lambda^{-1}(t)})L_{Q}^{*}P_{1}\|_{L^{2}})\\
\lesssim & -\frac{\lambda_{t}}{\lambda^{\frac{1}{2}}}\lesssim1.
\end{align*}
Thus, $t\mapsto\|(1-\chi_{R})L_{v}^{*}\mathbf{D}_{v}v(t)\|_{L^{2}}$
is integrable for $0\le t<T$, which means that $(1-\chi_{R})v(t)$
converges in $L^{2}$ in view of $i\partial_{t}(\chi_{R}v)=\chi_{R}L_{v}^{*}\mathbf{D}_{v}v$.
Therefore, there exists a function $v^{*}$ such that $(1-\chi_{R})v^{*}\in L^{2}$
and $(1-\chi_{R})v(t)\to(1-\chi_{R})v^{*}$ in $L^{2}$ for any $R>0$.

In view of 
\begin{align*}
v(t,x) & =[Q+P(b(t),\eta(t),\nu(t))+\eps]_{\lambda(t),\gamma(t),x(t)}\\
 & =[Q+\widehat{\eps}]_{\lambda(t),\gamma(t),x(t)},
\end{align*}
and from the asymptotics, $\frac{\ell(T-t)^{2}}{\lambda(t)}\to1$,
$\gamma\to\gamma^{*}$, and $x\to x^{*}$ as $t\to T$, we have 
\begin{align*}
[Q]_{\lambda(t),\gamma(t),x(t)}-[Q]_{\ell(T-t)^{2},\gamma^{*},x^{*}}\to0\quad\text{in}\quad L^{2}.
\end{align*}
Therefore, it suffices to show $[\widehat{\eps}]_{\lambda(t),\gamma(t),x(t)}\to v^{*}$
in $L^{2}$ as $t\to T$ and $v^{*}\in H^{1}$. We have 
\begin{align*}
\|[\widehat{\eps}]_{\lambda(t),\gamma(t),x(t)}\|_{\dot{H}^{1}}=\lambda^{-1}\|\widehat{\eps}\|_{\dot{H}^{1}}\lesssim1.
\end{align*}
Since $\|\widehat{\eps}\|_{L^{2}}<\frac{1}{20}\delta_{dec}$, $[\widehat{\eps}]_{\lambda(t),\gamma(t),x(t)}$
is bounded in $H^{1}$, and we have $[\widehat{\eps}]_{\lambda(t),\gamma(t),x(t)}\rightharpoonup v^{*}$
weakly in $H^{1}$ due to the uniqueness of the limit. Thus, we have $v^{*}\in H^{1}$.
By the Rellich--Kondrachov compactness theorem, $[\widehat{\eps}]_{\lambda(t),\gamma(t),x(t)}\to v^{*}$
in $L_{\textnormal{loc}}^{2}$. Thus, combining with the outer $L^{2}$-convergence,
we conclude that $[\widehat{\eps}]_{\lambda(t),\gamma(t),x(t)}\to v^{*}$
in $L^{2}$. This finishes the proof. 
\end{proof}

\subsection{\label{SubsectionChiralBlowup}Chiral blow-up solutions}

In this subsection, we prove the construction of chiral blow-up solutions,
Theorem~\ref{TheoremChiralBlowup}. The core analysis of the proof
is similar to that of Theorem~\ref{thm:precise statement codimensionone blowup}
in the sense that we look for a trapped solution. Keeping the bootstrap
proposition, we show that the trapped solution can be constructed
from a chiral initial data set, possibly of finite dimension. The
new ingredient is to construct a \emph{chiral} modified profile of
$\mathcal{R}$.

The first step is to construct a chiral modified profile $\mathcal{R}_{b,\eta,\nu}$
that well approximates the blow-up profile in the gauge transformed
side, i.e., $\mathcal{R}_{b,\eta,\nu}\approx\mathcal{R}$ and $\mathcal{-G}(\mathcal{R}_{b,\eta,\nu})\approx Q$
or $-\mathcal{G}(\mathcal{R}_{b,\eta,\nu})\approx Q+P$ depending
on topologies. We first modify the multiplication operator by $x$
to fit in the Hardy space. Since $C_{c}^{\infty}$ functions are not
in the Hardy space in general, simply taking smooth cut-offs is not
an appropriate localization. Let $R>0$ be a large constant to be
chosen later. Define $\omega_{R}$ by 
\begin{align*}
\omega_{R}(x)\coloneqq\frac{R}{i}(e^{\frac{ix}{R}}-1).
\end{align*}
It is easy to check that 
\[
\lim_{R\to+\infty}\omega_{R}(x)=x\text{ locally,}\quad\text{\ensuremath{\omega_{R}(x):L_{+}^{2}\to L_{+}^{2}},}\quad|\omega_{R}(y)|\lesssim\min(|y|,R).
\]
Define the chiral profiles $\mathcal{R}_{b,\eta,\nu}(x)$ by 
\begin{align*}
\mathcal{R}_{b,\eta,\nu}(x)\coloneqq e^{i\phi(b,\eta,\nu)}\mathcal{R}(x)(1-\eta\tfrac{1+\omega_{R}^{2}}{4}-ib\tfrac{\omega_{R}^{2}}{4}-i(\eta-\nu)\tfrac{\omega_{R}}{2})\in L_{+}^{2},
\end{align*}
where $\phi$ is a phase correction required for the gauge transform
$\mathcal{G}$, 
\begin{align*}
\phi(b,\eta,\nu)\coloneqq\frac{1}{2}\int_{-\infty}^{0}(|\mathcal{R}_{b,\eta,\nu}|^{2}-Q^{2}).
\end{align*}
We justify that these profiles $\mathcal{G}(\mathcal{R}_{b,\eta,\nu})(x)$
are appropriate blow-up profiles for \eqref{CMdnls-gauged} by establishing
the following scaling bounds. 
\begin{lem}
\label{lem:5.20}Let $b,\eta,\nu$ and $R$ satisfy $R>10$ and $(|b|+|\eta|)R^{\frac{3}{2}}\lesssim1$.
Then, we have 
\begin{align}
\|\mathcal{G}(\mathcal{R}_{b,\eta,\nu})+Q\|_{L^{2}} & \lesssim(|b|+|\eta|)R^{\frac{3}{2}}+|\eta-\nu|R^{\frac{1}{2}},\label{eq:HardyLocalizationL2LemmaEq}\\
\|\mathcal{G}(\mathcal{R}_{b,\eta,\nu})+Q\|_{\dot{\mathcal{H}}^{1}} & \lesssim(|b|+|\eta|)R^{\frac{1}{2}}+|\eta-\nu|,\label{eq:HardyLocalizationH1LemmaEq}\\
\|\mathcal{G}(\mathcal{R}_{b,\eta,\nu})+Q+P\|_{\dot{\mathcal{H}}^{2}} & \lesssim(|b|+|\eta|)R^{-\frac{1}{2}}+|\nu|R^{-1}.\label{eq:HardyLocalizationH2LemmaEq}
\end{align}
\end{lem}

\begin{proof}
We take the gauge transform $\mathcal{G}$ to obtain 
\begin{align*}
-\mathcal{G}(\mathcal{R}_{b,\eta,\nu})=Q(1-\eta\tfrac{1+\omega_{R}^{2}}{4}-ib\tfrac{\omega_{R}^{2}}{4}-i(\eta-\nu)\tfrac{\omega_{R}}{2})e^{-\frac{i}{2}\int_{0}^{x}(|\mathcal{R}_{b,\eta,\nu}|^{2}-Q^{2})dy},
\end{align*}
and 
\begin{align}
-(\mathcal{G}(\mathcal{R}_{b,\eta,\nu})+Q)= & Q(e^{-\frac{i}{2}\int_{0}^{x}(|\mathcal{R}_{b,\eta,\nu}|^{2}-Q^{2})dy}-1)\nonumber \\
 & +Q(-\eta\tfrac{1+\omega_{R}^{2}}{4}-ib\tfrac{\omega_{R}^{2}}{4}-i(\eta-\nu)\tfrac{\omega_{R}}{2})e^{-\frac{i}{2}\int_{0}^{x}(|\mathcal{R}_{b,\eta,\nu}|^{2}-Q^{2})dy}.\label{eq:HardyLocalization G Q decompse}
\end{align}
On the other hand, we have 
\begin{align*}
|\mathcal{R}_{b,\eta,\nu}|^{2}-Q^{2}= & Q^{2}(\Re(-\eta\tfrac{1+\omega_{R}^{2}}{2})+\Im(b\tfrac{\omega_{R}^{2}}{2}+(\eta-\nu)\omega_{R}))\\
 & +|\eta\tfrac{1+\omega_{R}^{2}}{4}+ib\tfrac{\omega_{R}^{2}}{4}+i(\eta-\nu)\tfrac{\omega_{R}}{2}|^{2})\\
\lesssim & Q^{2}(|\eta|\min(\langle x\rangle^{2},R^{2})+|b|\min(\tfrac{|x|^{3}}{R},R^{2})+|\eta-\nu|\min(\tfrac{|x|^{2}}{R},R)\\
 & +(|b|^{2}+|\eta|^{2})\min(\langle x\rangle^{4},R^{4})+|\nu|^{2}\min(\langle x\rangle^{2},R^{2})),
\end{align*}
where we used $|\Im(\omega_{R}^{2})|\lesssim\min(\tfrac{|x|^{3}}{R},R^{2})$.
Thus, we have 
\begin{align}
\int_{0}^{x}(|\mathcal{R}_{b,\eta,\nu}|^{2}-Q^{2})\lesssim & |\eta|\min(|x|,R)+|b|\min(\tfrac{|x|^{2}}{R},R)+|\eta-\nu|\min(\tfrac{|x|}{R},1)\nonumber \\
 & +(|b|^{2}+|\eta|^{2})\min(\langle x\rangle^{3},R^{3})+|\nu|^{2}\min(|x|,R).\label{eq:HardyLocalizationGaugeEstimate}
\end{align}
By \eqref{eq:HardyLocalization G Q decompse} and \eqref{eq:HardyLocalizationGaugeEstimate},
we deduce \eqref{eq:HardyLocalizationL2LemmaEq} and \eqref{eq:HardyLocalizationH1LemmaEq}.

Now we want to show \eqref{eq:HardyLocalizationH2LemmaEq}. We further
calculate to estimate $\dot{\mathcal{H}}^{2}$-level. 
\begin{align}
- & (\mathcal{G}(\mathcal{R}_{b,\eta,\nu})+Q+P)\nonumber \\
= & Q(e^{-\frac{i}{2}\int_{0}^{x}(|\mathcal{R}_{b,\eta,\nu}|^{2}-Q^{2})dy}-1-i\eta\frac{\omega_{R}}{2})\label{eq:HardyLocalizationH2FirstLine}\\
 & +Q(-\eta\tfrac{\omega_{R}^{2}-x^{2}}{4}-ib\tfrac{\omega_{R}^{2}-x^{2}}{4}+i\nu\tfrac{\omega_{R}-x}{2})\label{eq:HardyLocalizationH2SecondLine}\\
 & +Q(-\eta\tfrac{1+\omega_{R}^{2}}{4}-ib\tfrac{\omega_{R}^{2}}{4}-i(\eta-\nu)\tfrac{\omega_{R}}{2})(e^{-\frac{i}{2}\int_{0}^{x}(|\mathcal{R}_{b,\eta,\nu}|^{2}-Q^{2})dy}-1).\label{eq:HardyLocalizationH2ThirdLine}
\end{align}
In order to estimate \eqref{eq:HardyLocalizationH2FirstLine}, using \eqref{eq:HardyLocalizationGaugeEstimate}, we obtain  
\begin{align*}
\bigg|\int_{0}^{x}\big||\mathcal{R}_{b,\eta,\nu}|^{2}-Q^{2}\big|dy & -\eta\omega_{R}\bigg|\\
\lesssim & \eta R\textbf{1}_{|x|\gtrsim R}+|b|\min(\tfrac{|x|^{2}}{R},R)+|\eta-\nu|\min(\tfrac{|x|}{R},1)\\
 & +(|b|^{2}+|\eta|^{2})\min(\langle x\rangle^{3},R^{3})+|\nu|^{2}\min(|x|,R).
\end{align*}
Therefore, we deduce 
\begin{align*}
\|\eqref{eq:HardyLocalizationH2FirstLine}\|_{\dot{\mathcal{H}}^{2}}\lesssim(|b|+|\eta|)R^{-\frac{1}{2}}+|\eta-\nu|R^{-1}
\end{align*}
For \eqref{eq:HardyLocalizationH2SecondLine}, we use $|\omega_{R}^{2}-x^{2}|\lesssim\tfrac{|x|^{3}}{R^{2}}$
and $|\omega_{R}-x|\lesssim\tfrac{|x|^{2}}{R}$ for $\tfrac{|x|}{R}<\tfrac{1}{2}$
to estimate 
\begin{align*}
\|\eqref{eq:HardyLocalizationH2SecondLine}\|_{\dot{\mathcal{H}}^{2}}\lesssim(|b|+|\eta|)R^{-\frac{1}{2}}+|\nu|R^{-1}.
\end{align*}
From \eqref{eq:HardyLocalizationGaugeEstimate}, we also deduce 
\begin{align*}
\|\eqref{eq:HardyLocalizationH2ThirdLine}\|_{\dot{\mathcal{H}}^{2}}\lesssim(|b|^{2}+|\eta|^{2})R^{\frac{1}{2}}+|\eta-\nu|R^{-1}.
\end{align*}
This completes the proof. 
\end{proof}
We now fix a universal constant $\delta_{chi}$ such that $\delta_{chi}\ll\min(\epsilon,\delta_{dec})$,
say $\delta_{chi}=\frac{1}{100}\min(\epsilon,\delta_{dec})$, where
$\epsilon$ is given in Theorem~\ref{TheoremChiralBlowup}. Then,
we fix the parameter 
\[
R=\delta_{chi}^{\frac{2}{3}}\lambda^{-1}.
\]

As we did in Section~\ref{Subsection ExistenceTrappedSol}, we specify
the range of parameters for the initial data. Let $\mathcal{U}_{chi}\subset\mathbb{R}_{+}\times\mathbb{R}/2\pi\mathbb{Z}\times\mathbb{R}^{4}$
be a set of $(\mathring{\la}_{0},\mathring{\ga}_{0},\mathring{x}_{0},\mathring{b}_{0},\mathring{\eta}_{0},\mathring{\nu}_{0})$
satisfying 
\begin{align}
0<\mathring{b}_{0}<\tfrac{1}{3}b_{chi}^{*},\quad|\mathring{\eta}_{0}|<\tfrac{1}{3}\mathring{b}_{0}^{1+\kappa},\quad|\mathring{\nu}_{0}|<\tfrac{1}{3}\mathring{b}_{0}^{1-\kappa},\quad\mathring{b}_{0}^{2}=\mathring{\lambda}_{0}^{3},\label{eq:InitialDataUchi Chiral}
\end{align}
where $b_{chi}^{*}\ll1$ will be determined later. Note that we have
fixed the ratio of $\mathring{b}_{0}\text{ and }\mathring{\lambda}_{0}$
by $\mathring{b}_{0}^{2}=\mathring{\lambda}_{0}^{3}$. Under the conditions
\eqref{eq:InitialDataUchi Chiral}, Lemma~\ref{lem:5.20} reads 
\begin{align}
\begin{split}\|\mathcal{G}(\mathcal{R}_{\mathring{b}_{0},\mathring{\eta}_{0},\mathring{\nu}_{0}})+Q\|_{L^{2}} & \lesssim\delta_{chi},\\
\|\mathcal{G}(\mathcal{R}_{\mathring{b}_{0},\mathring{\eta}_{0},\mathring{\nu}_{0}})+Q\|_{\dot{\mathcal{H}}^{1}} & \lesssim\delta_{chi}^{\frac{1}{3}}\mathring{\la}_{0},\\
\|\mathcal{G}(\mathcal{R}_{\mathring{b}_{0},\mathring{\eta}_{0},\mathring{\nu}_{0}})+Q+P\|_{\dot{\mathcal{H}}^{2}} & \lesssim\delta_{chi}^{-\frac{1}{3}}\mathring{\la}_{0}^{2}.
\end{split}
\label{eq: GR+Q sim zero}
\end{align}

Although $\mathcal{R}_{b,\eta,\nu}\in H_{+}^{\infty}$ is smooth,
it does not have rapid spatial decay as $\mathcal{R}_{b,\eta,\nu}\notin\langle x\rangle^{-1}L^{2}$.
In order to find a blow-up solution with initial data in $\mathcal{S}_{+}$,
we carefully choose $\mathring{\eps}_{0}$ so that the initial data
$\mathcal{R}_{\mathring{b}_{0},\mathring{\eta}_{0},\mathring{\nu}_{0}}+\mathring{\eps}_{0}$
belongs to $\mathcal{S}_{+}$. 
\begin{lem}
\label{Lemma 5.9 schwartz initial} Suppose \eqref{eq:InitialDataUchi Chiral}.
There exists $\mathring{\eps}_{0}=\mathring{\eps}_{0}(\cdot;\mathring{b}_{0},\mathring{\eta}_{0},\mathring{\nu}_{0})\in H_{+}^{\infty}$
such that 
\begin{align*}
u_{0}\coloneqq\mathcal{R}_{\mathring{b}_{0},\mathring{\eta}_{0},\mathring{\nu}_{0}}+\mathring{\eps}_{0}\in\mathcal{S}_{+}\quad\text{with}\quad\|\mathring{\eps}_{0}\|_{H^{2}}\lesssim\mathring{b}_{0}^{2}.
\end{align*}
Moreover, the mapping $(\mathring{b}_{0},\mathring{\eta}_{0},\mathring{\nu}_{0})\mapsto u_{0}(\cdot;\mathring{b}_{0},\mathring{\eta}_{0},\mathring{\nu}_{0})$
is continuous. 
\end{lem}

In the sequel, we fix $u_{0}=u_{0}(\cdot;\mathring{b}_{0},\mathring{\eta}_{0},\mathring{\nu}_{0})\in\mathcal{S}_{+}$
for each $(\mathring{b}_{0},\mathring{\eta}_{0},\mathring{\nu}_{0})$
in order to find a smooth chiral blow-up solution. 
\begin{proof}[Proof of Lemma~\ref{Lemma 5.9 schwartz initial}]
For simplicity of notation, we write $\la,b,\eta,\nu$ for $\mathring{\la}_{0},\mathring{b}_{0},\mathring{\eta}_{0},\mathring{\nu}_{0}$.
We compute the Fourier transform of the soliton $\mathcal{R}$ as
\begin{align}
\mathcal{F}(\mathcal{R})(\xi)=-i2\sqrt{2}\pi e^{-\xi}\mathbf{1}_{\xi\geq0},\label{eq:R fourier}
\end{align}
and since \eqref{eq:R fourier} is exponentially decaying as $\xi\to\infty$,
it suffices to mollify the function near $\xi=0$. We can find a smooth
step function $\mathring{\chi}_{\lambda}$ supported on $(0,\infty)$
and satisfying 
\begin{align*}
\mathring{\chi}_{\lambda}(\xi)=1\quad\text{on}\quad(\lambda^{4},\infty),\quad\mathring{\chi}_{\lambda}(\xi)=\phi(\frac{x-\lambda^{4}}{\lambda^{4}})\quad\text{on}\quad(0,\lambda^{4}),
\end{align*}
where $\phi\in C_{c}^{\infty}$ is supported on $(-1,1)$. Then, we
have $\mathring{\chi}_{\lambda}(\xi)\wt{\mathcal{R}}(\xi)\in\mathcal{S}$
and $\text{supp}\mathring{\chi}_{\lambda}(\xi)\wt{\mathcal{R}}(\xi)\subset[0,\infty)$.
If $f\in H^{2}\cap\widehat{L}^{\infty}$, then we have 
\begin{align*}
\|\mathcal{F}^{-1}[(1-\mathring{\chi}_{\lambda})\mathcal{F}(f)]\|_{H^{2}}\lesssim\lambda^{5}\|\widehat{f}\|_{L^{\infty}}^{2}.
\end{align*}
We use a conventional notation for $f\in L^{2}$, 
\begin{align*}
f*\widehat{\mathring{\chi}_{\lambda}}\coloneqq\mathcal{F}^{-1}[\mathring{\chi}_{\lambda}\mathcal{F}(f)].
\end{align*}
From \eqref{eq:R fourier}, we have 
\begin{align*}
\|\mathcal{R}*\widehat{\mathring{\chi}_{\lambda}}-\mathcal{R}\|_{H^{2}}\lesssim\lambda^{4},\quad\mathcal{R}*\widehat{\mathring{\chi}_{\lambda}}\in\mathcal{S}\cap L_{+}^{2}.
\end{align*}
Now, we define $\mathring{\eps}$ by 
\begin{align}
\mathring{\eps}\coloneqq e^{i\phi(b,\eta,\nu)}(\mathcal{R}*\widehat{\mathring{\chi}_{\lambda}}-\mathcal{R})(x)(1-\eta\tfrac{1+\omega_{R}^{2}}{4}-ib\tfrac{\omega_{R}^{2}}{4}-i(\eta-\nu)\tfrac{\omega_{R}}{2}).\label{eq:def good epsilon}
\end{align}
Then, we have 
\begin{align*}
\mathcal{R}_{b,\eta,\nu}+\mathring{\eps}=e^{i\phi(b,\eta,\nu)}(\mathcal{R}*\widehat{\mathring{\chi}_{\lambda}})(x)(1-\eta\tfrac{1+\omega_{R}^{2}}{4}-ib\tfrac{\omega_{R}^{2}}{4}-i(\eta-\nu)\tfrac{\omega_{R}}{2})\in\mathcal{S}\cap L_{+}^{2}.
\end{align*}

Now, we prove $\|\mathring{\eps}\|_{H^{2}}\lesssim\mathring{b}^{2}$.
Using $\mathcal{F}(\omega_{R}f)=\tfrac{R}{i}(\widehat{f}(\xi-R^{-1})-\widehat{f}(\xi))$,
we have 
\begin{align*}
\mathcal{F}[\omega_{R}(\mathcal{R}*\widehat{\mathring{\chi}_{\lambda}})]=\tfrac{R}{i}[\widehat{\mathcal{R}}(\mathring{\chi}_{\lambda})(\xi-R^{-1})-\widehat{\mathcal{R}}(\mathring{\chi}_{\lambda})(\xi)],
\end{align*}
and then we have 
\begin{align*}
\|\omega_{R}(\mathcal{R}*\widehat{\mathring{\chi}_{\lambda}})\|_{H^{2}}\lesssim R\lambda^{4}\sim\lambda^{3}.
\end{align*}
Similarly, we have 
\begin{align*}
\|\omega_{R}^{2}(\mathcal{R}*\widehat{\mathring{\chi}_{\lambda}})\|_{H^{2}}\lesssim R^{2}\lambda^{4}\sim\lambda^{2},
\end{align*}
from which we conclude 
\begin{align*}
\|\mathring{\eps}\|_{H^{2}}\lesssim\lambda^{4}+|\nu|\lambda^{3}+b\lambda^{2}\lesssim\lambda^{3}.
\end{align*}
This finishes the proof. 
\end{proof}
We have the same scaling bound to \eqref{eq: GR+Q sim zero} for this
constructed initial data profile $u_{0}=\mathcal{R}_{\mathring{b}_{0},\mathring{\eta}_{0},\mathring{\nu}_{0}}+\mathring{\eps}_{0}$. 
\begin{lem}
Under \eqref{eq:InitialDataUchi Chiral} condition, we have 
\begin{align}
\|\mathcal{G}(u_{0})+Q\|_{L^{2}} & \lesssim\delta_{chi},\label{eq:HardyLocalizationL2Estimate}\\
\|\mathcal{G}(u_{0})+Q\|_{\dot{\mathcal{H}}^{1}} & \lesssim\delta_{chi}^{\frac{1}{3}}\mathring{\lambda}_{0},\label{eq:HardyLocalizationH1Estimate}\\
\|\mathcal{G}(u_{0})+Q+P\|_{\dot{\mathcal{H}}^{2}} & \lesssim\delta_{chi}^{-\frac{1}{3}}\mathring{\lambda}_{0}^{2}.\label{eq:HardyLocalizationH2Estimate}
\end{align}
\end{lem}

\begin{proof}
By \eqref{eq: GR+Q sim zero}, it suffices to show 
\begin{align*}
\|\mathcal{G}(\mathcal{R}_{\mathring{b}_{0},\mathring{\eta}_{0},\mathring{\nu}_{0}}+\mathring{\eps}_{0})-\mathcal{G}(\mathcal{R}_{\mathring{b}_{0},\mathring{\eta}_{0},\mathring{\nu}_{0}})\|_{H^{2}}\lesssim\mathring{b}_{0}^{2}.
\end{align*}
We define a polarization of the gauge transform $\mathcal{G}$ by $\mathcal{G}_{f}g$
\begin{align*}
\mathcal{G}_{f}g\coloneqq ge^{-i\frac{1}{2}\int_{-\infty}^{x}|f|^{2}dy}.
\end{align*}
Note that $\mathcal{G}_{f}f=\mathcal{G}f$. For the convenience of
notation, we simply denote $h=\mathcal{R}_{\mathring{b}_{0},\mathring{\eta}_{0},\mathring{\nu}_{0}}$.
We have 
\begin{align*}
\mathcal{G}_{h+\mathring{\eps}}(h+\mathring{\eps}_{0})-\mathcal{G}_{h}h=(\mathcal{G}_{h+\mathring{\eps}_{0}}-\mathcal{G}_{h})h+\mathcal{G}_{h+\mathring{\eps}_{0}}\mathring{\eps}_{0}.
\end{align*}
We first have 
\begin{align*}
\|\mathcal{G}_{h+\mathring{\eps}_{0}}\mathring{\eps}_{0}\|_{H^{2}}\lesssim(1+\|h\|_{H^{2}}^{C}+\|\mathring{\eps}_{0}\|_{H^{2}}^{C})\|\mathring{\eps}_{0}\|_{H^{2}}\lesssim\|\mathring{\eps}_{0}\|_{H^{2}}\lesssim\mathring{b}_{0}^{2},
\end{align*}
for some $C$. For $(\mathcal{G}_{h+\mathring{\eps}_{0}}-\mathcal{G}_{h})h$,
we have 
\begin{align*}
(\mathcal{G}_{h+\mathring{\eps}_{0}}-\mathcal{G}_{h})h=\mathcal{G}_{h}h(\exp{\textstyle \left(-\frac{i}{2}\int_{-\infty}^{x}|h+\mathring{\eps}_{0}|^{2}-|h|^{2}dy\right)-1).}
\end{align*}
Thus, we have 
\begin{align*}
\|(\mathcal{G}_{h+\mathring{\eps}_{0}}-\mathcal{G}_{h})h\|_{H^{2}}\lesssim\|\mathcal{G}_{h}h\|_{H^{2}}(\|\mathring{\eps}_{0}\|_{L^{2}}+\|\mathring{\eps}_{0}\|_{L^{2}}^{2}+\||h+\mathring{\eps}_{0}|^{2}-|h|^{2}\|_{H^{1}}).
\end{align*}
Using $h\in H^{\infty}$, we have 
\begin{align*}
\||h+\mathring{\eps}_{0}|^{2}-|h|^{2}\|_{H^{1}}\lesssim\|\mathring{\eps}_{0}\|_{H^{1}}.
\end{align*}
Therefore, we conclude 
\begin{align*}
\|\mathcal{G}_{h+\mathring{\eps}_{0}}(h+\mathring{\eps}_{0})-\mathcal{G}_{h}h\|_{H^{2}}\lesssim\|\mathring{\eps}_{0}\|_{H^{2}}\lesssim\mathring{b}_{0}^{2},
\end{align*}
and this finishes the proof. 
\end{proof}
Define the chiral initial data sets $\mathcal{O}_{chi}$ for \eqref{CMdnls}
and $\mathcal{O}_{chi}^{\mathcal{G}}$ for \eqref{CMdnls-gauged}
by 
\begin{align*}
\mathcal{O}_{chi} & \coloneqq\{[u_{0}]_{\mathring{\lambda}_{0},\mathring{\gamma}_{0},\mathring{x}_{0}}:(\mathring{\lambda}_{0},\mathring{\gamma}_{0},\mathring{x}_{0},\mathring{b}_{0},\mathring{\eta}_{0},\mathring{\nu}_{0})\in\mathcal{U}_{chi}\}\subset\mathcal{S}_{+}.\\
\mathcal{O}_{chi}^{\mathcal{G}} & \coloneqq\{-\mathcal{G}([u_{0}]_{\mathring{\lambda}_{0},\mathring{\gamma}_{0},\mathring{x}_{0}}):[u_{0}]_{\mathring{\lambda}_{0},\mathring{\gamma}_{0},\mathring{x}_{0}}\in\mathcal{O}_{chi}\}\subset H^{2}\cap\mathcal{O}_{dec}.
\end{align*}
We are going to apply the bootstrap proposition, Proposition~\ref{PropositionMainBootstrap}
for solutions to \eqref{CMdnls-gauged}. So, we start the analysis
for the gauge transformed data set $\mathcal{O}_{chi}^{\mathcal{G}}$.
Note that $\mathcal{O}_{chi}^{\mathcal{G}}$ is, in fact, a five-dimensional set in the chiral class, since there are five free parameters. On the gauge-transformed side, we denote the initial data by $v_{0}=[-\mathcal{G}(u_{0})]_{\mathring{\lambda}_{0},\mathring{\gamma}_{0},\mathring{x}_{0}}=-\mathcal{G}([u_{0}]_{\mathring{\lambda}_{0},\mathring{\gamma}_{0},\mathring{x}_{0}})\in\mathcal{O}_{chi}^{\mathcal{G}}$.

Let $v(t)$ be the solution to \eqref{CMdnls-gauged} with the initial
data $v_{0}\in\mathcal{O}_{chi}^{\mathcal{G}}$. Then, we have $v(t)\in\mathcal{O}_{dec}$,
and by the decomposition lemma, Lemma~\ref{LemmaDecomposition},
there exists a decomposition with $(\lambda,\gamma,x,b,\eta,\nu,\eps)\in\mathbb{R}_{+}\times\mathbb{R}/2\pi\mathbb{Z}\times\mathbb{R}^{4}\times\mathcal{Z}^{\perp}$
such that 
\begin{align}
v(t)=[Q+P(b,\eta,\nu)+\eps]_{\lambda,\gamma,x}.\label{eq:5.8 sol v decompose}
\end{align}
We denote the initial modulation parameters by 
\begin{align*}
(\lambda(0),\gamma(0),x(0),b(0),\eta(0),\nu(0),\eps(0))=(\lambda_{0},\gamma_{0},x_{0},b_{0},\eta_{0},\nu_{0},\eps_{0}),
\end{align*}
and then we will prepare to bootstrap.

Here, we have the transfer map of decompositions 
\[
(\mathring{\lambda}_{0},\mathring{\gamma}_{0},\mathring{x}_{0},\mathring{b}_{0},\mathring{\eta}_{0},\mathring{\nu}_{0})\mapsto(\lambda_{0},\gamma_{0},x_{0},b_{0},\eta_{0},\nu_{0},\eps_{0})
\]
defined by the relation 
\[
[-\mathcal{G}(u_{0}(\cdot;\mathring{b}_{0},\mathring{\eta}_{0},\mathring{\nu}_{0})]_{\mathring{\lambda}_{0},\mathring{\gamma}_{0},\mathring{x}_{0}}=v_{0}=[Q+P(b_{0},\eta_{0},\nu_{0})+\eps_{0}]_{\lambda_{0},\gamma_{0},x_{0}}.
\]
This transfer map is continuous by the continuity of $(\mathring{\lambda}_{0},\mathring{\gamma}_{0},\mathring{x}_{0},\mathring{b}_{0},\mathring{\eta}_{0},\mathring{\nu}_{0})\mapsto v_{0}$
and Lemma~\ref{LemmaDecomposition}. In particular, for fixed $\mathring{\lambda}_{0},\mathring{\gamma}_{0},\mathring{x}_{0},\mathring{b}_{0}$,
we will use a rescaled continuous map 
\begin{align}
\Phi_{chi}:\Big(\tfrac{10\mathring{\eta}_{0}}{\mathring{\lambda}_{0}^{\frac{3}{2}(1+\kappa)}},\tfrac{10\mathring{\nu}_{0}}{\mathring{\lambda}_{0}^{\frac{3}{2}(1-\kappa)}}\Big)\mapsto(\lambda_{0},\gamma_{0},x_{0},b_{0},\eta_{0},\nu_{0},\eps_{0}).\label{eq:DecompositionMapChiral}
\end{align}
Next, we will apply the main bootstrap proposition. For this purpose,
we set the initial data condition for chiral solution \eqref{eq:InitialDataUchi Chiral}
so that these initial modulation parameters given by \eqref{eq:DecompositionMapChiral}
satisfy the initial bound condition \eqref{eq:5.2 ini boots assump}.

To obtain a chiral blow-up solution with initial data arbitrarily
close to $\mathcal{R}$ in $L^{2}$, we will further refine $b_{chi}^{*}$ to incorporate the smallness of $\epsilon$ with keeping
Proposition~\ref{PropositionMainBootstrap} on hand (and \eqref{eq:5.2 rem para dependecy}).
We will have the following parameter dependency: 
\begin{align}
0<b_{chi}^{\ast}\ll\delta_{chi}^{\frac{100}{\kappa}\cdot\frac{1}{3}}\ll\min(\epsilon,\delta_{dec})^{\frac{100}{\kappa}\cdot\frac{1}{4}}\ll1.\label{eq:5.8 chiral para dependecy}
\end{align}
We note that since we have already fixed $\delta_{chi}=\frac{1}{100}\min(\epsilon,\delta_{dec})\ll1$,
the last two dependences are established naturally. 
\begin{lem}
\label{LemmaHardyLocalBootstrapAssumption} For arbitrary small $\epsilon>0$,
there exist constants $b_{chi}^{*}>0$ which satisfy
the following: (i) $b_{chi}^{*}$ satisfies the parameter
dependence \eqref{eq:5.2 rem para dependecy}. (ii) $(\lambda_{0},\gamma_{0},x_{0},b_{0},\eta_{0},\nu_{0},\eps_{0})$
satisfy the initial bound condition \eqref{eq:5.2 ini boots assump}
with $b^{*}=b_{chi}^{*}$. 
\end{lem}

\begin{proof}
Fix $b_{chi}^{*}$ satisfying \eqref{eq:5.8 chiral para dependecy},
and then $b_{chi}^{*}$ also satisfies \eqref{eq:5.2 rem para dependecy}.
For this $b_{chi}^{*}$, it suffices to show (ii).

We note that $\lambda\in\bbR_{+}$, $\gamma\in\bbR/2\pi\bbZ$, and
$(x,b,\eta,\nu)\in\bbR^{4}$. Here, we equip $\bbR^{4}$ with the
Euclidean metric. We also equip $\bbR_{+}$ with the metric $\text{dist}(\lambda_{1},\lambda_{2})=|\log(\frac{\lambda_{1}}{\lambda_{2}})|$,
and equip $\bbR/2\pi\bbZ$ with the induced metric from $\bbR$. We
first claim that 
\begin{align}
\text{dist}((\mathring{\lambda}_{0},\mathring{\gamma}_{0},\mathring{x}_{0},\mathring{b}_{0},\mathring{\eta}_{0},\mathring{\nu}_{0}),(\lambda_{0},\gamma_{0},x_{0},b_{0},\eta_{0},\nu_{0}))\lesssim\delta_{chi}^{-\frac{1}{3}}\mathring{\lambda}_{0}^{2},\label{eq:ChiralityInitialModulationDifference}
\end{align}
where the parameters are defined in \eqref{eq:InitialDataUchi Chiral}
and \eqref{eq:5.8 sol v decompose}. We define the function $\mathbf{F}$
by 
\begin{align*}
\mathbf{F}(\lambda,\gamma,x;v)=((\widehat{\eps},\mathcal{Z}_{1})_{r},(\widehat{\eps},\mathcal{Z}_{2})_{r},(\widehat{\eps},\mathcal{Z}_{3})_{r})^{T},\quad\text{with}\quad\widehat{\eps}=[v]_{\lambda^{-1},-\gamma,-x}-Q.
\end{align*}
By the proof of Lemma~\ref{LemmaDecomposition} and the proof of
the implicit function theorem, we obtain the difference estimate using
$\mathbf{F}(\lambda_{0},\gamma_{0},x_{0};v_{0})=0$, 
\begin{align*}
\text{dist}((\mathring{\lambda}_{0},\mathring{\gamma}_{0},\mathring{x}_{0}),(\lambda_{0},\gamma_{0},x_{0})) & <|\mathbf{F}(\mathring{\lambda}_{0},\mathring{\gamma}_{0},\mathring{x}_{0};v_{0})-\mathbf{F}(\lambda_{0},\gamma_{0},x_{0};v_{0})|\\
 & =|\mathbf{F}(\mathring{\lambda}_{0},\mathring{\gamma}_{0},\mathring{x}_{0};v_{0})|,
\end{align*}
where $v_{0}=[-\mathcal{G}(u_{0})]_{\mathring{\lambda}_{0},\mathring{\gamma}_{0},\mathring{x}_{0}}$.
Thanks to the transversality \eqref{eq:transversality} and the scaling
bound \eqref{eq:HardyLocalizationH2Estimate}, we have 
\begin{align*}
|\mathbf{F}(\mathring{\lambda}_{0},\mathring{\gamma}_{0},\mathring{x}_{0};v_{0})|={\textstyle \sum_{k=1}^{3}}|(-\mathcal{G}(u_{0})-(Q+P(\mathring{b}_{0},\mathring{\eta}_{0},\mathring{\nu}_{0})),\mathcal{Z}_{k})_{r}|\lesssim\delta_{chi}^{-\frac{1}{3}}\mathring{\lambda}_{0}^{2},
\end{align*}
and hence we deduce 
\begin{align*}
\text{dist}((\mathring{\lambda}_{0},\mathring{\gamma}_{0},\mathring{x}_{0}),(\lambda_{0},\gamma_{0},x_{0}))\lesssim\delta_{chi}^{-\frac{1}{3}}\mathring{\lambda}_{0}^{2}.
\end{align*}
Now, we want to estimate the difference for $b,\eta,\nu$. By \eqref{eq:DecompositionAppendix b eta nu definition}
and transversality \eqref{eq:transversality}, we have 
\begin{align*}
b_{0}=-\frac{(\widehat{\eps_{0}},\mathcal{Z}_{4})_{r}}{(i\frac{y^{2}}{4}Q,\mathcal{Z}_{4})_{r}}=\mathring{b}_{0}-\frac{(\eps_{0},\mathcal{Z}_{4})_{r}}{(i\frac{y^{2}}{4}Q,\mathcal{Z}_{4})_{r}}
\end{align*}
where $\widehat{\eps}_{0}$ and $\eps_{0}$ are given by 
\begin{align*}
\widehat{\eps}_{0}=[v_{0}]_{\lambda_{0}^{-1},-\gamma_{0},-x_{0}}-Q,\quad\text{and}\quad\eps_{0}=\widehat{\eps}_{0}-P(b_{0},\eta_{0},\nu_{0}).
\end{align*}
Thus, we have $|b_{0}-\mathring{b}_{0}|\lesssim\delta_{chi}^{-\frac{1}{3}}\mathring{\lambda}_{0}^{2}$.
Using a similar argument, we also have $|\eta_{0}-\mathring{\eta}_{0}|+|\nu_{0}-\mathring{\nu}_{0}|\lesssim\delta_{chi}^{-\frac{1}{3}}\mathring{\lambda}_{0}^{2}$.
Therefore, we conclude the claim \eqref{eq:ChiralityInitialModulationDifference}.

Now, we show \eqref{eq:5.2 ini boots assump} with $b_{chi}^{*}$. By \eqref{eq:InitialDataUchi Chiral} and \eqref{eq:ChiralityInitialModulationDifference},
we deduce 
\begin{align*}
0<b_{0}<\tfrac{1}{2}b_{chi}^{*},\quad|\eta_{0}|<\tfrac{1}{2}b_{0}^{1+\kappa},\quad|\nu_{0}|<\tfrac{1}{2}b_{0}^{1-\kappa},\quad\tfrac{1}{1.1}\leq\tfrac{b_{0}^{2}}{\lambda_{0}^{3}}\leq1.1.
\end{align*}
Recall that $w=[v]_{\lambda^{-1},-\gamma,-x}$ and 
\begin{align*}
w-Q=\widehat{\eps},\quad w_{1}=\mathbf{D}_{w}w,\quad A_{Q}=\partial_{y}(y-\mathcal{H})\langle y\rangle^{-1}.
\end{align*}
We will show that 
\begin{align*}
\|\widehat{\eps}_{0}\|_{L^{2}}\leq\tfrac{1}{40}\delta_{dec}.
\end{align*}
Indeed, we can write 
\begin{align*}
\widehat{\eps}_{0}= & [-\mathcal{G}(u_{0})]_{\frac{\mathring{\lambda}_{0}}{\lambda_{0}},\mathring{\gamma}_{0}-\gamma_{0},\mathring{x}_{0}-x_{0}}-Q\\
= & ([-\mathcal{G}(u_{0})]_{\frac{\mathring{\lambda}_{0}}{\lambda},\mathring{\gamma}_{0}-\gamma,\mathring{x}_{0}-x}+\mathcal{G}(u_{0}))-(\mathcal{G}(u_{0})+Q),
\end{align*}
 then it follows that $\|\widehat{\eps}_{0}\|_{L^{2}}\leq\frac{1}{40}\delta_{dec}$
directly from \eqref{eq:HardyLocalizationL2Estimate} and \eqref{eq:ChiralityInitialModulationDifference}.
Similarly, by \eqref{eq:HardyLocalizationH1Estimate} and \eqref{eq:ChiralityInitialModulationDifference},
\begin{align*}
\|\widehat{\eps}_{0}\|_{\dot{\mathcal{H}}^{1}}\leq & \|[-\mathcal{G}(u_{0})]_{\frac{\mathring{\lambda}_{0}}{\lambda},\mathring{\gamma}_{0}-\gamma,\mathring{x}_{0}-x}+\mathcal{G}(u_{0})\|_{\dot{\mathcal{H}}^{1}}+\|\mathcal{G}(u_{0})+Q\|_{\dot{\mathcal{H}}^{1}}\\
\leq & C\delta_{chi}^{-\frac{1}{3}}\lambda_{0}^{2}+C\delta_{chi}^{\frac{1}{3}}\lambda_{0}.
\end{align*}
Therefore, we obtain 
\begin{equation}
\|\widehat{\eps}_{0}\|_{\dot{\mathcal{H}}^{1}}\leq\tfrac{1}{100}\lambda_{0}.\label{eq:hat epsilon_0}
\end{equation}
This finishes the proof. 
\end{proof}

Now, we are ready to prove Theorem~\ref{TheoremChiralBlowup}. 

\begin{proof}[Proof of Theorem~\ref{TheoremChiralBlowup}]
The proof goes along with the proof of Theorem~\ref{thm:precise statement codimensionone blowup}
in Section~\ref{Subsection ExistenceTrappedSol}. We will apply the
bootstrap proposition and topological lemmas to obtain a trapped solution
to the gauged equation \eqref{CMdnls-gauged}. Then, we go back to
\eqref{CMdnls} and show the convergence of the radiation part in \eqref{CMdnls}.

In a similar spirit to the proof of Theorem~\ref{thm:precise statement codimensionone blowup},
we now prepare a two-dimensional set of \emph{chiral} initial data.
Fix $\mathring{\lambda}_{0},\mathring{\gamma}_{0},\mathring{x}_{0},\mathring{b}_{0}$
such that $0<\mathring{b}_{0}<\frac{1}{3}b_{chi}^{\ast}$ and $\mathring{b}_{0}^{2}=\mathring{\lambda}_{0}^{3}$.
Consider $u_{0}=u_{0}(\cdot;\mathring{b}_{0},\mathring{\eta}_{0},\mathring{\nu}_{0})\in\mathcal{S}_{+}$
with $(\mathring{\eta}_{0},\mathring{\nu}_{0})$ varying in $[-\frac{1}{4}\mathring{b}_{0}^{1+\kappa},\frac{1}{4}\mathring{b}_{0}^{1+\kappa}]\times[-\frac{1}{4}\mathring{b}_{0}^{1-\kappa},\frac{1}{4}\mathring{b}_{0}^{1-\kappa}]$.
Set $v_{0}=-[\G(u_{0}(\mathring{b}_{0},\mathring{\eta}_{0},\mathring{\nu}_{0}))]_{\mathring{\lambda}_{0},\mathring{\gamma}_{0},\mathring{x}_{0}}\in\mathcal{O}_{chi}^{\mathcal{G}}$.
Let $v(t)$ be the (forward) maximal solution to \eqref{CMdnls-gauged}
with initial data $v_{0}\in\mathcal{O}_{chi}^{\mathcal{G}}$ and lifespan
$[0,T)$. We hope to find special special $(\mathring{\eta}_{0},\mathring{\nu}_{0})$
such that $v(t)$ becomes a trapped solution.

Decompose $v(t)$ as \eqref{eq:5.8 sol v decompose} with the modulation
parameters $(\lambda,\gamma,x,b,\eta,\nu,\eps)(t)$ and denote by
$T_{exit}$ the exit time of the bootstrap hypotheses \eqref{eq:BootstapAssumption}
as in \eqref{eq:DefinitionTexit}. Thanks to Lemma~\ref{LemmaHardyLocalBootstrapAssumption},
the initial modulation parameters $(\lambda,\gamma,x,b,\eta,\nu,\eps)(0)$
satisfy \eqref{eq:5.2 ini boots assump} with $b^{\ast}=b_{chi}^{\ast}$. In particular, they satisfy \eqref{eq:BootstapAssumption}
and we have $T_{exit}>0$. Moreover, the bootstrap Proposition~\ref{PropositionMainBootstrap}
can be applied.

We show that there exists $(\mathring{\eta}_{0},\mathring{\nu}_{0})\in[-\frac{1}{4}\mathring{b}_{0}^{1+\kappa},\frac{1}{4}\mathring{b}_{0}^{1+\kappa}]\times[-\frac{1}{4}\mathring{b}_{0}^{1-\kappa},\frac{1}{4}\mathring{b}_{0}^{1-\kappa}]$
such that $v(t)$ becomes a trapped solution. As before, we proceed
with a contradiction argument; assume $T_{exit}<T$ for all such $(\mathring{\eta}_{0},\mathring{\nu}_{0})$.
By Proposition~\ref{PropositionMainBootstrap}, we have \eqref{eq:TexitBootstrapfail}
at the exit time $T_{exit}$. Then, $\widetilde{T}_{exit}$ given
by \eqref{eq:DefinitionTexitTilde} is well defined for any $(\mathring{\eta}_{0},\mathring{\nu}_{0})\in[-0.1\mathring{b}_{0}^{1+\kappa},0.1\mathring{b}_{0}^{1+\kappa}]\times[-0.1\mathring{b}_{0}^{1-\kappa},0.1\mathring{b}_{0}^{1-\kappa}]$.
Now, we define a function $\Psi_{chi}:[-1,1]^{2}\to[-1.2,1.2]^{2}\setminus[-0.9,0.9]^{2}$
with 
\begin{align}
\Psi_{chi}\Big(\tfrac{10\mathring{\eta}_{0}}{\mathring{\lambda}_{0}^{\tfrac{3}{2}+\kappa}},\tfrac{10\mathring{\nu}_{0}}{\mathring{\lambda}_{0}^{\frac{3}{2}-\kappa}}\Big)=\left(\left(\tfrac{10\widetilde{\eta}}{\lambda^{\frac{3}{2}+\kappa}}\right)(\widetilde{T}_{exit}),\left(\tfrac{10\widetilde{\nu}}{\lambda^{\frac{3}{2}-\kappa}}\right)(\widetilde{T}_{exit})\right),\label{eq:def Psi chi}
\end{align}
where $\widetilde{\eta}$, $\widetilde{\nu}$ are the refined modulation
parameters given by \eqref{eq:RefinedModulationDefinition}. First,
we show that $\Psi_{chi}$ is continuous. To see this fact more precisely,
we observe that $\Psi_{chi}$ is a composition of continuous maps,
\begin{equation}
\Phi\circ\Phi_{chi}=\Psi_{chi}.\label{eq:Psi_chi}
\end{equation}
Here, we recall the map $\Phi:\mathcal{I}\to[-1.2,1.2]^{2}\setminus[-0.9,0.9]^{2}$
defined by \eqref{eq:5.2 def Phi} where $\mathcal{I}$ is defined
by \eqref{eq:5.2. def set I}. Also, $\Phi_{chi}$ is the transfer
map of parameters, given in \eqref{eq:DecompositionMapChiral}. Thanks
to Lemma~\ref{LemmaHardyLocalBootstrapAssumption} with \eqref{eq:ChiralityInitialModulationDifference},
we have $\Phi_{chi}\big(\frac{10\mathring{\eta}_{0}}{\mathring{\lambda}_{0}^{\frac{3}{2}(1+\kappa)}},\frac{10\mathring{\nu}_{0}}{\mathring{\lambda}_{0}^{\frac{3}{2}(1-\kappa)}}\big)\in\mathcal{I}$
for any $(\mathring{\eta}_{0},\mathring{\nu}_{0})\in[-0.1\mathring{b}_{0}^{1+\kappa},0.1\mathring{b}_{0}^{1+\kappa}]\times[-0.1\mathring{b}_{0}^{1-\kappa},0.1\mathring{b}_{0}^{1-\kappa}]$,
so the composition map $\Phi\circ\Phi_{chi}$ from $[-1,1]^{2}$ is
well defined, verifying \eqref{eq:Psi_chi}. From \eqref{eq:DecompositionMapChiral},
$\Phi_{chi}$ is continuous, and $\Phi$ is also continuous by Lemma~\ref{LemmaTopologicalPhiContinuous}.
Hence $\Psi_{chi}$ is continuous. By the definition of $\Psi_{chi}$
\eqref{eq:def Psi chi}, \eqref{eq:ChiralityInitialModulationDifference},
and \eqref{eq:RefinedModulationDiffer}, we verify (2) of Lemma~\ref{LemmaPsiNonExistenceBrouwer}.
By Lemma~\ref{LemmaPsiNonExistenceBrouwer}, such a function $\Psi_{chi}$
cannot exist, which is a contradiction. Therefore, there exists a
$(\mathring{\eta}_{0},\mathring{\nu}_{0})\in[-\frac{1}{4}\mathring{b}_{0}^{1+\kappa},\frac{1}{4}\mathring{b}_{0}^{1+\kappa}]\times[-\frac{1}{4}\mathring{b}_{0}^{1-\kappa},\frac{1}{4}\mathring{b}_{0}^{1-\kappa}]$
such that $v(t)$ with the initial data $v_{0}\in\mathcal{O}_{chi}^{\mathcal{G}}$
is a trapped solution to \eqref{CMdnls-gauged}.

By Proposition~\ref{PropSharpDescription}, the trapped solution
$v(t)$ found in the previous paragraph blows up in finite time (together
with sharp descriptions). Applying the gauge transform $\mathcal{-G}^{-1}$
to $v(t)$, we have a blow-up solution $u(t)$ with the initial data
$[u_{0}(\cdot;\mathring{b}_{0},\mathring{\eta}_{0},\mathring{\nu}_{0})]_{\mathring{\lambda}_{0},\mathring{\gamma}_{0},\mathring{x}_{0}}\in\mathcal{O}_{chi}$,
the Schwartz chiral initial data.

Now, it remains to show the convergence of the remainder part of the
solution $u=-\mathcal{G}^{-1}(v)$ to \eqref{CMdnls}. We know by
Proposition~\ref{PropSharpDescription}, the reminder of $v$, $[\wt{\eps}]_{\la,\ga,x}$
converges to $v^{*}$. We will show that for some constant $c$ and
asymptotic profile $u^{*}$, 
\begin{align*}
u(t,x)-[\mathcal{R}]_{\lambda,\gamma^{*},x^{*}}e^{ic}\to u^{*}\quad\text{as}\quad t\to T\quad\text{in}\quad L^{2},
\end{align*}
where $\ga^{*}$ and $x^{*}$ are given in Proposition~\ref{PropSharpDescription}.

Using the inverse of gauge transform 
\begin{align*}
\mathcal{G}^{-1}(f)(x)=f(x)e^{\frac{i}{2}\int_{-\infty}^{x}|f(y)|^{2}dy},
\end{align*}
and $\mathcal{G}^{-1}([f]_{\lambda,\gamma,x})=[\mathcal{G}^{-1}(f)]_{\lambda,\gamma,x}$,
we write 
\begin{align}
-u(t,x)-\mathcal{G}^{-1}([Q]_{\lambda,\gamma,x})e^{ic}=\mathcal{G}^{-1}([Q+\widehat{\eps}]_{\lambda,\gamma,x})-\mathcal{G}^{-1}([Q]_{\lambda,\gamma,x})e^{ic},\label{eq:DescriptionUnGaugedU}
\end{align}
where $c=c(x^{*},v^{*})$ is a constant given by 
\begin{align}
c\coloneqq-\frac{1}{2}\int_{-\infty}^{x^{*}}|v^{*}|^{2}dy.\label{eq:Definitionconstantxstarvstar}
\end{align}
Since $\gamma(t)\to\gamma^{*}$ and $x(t)\to x^{*}$ as $t\to T$,
we will omit $\ga$ and $x$ parameters subscripts, $[f]_{\lambda}\coloneqq[f]_{\lambda,\gamma^{*},x^{*}}$.
We have 
\begin{align}
\eqref{eq:DescriptionUnGaugedU}= & \mathcal{G}^{-1}([\widehat{\eps}]_{\lambda}){\textstyle \exp\left(\frac{i}{2}\int_{-\infty}^{x}|[Q+\widehat{\eps}]_{\lambda}|^{2}-|[\widehat{\eps}]_{\lambda}|^{2}dy\right)}\label{eq:description u^*1}\\
 & +{\textstyle \left[\mathcal{G}^{-1}(Q)\left\{ \exp\left(\frac{i}{2}\int_{-\infty}^{x}|Q+\widehat{\eps}|^{2}-|Q|^{2}dy\right)-e^{ic}\right\} \right]_{\lambda}.}\label{eq:description u^*2}
\end{align}
Since $[\widehat{\eps}]_{\lambda(t),\gamma(t),x(t)}\to v^{*}$ in
$L^{2}$ and $\lambda(t)\to0$, we have $\widehat{\eps}\rightharpoonup0$
weakly in $L^{2}$. For \eqref{eq:description u^*1}, we have $\mathcal{G}^{-1}([\widehat{\eps}]_{\lambda})\to\mathcal{G}^{-1}(v^{*})$
in $L^{2}$. We compute 
\begin{align*}
{\textstyle \int_{-\infty}^{x}|[Q+\widehat{\eps}]_{\lambda}|^{2}-|[\widehat{\eps}]_{\lambda}|^{2}dy=\Re\int_{-\infty}^{\frac{x-x^{*}}{\lambda}}Q^{2}+2Q\widehat{\eps}dy.}
\end{align*}
We have the pointwise convergence $\int_{-\infty}^{\frac{x-x^{*}}{\lambda}}Q^{2}dy\to0$
for $x<x^{*}$ and $\int_{-\infty}^{\frac{x-x^{*}}{\lambda}}Q^{2}dy\to M(Q)=2\pi$
for $x>x^{*}$. We have $\Re\int_{-\infty}^{\frac{x-x^{*}}{\lambda}}Q\widehat{\eps}dy\to0$
since $\widehat{\eps}\rightharpoonup0$ weakly in $L^{2}$. Therefore,
using DCT with $\mathcal{G}^{-1}(v^{*})\in H^{1}$, we have 
\[
\eqref{eq:description u^*1}\to\mathcal{G}^{-1}(v^{*})\mathbf{1}_{x<x^{*}}-\mathcal{G}^{-1}(v^{*})\mathbf{1}_{x>x^{*}}\text{ in }L^{2}.
\]
Now, we want to show \eqref{eq:description u^*2} $\to0$ in $L^{2}$.
We have 
\begin{align*}
{\textstyle \int_{-\infty}^{x}|Q+\widehat{\eps}|^{2}-|Q|^{2}dy=\Re\int_{-\infty}^{x}2Q\widehat{\eps}+|\widehat{\eps}|^{2}dy.}
\end{align*}
As in \eqref{eq:description u^*1}, we also have $\Re\int_{-\infty}^{x}2Q\widehat{\eps}dy\to0$
as $t\to T$. By \eqref{eq:Definitionconstantxstarvstar}, it suffices
to show the pointwise convergence, 
\begin{align*}
{\textstyle \int_{-\infty}^{x}|\widehat{\eps}|^{2}dy-\int_{-\infty}^{x^{*}}|v^{*}|^{2}dy\to0.}
\end{align*}
Since $[\widehat{\eps}]_{\lambda}\to v^{*}$ in $L^{2}$ and $\lambda\to0$,
we have 
\begin{align*}
{\textstyle \int_{-\infty}^{x}|\widehat{\eps}|^{2}dy=\int_{-\infty}^{\lambda x+x^{*}}|[\widehat{\eps}]_{\lambda}|^{2}dy\to\int_{-\infty}^{x^{*}}|v^{*}|^{2}dy.}
\end{align*}
Again using DCT with $e^{ic}\mathcal{G}^{-1}(Q)\in H^{1}$, we deduce
\begin{align*}
\|\eqref{eq:description u^*2}\|_{L^{2}}\to0\quad\text{as}\quad t\to T.
\end{align*}
Combining these results, we have 
\begin{align*}
u(t,x)-[\mathcal{R}]_{\lambda,\gamma,x}e^{ic}\to u^{*}\quad\text{as}\quad t\to T\quad\text{in}\quad L^{2},
\end{align*}
where $u^{*}$ is given by 
\begin{align}
u^{*}=-\mathcal{G}^{-1}(v^{*})\mathbf{1}_{x<x^{*}}+\mathcal{G}^{-1}(v^{*})\mathbf{1}_{x>x^{*}},\label{eq:u-ast-formula}
\end{align}
with $v^{*}\in H^{1}$. Since $u(t),\mathcal{R}\in L_{+}^{2}$, we
have $u^{*}\in L_{+}^{2}$. Redenoting $\gamma^{*}+c$ by $\gamma^{*}$,
we show the sharp description for a chiral blow-up solution. This finishes
the proof of Theorem~\ref{TheoremChiralBlowup}.
\end{proof}
\begin{rem}
In view of \eqref{eq:u-ast-formula}, we deduce only $u^{\ast}\in L^{2}$
although $v^{\ast}\in H^{1}$. If $v^{\ast}$ degenerates at the blow-up
point $x^{\ast}$ (say $v^{\ast}(x^{\ast})=0$), then the regularity
of $u^{\ast}$ can be further improved. But we do not pursue this
direction here. 
\end{rem}

\section{Proof of the nonlinear estimate, Lemma~\ref{LemmaNonlinearEstimate1}}\label{Section Proof of nonlinear estimates} 

In this section, we provide the proof of nonlinear estimate, Lemma~\ref{LemmaNonlinearEstimate1}. Before proving this
lemma, we record that the pointwise bound for modified profiles $P$
and $P_{1}$:
\begin{align}
\begin{split}
	|\partial_{y}^{k}P_{e}|(y)\lesssim bQ^{k-1},\quad 
	|\partial_{y}^{k}P_{o}|(y)\lesssim b^{1-\kappa}Q^{k} \text{ for }k\geq0,\\
	|\partial_{y}^{k}P_{1,e}|(y)\lesssim \lambda Q^{k+1},\quad 
	|\partial_{y}^{k}P_{1,o}|(y)\lesssim bQ^{k} \text{ for }k\geq0.
\end{split}
\label{eq:ProfilePointwiseBound}
\end{align}

\subsection{Comments on Lemma~\ref{LemmaNonlinearEstimate1}}\label{sec:proof of lemma 5.17}

We aim to prove Lemma~\ref{LemmaNonlinearEstimate1} by extracting from the nonlinear part the profile terms that are strictly larger than $b^2$.
Recall that 
\begin{align*}
    \textnormal{NL}_{1}=-\tfrac{1}{4}(|w|^{4}-Q^{4})w_{1}+\widehat{\eps}|D|\overline{\widehat{\eps}}w_{1}+Q|D|\overline{\widehat{\eps}}w_{1}+\widehat{\eps}|D|Qw_{1}.
\end{align*}
We need to identify terms in $\textnormal{NL}_{1}$ that are quadratic in $b,\eta,\nu$. Formally, they are written as 
\begin{align}
    -\Re(Q^{3}P)P_{1}+Q|D|(\overline{P}P_{1})+P|D|(QP_{1}).\label{eq:NL1 Quadratic Profile Formal}
\end{align}
From the derivation of the modulation laws \eqref{eq:FormalBlowupLawApproxEquModified},
we have \emph{formally} 
\begin{align}
    \eqref{eq:NL1 Quadratic Profile Formal}=-i(b\Lambda_{-1}+\tfrac{\eta}{2}i-\nu\partial_{y})P_{1}-i(\tfrac{3b^{2}}{2}+\tfrac{\eta^{2}}{2})\tfrac{i}{2}yQ-ib\eta\tfrac{1}{2}yQ+ib\nu\tfrac{i}{2}Q+ib\mu\tfrac{1}{2}Q.\label{eq:NL1 Quadratic Profile another form}
\end{align}
More precisely, in \eqref{eq:NL1 Quadratic Profile Formal} we cannot
define $Q|D|(\overline{P}P_{1})$ since it contains $|D|(y^{3}Q^{2})$, which formally produces $\mathcal{H}(1)$ as seen in \eqref{eq:QuadraticModulation calculation nonlinearpart}. In particular, this yields the term $-(b^{2}+\eta^{2})\tfrac{i}{4}Q\mathcal{H}(1)$ in \eqref{eq:QuadraticModulation calculation nonlinearpart}.
To avoid this issue, we simply drop the harmful term $Q|D|(\overline{P}_{e}P_{1,o})$ from $Q|D|(\overline{P}P_{1})$. Accordingly, we define
\begin{align}
    \mathcal{P} \coloneqq Q\mathcal{H}\partial_{y}(\overline{P}P_{1} - \overline{P}_{e}P_{1,o}). \label{eq:def calP}
\end{align}
Using this, we identify the profile part of $\textnormal{NL}_{1}$ as
\begin{align}
    -\Re(Q^{3}P)P_{1}+ \mathcal{P} + P|D|(Q P_{1}). \label{eq:profile redefine}
\end{align}
Thus, our goal is to show 
\begin{align}
	\|\textnormal{NL}_{1}-(-\Re(Q^{3}P)P_{1}+ \mathcal{P} + P|D|(Q P_{1}))\|_{L^{2}}\lesssim b^{2}.\label{eq:NL-1 estimate 6.1 2}
\end{align}
Since $Q|D|(\overline{P}_{e}P_{1,o})$ is formally odd, we observe that the even part of \eqref{eq:NL-1 estimate 6.1 2} coincides with \eqref{eq:NL1eWithoutProfileEstimate}. Moreover, a direct computation shows that the odd part of \eqref{eq:profile redefine} is given by $-\frac{1}{2}(\nu^2 - i\nu\mu) Q_y$, which yields \eqref{eq:NL1Estimate}.

\begin{rem}[Commuting $\mathcal{H}$ and $\partial_{y}$]
\label{Remark commuting Hilbert derivative} Even if we restrict ourselves
to even parts and exclude $|D|(y^{3}Q^{2})$, we still have an issue
of handling $|D|$ in nonlinear terms. In view of the Fourier symbol,
we see $|D|=\mathcal{H}\partial_{y}$ in $H^{1}$. If the applied function
$f$ has less decay, then the commutation must be used with caution.
For example, one cannot directly have $|D|(y^{2}Q^{2})=\mathcal{H}\partial_{y}(y^{2}Q^{2})$.
Such a situation may occur when we handle a nonlinear term as like $|D|([\ol{\wt{\eps}}w_{1}]_{e})$.
Since $\overline{\widehat{\eps}}w_{1}$ contains $[\ol PP_{1}]_{e}$,
which contains $y^{2}Q^{2}$. However, this issue is avoidable since
$\overline{\widehat{\eps}}w_{1}\in H^{1}$. We apply $|D|=\mathcal{H}\partial_{y}$
before decomposing $\overline{\widehat{\eps}}w_{1}$ into several
terms. Indeed, we can proceed as 
\begin{align*}
Q|D|(\overline{\widehat{\eps}}w_{1})=Q\mathcal{H}\partial_{y}(\overline{\widehat{\eps}}w_{1})= & Q\mathcal{H}\partial_{y}(\overline{\widehat{\eps}}w_{1})_{o}+Q\mathcal{H}\partial_{y}(\overline{\widehat{\eps}}w_{1})_{e}\\
= & Q\mathcal{H}\partial_{y}(\overline{\widehat{\eps}}w_{1})_{o}+Q\mathcal{H}\partial_{y}(\overline{\widehat{\eps}}w_{1}-\overline{P}P_{1})_{e}+Q\mathcal{H}\partial_{y}(\overline{P}P_{1})_{e}.
\end{align*}
This makes sense since each term is well defined, (e.g., $Q\mathcal{H}\partial_{y}(\overline{P}P_{1})_{e}=\tfrac{b\nu}{8}Q\mathcal{H}\partial_{y}(y^{2}Q^{2})+\cdots$
and $\partial_{y}(y^{2}Q^{2})=-\partial_{y}(Q^{2})\in L^{2}$). In
this way, we can justify $|D|=\mathcal{H}\partial_{y}=\partial_{y}\mathcal{H}$
for most cases, which has no growth like $y^{2}Q^{2}$. However,
this argument is not applicable to growing functions like $y^{3}Q^{2}$.
So, as explained above, we have to set aside by taking even part.
In the sequel, we will abuse the notation, if applicable, $|D|$ as
$\mathcal{H}\partial_{y}$ or $\partial_{y}\mathcal{H}$ for simplicity
of presentation. 
\end{rem}

\subsection{Proof of Lemma~\ref{LemmaNonlinearEstimate1}}

We will first show \eqref{eq:NL-1 estimate 6.1 2}. Then, taking odd and even parts, we conclude \eqref{eq:NL1Estimate} and \eqref{eq:NL1eWithoutProfileEstimate}.

\textbf{Step 1. Estimate of quintic term $-\tfrac{1}{4}(|w|^{4}-Q^{4})w_{1}$.}

Rewrite by expanding $w=Q+\wt{\eps}$, 
\begin{align}
-\tfrac{1}{4}(|w|^{4}-Q^{4})w_{1}=-(Q^{2}\Re(Q\widehat{\eps})+\tfrac{1}{4}M_{Q}(\widehat{\eps}))w_{1}\label{eq:NonlinearEstimateQuintic}
\end{align}
where the quadratic and higher order in $\wt{\eps}$ is 
\begin{align}
M_{Q}(\widehat{\eps})=2Q^{2}|\widehat{\eps}|^{2}+4(\Re(Q\widehat{\eps}))^{2}+4|\widehat{\eps}|^{2}\Re(Q\widehat{\eps})+|\widehat{\eps}|^{4}.\label{eq:DefinitionMQ}
\end{align}
To show \eqref{eq:NL-1 estimate 6.1 2}, we further decompose $\widehat{\eps}=P+\eps$
to have 
\begin{align*}
\eqref{eq:NonlinearEstimateQuintic}=-\tfrac{1}{4}\big(4\Re(Q^{3}P)+4\Re(Q^{3}\eps)+M_{Q}(\widehat{\eps})\big)w_{1}.
\end{align*}
As $-\Re(Q^{3}P)P_{1}$ is the profile term, we will estimate 
\begin{equation}
\left\Vert -\tfrac{1}{4}(|w|^{4}-Q^{4})w_{1}+\Re(Q^{3}P)P_{1}\right\Vert _{L^{2}}\lesssim b^{2}.\label{eq:NL1 6.1 step 1}
\end{equation}

By \eqref{eq:CoercivityEstimateH2epsilon}, \eqref{eq:InterpolationEstimateW1,inf},
\eqref{eq:CoercivityEstimateH1}, and \eqref{eq:InterpolationEstimateW,inf},
we deduce 
\begin{align}
\|(4\Re(Q^{3}\eps)+M_{Q}(\widehat{\eps}))w_{1}\|_{L^{2}} & \lesssim\|(Q^{3}|\eps|+Q^{2}|\widehat{\eps}|^{2}+|\widehat{\eps}|^{4})w_{1}\|_{L^{2}}\nonumber \\
 & \lesssim\|Q^{2}\eps\|_{L^{2}}\|Qw_{1}\|_{L^{\infty}}+\|Q\widehat{\eps}\|_{L^{\infty}}^{2}\|w_{1}\|_{L^{2}}+\|\widehat{\eps}\|_{L^{\infty}}^{4}\|w_{1}\|_{L^{2}}\nonumber \\
 & \lesssim\lambda^{2}\|Qw_{1}\|_{L^{\infty}}+b^{2-}\lambda+b^{2}\nonumber \\
 & \lesssim b^{2}+\lambda^{2}\|Qw_{1}\|_{L^{\infty}}.\label{eq:NL1Estimate Quintic 1}
\end{align}
Using $w_{1}=P_{1}+\eps_{1}$ and \eqref{eq:InterpolationEstimateW1,inf +} $\|f\|_{L^{\infty}}^{2}\lesssim\|f\|_{L^{2}}\|\partial_{y}f\|_{L^{2}}$,
we have 
\begin{align}
\|Qw_{1}\|_{L^{\infty}}  \leq\|QP_{1}\|_{L^{\infty}}+\|Q\eps_{1}\|_{L^{\infty}}
 \lesssim \lambda. \label{eq:NL1Estimate Quintic 2}
\end{align}
By \eqref{eq:NL1Estimate Quintic 1} and \eqref{eq:NL1Estimate Quintic 2},
we deduce 
\begin{align}
\|(4\Re(Q^{3}\eps)+M_{Q}(\widehat{\eps}))w_{1}\|_{L^{2}}\lesssim b^{2}.\label{eq:NL1Estimate Quintic 4}
\end{align}
For $4\Re(Q^{3}P)w_{1}$, we decompose 
\begin{align}
4\Re(Q^{3}P)w_{1}=4\Re(Q^{3}P)P_{1}+4\Re(Q^{3}P)\eps_{1},\label{eq:NL1Estimate Quintic 5}
\end{align}
and, from \eqref{eq:CoercivityEstimateH2epsilon} and the definition
of $P$, 
\begin{align}
\|Q^{3}P\eps_{1}\|_{L^{2}}\lesssim b^{1-}\lambda^{2}\lesssim b^{2},\label{eq:NL1Estimate Quintic 6}
\end{align}
Therefore, \eqref{eq:NL1Estimate Quintic 4}, \eqref{eq:NL1Estimate Quintic 5},
and \eqref{eq:NL1Estimate Quintic 6} implies \eqref{eq:NL1 6.1 step 1}.

\textbf{Step 2.} \textbf{Estimate of cubic term $\widehat{\eps}|D|\overline{\widehat{\eps}}w_{1}$.}
In this step, we claim that 
\begin{align}
\|\widehat{\eps}|D|\overline{\widehat{\eps}}w_{1}\|_{L^{2}}\lesssim b^{2}.\label{eq:NL1 6.1 step 2}
\end{align}
In this term, there are no quadratic profile terms.

Using the decompositions $\wt\eps=P+\eps$ and $w_{1}=P_{1}+\eps_{1}$ with \eqref{eq:InterpolationEstimateW,inf}, we have
\begin{align*}
	\|\widehat{\eps}|D|\overline{\widehat{\eps}}w_{1}\|_{L^{2}}
	\leq
	\|\widehat{\eps}|D|(\overline{\widehat{\eps}}w_{1}-\ol{P}P_{1,e} )\|_{L^{2}}
	+
	\|\widehat{\eps}|D|(\ol{P}P_{1,e} )\|_{L^{2}}.
\end{align*}
For the second term, we have
\begin{align*}
	|D|(\ol{P}P_{1,e} )=\partial_y(Q^2) \cdot O((|\nu|+|\mu|)|\nu|)+ \partial_y(yQ^2)\cdot O((|\nu|+|\mu|)b),
\end{align*}
which yields
\begin{align*}
		\|\widehat{\eps}|D|(\ol{P}P_{1,e} )\|_{L^{2}}
		\lesssim b^{1-}\lambda\|Q^3\wt\eps\|_{L^2}+b\lambda \|Q^2\wt\eps\|_{L^2}
		\lesssim b^{-1}\lambda^2\lesssim b^2.
\end{align*}
For the first term, we note that
\begin{align*}
	(\overline{\partial_y\widehat{\eps}})w_{1}-(\ol{\partial_yP})P_{1,e}
	=&
	(\overline{\partial_y\eps})w_{1}+(\ol{\partial_yP})(w_1-P_{1,e}),
	\\
	\overline{\widehat{\eps}}\partial_yw_{1}-\ol{P}\partial_yP_{1,e}
	=&
	\overline{\eps}\partial_yP_{1,e}+\overline{\widehat{\eps}}\partial_yP_{1,o}+\overline{\widehat{\eps}}\partial_y\eps_{1}.
\end{align*}
By \eqref{eq:CoercivityEstimateH1}, \eqref{eq:InterpolationEstimateW,inf}, \eqref{eq:CoercivityEstimateH2epsilon}, and \eqref{eq:InterpolationEstimateW1,inf +}, we have
\begin{align}
	\|\widehat{\eps}|D|(\overline{\widehat{\eps}}w_{1}-\ol{P}P_{1,e} )\|_{L^{2}}
	\lesssim& \|\widehat{\eps}\|_{L^{\infty}}
	(\|(\overline{\partial_y\eps})w_{1}\|_{L^2}
	+\|(\ol{\partial_yP})(w_1-P_{1,e})\|_{L^2} \nonumber
	\\
	&\qquad\quad +\|\overline{\eps}\partial_yP_{1,e}\|_{L^2}
	+\|\overline{\widehat{\eps}}\partial_yP_{1,o}\|_{L^2}
	+\|\overline{\widehat{\eps}}\partial_y\eps_{1}\|_{L^2}) \nonumber
	\\
	\lesssim &b^2+\lambda^{\frac12}\|(\ol{\partial_yP})(w_1-P_{1,e})\|_{L^2}. \label{eq:NL1Estimate eDew1 1}
\end{align}
We control $\lambda^{\frac12}\|(\ol{\partial_yP})(w_1-P_{1,e})\|_{L^2}$ as follows:
\begin{align}
\lambda^{\frac12}\|(\ol{\partial_yP})(w_1-P_{1,e})\|_{L^2}
\leq& \lambda^{\frac12}\|\partial_yP_e\|_{L^\infty}(\|\eps_{1,e}\|_{L^2}+\|w_{1,o}\|_{L^2}) \nonumber
\\
&+ \lambda^{\frac12}(\|(\partial_yP_o)\eps_{1,e}\|_{L^2}+\|(\partial_yP_o)w_{1,o}\|_{L^2}) \nonumber
\\
\lesssim& b^2+b^{1-}\lambda^{\frac52}+\lambda^{\frac12}\|(\partial_yP_o)w_{1,o}\|_{L^2} \nonumber
\\
\lesssim & b^2+\lambda^{\frac12}\|(\partial_yP_o)w_{1,o}\|_{L^2}. \label{eq:NL1Estimate eDew1 2}
\end{align}
For $\lambda^{\frac12}\|(\partial_yP_o)w_{1,o}\|_{L^2}$,
using the pointwise bound \eqref{eq:ProfilePointwiseBound} and the decomposition
$w_{1,o}=P_{1,o}+\eps_{1,o}$,
\begin{align}
\lambda^{\frac12}\|(\partial_yP_o)w_{1,o}\|_{L^2} & \lesssim\lambda^{\frac{1}{2}}|\nu|(\|QP_{1,o}\|_{L^{2}}+\|Q\eps_1\|_{L^{2}})\nonumber \\
 & \lesssim 
 b^{1-}\lambda^{2} \lesssim b^{2}.\label{eq:NL1Estimate eDew1 3}
\end{align}
By \eqref{eq:NL1Estimate eDew1 1}--\eqref{eq:NL1Estimate eDew1 3}, we have \eqref{eq:NL1 6.1 step 2}.

\textbf{Step 3.} \textbf{Estimate of $\widehat{\eps}|D|Qw_{1}$}.
We claim that 
\begin{align}
\|\widehat{\eps}|D|Qw_{1}-P|D|QP_{1}\|_{L^{2}}\lesssim b^{2},\label{eq:NL1 6.1 step 3}
\end{align}
where the quadratic profile term is $P|D|QP_{1}$.

Using $(1+y^{2})\langle y\rangle^{-2}=1$ and \eqref{eq:CommuteHilbert},
we write 
\begin{align*}
\widehat{\eps}|D|Qw_{1} & =\widehat{\eps}\langle y\rangle^{-2}|D|Qw_{1}+\widehat{\eps}\langle y\rangle^{-2}y^{2}|D|Qw_{1}\\
 & =\widehat{\eps}\langle y\rangle^{-2}|D|Qw_{1}+\widehat{\eps}\langle y\rangle^{-2}y\mathcal{H}[y\partial_{y}(Qw_{1})].
\end{align*}
From \eqref{eq:CommuteHilbertDerivative} and the decompositions $\widehat{\eps}=P+\eps$,
$w_{1}=P_{1}+\eps_{1}$, we have 
\begin{align}
\widehat{\eps}|D|Qw_{1}-P|D|QP_{1}= & \widehat{\eps}\langle y\rangle^{-2}|D|Q\eps_{1}+\widehat{\eps}\langle y\rangle^{-2}y\mathcal{H}[y\partial_{y}(Q\eps_{1})]\nonumber \\
 & +\eps\langle y\rangle^{-2}|D|QP_{1}+\eps\langle y\rangle^{-2}y\mathcal{H}[y\partial_{y}(QP_{1})].\label{eq:NL1Estimate eDQw1 1}
\end{align}
Using the pointwise bound, $||D|QP_{1}|+|y\mathcal{H}[y\partial_{y}(QP_{1})]|\lesssim \lambda$,
we obtain 
\begin{align*}
\|\eqref{eq:NL1Estimate eDQw1 1}\|_{L^{2}}\lesssim\|\langle y\rangle^{-1}\widehat{\eps}\|_{L^{\infty}}(\|\partial_{y}(Q\eps_{1})\|_{L^{2}}+\|y\partial_{y}(Q\eps_{1})\|_{L^{2}})+\lambda\|Q^{2}\eps\|_{L^{2}},
\end{align*}
and by Lemma~\ref{LemmaNonlinearCoercivity}, we conclude \eqref{eq:NL1 6.1 step 3}.

\textbf{Step 4.} \textbf{Estimate of $Q|D|(\overline{\widehat{\eps}}w_{1})$}.
In this step, we claim that 
\begin{align}
\|Q|D|(\overline{\widehat{\eps}}w_{1})-\mathcal{P}\|_{L^{2}}\lesssim b^{2}.\label{eq:NL1 6.1 step 4}
\end{align}
We decompose $Q|D|\overline{\widehat{\eps}}w_{1}-\mathcal{P}$ into
three parts, 
\begin{align}
2(Q|D|\overline{\widehat{\eps}}w_{1}-\mathcal{P})= & Q\mathcal{H}\partial_{y}(\overline{\widehat{\eps}}_{e}w_{1,o})\label{eq:NL1Estimate QDew1 decompose 1}\\
 & +Q\mathcal{H}\partial_{y}(\overline{\widehat{\eps}}_{o}w_{1,e}-\ol P_{o}P_{1,e})\label{eq:NL1Estimate QDew1 decompose 2}\\
 & +Q\mathcal{H}\partial_{y}(\overline{\widehat{\eps}}w_{1}-\ol PP_{1})_{e}.\label{eq:NL1Estimate QDew1 decompose 3}
\end{align}

To estimate \eqref{eq:NL1Estimate QDew1 decompose 1}, we use $(1+y^{2})\langle y\rangle^{-2}=1$
to further decompose and locate extra decay in the Hilbert transform.
Using \eqref{eq:CommuteHilbert} with $\int_{\R}f_{o}=0$, we have
\begin{align}
\eqref{eq:NL1Estimate QDew1 decompose 1}= & Q\mathcal{H}[\langle y\rangle^{-2}\partial_{y}(\overline{\widehat{\eps}}_{e}w_{1,o})]+Q\mathcal{H}[y^{2}\langle y\rangle^{-2}\partial_{y}(\overline{\widehat{\eps}}_{e}w_{1,o})]\nonumber \\
= & Q\mathcal{H}[\langle y\rangle^{-2}\partial_{y}(\overline{\widehat{\eps}}_{e}w_{1,o})]+yQ\mathcal{H}[y\langle y\rangle^{-2}\partial_{y}(\overline{\widehat{\eps}}_{e}w_{1,o})]\label{eq:NL1Estimate QDew1 decompose 1-1}
\end{align}
For the second term of \eqref{eq:NL1Estimate QDew1 decompose 1-1},
we have 
\begin{align}
\|yQ\mathcal{H}[y\langle y\rangle^{-2}\partial_{y}(\overline{\widehat{\eps}}_{e}w_{1,o})]\|_{L^{2}} & \lesssim\|y\langle y\rangle^{-2}\partial_{y}(\overline{\widehat{\eps}}_{e})w_{1,o})\|_{L^{2}}+\|y\langle y\rangle^{-2}\overline{\widehat{\eps}}_{e}(\partial_{y}w_{1,o})\|_{L^{2}}\nonumber \\
 & \lesssim\|\partial_{y}\overline{\widehat{\eps}}_{e}\|_{L^{\infty}}\|Qw_{1,o}\|_{L^{2}}+\|Q\overline{\widehat{\eps}}_{e}\|_{L^{\infty}}\|\partial_{y}w_{1,o}\|_{L^{2}}\nonumber \\
 & \lesssim b^{2}.\label{eq:NL1Estimate QDew1 decompose 1-2}
\end{align}
Applying a similar argument to the first term of \eqref{eq:NL1Estimate QDew1 decompose 1-1},
we conclude 
\begin{align}
\|\eqref{eq:NL1Estimate QDew1 decompose 1}\|_{L^{2}}\lesssim b^{2}.\label{eq:NL1Estimate QDew1 decompose 1 goal}
\end{align}

For \eqref{eq:NL1Estimate QDew1 decompose 2} and \eqref{eq:NL1Estimate QDew1 decompose 3},
we first note that 
\begin{align}
\overline{\widehat{\eps}}_{o}w_{1,e}-\overline{P}_{o}P_{1,e}=\overline{\widehat{\eps}}_{o}\eps_{1,e}+\overline{\eps}_{o}P_{1,e},\quad\overline{\widehat{\eps}}w_{1}-\overline{P}P_{1}=\overline{\widehat{\eps}}\eps_{1}+\overline{\eps}P_{1}.\label{eq:NL1Estimate QDew1 ewpp1 equal eeep}
\end{align}
By \eqref{eq:NL1Estimate QDew1 ewpp1 equal eeep}, we rewrite \eqref{eq:NL1Estimate QDew1 decompose 2}
as 
\begin{align}
\eqref{eq:NL1Estimate QDew1 decompose 2}=Q\mathcal{H}\partial_{y}(\overline{\widehat{\eps}}_{o}\eps_{1,e})+Q\mathcal{H}\partial_{y}(\overline{\eps}_{o}P_{1,e}).\label{eq:NL1Estimate QDew1 decompose 2-1}
\end{align}
For the first term of \eqref{eq:NL1Estimate QDew1 decompose 2-1},
we apply a similar argument to \eqref{eq:NL1Estimate QDew1 decompose 1}.
Using $(1+y^{2})\langle y\rangle^{-2}=1$ and \eqref{eq:CommuteHilbert}
with $\int_{\R}f_{o}=0$, we have 
\begin{align}
Q\mathcal{H}\partial_{y}(\overline{\widehat{\eps}}_{o}\eps_{1,e})=Q\mathcal{H}[\langle y\rangle^{-2}\partial_{y}(\overline{\widehat{\eps}}_{o}\eps_{1,e})]+yQ\mathcal{H}[y\langle y\rangle^{-2}\partial_{y}(\overline{\widehat{\eps}}_{o}\eps_{1,e})].\label{eq:NL1Estimate QDew1 decompose 2-2}
\end{align}
By a similar argument to \eqref{eq:NL1Estimate QDew1 decompose 1-2}
with substituting $\eps_{1}$ for $w_{1}$ and $\ol{\wt{\eps}}_{o}$
for $\ol{\wt{\eps}}_{e}$, we have 
\begin{align*}
\|\eqref{eq:NL1Estimate QDew1 decompose 2-2}\|_{L^{2}}\lesssim b^{2}\lambda^{\frac{1}{2}-}.
\end{align*}
For the second term of \eqref{eq:NL1Estimate QDew1 decompose 2-1},
we have 
\begin{align*}
\|Q\mathcal{H}\partial_{y}(\overline{\eps}_{o}P_{1,e})\|_{L^{2}}\lesssim\|\partial_{y}(\overline{\eps}_{o}P_{1,e})\|_{L^{2}}\lesssim (|v|+|\mu|)\|\eps\|_{\dot{\mathcal{H}}^{2}}\lesssim b^{2},
\end{align*}
which give a better bound than to prove. Hence, we conclude 
\begin{align}
\|\eqref{eq:NL1Estimate QDew1 decompose 2}\|_{L^{2}}\lesssim b^{2}.\label{eq:NL1Estimate QDew1 decompose 2 goal}
\end{align}

Now, we estimate \eqref{eq:NL1Estimate QDew1 decompose 3}. Again,
using $(1+y^{2})\langle y\rangle^{-2}=1$ and \eqref{eq:CommuteHilbertDerivative},
we have 
\begin{align}
\eqref{eq:NL1Estimate QDew1 decompose 3}= & Q\mathcal{H}\partial_{y}[\langle y\rangle^{-2}(\overline{\widehat{\eps}}w_{1}-\ol PP_{1})_{e}]+Q\mathcal{H}\partial_{y}[y^{2}\langle y\rangle^{-2}(\overline{\widehat{\eps}}w_{1}-\ol PP_{1})_{e}]\nonumber \\
= & Q\mathcal{H}\partial_{y}[\langle y\rangle^{-2}(\overline{\widehat{\eps}}w_{1}-\ol PP_{1})_{e}]+yQ\mathcal{H}\partial_{y}[y\langle y\rangle^{-2}(\overline{\widehat{\eps}}w_{1}-\ol PP_{1})_{e}]\label{eq:NL1Estimate QDew1 decompose 3-1}\\
 & +Q\mathcal{H}[y\langle y\rangle^{-2}(\overline{\widehat{\eps}}w_{1}-\ol PP_{1})_{e}],\label{eq:NL1Estimate QDew1 decompose 3-2}
\end{align}
where the second line is well justified since 
\begin{align*}
	y\langle y\rangle^{-2}(\overline{\widehat{\eps}}w_{1}-\ol PP_{1})_{e}\in H^{1},\text{ and }
	\partial_{y}\big((y\langle y\rangle^{-2}(\overline{\widehat{\eps}}w_{1}-\ol PP_{1})_{e}\big)\in L^{2}(xdx).
\end{align*}
For \eqref{eq:NL1Estimate QDew1 decompose 3-1}, using \eqref{eq:NL1Estimate QDew1 ewpp1 equal eeep}
and applying a similar argument to \eqref{eq:NL1Estimate QDew1 decompose 1-2},
we have 
\begin{align*}
\|\eqref{eq:NL1Estimate QDew1 decompose 3-1}\|_{L^{2}}\lesssim b^{2}.
\end{align*}
Now, we estimate the trickiest case \eqref{eq:NL1Estimate QDew1 decompose 3-2}.
We have 
\begin{align}
\|\eqref{eq:NL1Estimate QDew1 decompose 3-2}\|_{L^{2}}\lesssim\|Q\|_{L^{2}}\|\mathcal{H}[y\langle y\rangle^{-2}(\overline{\widehat{\eps}}w_{1}-\overline{P}P_{1})_{e}]\|_{L^{\infty}}.\label{eq:NL1Estimate QDew1 decompse lastterm evenpart}
\end{align}
We further decompose 
\begin{align*}
(\overline{\widehat{\eps}}w_{1}-\overline{P}P_{1})_{e}= & \overline{\widehat{\eps}}_{e}\eps_{1,e}+\overline{\eps}_{o}w_{1,o}+\overline{P_{o}}\eps_{1,o}+\overline{\eps}_{e}P_{1,e}\\
= & E_{1}+E_{2},
\end{align*}
where 
\begin{align*}
E_{1} & =\overline{\widehat{\eps}}_{e}\eps_{1,e}+\overline{\eps}_{o}w_{1,o},\\
E_{2} & =\overline{P_{o}}\eps_{1,o}+\overline{\eps}_{e}P_{1,e}.
\end{align*}
Thus, we decompose \eqref{eq:NL1Estimate QDew1 decompse lastterm evenpart}
by 
\begin{align}
\eqref{eq:NL1Estimate QDew1 decompse lastterm evenpart}\lesssim & \|\mathcal{H}[y\langle y\rangle^{-2}E_{1}]\|_{L^{\infty}}\label{eq:NL1Estimate QDew1 decompse lastterm evenpartFurther1}\\
&+  \|\mathcal{H}[y\langle y\rangle^{-2}E_{2}]\|_{L^{\infty}}.\label{eq:NL1Estimate QDew1 decompse lastterm evenpartFurther2}
\end{align}
For \eqref{eq:NL1Estimate QDew1 decompse lastterm evenpartFurther1},
thanks to the Sobolev embedding, we have 
\begin{align}
\eqref{eq:NL1Estimate QDew1 decompse lastterm evenpartFurther1}\lesssim\|y\langle y\rangle^{-2}E_{1}\|_{L^{2}}^{\frac{1}{2}}\|\partial_{y}[y\langle y\rangle^{-2}E_{1}]\|_{L^{2}}^{\frac{1}{2}}.\label{eq:NL1Estimate QDew1 decompose 3 goal1-1}
\end{align}
For the first term of RHS of \eqref{eq:NL1Estimate QDew1 decompose 3 goal1-1},
by Lemma~\ref{LemmaNonlinearCoercivity} we have 
\begin{align*}
\|y\langle y\rangle^{-2}E_{1}\|_{L^{2}}\leq & \|y\langle y\rangle^{-2}\overline{\widehat{\eps}}_{e}\eps_{1,e}\|_{L^{2}}+\|y\langle y\rangle^{-2}\overline{\eps}_{o}w_{1,o}\|_{L^{2}}\\
\lesssim & \|\langle y\rangle^{-1}\overline{\widehat{\eps}}_{e}\|_{L^{\infty}}\|\eps_{1,e}\|_{L^{2}}+\|\langle y\rangle^{-1}\overline{\eps}_{o}\|_{L^{\infty}}\|w_{1,o}\|_{L^{2}}\\
\lesssim & \lambda b+\la b\sim\la^{2+\frac{1}{2}},
\end{align*}
and for the second term of the RHS of \eqref{eq:NL1Estimate QDew1 decompose 3 goal1-1},
we have 
\begin{align*}
\|\partial_{y}[y\langle y\rangle^{-2}E_{1}]\|_{L^{2}}\lesssim & \|\partial_{y}[y\langle y\rangle^{-2}\overline{\widehat{\eps}}_{e}\eps_{1,e}]\|_{L^{2}}+\|\partial_{y}[y\langle y\rangle^{-2}\overline{\eps}_{o}w_{1,o}]\|_{L^{2}}\\
\lesssim & \|\partial_{y}\overline{\widehat{\eps}}_{e}\|_{L^{\infty}}\|\langle y\rangle^{-1}\eps_{1,e}\|_{L^{2}}+\|\langle y\rangle^{-1}\overline{\widehat{\eps}}_{e}\|_{L^{\infty}}\|\eps_{1,e}\|_{\dot{\mathcal{H}}^{1}}\\
 & +\|\overline{\eps}_{o}\|_{\dot{\mathcal{H}}^{2}}\|w_{1,o}\|_{L^{\infty}}+\|\langle y\rangle^{-1}\overline{\eps}_{o}\|_{L^{\infty}}\|\partial_{y}\eps_{1,o}\|_{L^{2}}+b\|\overline{\eps}_{o}\|_{\dot{\mathcal{H}}^{2}}\\
\lesssim & \lambda^{3+\frac{1}{2}}.
\end{align*}
Here we used $\partial_{y}w_{1,o}=\partial_{y}\eps_{1,o}+\partial_{y}P_{1,o}\sim\partial_{y}\eps_{1,o}+(|b|+|\eta|)Q$.
Therefore, we have 
\begin{align}
\eqref{eq:NL1Estimate QDew1 decompse lastterm evenpartFurther1}\lesssim\lambda^{3}\sim b^{2}.\label{eq:NL1Estimate QDew1 decompose 3 goal1}
\end{align}
For \eqref{eq:NL1Estimate QDew1 decompse lastterm evenpartFurther2},
again by Sobolev embedding, we have 
\begin{align}
\eqref{eq:NL1Estimate QDew1 decompse lastterm evenpartFurther2}\lesssim\|y\langle y\rangle^{-2}E_{2}\|_{L^{2}}^{\frac{1}{2}}\|\partial_{y}[y\langle y\rangle^{-2}E_{2}]\|_{L^{2}}^{\frac{1}{2}}.\label{eq:NL1Estimate QDew1 decompose 3 goal2-1}
\end{align}
We have 
\begin{align*}
\|y\langle y\rangle^{-2}E_{2}\|_{L^{2}}\lesssim|\nu|\|\eps_{1,o}\|_{\dot{\mathcal{H}}^{1}}+(|\nu|+|\mu|)\|\eps_{e}\|_{\dot{\mathcal{H}}^{2}} & \lesssim b^2,\\
\|\partial_{y}[y\langle y\rangle^{-2}E_{2}]\|_{L^{2}}\lesssim|\nu|\|\eps_{1,o}\|_{\dot{\mathcal{H}}^{1}}+(|\nu|+|\mu|)\|\eps_{e}\|_{\dot{\mathcal{H}}^{2}} & \lesssim b^2,
\end{align*}
and 
\begin{align*}
\|y\langle y\rangle^{-2}E_{2}\|_{L^{2}}^{\frac{1}{2}}\|\partial_{y}[y\langle y\rangle^{-2}E_{2}]\|_{L^{2}}^{\frac{1}{2}}\lesssim b^2.
\end{align*}
Thus, we have 
\begin{align}
\eqref{eq:NL1Estimate QDew1 decompse lastterm evenpartFurther2}\lesssim b^{2}.\label{eq:NL1Estimate QDew1 decompose 3 goal2}
\end{align}
Therefore, collecting \eqref{eq:NL1Estimate QDew1 decompose 1 goal},
\eqref{eq:NL1Estimate QDew1 decompose 2 goal}, \eqref{eq:NL1Estimate QDew1 decompose 3 goal1},
and \eqref{eq:NL1Estimate QDew1 decompose 3 goal2}, we conclude \eqref{eq:NL1 6.1 step 4}.

\textbf{Step 5.} By the previous steps, \eqref{eq:NL1 6.1 step 1},
\eqref{eq:NL1 6.1 step 2}, \eqref{eq:NL1 6.1 step 3}, and \eqref{eq:NL1 6.1 step 4},
we conclude 
\begin{align}
\|\textnormal{NL}_{1}-(-\Re(Q^{3}P)P_{1}+P|D|QP_{1}+\mathcal{P})\|_{L^{2}}\lesssim b^{2}.\label{eq:NL1Estimate  without Profiles}
\end{align}
From the definition of $\mathcal{P}$ \eqref{eq:def calP}, we have that $\mathcal{P}+Q|D|(\overline{P}_{e}P_{1,o})=Q|D|\overline{P}P_{1}$
formally. In view of \eqref{eq:QDPP1}, the term $Q|D|(\overline{P}_{e}P_{1,o})$ corresponds to
\begin{align*}
	-(b^{2}+\eta^{2})\tfrac{1}{4}Q\mathcal{H}(1)-i(i\tfrac{b^{2}}{2}+\tfrac{b\eta}{2})\tfrac{1}{(1+y^{2})^{2}}yQ.
\end{align*}
Combining this and \eqref{eq:NL1 Quadratic Profile another form}, we observe that
\begin{equation}
	\begin{aligned}
		&-\Re(Q^{3}P)P_{1}+P|D|QP_{1}+\mathcal{P}-i(i\tfrac{b^{2}}{2}+\tfrac{b\eta}{2})\tfrac{1}{(1+y^{2})^{2}}yQ
		\\
		&=
		-i(b\Lambda_{-1}+\tfrac{\eta}{2}i-\nu\partial_{y})P_{1}-i(\tfrac{3b^{2}}{2}+\tfrac{\eta^{2}}{2})\tfrac{i}{2}yQ-ib\eta\tfrac{1}{2}yQ+ib\nu\tfrac{i}{2}Q+ib\mu\tfrac{1}{2}Q.
		\\
		&=
		-[i(b\Lambda_{-1}+\tfrac{\eta}{2}i-\nu\partial_{y})P_{1}+b\nu\tfrac{1}{2}Q-b\mu\tfrac{i}{2}Q]_{e}
		\\
		&\quad -i(i\tfrac{b^{2}}{2}+\tfrac{b\eta}{2})\tfrac{1}{(1+y^{2})^{2}}yQ
		-\tfrac{1}{2}(\nu^2-i\nu\mu)Q_y.
	\end{aligned}\label{eq:NL1Estimate Profiles collection}
\end{equation}
Here, the last equality comes from a direct computation. 
Thus, by \eqref{eq:NL1Estimate  without Profiles} and
\eqref{eq:NL1Estimate Profiles collection}, we derive \eqref{eq:NL1Estimate} and \eqref{eq:NL1eWithoutProfileEstimate}. This finishes the proof of Lemma~\ref{LemmaNonlinearEstimate1}. \qed


\appendix

\section{Subcoercivity}

\label{AppendixSubcoercivity} In this subsection, we prove the subcoercivity
and coercivity of adapted derivatives introduced in Section~\ref{Section 3 Linearization}.
The adapted function spaces $\dot{\mathcal{H}}^{1}$ and $\dot{\mathcal{H}}^{2}$
are naturally derived from the coercivity properties.

We mainly use the following weighted Hardy inequalities: 
\begin{lem}[{{\cite[Lemma A.1]{KimKwon2020blowup}}}]
\label{WeightedHardyIneq} Let $0<r_{1}<r_{2}<\infty$. In addition,
let $\varphi:[r_{1},r_{2}]\to\mathbb{R}_{+}$ be a $C^{1}$ weight
function such that $\partial_{r}\varphi$ is non-vanishing and $\varphi\lesssim|r\partial_{r}\varphi|$.
Then, for smooth $f:[r_{1},r_{2}]\to\mathbb{C}$, we have 
\begin{align*}
\int_{r_{1}}^{r_{2}}\left|\frac{f}{r}\right|^{2}|r^{2}\partial_{r}\varphi|dr\lesssim\int_{r_{1}}^{r_{2}}|\partial_{r}f|^{2}r\varphi dr+\begin{cases}
\varphi(r_{2})|f(r_{2})|^{2} & if\ \partial_{r}\varphi>0,\\
\varphi(r_{1})|f(r_{1})|^{2} & if\ \partial_{r}\varphi<0.
\end{cases}
\end{align*}
\end{lem}

\begin{lem}[{{\cite[Corollary A.3]{KimKwon2020blowup}}}]
\label{WeightedLogHardyIneq} Let $k,l\in\mathbb{R}$ and $k\neq l$.
Then, we have 
\begin{align*}
\int_{r_{1}}^{r_{2}}\left|\frac{f}{r^{k+\frac{1}{2}}}\right|^{2}dr\lesssim\int_{r_{1}}^{r_{2}}\bigg|\frac{(\partial_{r}-\frac{l}{r})f}{r^{k-\frac{1}{2}}}\bigg|^{2}dr+\begin{cases}
|(r_{2})^{-k}f(r_{2})|^{2} & if\ l>k,\\
|(r_{1})^{-k}f(r_{1})|^{2} & if\ l<k.
\end{cases}
\end{align*}
\end{lem}

In the following, we will argue for $x\in\R_{+}$. Then, the case $x\in\R_{-}$
follows in a similar way. So, when we restrict ourselves to $\R_{+}$, we use the variable
$r$ instead of $x.$ 
\begin{lem}[Subcoercivity for $L_{Q}$ on $\dot{\mathcal{H}}^{1}$]
\label{LemmaAppendix LQ subcoer} For $v\in\dot{\mathcal{H}}^{1}$,
we have 
\begin{align*}
\|L_{Q}v\|_{L^{2}}+\|{\bf 1}_{|x|\lesssim1}v\|_{L^{2}}\sim\|v\|_{\dot{\mathcal{H}}^{1}}.
\end{align*}
\end{lem}

\begin{proof}
We recall 
\begin{align*}
L_{Q}v=\partial_{x}v+\frac{1}{2}\mathcal{H}(Q^{2})v+\mathcal{H}(\text{Re}(Qv))Q=\mathbf{D}_{Q}v+\mathcal{H}(\text{Re}(Qv))Q.
\end{align*}
Using \eqref{eq:HilbertDecayExchange1}, \eqref{eq:CommuteHilbert},
and the Hölder inequality with the fact $Q^{\frac{1}{2}+}\in L^{2}$,
we observe that 
\begin{align*}
\|\mathcal{H}(\text{Re}(Qv))Q\|_{L^{2}}\lesssim & \|xQ\mathcal{H}(\text{Re}(xQ^{3}v))\|_{L^{2}}+\|Q\mathcal{H}(\text{Re}(Q^{3}v))\|_{L^{2}}\\
 & +{\textstyle \left\Vert Q\int_{\bbR}\text{Re}(xQ^{3}v)dx\right\Vert _{L^{2}}}\\
\lesssim & \|{\bf 1}_{|x|\leq r_{0}}v\|_{L^{2}}+r_{0}^{-\frac{1}{2}+}\|v\|_{\dot{\mathcal{H}}^{1}}.
\end{align*}
Thus, taking sufficiently large $r_{0}>0$, it suffices to show 
\begin{align}
\|\mathbf{D}_{Q}v\|_{L^{2}}+\|{\bf 1}_{|x|\lesssim1}v\|_{L^{2}}\sim\|v\|_{\dot{\mathcal{H}}^{1}}.\label{eq:A. Lemma LQ 1}
\end{align}
We will apply a similar argument to that in \cite[Lemma A.3]{KimKwonOh2020blowup}.
By symmetry, it suffices to prove \eqref{eq:A. Lemma LQ 1} on $\mathbb{R}_{+}$.
Indeed, we claim 
\begin{align}
\|\mathbf{D}_{Q}v\|_{L^{2}(\mathbb{R}_{+})}+\|{\bf 1}_{r\lesssim1}v\|_{L^{2}(\mathbb{R}_{+})}\sim\|v\|_{\dot{\mathcal{H}}^{1}(\mathbb{R}_{+})}.\label{eq:AppendixLQSubcoercivityReduce}
\end{align}
For simplicity, we omit $\mathbb{R}_{+}$. $(\lesssim)$ is trivial.

For $(\gtrsim)$, we divide the proof into two regions, $r\lesssim1$
and $r\gtrsim1$. For $r\lesssim1$, using $\|\mathbf{D}_{Q}v\|_{L^{2}}=\|\partial_{r}v\|_{L^{2}}+O(\left\Vert Qv\right\Vert _{L^{2}})$,
we directly deduce \eqref{eq:AppendixLQSubcoercivityReduce} on $r\lesssim1$.

Now, we prove \eqref{eq:AppendixLQSubcoercivityReduce} on $r\gtrsim1$.
We note that $\mathbf{D}_{Q}=Q\partial_{x}Q^{-1}$. Applying Lemma~\ref{WeightedHardyIneq}
for $f=Q^{-1}v$ with $\varphi=Q^{2}$ in $r\geq10$, we have 
\begin{align}
\begin{split}\left\Vert {\bf 1}_{[r_{0},\infty)}r^{-1}v\right\Vert _{L^{2}}^{2}={\textstyle \int_{r_{0}}^{\infty}r^{-2}(Q^{-1}v)^{2}Q^{2}dr\lesssim} & {\textstyle \int_{r_{0}}^{\infty}r^{-2}(Q^{-1}v)^{2}|r\partial_{r}(Q^{2})|rdr}\\
\lesssim & \|{\bf 1}_{[r_{0},\infty)}\mathbf{D}_{Q}v\|_{L^{2}}^{2}+|v(r_{0})|^{2},
\end{split}
\label{eq:A. Lemma LQ 2}
\end{align}
for $r_{0}\geq10$. Here, we used $|r\partial_{r}(Q^{2})|\sim r^{-2}$
for $r>10$. Averaging the boundary term of \eqref{eq:A. Lemma LQ 2},
we deduce 
\begin{align*}
\left\Vert {\bf 1}_{[20,\infty)}r^{-1}v\right\Vert _{L^{2}}^{2}\lesssim\|{\bf 1}_{[10,\infty)}\mathbf{D}_{Q}v\|_{L^{2}}^{2}+\|{\bf 1}_{[10,20]}v\|_{L^{2}}^{2}.
\end{align*}
We also have 
\begin{align*}
\|{\bf 1}_{r\gtrsim1}\partial_{r}v\|_{L^{2}}^{2}\leq\|{\bf 1}_{r\gtrsim1}\mathbf{D}_{Q}v\|_{L^{2}}^{2}+\|{\bf 1}_{r\gtrsim1}r^{-1}v\|_{L^{2}}^{2}\lesssim\|{\bf 1}_{r\gtrsim1}\mathbf{D}_{Q}v\|_{L^{2}}^{2}+\|{\bf 1}_{r\sim1}v\|_{L^{2}}^{2},
\end{align*}
and we conclude the proof. 
\end{proof}
Using this subcoercivity and the kernel characterization of Proposition~\ref{PropKernel mathcalLQ},
we prove the coercivity for $L_{Q}$ with Proposition~\ref{PropKernel mathcalLQ}. 
\begin{proof}[Proof of Proposition~\ref{PropCoercivityLQ}]
We note that $(\lesssim)$ immediately comes from the subcoercivity
Lemma~\ref{LemmaAppendix LQ subcoer}.

For $(\gtrsim)$, we suppose that it is not true. Then, we can choose
a sequence $\{v_{n}\}_{n\in\mathbb{N}}\subset\dot{\mathcal{H}}^{1}$
such that $\|L_{Q}v_{n}\|_{L^{2}}=\frac{1}{n}$, $\|v_{n}\|_{\dot{\mathcal{H}}^{1}}=1$,
and $(\psi_{j},v_{n})_{r}=0$ for $j=1,2,3$. After passing to a subsequence,
we can find $v_{\infty}\in\dot{\mathcal{H}}^{1}$ such that $v_{n}\rightharpoonup v_{\infty}$
weakly in $\dot{\mathcal{H}}^{1}$ and $v_{n}\to v_{\infty}$ strongly
in $L_{loc}^{2}$. By weak convergence, we have $L_{Q}v_{\infty}=0$
and $(\psi_{j},v_{\infty})_{r}=0$ for $j=1,2,3$. By the kernel of
$L_{Q}$, Proposition~\ref{PropKernel mathcalLQ}, and the orthogonality
conditions, we deduce that $v_{\infty}=0$. From the subcoercivity,
Lemma~\ref{LemmaAppendix LQ subcoer}, we deduce that $\liminf_{n\to\infty}\|{\bf 1}_{|x|\lesssim1}v_{n}\|_{L^{2}}\gtrsim1$.
However, this is a contradiction with $v_{\infty}=0$ and $v_{n}\to v_{\infty}$
strongly in $L_{loc}^{2}$. 
\end{proof}
\begin{lem}[Subcoercivity for $A_{Q}$ on $\dot{\mathcal{H}}^{1}$]
\label{LemmaAppendix AQ Subcoercivity} For $v\in\mathcal{S}$, we
have 
\begin{align*}
\|A_{Q}v\|_{L^{2}}+\||D|^{\frac{1}{2}}Qv\|_{L^{2}}\sim\|v\|_{\dot{\mathcal{H}}^{1}}.
\end{align*}
\end{lem}

\begin{proof}
We recall that 
\begin{align*}
A_{Q}^{*}A_{Q}f=-\partial_{xx}f+\tfrac{1}{4}Q^{4}f-Q|D|Qf.
\end{align*}
We have 
\begin{align*}
\|A_{Q}v\|_{L^{2}}^{2}=(A_{Q}^{*}A_{Q}v,v)_{r}=\|\partial_{x}v\|_{L^{2}}^{2}+\tfrac{1}{4}\|Q^{2}v\|_{L^{2}}^{2}-(Q|D|Qv,v)_{r},
\end{align*}
and 
\begin{align*}
(Q|D|Qv,v)_{r}=\||D|^{\frac{1}{2}}Qv\|_{L^{2}}^{2}.
\end{align*}
That is, we have 
\begin{align}
\|A_{Q}v\|_{L^{2}}+\||D|^{\frac{1}{2}}Qv\|_{L^{2}}=\|\partial_{x}v\|_{L^{2}}+\tfrac{1}{4}\|Q^{2}v\|_{L^{2}}^{2}.\label{eq:AppendixAQSubcoercivityEqual}
\end{align}
By \eqref{eq:AppendixAQSubcoercivityEqual}, we can deduce $(\lesssim)$.

For $(\gtrsim)$, again using \eqref{eq:AppendixAQSubcoercivityEqual},
it suffices to show 
\begin{align*}
\|\partial_{x}v\|_{L^{2}}+\|Q^{2}v\|_{L^{2}}\gtrsim\|Qv\|_{L^{2}}.
\end{align*}
We also have $\|Qv\|_{L^{2}}^{2}=\|Qv{\bf 1}_{|x|\lesssim1}\|_{L^{2}}^{2}+\|Qv{\bf 1}_{|x|\gtrsim1}\|_{L^{2}}^{2}$.
For $\|Qv{\bf 1}_{|x|\lesssim1}\|_{L^{2}}$, we have $\|Qv{\bf 1}_{|x|\lesssim1}\|_{L^{2}}\lesssim\|Q^{2}v{\bf 1}_{|x|\lesssim1}\|_{L^{2}}$.
For $\|Qv{\bf 1}_{|x|\gtrsim1}\|_{L^{2}}\sim\||x|^{-1}v{\bf 1}_{|x|\gtrsim1}\|_{L^{2}}$,
we apply Lemma~\ref{WeightedHardyIneq} with $\varphi=r^{-1}$. Then,
we have 
\begin{align*}
\||x|^{-1}v{\bf 1}_{|x|\gtrsim1}\|_{L^{2}}\lesssim\|\partial_{x}v{\bf 1}_{|x|\gtrsim1}\|_{L^{2}}+\|v{\bf 1}_{|x|\sim1}\|_{L^{2}},
\end{align*}
and this proves $(\gtrsim)$. 
\end{proof}
We note that we do not use the coercivity for $A_{Q}$. It is possible
to show a coercivity of $A_{Q}$ under a suitable orthogonality, but
we do not prove it here. Next, we show an equivalence of adapted Sobolev
norm. 
\begin{lem}[Comparison of $\dot{\mathcal{H}}^{2}$ and $\dot{H}^{2}$]
\label{LemmaH2comparison} We have 
\begin{align*}
\|\partial_{xx}v\|_{L^{2}}+\|{\bf 1}_{|x|\sim1}(|\partial_{x}v|+|v|)\|_{L^{2}}\sim\|v\|_{\dot{\mathcal{H}}^{2}}.
\end{align*}
\end{lem}

\begin{proof}
It suffices to show $(\gtrsim)$ from its definition. Indeed, we claim
that 
\begin{align}
\left\Vert \langle x\rangle^{-1}\langle v\rangle_{-1}\right\Vert _{L^{2}}^{2}\lesssim\|\partial_{xx}v\|_{L^{2}}^{2}+\|{\bf 1}_{|x|\lesssim1}\partial_{x}v\|_{L^{2}}^{2}+\|{\bf 1}_{|x|\lesssim1}v\|_{L^{2}}^{2}.\label{eq:Appendixs cal H2 H2 goal}
\end{align}
We decompose the LHS of \eqref{eq:Appendixs cal H2 H2 goal} by $|x|\geq1$
and $|x|\leq1$, 
\begin{align*}
\left\Vert \langle x\rangle^{-1}\langle v\rangle_{-1}\right\Vert _{L^{2}}^{2} & =\left\Vert {\bf 1}_{|x|\lesssim1}\langle x\rangle^{-1}\langle v\rangle_{-1}\right\Vert _{L^{2}}^{2}+\left\Vert {\bf 1}_{|x|\gtrsim1}\langle x\rangle^{-1}\langle v\rangle_{-1}\right\Vert _{L^{2}}^{2}\\
 & \sim\left\Vert {\bf 1}_{|x|\lesssim1}(|v|+|\partial_{x}v|)\right\Vert _{L^{2}}^{2}+\left\Vert {\bf 1}_{|x|\gtrsim1}|x|^{-2}(||v|+|x\partial_{x}v||)\right\Vert _{L^{2}}^{2}.
\end{align*}
By symmetry, we restrict ourselves to $0<x\leq1$ and $x\geq1$. For $r=x\gtrsim1$,
we claim that 
\begin{align}
\left\Vert {\bf 1}_{r\gtrsim1}|r|^{-2}(||v|+|r\partial_{r}v||)\right\Vert _{L^{2}}\lesssim\|\partial_{rr}v\|_{L^{2}}+\|{\bf 1}_{r\sim1}\partial_{r}v\|_{L^{2}}+\|{\bf 1}_{r\sim1}v\|_{L^{2}}.\label{eq:Appendixs cal H2 H2 r geq1 goal}
\end{align}
By Lemma~\ref{WeightedLogHardyIneq} with $k=\frac{3}{2}$ and $l=1$,
we have 
\begin{align}
\left\Vert {\bf 1}_{r\geq r_{0}}r^{-2}v\right\Vert _{L^{2}}\lesssim\left\Vert {\bf 1}_{r\geq r_{0}}r^{-1}\left(\partial_{r}-\tfrac{1}{r}\right)v\right\Vert _{L^{2}}+|v(r_{0})|.\label{eq:A calH H compare lemma 1}
\end{align}
Averaging the boundary term of \eqref{eq:A calH H compare lemma 1},
we have 
\begin{align}
\left\Vert {\bf 1}_{r\gtrsim1}r^{-2}v\right\Vert _{L^{2}}\lesssim\left\Vert {\bf 1}_{r\gtrsim1}r^{-1}\left(\partial_{r}-\tfrac{1}{r}\right)v\right\Vert _{L^{2}}+\left\Vert {\bf 1}_{r\sim1}v\right\Vert _{L^{2}}.\label{eq:Appendixs cal H2 H2 r geq1 r-2v auxiliary1}
\end{align}
Applying Lemma~\ref{WeightedHardyIneq} with $\varphi=r^{-3}$ and
$f=r(\partial_{r}-\frac{1}{r})v$, we have 
\begin{align}
\left\Vert {\bf 1}_{r\geq1}r^{-1}\left(\partial_{r}-\tfrac{1}{r}\right)v\right\Vert _{L^{2}} & \lesssim\left\Vert r^{-1}\partial_{r}[r\left(\partial_{r}-\tfrac{1}{r}\right)v]\right\Vert _{L^{2}}+\|{\bf 1}_{r\sim1}(|\partial_{r}v|+|v|)\|_{L^{2}}\nonumber \\
 & =\|\partial_{rr}v\|_{L^{2}}+\|{\bf 1}_{r\sim1}(|\partial_{r}v|+|v|)\|_{L^{2}}.\label{eq:Appendixs cal H2 H2 r geq1 r-2v auxiliary2}
\end{align}
We used $r^{-1}\partial_{r}r=\partial_{r}+\frac{1}{r}$ and $(\partial_{r}+\frac{1}{r})(\partial_{r}-\frac{1}{r})=\partial_{rr}$.
Since $r^{-1}\partial_{r}v=\partial_{r}(r^{-1}v)+r^{-2}v$, $\partial_{r}=(\partial_{r}+\frac{1}{r})-\frac{1}{r}$,
and $(\partial_{r}+\frac{1}{r})(\partial_{r}-\frac{1}{r})=\partial_{rr}$,
we deduce 
\begin{align}
\left\Vert {\bf 1}_{r\gtrsim1}r^{-1}\partial_{r}v\right\Vert _{L^{2}} & \lesssim\left\Vert \partial_{rr}v\right\Vert _{L^{2}}+\left\Vert {\bf 1}_{r\gtrsim1}r^{-2}v\right\Vert _{L^{2}}+\left\Vert {\bf 1}_{r\gtrsim1}\partial_{r}\left[\left(\partial_{r}-\tfrac{1}{r}\right)v\right]\right\Vert _{L^{2}}\nonumber \\
 & \lesssim\left\Vert \partial_{rr}v\right\Vert _{L^{2}}+\left\Vert {\bf 1}_{r\gtrsim1}r^{-2}v\right\Vert _{L^{2}}+\left\Vert {\bf 1}_{r\gtrsim1}r^{-1}\left(\partial_{r}-\tfrac{1}{r}\right)v\right\Vert _{L^{2}}.\label{eq:Appendixs cal H2 H2 r geq1 r-1drv}
\end{align}
By \eqref{eq:Appendixs cal H2 H2 r geq1 r-2v auxiliary1}, \eqref{eq:Appendixs cal H2 H2 r geq1 r-2v auxiliary2},
and \eqref{eq:Appendixs cal H2 H2 r geq1 r-1drv}, we conclude the
claim \eqref{eq:Appendixs cal H2 H2 r geq1 goal}.

For $0<r=x\lesssim1$, we claim that 
\begin{align}
\left\Vert {\bf 1}_{r\lesssim1}(||v|+|\partial_{r}v||)\right\Vert _{L^{2}}\lesssim\|\partial_{rr}v\|_{L^{2}}+\|{\bf 1}_{r\sim1}\partial_{r}v\|_{L^{2}}+\|{\bf 1}_{r\sim1}v\|_{L^{2}}.\label{eq:Appendix cal H2 H2 r leq1 goal}
\end{align}
and this immediately follows from Lemma~\ref{WeightedLogHardyIneq}
with $k=-\frac{1}{2}$ and $l=0$. 
\end{proof}
Next, we discuss the subcoercivity of $A_{Q}\tilde L_{Q}$. 
\begin{lem}[Subcoercivity for $A_{Q}\widetilde{L}_{Q}$ on $\dot{\mathcal{H}}^{2}$]
\label{lem:subcoercivity AQLQ}We have 
\begin{align*}
\|A_{Q}\widetilde{L}_{Q}v\|_{L^{2}}+\|{\bf 1}_{|x|\lesssim1}(|\partial_{x}v|+|v|)\|_{L^{2}}\sim\|v\|_{\dot{\mathcal{H}}^{2}}.
\end{align*}
\end{lem}

\begin{proof}
We recall that 
\begin{align}
A_{Q}\widetilde{L}_{Q}v=A_{Q}[\mathbf{D}_{Q}v+Q^{-1}\mathcal{H}(\text{Re}(Q^{3}v))].\label{eq:AQLQCoercivityRedefine}
\end{align}
We also note that by subcoercivity for $A_{Q}$ with \eqref{eq:CommuteHilbertDerivative}
and $(1+y^{2})Q^{2}=2$. We first observe that the second term is
perturbative; 
\begin{align*}
\|A_{Q}[Q^{-1}\mathcal{H}(\text{Re}(Q^{3}v))]\|_{L^{2}}\lesssim\|Q^{3}v\|_{L^{2}}+\|Q^{2}\partial_{x}v\|_{L^{2}}.
\end{align*}
Hence, it suffices to prove 
\begin{align*}
\|A_{Q}\mathbf{D}_{Q}v\|_{L^{2}}+\|{\bf 1}_{|x|\lesssim1}(|\partial_{x}v|+|v|)\|_{L^{2}}\sim\|v\|_{\dot{\mathcal{H}}^{2}}.
\end{align*}

We obtain ($\lesssim$) using the subcoercivity of $A_{Q}$ and
the definition of $\mathbf{D}_{Q}$. For ($\gtrsim$), we use Lemma~\ref{LemmaH2comparison}
to reduce to show 
\begin{align*}
\|\partial_{xx}v\|_{L^{2}}\lesssim\|A_{Q}\mathbf{D}_{Q}v\|_{L^{2}}+\|{\bf 1}_{|x|\lesssim1}(|\partial_{x}v|+|v|)\|_{L^{2}}.
\end{align*}
From the subcoercivity of $A_{Q}$, Lemma~\ref{LemmaAppendix AQ Subcoercivity},
we have 
\begin{align*}
\|A_{Q}\mathbf{D}_{Q}v\|_{L^{2}} & \sim\|\mathbf{D}_{Q}v\|_{\dot{\mathcal{H}}1}-\||D|^{\frac{1}{2}}Q\mathbf{D}_{Q}v\|_{L^{2}}\geq\|\partial_{x}\mathbf{D}_{Q}v\|_{L^{2}}-\||D|^{\frac{1}{2}}Q\mathbf{D}_{Q}v\|_{L^{2}}.
\end{align*}
By \eqref{eq:CommuteHilbertDerivative} with $(1+x^{2})\langle x\rangle^{2}=1$,
we have 
\begin{align*}
\||D|^{\frac{1}{2}}Q\mathbf{D}_{Q}v\|_{L^{2}}^{2} & =(|D|Q\mathbf{D}_{Q}v,Q\mathbf{D}_{Q}v)_{r}\\
 & =(\mathcal{H}[x\partial_{x}(Q\mathbf{D}_{Q}v)],x\langle x\rangle^{-2}Q\mathbf{D}_{Q}v)_{r}+(|D|Q\mathbf{D}_{Q}v,\langle x\rangle^{-2}Q\mathbf{D}_{Q}v)_{r}\\
 & \lesssim\|Q^{2}\mathbf{D}_{Q}v\|_{L^{2}}\|v\|_{\dot{\mathcal{H}}^{2}}\\
 & \lesssim(\|{\bf 1}_{|x|\leq r_{0}}(|\partial_{x}v|+|v|)\|_{L^{2}}+r_{0}^{-1}\|v\|_{\dot{\mathcal{H}}^{2}})\|v\|_{\dot{\mathcal{H}}^{2}}.
\end{align*}
Furthermore, we have 
\begin{align*}
\|{\bf 1}_{|x|\leq r_{0}}(|\partial_{x}v|+|v|)\|_{L^{2}}\|v\|_{\dot{\mathcal{H}}^{2}}\leq r_{0}\|{\bf 1}_{|x|\leq r_{0}}(|\partial_{x}v|+|v|)\|_{L^{2}}^{2}+r_{0}^{-1}\|v\|_{\dot{\mathcal{H}}^{2}}^{2}.
\end{align*}
Thus, taking $r_{0}>1$ sufficiently large, we reduce to prove 
\begin{align*}
\|\partial_{xx}v\|_{L^{2}}\lesssim\|\partial_{x}\mathbf{D}_{Q}v\|_{L^{2}}+\|{\bf 1}_{|x|\lesssim1}(|\partial_{x}v|+|v|)\|_{L^{2}}.
\end{align*}
By definition of $\mathbf{D}_{Q}$, we have 
\begin{align*}
\|\partial_{xx}v\|_{L^{2}}\leq\|\partial_{x}\mathbf{D}_{Q}v\|_{L^{2}}+\left\Vert \partial_{x}\left(x\langle x\rangle^{-2}v\right)\right\Vert _{L^{2}}.
\end{align*}
Hence, it suffices to show 
\begin{align*}
\left\Vert \partial_{x}\left(x\langle x\rangle^{-2}v\right)\right\Vert _{L^{2}}\lesssim\|\partial_{x}\mathbf{D}_{Q}v\|_{L^{2}}+\|{\bf 1}_{|x|\lesssim1}(|\partial_{x}v|+|v|)\|_{L^{2}}.
\end{align*}
We check that 
\begin{align*}
\partial_{x}\left(x\langle x\rangle^{-2}v\right)=x\langle x\rangle^{-2}\mathbf{D}_{Q}v+(1-2x^{2})\langle x\rangle^{-4}v.
\end{align*}
For $|x|\lesssim1$, we have $\Vert{\bf 1}_{|x|\lesssim1}\partial_{x}(x\langle x\rangle^{-2}v)\Vert_{L^{2}}\lesssim\|{\bf 1}_{|x|\lesssim1}(|\partial_{x}v|+|v|)\|_{L^{2}}$.
Therefore, we reduce the proof to show 
\begin{align}
\||x|^{-1}\mathbf{D}_{Q}v{\bf 1}_{|x|\gtrsim1}\|_{L^{2}}+\||x|^{-2}v{\bf 1}_{|x|\gtrsim1}\|_{L^{2}}\lesssim\|\partial_{x}\mathbf{D}_{Q}v\|_{L^{2}}+\|{\bf 1}_{|x|\lesssim1}(|\partial_{x}v|+|v|)\|_{L^{2}}.\label{eq:SubcoercivityAQLQgoal}
\end{align}
By symmetry, we restrict the domain to $\mathbb{R}^{+}$ with $0<r=x$.
Applying Lemma~\ref{WeightedHardyIneq} with $f=\mathbf{D}_{Q}v$
and $\varphi=r^{-1}$, we have 
\begin{align}
{\textstyle \int_{20}^{\infty}|r^{-1}\mathbf{D}_{Q}v|^{2}dr}\lesssim{\textstyle \int_{20}^{\infty}|\partial_{r}\mathbf{D}_{Q}v|^{2}dr}+\|{\bf 1}_{r\sim1}|v|_{-1}\|_{L^{2}}.\label{eq:SubcoercivityAQLQ r-1DQv}
\end{align}
Using Lemma~\ref{WeightedHardyIneq} with $f=\langle r\rangle v$
and $\varphi=r^{-5}$ and $\mathbf{D}_{Q}v=\langle r\rangle^{-1}\partial_{r}(\langle r\rangle v)$,
we deduce 
\begin{align}
{\textstyle \int_{20}^{\infty}|r^{-2}v|^{2}dr\sim\int_{20}^{\infty}|r^{-1}\langle r\rangle v|^{2}r^{-4}dr} & \lesssim{\textstyle \int_{20}^{\infty}|\partial_{r}(\langle r\rangle v)|^{2}r^{-4}dr}+\|{\bf 1}_{r\sim1}v\|_{L^{2}}\nonumber \\
 & \sim{\textstyle \int_{20}^{\infty}|r^{-1}\mathbf{D}_{Q}v|^{2}dr}+\|{\bf 1}_{r\sim1}v\|_{L^{2}}.\label{eq:SubcoercivityAQLQ r-2v}
\end{align}

By plugging \eqref{eq:SubcoercivityAQLQ r-2v} and \eqref{eq:SubcoercivityAQLQ r-1DQv}
into \eqref{eq:SubcoercivityAQLQgoal}, we complete the proof. 
\end{proof}

\section{Decomposition lemma}

\label{AppendixDecomposition} In this subsection, we provide the
proof of the decomposition lemma (Lemma~\ref{LemmaDecomposition}). 
\begin{proof}[Proof of Lemma~\ref{LemmaDecomposition}]
\textbf{Step 1.} We equip $\bbR_{+}$ with the metric $\text{dist}(\lambda_{1},\lambda_{2})=|\log(\frac{\lambda_{1}}{\lambda_{2}})|$,
and equip $\bbR/2\pi\bbZ$ with the induced metric of $\bbR$. We
first decompose $v$ by $[Q+\widehat{\eps}]_{\lambda,\gamma,x}$.
We define 
\begin{align*}
\mathbf{F}(\lambda,\gamma,x;v)=((\widehat{\eps},\mathcal{Z}_{1})_{r},(\widehat{\eps},\mathcal{Z}_{2})_{r},(\widehat{\eps},\mathcal{Z}_{3})_{r})^{T},
\end{align*}
where 
\begin{align*}
w\coloneqq\lambda^{1/2}e^{-i\gamma}v(\lambda\cdot+x),\quad\widehat{\eps}\coloneqq w-Q.
\end{align*}
We have 
\begin{align*}
\partial_{\lambda}\mathbf{F}_{k} & =\lambda^{-1}(\Lambda Q,[\mathcal{Z}_{k}]_{\lambda,\gamma,x})_{r}-\lambda^{-1}(v-Q,[\Lambda\mathcal{Z}_{k}]_{\lambda,\gamma,x})_{r},\\
\partial_{\gamma}\mathbf{F}_{k} & =(-iQ,[\mathcal{Z}_{k}]_{\lambda,\gamma,x})_{r}+(v-Q,[i\mathcal{Z}_{k}]_{\lambda,\gamma,x})_{r},\\
\partial_{x}\mathbf{F}_{k} & =\lambda^{-1}(Q_{y},[\mathcal{Z}_{k}]_{\lambda,\gamma,x})_{r}-\lambda^{-1}(v-Q,[\partial_{y}\mathcal{Z}_{k}]_{\lambda,\gamma,x})_{r}.
\end{align*}
We note that $\partial_{\lambda,\gamma,x}\mathbf{F}(1,0,0,;v)=\partial_{\lambda,\gamma,x}\mathbf{F}(1,0,0;Q)+M^{C}O(\|v-Q\|_{L^{2}})$.
\begin{align*}
\partial_{\lambda,\gamma,x}\mathbf{F}(1,0,0;Q)=\begin{pmatrix}(\Lambda Q,\mathcal{Z}_{1})_{r} & 0 & 0\\
0 & -(iQ,\mathcal{Z}_{2})_{r} & 0\\
0 & 0 & (Q_{y},\mathcal{Z}_{3})_{r}
\end{pmatrix}.
\end{align*}
In the neighborhood at $(1,0,0)$, we can check that all components
of $\partial_{\lambda,\gamma,x}\mathbf{F}(\lambda,\gamma,x;v)$ are
well defined and continuous for $\lambda,\gamma,x$ when $\|v-Q\|_{L^{2}}<\infty$.
Thus, $\mathbf{F}$ is $C^{1}$ for $(\lambda,\gamma,x)$ and is invertible
at $(1,0,0,Q)$. For given $M\gg1$, there exist $0<\delta_{dec}\ll\delta_{1}$
and $C^{1}$ maps $\textbf{G}_{1,0,0}:B_{\delta_{dec}}(Q)\to B_{\delta_{1}}(1,0,0)$
such that for given $v\in B_{\delta_{dec}}(Q)\subset L^{2}$, $\textbf{G}_{1,0,0}(v)$
are unique solutions to $\textbf{F}(\textbf{G}_{1,0,0}(v),v)=0$ in
$B_{\delta_{1}}(1,0,0)$.

Now, we want to find a unique $\textbf{G}:B_{\delta_{dec}}(Q)\to\mathbb{R}\times\mathbb{R}/2\pi\mathbb{Z}\times\mathbb{R}$
such that $\textbf{F}(\textbf{G}(v),v)=0$. Using scale, phase,
and translation invariances, we define $\textbf{G}_{\lambda,\gamma,x}$
for $(\lambda,\gamma,x)\in\mathbb{R}\times\mathbb{R}/2\pi\mathbb{Z}\times\mathbb{R}$
in the obvious way. We note that $\textbf{G}_{\lambda,\gamma,x}:B_{\delta_{dec}}(Q)_{\lambda,\gamma,x}\to B_{\delta_{1}}(\lambda,\gamma,x)$
where $B_{\delta_{dec}}(Q)_{\lambda,\gamma,x}=\{v:\|[Q]_{\lambda,\gamma,x}-v\|_{L^{2}}<\delta_{dec}\}$.
We define 
\begin{align*}
\textbf{G}\coloneqq\bigcup_{\lambda,\gamma,x}\textbf{G}_{\lambda,\gamma,x}:\bigcup_{\lambda,\gamma,x}B_{\delta_{dec}}(Q)_{\lambda,\gamma,x}\to\mathbb{R}\times\mathbb{R}/2\pi\mathbb{Z}\times\mathbb{R}.
\end{align*}

We first show that $\textbf{G}$ is well defined. If $v\in B_{\delta_{dec}}(Q)_{\lambda_{1},\gamma_{1},x_{1}}\cap B_{\delta_{dec}}(Q)_{\lambda_{2},\gamma_{2},x_{2}}$,
we have $\text{dist}((\lambda_{1},\gamma_{1},x_{1}),(\lambda_{2},\gamma_{2},x_{2}))\lesssim\delta_{dec}$.
Therefore, we have 
\begin{align*}
\text{dist}(\textbf{G}_{\lambda_{2},\gamma_{2},x_{2}}(v),(\lambda_{1},\gamma_{1},x_{1}))\lesssim\delta_{dec}\ll\delta_{1}.
\end{align*}
Since $\textbf{F}(\textbf{G}_{\lambda_{2},\gamma_{2},x_{2}}(v),v)=0$,
we deduce $\textbf{G}_{\lambda_{1},\gamma_{1},x_{1}}=\textbf{G}_{\lambda_{2},\gamma_{2},x_{2}}$
from the uniqueness of $\textbf{G}_{\lambda_{1},\gamma_{1},x_{1}}$
in $B_{\delta_{1}}(\lambda_{1},\gamma_{1},x_{1})$. Thus, $\textbf{G}$
is well defined.

Now, we turn to the uniqueness of $\textbf{G}$. We define the maps
\begin{align*}
\widehat{\eps}:\bigcup_{\lambda,\gamma,x}B_{\delta_{dec}}(Q)_{\lambda,\gamma,x}\to B_{\delta^{\prime}}(0)\subset L^{2}
\end{align*}
by $\widehat{\eps}(v)=v_{\lambda^{-1},-\gamma,-x}-Q$, where $(\lambda,\gamma,x)=\textbf{G}(v)$.
We claim that given $v\in\bigcup_{\lambda,\gamma,x}B_{\delta_{dec}}(Q)_{\lambda,\gamma,x}$,
$\textbf{G}(v)\in\mathbb{R}\times\mathbb{R}/2\pi\mathbb{Z}\times\mathbb{R}$
is the unique solution to $\textbf{F}(\textbf{G}(v),v)=0$ such that
$\|\widehat{\eps}\|_{L^{2}}<\delta^{\prime}$. Let $\textbf{G}^{\prime}=(\lambda^{\prime},\gamma^{\prime},x^{\prime})$
be a solution to $\textbf{F}(\textbf{G}^{\prime}(v),v)=0$ with $\widehat{\eps}^{\prime}(v)=v_{(\lambda^{\prime})^{-1},-\gamma^{\prime},-x^{\prime}}$
and $\|\widehat{\eps}^{\prime}\|_{L^{2}}<\delta^{\prime}$. If $\text{dist}(\textbf{G}(v),\textbf{G}^{\prime}(v))<\delta_{1}$,
then $\textbf{G}=\textbf{G}^{\prime}$ by the uniqueness of $\textbf{G}$.
If $\text{dist}(\textbf{G}(v),\textbf{G}^{\prime}(v))\geq\delta_{1}$,
then 
\begin{align*}
\|[Q]_{\lambda,\gamma,x}-[Q]_{\lambda^{\prime},\gamma^{\prime},x^{\prime}}\|_{L^{2}}\gtrsim\delta_{1}.
\end{align*}
However, this is a contradiction with $v=[Q+\widehat{\eps}]_{\lambda,\gamma,x}=[Q+\widehat{\eps}^{\prime}]_{\lambda^{\prime},\gamma^{\prime},x^{\prime}}$
and $\|\widehat{\eps}\|_{L^{2}},\|\widehat{\eps}^{\prime}\|_{L^{2}}<\delta^{\prime}\ll\delta_{1}$.\footnote{Here we choose $\delta^{\prime}$ such that $\delta^{\prime}\ll\delta_{1}$,
and then shrink $\delta_{dec}$ to define the map $\widehat{\eps}$.}

\textbf{Step 2.} We further decompose $\widehat{\eps}=P+\eps$. Since
$\widehat{\eps}(v)=v_{\lambda^{-1},-\gamma,-x}-Q$, where $(\lambda,\gamma,x)=\textbf{G}(v)$,
we have $(\widehat{\eps},\mathcal{Z}_{k})_{r}=0$ for $k=1,2,3$.
we claim that there is a unique $(b,\eta,\nu)$ such that 
\begin{align*}
(\eps,\mathcal{Z}_{k})_{r}=0\quad\text{for }k=1,2,\cdots,6.
\end{align*}
For $k=1,2,3$, we conclude the claim by the transversality and
orthogonality condition for $\widehat{\eps}$. For $k=4,5,6$, we
have 
\begin{align*}
0=(\eps,\mathcal{Z}_{k})_{r}=(\widehat{\eps},\mathcal{Z}_{k})_{r}+b(i\tfrac{y^{2}}{4}Q,\mathcal{Z}_{k})_{r}+\eta(\tfrac{1+y^{2}}{4}Q,\mathcal{Z}_{k})_{r}-\nu(i\tfrac{y}{2}Q,\mathcal{Z}_{k})_{r}.
\end{align*}
Again by the transversality, $\textbf{H}=(b,\eta,\nu)$ is uniquely
defined by 
\begin{align}
b=-\frac{(\widehat{\eps},\mathcal{Z}_{4})_{r}}{(i\frac{y^{2}}{4}Q,\mathcal{Z}_{4})_{r}},\quad\eta=-\frac{(\widehat{\eps},\mathcal{Z}_{5})_{r}}{(\frac{1+y^{2}}{4}Q,\mathcal{Z}_{5})_{r}},\quad\nu=\frac{(\widehat{\eps},\mathcal{Z}_{6})_{r}}{(i\frac{y}{2}Q,\mathcal{Z}_{6})_{r}}.\label{eq:DecompositionAppendix b eta nu definition}
\end{align}
We also have $b,\eta,\nu$ is $C^{1}$ from $\widehat{\eps}(v)=[v]_{\lambda^{-1},-\gamma,-x}-Q$,
where $(\lambda,\gamma,x)=\textbf{G}(v)$. Moreover, we have 
\begin{align*}
\|\langle y\rangle^{-2}\eps\|_{L^{2}}\lesssim\|\langle y\rangle^{-2}\widehat{\eps}\|_{L^{2}}+|b|+|\eta|+|\nu|\lesssim\delta^{\prime}
\end{align*}
by $\|\widehat{\eps}\|_{L^{2}}<\delta^{\prime}$. Thus, we can find
$\delta^{\prime\prime}>0$ so that $\|\langle y\rangle^{-2}\eps\|_{L^{2}}+|b|+|\eta|+|\nu|<\delta^{\prime\prime}$

\textbf{Step 3.} We note that the map $v\mapsto\eps(v)=[v]_{\lambda^{-1},-\gamma,-x}-Q-P(b,\eta,\nu)$
is continuous on the $H^{2}$-topology, that is, 
\begin{align*}
\eps:\bigcup_{\lambda,\gamma,x}B_{\delta_{dec}}^{H^{2}}(Q)_{\lambda,\gamma,x}\to B_{\delta^{\prime\prime}}^{\mathcal{Z}^{\perp}}(0),
\end{align*}
where 
\begin{align*}
B_{\delta^{\prime\prime}}^{\mathcal{Z}^{\perp}}(0)=\{\eps\in\langle y\rangle^{2}L^{2}:\|\langle y\rangle^{-2}\eps\|_{L^{2}}\leq\delta^{\prime\prime},\eps\in\mathcal{Z}^{\perp}\}.
\end{align*}
We note that the domain of $\eps(v)$ is $H^{2}$, and the range is
$\langle y\rangle^{2}L^{2}$. By the definition of $\eps$, the map
$v\mapsto(\mathbf{G}(v),\mathbf{H}(v),\eps(v))$ with $v\in\bigcup_{\lambda,\gamma,x}B_{\delta_{dec}}^{H^{2}}(Q)_{\lambda,\gamma,x}$
has a continuous left inverse 
\begin{align*}
\Phi:\mathbb{R}\times\mathbb{R}/2\pi\mathbb{Z}\times\mathbb{R}\times B_{\delta^{\prime\prime}}(0)^{3}\times B_{\delta^{\prime\prime}}^{\mathcal{Z}^{\perp}}(0)\to\langle y\rangle^{2}L^{2},\\
(\lambda,\gamma,x,b,\eta,\nu,\eps)\mapsto[Q+P(b,\eta,\nu)+\eps]_{\lambda,\gamma,x}.
\end{align*}
Furthermore, we have the uniqueness of $(\mathbf{G},\mathbf{H})$
since $\mathbf{G}$ and $\mathbf{H}$ are unique. This implies that
$\Im(\mathbf{G},\mathbf{H},\eps)=\Phi^{-1}(\bigcup_{\lambda,\gamma,x}B_{\delta_{dec}}^{H^{2}}(Q)_{\lambda,\gamma,x})$,
and it is open. Thus, $\Phi|_{\Im(\mathbf{G},\mathbf{H},\eps)}$ is
a right inverse of $(\mathbf{G},\mathbf{H},\eps)$. Therefore, the
restriction 
\begin{align*}
\Phi|_{\Im(\mathbf{G},\mathbf{H},\eps)}:\Im(\mathbf{G},\mathbf{H},\eps)\to\bigcup_{\lambda,\gamma,x}B_{\delta_{dec}}^{H^{2}}(Q)_{\lambda,\gamma,x}
\end{align*}
is a homeomorphism with the inverse $(\mathbf{G},\mathbf{H},\eps)$.

\textbf{Step 4.} We finish the proof.

(1) $\sim$ (3): We almost proved these in step $1$ and step $2$.
$C^{1}$ property on the $H^{2}$ topology is immediate from the embedding
$H^{2}\hookrightarrow L^{2}\hookrightarrow\langle y\rangle^{2}L^{2}$.

(4): Taking $b^{*}\ll\delta_{dec}$, we have $\overline{\mathcal{O}}_{init}\subseteq\bigcup_{\lambda,\gamma,x}B_{\delta_{dec}}^{H^{2}}(Q)_{\lambda,\gamma,x}$.
Due to the uniqueness of $(\mathbf{G},\mathbf{H})$, we have $\overline{\mathcal{U}}_{init}\subseteq\Im(\mathbf{G},\mathbf{H},\eps)$.
By Step 3, we conclude $\Phi|_{\overline{\mathcal{U}}_{init}}$ is
a homeomorphism from $\overline{\mathcal{U}}_{init}$ to $\overline{\mathcal{O}}_{init}$.
In particular, $\mathcal{O}_{init}$ is open since $\mathcal{U}_{init}$
is open. 
\end{proof}

\section{Proofs of some preliminary lemmas}

\label{AppendixPreliminaryProof} In this section, we prove some preliminary
lemmas given without proof in Section~\ref{sec:notation and preliminaries}.
The proofs are fairly standard. 
\begin{proof}[Proof of Lemma~\ref{LemmaHilbertUsefulEquation}]
We first assume that $f$ and $g$ are real-valued. We have for any (complex-valued)
function $h\in C_{c}^{\infty}$, $\Re{\textstyle \int_{\R}fg\overline{h}={\textstyle \int_{\R}fg\Re\overline{h}.}}$
We denote $\Re\overline{h}=\tilde h$, then $\tilde h$ is also real-valued.
Thus, we have 
\begin{align*}
{\textstyle \int_{\R}fg\tilde h} & ={\textstyle \int_{\R}(\Pi_{+}+\Pi_{-})f(\Pi_{+}+\Pi_{-})g(\Pi_{+}+\Pi_{-})\tilde h}\\
 & ={\textstyle \int_{\R}(\Pi_{+}f\Pi_{-}g+\Pi_{-}f\Pi_{+}g)\tilde h+\Pi_{+}f\Pi_{+}g\Pi_{-}\tilde h+\Pi_{-}f\Pi_{-}g\Pi_{+}\tilde h.}
\end{align*}
Using $\Pi_{+}=\tfrac{1}{2}(1+i\mathcal{H})$, we deduce 
\begin{align}
\Re{\textstyle \int_{\R}fg\overline{h}} & ={\textstyle \int_{\R}(\mathcal{H}f\cdot\mathcal{H}g-\mathcal{H}(f\cdot\mathcal{H}g+\mathcal{H}f\cdot g))\tilde h\nonumber}\\
 & =\Re{\textstyle \int_{\R}(\mathcal{H}f\cdot\mathcal{H}g-\mathcal{H}(f\cdot\mathcal{H}g+\mathcal{H}f\cdot g))\overline{h}},\label{eq:HilbertProductRuleWeakProve}
\end{align}
and this leads to \eqref{eq:HilbertProductRule} in a weak sense when
$f$ and $g$ are real-valued. We know that $\mathcal{H}$ commutes
with $i$, and \eqref{eq:HilbertProductRule} is linear for $f$ and
$g$. Thus, we conclude \eqref{eq:HilbertProductRule} for any complex-valued
$f$ and $g$ in a weak sense. Now, in \eqref{eq:HilbertProductRuleWeakProve},
we take $h$ by 
\begin{align*}
h=\mathcal{H}f\cdot\mathcal{H}g-\mathcal{H}(f\cdot\mathcal{H}g+\mathcal{H}f\cdot g)-fg.
\end{align*}
Then we obtain $h\in H^{1}$ for $f,g\in H^{1}$, and \eqref{eq:HilbertProductRule}
is satisfied in the $H^{1}$ sense. Thanks to the Sobolev embedding, we conclude
\eqref{eq:HilbertProductRule} in the pointwise sense. 
\end{proof}
\begin{proof}[Proof of Lemma~\ref{LemmaCommuteHilbert}]
We assume $f\in C_{c}^{\infty}$. Using \eqref{eq:HilbertIntegralFormula},
we have 
\begin{align*}
\mathcal{H}(xf)(x) & ={\textstyle \frac{1}{\pi}\int_{0+}^{\infty}\frac{(x-y)f(x-y)-(x+y)f(x+y)}{y}dy}\\
 & =x\cdot{\textstyle \frac{1}{\pi}\int_{0+}^{\infty}\frac{f(x-y)-f(x+y)}{y}dy+{\textstyle \frac{1}{\pi}\int_{0+}^{\infty}(-f(x-y)-f(x+y))dy}}\\
 & =x\mathcal{H}f(x)-{\textstyle \frac{1}{\pi}\int_{\mathbb{R}}f(x+y)dy.}
\end{align*}
Thus, we have \eqref{eq:CommuteHilbert} for $f\in C_{c}^{\infty}$.
Taking $f_{n}\to f$ in $\langle x\rangle^{-1}L^{2}$ with $f_{n}\in C_{c}^{\infty}$,
we conclude \eqref{eq:CommuteHilbert} in almost everywhere sense.

Assuming $f\in C_{c}^{\infty}$ and using \eqref{eq:CommuteHilbert},
we have 
\begin{align*}
[x,\mathcal{H}]\partial_{x}f(x)=\tfrac{1}{\pi}{\textstyle \int_{\mathbb{R}}\partial_{x}f(y)dy=\tfrac{1}{\pi}(f(+\infty)-f(-\infty))=0.}
\end{align*}
Again using the density argument, we deduce \eqref{eq:CommuteHilbertDerivative}
in almost everywhere sense. 
\end{proof}

\section{Unconditional Lax equation}

\label{AppendixUnconditionalLaxProof} In this subsection, we prove
the Lax pair equation for \eqref{CMdnls-gauged}, Proposition~\ref{PropositionUnconditionalLax}.
In fact, we show that this Lax structure holds for general $H^{1}$
solutions and then we deduce a Lax pair structure \eqref{eq:LaxPair}
for \eqref{CMdnls} \emph{without the chiral condition}. Here, we compute the Lax pair equation in \eqref{CMdnls-gauged}, but this is essentially the same as in \cite{GerardLenzmann2022}. 
With a suitable change
\begin{align*}
    -i\partial_{x}-\Pi_{+}u\Pi_{+}\overline{u}
    &\;\to\; {\mathcal{L}}_{\textnormal{Lax}}, 
    \\
    \Pi_{+}u \Pi_{+}\partial_x \overline{u}-\Pi_{+}\partial_x u\Pi_{+}\overline{u}
    +i(\Pi_{+} u \Pi_{+}{\overline{u}})^2
    &\;\to\; \td{\mathcal{P}}_{\textnormal{Lax}},
\end{align*}
where
\begin{align*}
    {\mathcal{L}}_{\textnormal{Lax}}&=-i\partial_{x}-u\Pi_{+}\overline{u},
    \\
    \td{\mathcal{P}}_{\textnormal{Lax}}
    &\coloneqq u\Pi_{+}{\partial_x \overline{u}}
    -\partial_x u\Pi_{+}{\overline{u}}
    +i(u \Pi_{+}{\overline{u}})^2,
\end{align*}
one can then verify that the computation in \cite{GerardLenzmann2022} remains valid without the chiral condition. This Lax pair $({\mathcal{L}}_{\textnormal{Lax}},\td{\mathcal{P}}_{\textnormal{Lax}})$ is equivalent to \eqref{eq:Appendix Lax structure equ}, since we have 
\begin{align*}
    \td{\mathcal{P}}_{\textnormal{Lax}}={\mathcal{P}}_{\textnormal{Lax}}+i{\mathcal{L}}_{\textnormal{Lax}}^2.
\end{align*}

\begin{proof}[Proof of Proposition~\ref{PropositionUnconditionalLax}]
Recall that 
\begin{align*}
\widetilde{\bfD}_{v}=\partial_{x}+\frac{1}{2}v\mathcal{H}(\overline{v}\ \cdot\ ),\quad H_{v}=-\partial_{xx}+\frac{1}{4}|v|^{4}-v|D|\overline{v}.
\end{align*}
For convenience, we use $2\widetilde{\bfD}_{v}$. We want to calculate
$[2\widetilde{\bfD}_{v},H_{v}]$. We have 
\begin{align*}
2\widetilde{\bfD}_{v}H_{v}f= & 2\partial_{x}H_{v}f+v\mathcal{H}(\overline{v}H_{v}f)\\
= & -2\partial_{xxx}f+\tfrac{1}{2}\partial_{x}(|v|^{4})f+\tfrac{1}{2}|v|^{4}\partial_{x}f-2\partial_{x}(v|D|\overline{v}f)+v\mathcal{H}(\overline{v}H_{v}f),
\end{align*}
and 
\begin{align*}
2H_{v}\widetilde{\bfD}_{v}f= & H_{v}[2\partial_{x}f+v\mathcal{H}(\overline{v}f)]\\
= & -2\partial_{xxx}f+\tfrac{1}{2}|v|^{4}\partial_{x}f-2v|D|(\overline{v}\partial_{x}f)+H_{v}[v\mathcal{H}(\overline{v}f)].
\end{align*}
Thus, 
\begin{align}
[2\widetilde{\bfD}_{v},H_{v}]f= & \tfrac{1}{2}\partial_{x}(|v|^{4})f-2(\partial_{x}v)|D|(\overline{v}f)-2v|D|(\overline{\partial_{x}v}f)\nonumber \\
 & +v\mathcal{H}(\overline{v}H_{v}f)-H_{v}[v\mathcal{H}(\overline{v}f)]\label{eq:Laxpair-proof-1}
\end{align}
Now, we claim that 
\begin{align}
v\mathcal{H}(\overline{v}H_{v}f)-H_{v}[v\mathcal{H}(\overline{v}f)]= & -\tfrac{1}{2}\partial_{x}(|v|^{4})f+2(\partial_{x}v)|D|(\overline{v}f)+2v|D|(\overline{\partial_{x}v}f)\label{eq:LaxpairClaimType1}\\
 & -(H_{v}v)\mathcal{H}(\overline{v}f)\label{eq:LaxpairClaimType2}\\
 & +v\mathcal{H}(\overline{H_{v}v}f),\label{eq:LaxpairClaimType3}
\end{align}
and thus $\eqref{eq:Laxpair-proof-1}=-(H_{v}v)\mathcal{H}(\overline{v}f)+v\mathcal{H}(\overline{H_{v}v}f)$.

We divide $v\mathcal{H}(\overline{v}H_{v}f)$ and $-H_{v}[v\mathcal{H}(\overline{v}f)]$
into three types: \eqref{eq:LaxpairClaimType1} type, \eqref{eq:LaxpairClaimType2}
type, and \eqref{eq:LaxpairClaimType3} type.

Using $H_{v}=-\partial_{xx}+\frac{1}{4}|v|^{4}-v|D|\overline{v}$,
we break $v\mathcal{H}(\overline{v}H_{v}f)$ into 
\begin{align}
v\mathcal{H}(\overline{v}H_{v}f)= & -v\mathcal{H}(\overline{v}\partial_{xx}f)-v\mathcal{H}[|v|^{2}|D|(\overline{v}f)]\tag{\eqref{eq:LaxpairClaimType1} \text{type}}\\
 & +v\mathcal{H}(\tfrac{1}{4}|v|^{4}\overline{v}f).\tag{\eqref{eq:LaxpairClaimType3} \text{type}}
\end{align}
Also, by $H_{v}=-\partial_{xx}+\frac{1}{4}|v|^{4}-v|D|\overline{v}$,
we calculate $-H_{v}[v\mathcal{H}(\overline{v}f)]$ as 
\begin{align}
-H_{v}[v\mathcal{H}(\overline{v}f)]= & 2(\partial_{x}v)|D|(\overline{v}f)\tag{\eqref{eq:LaxpairClaimType1} \text{type}}\\
 & +(\partial_{xx}v)\mathcal{H}(\overline{v}f)-(\tfrac{1}{4}|v|^{4}v)\mathcal{H}(\overline{v}f)\tag{\eqref{eq:LaxpairClaimType2} \text{type}}\\
 & +v|D|\{|v|^{2}\mathcal{H}(\overline{v}f)\}\label{eq:LaxpairClaimFurtherDecompose1}\\
 & +v|D|\partial_{x}(\overline{v}f).\label{eq:LaxpairClaimFurtherDecompose2}
\end{align}
For \eqref{eq:LaxpairClaimFurtherDecompose1}, we have 
\begin{align}
\eqref{eq:LaxpairClaimFurtherDecompose1}= & v\mathcal{H}[|v|^{2}|D|(\overline{v}f)]\tag{\eqref{eq:LaxpairClaimType1} \text{type}}\\
 & +v\mathcal{H}[\partial_{x}(|v|^{2})\mathcal{H}(\overline{v}f)].\label{eq:LaxpairClaimFurtherDecompose1-1}
\end{align}
Using the product rule for $\mathcal{H}$ \eqref{eq:HilbertProductRule}
with $\partial_{x}(|v|^{2})$ and $\overline{v}f$, we have 
\begin{align*}
\eqref{eq:LaxpairClaimFurtherDecompose1-1}= & -v\mathcal{H}[|D|(|v|^{2})(\overline{v}f)]-v\partial_{x}(|v|^{2})(\overline{v}f)+v|D|(|v|^{2})\cdot\mathcal{H}(\overline{v}f)\\
= & -\tfrac{1}{2}\partial_{x}(|v|^{4})f\tag{\eqref{eq:LaxpairClaimType1} \text{type}}\\
 & +v|D|(\overline{v}\cdot v)\cdot\mathcal{H}(\overline{v}f)\tag{\eqref{eq:LaxpairClaimType2} \text{type}}\\
 & +v\mathcal{H}[\overline{-v|D|(\overline{v}\cdot v)}f].\tag{\eqref{eq:LaxpairClaimType3} \text{type}}
\end{align*}
Our remaining term is \eqref{eq:LaxpairClaimFurtherDecompose2}. We
have 
\begin{align}
\eqref{eq:LaxpairClaimFurtherDecompose2}=v|D|(\overline{\partial_{x}v}f)+v|D|(\overline{v}\partial_{x}f)= & v|D|(\overline{\partial_{x}v}f)+v\mathcal{H}(\overline{v}\partial_{xx}f)\tag{\eqref{eq:LaxpairClaimType1} \text{type}}\\
 & +v\mathcal{H}(\overline{\partial_{x}v}\partial_{x}f).\label{eq:LaxpairClaimFurtherDecompose3}
\end{align}
We also have 
\begin{align*}
\eqref{eq:LaxpairClaimFurtherDecompose3}= & v|D|(\overline{\partial_{x}v}f)\tag{\eqref{eq:LaxpairClaimType1} \text{type}}\\
 & -v\mathcal{H}(\overline{\partial_{xx}v}f).\tag{\eqref{eq:LaxpairClaimType3} \text{type}}
\end{align*}
Now, collecting each type of terms, we conclude the claim.

Applying our claim to \eqref{eq:Laxpair-proof-1}, we deduce 
\begin{align*}
[2\widetilde{\bfD}_{v},H_{v}]f=-(H_{v}v)\mathcal{H}(\overline{v}f)+v\mathcal{H}(\overline{H_{v}v}f).
\end{align*}
Since \eqref{CMdnls-gauged} can be written as $-i\partial_{t}v=-H_{v}v$,
we have 
\begin{align*}
-i\partial_{t}2\widetilde{\bfD}_{v}f=-(H_{v}v)\mathcal{H}(\overline{v}f)+v\mathcal{H}(\overline{H_{v}v}f).
\end{align*}
Therefore, we conclude 
\begin{align*}
-i\partial_{t}2\widetilde{\bfD}_{v}f=2\widetilde{\bfD}_{v}H_{v}f-2H_{v}\widetilde{\bfD}_{v}f,
\end{align*}
and this finishes the proof. 
\end{proof}
Now, we prove that \eqref{eq:LaxPair} is a Lax pair on $H^{1}$.
Denote $\theta=\theta(t,x)\coloneqq-\frac{1}{2}\int_{-\infty}^{x}|u(t,y)|^{2}dy$.
Then we have $v(t,x)=\mathcal{G}(u)(t,x)=-u(t,x)e^{i\theta(t,x)}$. 
\begin{prop}
\label{Appendix Proposition Laxpair upto gauge} Let $u$ be the solution
to \eqref{CMdnls} in $H^{1}$, and $v=-\mathcal{G}(u)$. Then, we
have 
\begin{align}
e^{-i\theta}(-i\widetilde{\D}_{v})e^{i\theta} & =\mathcal{L}_{\textnormal{Lax}}=-i\partial_{x}-u\Pi_{+}\overline{u},\label{eq:Appendix Lax L equ}\\
e^{-i\theta}(-iH_{v})e^{i\theta} & =\mathcal{P}_{\textnormal{Lax}}+i\partial_{t}\theta=i\partial_{xx}+i2uD_{+}\overline{u}+i\partial_{t}\theta.\label{eq:Appendix Lax P equ}
\end{align}
Moreover, we have 
\begin{align}
\partial_{t}(\mathcal{L}_{\textnormal{Lax}})=[\mathcal{P}_{\textnormal{Lax}},\mathcal{L}_{\textnormal{Lax}}]\label{eq:Appendix Lax structure equ}
\end{align}
without the chiral condition, that is, \eqref{eq:Appendix Lax structure equ}
holds true for any $u\in H^{1}$. 
\end{prop}

\begin{proof}
By direct calculation, we have 
\begin{align}
e^{-i\theta}\widetilde{\D}_{v}e^{i\theta} & =\widetilde{\D}_{u}+i\theta_{x},\label{eq:Appendix Lax Dtilde equ}\\
e^{-i\theta}H_{v}e^{i\theta} & =H_{u}-2i\theta_{x}\partial_{x}+\theta_{x}^{2}-i\theta_{xx}.\label{eq:Appendix Lax H equ}
\end{align}
Putting $\theta=-\frac{1}{2}\int_{-\infty}^{x}|u|^{2}dy$ in \eqref{eq:Appendix Lax Dtilde equ}
with $\Pi_{+}=\frac{1}{2}(1+i\mathcal{H})$, we deduce \eqref{eq:Appendix Lax L equ}.
In addition, we compute 
\begin{align*}
e^{-i\theta}(-iH_{v})e^{i\theta}= & -iH_{u}-2\theta_{x}\partial_{x}-i\theta_{x}^{2}-\theta_{xx}\\
= & i\partial_{xx}+iu|D|\overline{u}+|u|^{2}\partial_{x}+\tfrac{1}{2}(u\overline{\partial_{x}u}+\overline{u}\partial_{x}u)-\tfrac{i}{2}|u|^{4}\\
= & i\partial_{xx}+iu|D|\overline{u}+u\partial_{x}(\overline{u}\ \cdot\ )+\tfrac{1}{2}(-u\overline{\partial_{x}u}+\overline{u}\partial_{x}u)-\tfrac{i}{2}|u|^{4}\\
= & i\partial_{xx}+i2uD_{+}\overline{u}+\tfrac{1}{2}(-u\overline{u_{x}}+\overline{u}u_{x})-\tfrac{i}{2}|u|^{4}.
\end{align*}
Also, we compute 
\begin{align*}
i\theta_{t}=-i\Re(i\overline{u}u_{x})-\tfrac{i}{2}|u|^{4}=\tfrac{1}{2}(-u\overline{u_{x}}+\overline{u}u_{x})-\tfrac{i}{2}|u|^{4},
\end{align*}
and then we obtain \eqref{eq:Appendix Lax P equ}.

Now, we check \eqref{eq:Appendix Lax structure equ}. By \eqref{eq:Appendix Lax L equ},
\eqref{eq:Appendix Lax P equ}, and \eqref{eq:LaxEqu Unconditional},
we have 
\begin{align*}
\partial_{t}(\mathcal{L}_{\textnormal{Lax}})= & \partial_{t}(e^{-i\theta}(-i\widetilde{\D}_{v})e^{i\theta})\\
= & -i\theta_{t}\mathcal{L}_{\textnormal{Lax}}+\mathcal{L}_{\textnormal{Lax}}i\theta_{t}+e^{-i\theta}\partial_{t}(-i\widetilde{\D}_{v})e^{i\theta}\\
= & [-i\theta_{t},\mathcal{L}_{\textnormal{Lax}}]+e^{-i\theta}[-iH_{v},-i\widetilde{\mathbf{D}}_{v}]e^{i\theta}\\
= & [-i\theta_{t},\mathcal{L}_{\textnormal{Lax}}]+[\mathcal{P}_{\textnormal{Lax}}+i\theta_{t},\mathcal{L}_{\textnormal{Lax}}]=[\mathcal{P}_{\textnormal{Lax}},\mathcal{L}_{\textnormal{Lax}}].
\end{align*}
This finishes the proof. 
\end{proof}

\section{Some algebraic identities}

\label{AppendixAlgebraicIdentity}

In this section, we collect various algebraic identities for the convenience
of the readers. Some identities may not be used in this paper.

We note some formulae related to the soliton $Q$. We compute 
\[
\mathcal{G}(\mathcal{R})=\frac{\sqrt{2}}{x+i}e^{-\frac{i}{2}\int_{-\infty}^{x}|\mathcal{R}(y)|^{2}dy}=\frac{\sqrt{2}}{x+i}\frac{\sqrt{1+x^{2}}}{i-x}=-Q.
\]
We have 
\begin{align*}
\Lambda Q & =\tfrac{1-y^{2}}{\sqrt{2}(1+y^{2})^{3/2}},\quad Q_{y}=-\tfrac{\sqrt{2}y}{(1+y^{2})^{3/2}}=-\tfrac{1}{2}yQ^{3},\\
\Lambda(yQ) & =\tfrac{3y+y^{3}}{\sqrt{2}(1+y^{2})^{3/2}}=\tfrac{1}{2}yQ+\tfrac{1}{2}yQ^{3}.
\end{align*}
We note the Fourier transforms of some functions. 
\begin{align*}
\mathcal{F}(\langle y\rangle^{-1}Q)(\xi) & =\sqrt{2}\pi e^{-|\xi|},\\
\begin{split}\mathcal{F}(\langle y\rangle^{-3}Q)(\xi) & =\tfrac{\pi}{\sqrt{2}}(1+|\xi|)e^{-|\xi|},\\
\mathcal{F}(\langle y\rangle^{-3}\Lambda Q)(\xi) & =\tfrac{\pi}{4\sqrt{2}}e^{-|\xi|}(\xi^{2}+|\xi|+1),
\end{split}
\begin{split}\mathcal{F}(\langle y\rangle^{-3}yQ)(\xi) & =-i\tfrac{\pi}{\sqrt{2}}\xi e^{-|\xi|},\\
\mathcal{F}(\langle y\rangle^{-3}Q_{y})(\xi) & =i\tfrac{\pi}{4\sqrt{2}}\xi e^{-|\xi|}(1+|\xi|).
\end{split}
\end{align*}
We note some identities with respect to the Hilbert transform $\mathcal{H}$.
Using $\mathbf{D}_{Q}Q=0$. we have 
\begin{align}
\mathcal{H}(\tfrac{1}{1+y^{2}})=\tfrac{y}{1+y^{2}},\quad\text{or}\quad\mathcal{H}(Q^{2})=yQ^{2}.\label{eq:AppendixAlgebraic DQQequal0}
\end{align}
By differentiating \eqref{eq:AppendixAlgebraic DQQequal0}, we deduce
\begin{align*}
\begin{split}\mathcal{H}(\tfrac{2}{(1+y^{2})^{2}}) & =\tfrac{3y+y^{3}}{(1+y^{2})^{2}},\\
\mathcal{H}(\tfrac{2y^{2}}{(1+y^{2})^{2}}) & =\tfrac{y^{3}-y}{(1+y^{2})^{2}},
\end{split}
\begin{split}\mathcal{H}(\tfrac{2y}{(1+y^{2})^{2}}) & =\tfrac{y^{2}-1}{(1+y^{2})^{2}},\\
\mathcal{H}(\tfrac{2y^{3}}{(1+y^{2})^{2}}) & =-\tfrac{3y^{2}+1}{(1+y^{2})^{2}}.
\end{split}
\end{align*}
We have some (\emph{formal}) identities related to $|D|=\partial_{x}\mathcal{H}=\mathcal{H}\partial_{x}$.
\begin{align}
|D|(yQ^{2})=yQ^{4},\quad|D|(y^{2}Q^{2})=2\tfrac{y^{2}-1}{(1+y^{2})^{2}},\quad|D|(y^{3}Q^{2})=2\mathcal{H}(1)-yQ^{4}.\label{eq:AppendixAlgebraic |D| equ}
\end{align}
In fact, \eqref{eq:AppendixAlgebraic |D| equ} follows from 
\begin{align*}
|D|(yQ^{2}) & =\partial_{y}\mathcal{H}(yQ^{2})=-\partial_{y}(Q^{2})=yQ^{4},\\
|D|(y^{2}Q^{2}) & =|D|(y^{2}Q^{2}-2)=-\partial_{y}(yQ^{2})=2\tfrac{y^{2}-1}{(1+y^{2})^{2}}\\
|D|(y^{3}Q^{2}) & =2\mathcal{H}\partial_{y}(y-\tfrac{y}{1+y^{2}})=2\mathcal{H}(1)-|D|yQ^{2}=2\mathcal{H}(1)-yQ^{4}.
\end{align*}
From this, we deduce 
\begin{align*}
Q|D|(y^{3}Q^{2})-y^{2}Q|D|(yQ^{2})=2Q\mathcal{H}(1)-2yQ^{3}.
\end{align*}
Using \eqref{eq:AppendixAlgebraic |D| equ}, we also deduce 
\begin{align*}
|D|QP_{1}=\tfrac{ib+\eta}{2}\partial_{y}(Q^{2})+\tfrac{i\nu+\mu}{2}\partial_{y}(yQ^{2})\lesssim|b|Q^{3}+|b|^{1-\kappa}Q^{2}.
\end{align*}
We note some identities for linear operators. 
\begin{itemize}
\item $\mathbf{D}_{Q}$: 
\begin{align*}
\begin{split}\mathbf{D}_{Q}(Q) & =0,\\
\mathbf{D}_{Q}(yQ) & =Q,
\end{split}
\begin{split}\mathbf{D}_{Q}(Q_{y}) & =Q\cdot\tfrac{y^{2}-1}{(1+y^{2})^{2}},\\
\mathbf{D}_{Q}(y^{2}Q) & =2yQ,
\end{split}
\quad\begin{split}\mathbf{D}_{Q}(\Lambda Q) & =Q^{2}Q_{y},\\
\mathbf{D}_{Q}(Q^{-1}) & =yQ.
\end{split}
\end{align*}
\item $L_{Q}$: 
\begin{align*}
\begin{gathered}L_{Q}iQ=L_{Q}Q_{y}=L_{Q}\Lambda Q=0,\\
L_{Q}Q=yQ^{3},\quad L_{Q}Q^{-1}=yQ+Q\mathcal{H}(1),\quad L_{Q}iy^{2}Q=2L_{Q}iQ^{-1}=2iyQ.
\end{gathered}
\end{align*}
\item $L_{Q}^{*}$: 
\begin{align*}
L_{Q}^{*}Q=L_{Q}^{*}iQ^{-1}=0,\quad L_{Q}^{*}Q^{-1}=-Q\mathcal{H}(1),\\
L_{Q}^{*}iQ=-2iQ_{y},\quad L_{Q}^{*}yQ=Q,\quad L_{Q}^{*}iyQ=-2i\Lambda Q.
\end{align*}
Moreover, we have 
\begin{align*}
L_{Q}^{*}P_{1}=bi\Lambda Q-\tfrac{\eta}{2}Q-\nu iQ_{y}.
\end{align*}
\item $\widetilde{L}_{Q}$: 
\begin{align}
\widetilde{L}_{Q}\Lambda Q=\tfrac{1}{4}yQ,\quad\widetilde{L}_{Q}Q_{y}=-\tfrac{1}{4}Q,\quad\widetilde{L}_{Q}Q^{-1}=2yQ.\label{eq:Appendix LQtilde identity}
\end{align}
Using \eqref{eq:LQ LQtilde Difference}, we deduce the first two identities
of \eqref{eq:Appendix LQtilde identity}. The last identity is obtained
from the definition of $\widetilde{L}_{Q}$. 
\end{itemize}
 \bibliographystyle{abbrv}
\bibliography{referenceCM}

\end{document}